\documentclass[11pt,leqno]{article}

\usepackage{amsmath}
\usepackage{amssymb}
\usepackage{dsfont}
\usepackage{mathrsfs}
\usepackage{epsfig}
\usepackage{float}

\numberwithin{equation}{section}

\newtheorem{Theorem}{Theorem}
\newtheorem{Corollary}[Theorem]{Corollary}
\newtheorem{Lemma}[Theorem]{Lemma}
\newtheorem{Proposition}[Theorem]{Proposition}
\newtheorem{Definition}[Theorem]{Definition}
\newtheorem{remark}[Theorem]{Remark}

\numberwithin{Theorem}{section}

\bibstyle{plain}

\newcommand{\al}{\alpha}
\newcommand{\R}{{\mathbf R}}
\newcommand{\ds}{\displaystyle}
\newcommand{\e}{\varepsilon}
\newcommand{\Om}{\Omega}
\newcommand{\lra}{\longrightarrow}
\newcommand{\ra}{\rightarrow}
\newcommand{\p}{\partial}
\newcommand{\la}{\lambda}
\newcommand{\g}{\gamma}
\newcommand{\ov}{\overline}
\newcommand{\C}{{\mathbf C}}

\newcommand{\N}{{\mathbf N}}
\newcommand{\Z}{{\mathbf Z}}
\newcommand{\Q}{{\mathbf Q}}
\newcommand{\E}{{\mathbf E}}
\newcommand{\1}{{\mathds 1}}

\newcommand{\si}{\sigma}

\newcommand{\de}{\delta}

\newcommand{\ph}{\varphi}

\newcommand{\Po}{{\mathcal P}}
\newcommand{\Xo}{{\mathcal X}}
\newcommand{\Yo}{{\mathcal Y}}
\newcommand{\So}{{\mathcal S}}

\newcommand{\Lo}{{\mathcal L}}
\newcommand{\Do}{{\mathcal D}}
\newcommand{\Co}{{\mathcal C}}
\newcommand{\Eo}{{\mathcal E}}
\newcommand{\Oo}{{\mathcal O}}
\newcommand{\Sr}{{\mathscr{S}}}
\newcommand{\DR}{{\mathcal D}^{1,2}(\R^N)}

\newcommand{\wh}{\widehat}

\newcommand{\ii}{\int}

\newcommand{\beq}{\begin{equation}}
\newcommand{\eeq}{\end{equation}}

\newcommand{\references}[1]{\theinstitutions 
\begin{thebibliography}{#1}}

\newcommand{\Endrefs}{\end{thebibliography}}

\title{ Traveling waves for nonlinear Schr\"odinger equations with 
nonzero conditions at infinity, II}

\date{  }

\author{ David CHIRON\footnote{Laboratoire J.A. Dieudonn\'e UMR 6621, Universit{\'e} 
de Nice-Sophia Antipolis, Parc Valrose, 06108 Nice Cedex 02, France. 
{\sf e-mail}: chiron@unice.fr.} 
\quad \quad  and \quad \quad  
Mihai MARI\c S \footnote{Institut de Math\'ematiques de 
Toulouse UMR 5219, Universit\'e Paul Sabatier, 118 route de Narbonne, 31062 Toulouse 
Cedex 9, France, and Institut Universitaire de France. {\sf e-mail}: mihai.maris@math.univ-toulouse.fr.}}

\begin{document}

\maketitle

\vspace{-10pt}

\begin{abstract}
We prove the existence of nontrivial finite energy traveling waves for a large class of 
nonlinear Schr\"odinger equations with  nonzero conditions at infinity (includindg the Gross-Pitaevskii 
and the so-called "cubic-quintic" equations)
in space dimension $ N \geq 2$.
We show that minimization of the energy at fixed momentum can be used whenever the 
associated nonlinear potential  is nonnegative
and it gives a set of orbitally stable traveling waves, while  minimization of 
the action at constant kinetic energy can be used in all cases.
We also explore the relationship between the families of traveling waves obtained by different methods and 
 we prove a sharp nonexistence result for traveling waves with small energy.

\smallskip


\smallskip


\end{abstract}

{
\rm 

\tableofcontents
}

\section{Introduction} 
\label{intro}

We consider  the 
nonlinear Schr\"odinger equation
\beq
\label{1.1}
i \frac{ \p \Phi}{\p t} + \Delta \Phi + F( | \Phi |^2) \Phi = 0
\qquad \mbox{ in } \R^N,
\eeq
where $ \Phi$ is a complex-valued function on $\R^N$  satisfying the "boundary condition" 
$|\Phi | \lra r_0 $ as $ |x| \lra \infty $, 
$r_0 >0$ and $F$ is a real-valued function on $\R_+$ such that  $ F(r_0^2) = 0$.

Equation  (\ref{1.1}),  
with the considered non-zero conditions at infinity,  arises in the modeling of a 
great variety of physical phenomena  such as  superconductivity, superfluidity 
in Helium II, phase transitions and Bose-Einstein condensate
(\cite{AHMNPTB}, \cite{barashenkov2}, \cite{barashenkov1}, \cite{berloff}, \cite{coste}, \cite{GR}, 
\cite{gross}, \cite{IS}, \cite{JR}, \cite{JPR}, \cite{RB}).
In nonlinear optics, it appears in the context of dark solitons 
(\cite{KL}, \cite{KPS}), which are localized  nonlinear waves  moving on a stable, nonzero background  at rest at infinity. 
Two important model cases for (1.1) have been extensively studied 
both in the physical and  mathematical literature:
the Gross-Pitaevskii equation (where $F(s)=  1-s$)
and the so-called "cubic-quintic"  Schr\"odinger equation (where
$F(s) = - \al _1 + \al _3 s - \al _5 s ^2 $, $\; \al_1, \, \al _3, \, \al _5$
are positive and
$F$ has two positive roots).

In contrast to the case of zero boundary conditions at infinity 
(when the dynamics associated to (\ref{1.1}) is essentially governed 
by dispersion and scattering), the non-zero boundary conditions  allow a 
much richer dynamics and 
give rise to a remarkable variety of special solutions, such as 
traveling waves, standing waves or vortex solutions.  

\medskip

Using the Madelung transformation
$\Phi(x,t) = \sqrt{\rho(x,t)} e^{i\theta (x,t )}$
(which is well-defined in any  domain  where $ \Phi \neq 0$), 
equation (\ref{1.1}) is equivalent to the system
$$ 
\left\{\begin{array}{ll}
\displaystyle{ \p_t \rho + 2 {\rm div} ( \rho \nabla \theta ) } = 0 \\ \ \\ 
\displaystyle{ \p_{t} \nabla \theta + 2 ( \nabla \theta \cdot \nabla) \nabla \theta 
- \nabla ( F ( \rho) ) = \nabla \left( \frac{\Delta \rho}{ 2 \rho} - \frac{|\nabla \rho|^2}{ 4 \rho} \right) } .
\end{array}\right. 
$$
This is the system of Euler's equations for a compressible inviscid fluid of density $ \rho $ 
and velocity $ 2 \nabla \theta$ with an additional dispersive term usually  called {\it quantum pressure.} 
It has been shown that, if $F$ is $C^1$ near $ r_0^2$, 
 $ \; F'(r_0^2) < 0$ and the density varies slowly at infinity, the sound velocity at infinity associated to (\ref{1.1})  is 
$ v_s = r_0 \sqrt{ - 2 F'(r_0^2) } $ (see the introduction of \cite{M8}).

If $ F'( r_0^2) <0 $ (which means that (\ref{1.1}) is defocusing), 
a simple scaling enables us to assume that $ r_0 = 1 $ and $ F'(r_0 ^2) = -1$ (see \cite{M10}, p. 108); 
we will do so throughout the rest of this paper. The sound velocity at infinity is then $ v_s = \sqrt{2}$.

Equation (\ref{1.1}) has a  Hamiltonian structure.  
Indeed, let  $ V( s) =  \ii_s^{1} F( \tau ) \, d \tau$. 
It is then easy to see that, at least formally, the "energy"
\beq
\label{1.2}
E(\Phi ) = \int_{\R^N} | \nabla \Phi |^2 \, dx
 + \int_{\R^N} V( |\Phi |^2) \, dx
\eeq
is  conserved.
Another quantity   conserved  by the flow of (\ref{1.1}) is the momentum, 
$ \mathbf{P} ( \Phi) = ( P_1( \Phi), \dots , P_N( \Phi ))$.
A rigorous definition of the momentum will be given in the next section. 
If $\Phi$ is a  function sufficiently localized in space, 
we have $P_k( \Phi ) = \ii _{\R^N} \langle i \Phi_{x_k}, \Phi \rangle\, dx $, 
where $ \langle \cdot , \cdot \rangle $ is the usual scalar product in 
$\C \simeq \R^2$.

In a series of papers   (see, e.g., \cite{barashenkov2},  \cite{barashenkov1}, 
\cite{GR}, \cite{JR}, \cite{JPR}), 
particular attention has been paid to the traveling  waves of (\ref{1.1}). 
These are solutions of the form $ \Phi (x, t) = \psi (x + ct \omega)$, where 
$ \omega \in S^{N-1}$ is the direction of propagation and $ c \in \R^*$ is 
the speed of the traveling wave. 
They are supposed to play an important role in the dynamics of (\ref{1.1}).
We say that $\psi $ has finite energy 
if $ \nabla \psi \in L^2(\R^N)$ and $ V(|\psi |^2) \in L^1(  \R^N)$. 
Since the equation (\ref{1.1}) is rotation invariant, we may assume that 
$ \omega = (1, 0, \dots, 0)$. Then a traveling wave of speed $c $ satisfies 
the equation 
\beq
\label{1.3}
i c \frac{\p \psi}{\p x_1} + \Delta \psi + F(|\psi|^2) \psi  = 0 
\qquad \mbox{ in  } \R^N.
\eeq
It is obvious that a function $ \psi $ satisfies (\ref{1.3}) for some velocity 
$c$ if and only if $\psi( - x_1, x')$ satisfies (\ref{1.3}) with $c$ replaced 
by $-c$. Hence it suffices to consider the case $ c \geq 0$.  

\medskip

\begin{figure}[H]
\label{diaggrospit}
\includegraphics[width=7cm, height=4.5cm]{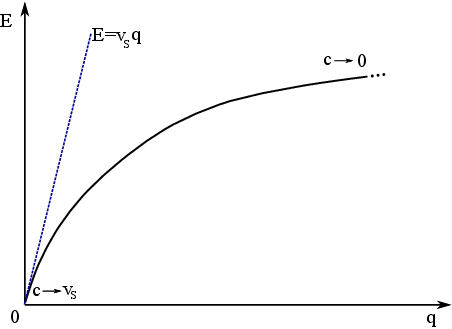}
\qquad \quad
\includegraphics[width=7cm, height=4.5cm]{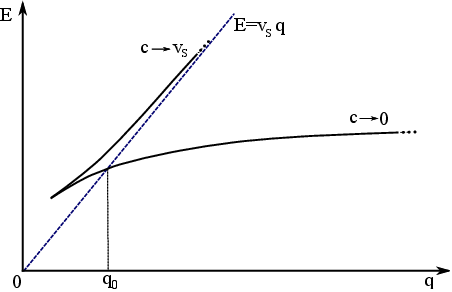}

\caption{{\small Energy ($E$) - momentum ($q$)  diagrams for (GP): (a) in dimension $2$; (b) in dimension $3$.}}
\end{figure}

The formal computations and numerical experiments led to a list of conjectures, often called 
{\it the Roberts programme, }
about the existence, the qualitative properties, the stability/instability and  the role of 
traveling waves into the dynamics of (\ref{1.1}); 
see \cite{BS2} or the introduction of \cite{M10} for a brief presentation.  
 It has been  conjectured 
that finite energy traveling waves of speed $c$ exist only for subsonic speeds: 
$ c < v_s$. The nonexistence of traveling waves for supersonic speeds ($c > v_s$)
has been rigorously proven  (see \cite{M8} and references therein). 
The numerical investigation of the traveling waves 
of the Gross-Pitaevskii equation ($F(s) = 1-s$) has been carried out 
in \cite{JR}. The method  used there  was a continuation argument with 
respect to the speed, solving \eqref{1.3} by the Newton  algorithm. 
Denoting $Q(\psi ) = P_1( \psi)= q$ the momentum of $\psi $ with respect 
to the $x_1-$direction, the representation of  solutions in the 
energy vs. momentum diagram gives the curves  in figure \ref{diaggrospit} below
(the straight line is the line $E = v_s q$). 
Notice that in dimension $N=2$ the curve is concave, while in dimension $N=3$ it consists of two branches, 
the lower branch being concave and the upper branch being convex. 

The rigorous proof of the existence of traveling waves has been a 
long lasting problem and was considered in a series of papers, 
see \cite{BS}, \cite{BOS}, \cite{chiron}, \cite{BGS}, \cite{M10}.  
At least formally, traveling waves are critical points of the functional 
$E - cQ$. Therefore, it is  a natural idea to look for such solutions 
as minimizers of the energy at fixed momentum, the speed $c$ being  then the 
Lagrange multiplier associated to the minimization problem. In the case of 
the Gross-Pitaevskii equation, in view of the above diagrams, this method 
is expected to give the full curve of traveling waves if $N=2$ and 
only the lower part that lies under the line $E = v_s q$ if $N=3$ (it is clear that  minimizers of $E$ 
at fixed $Q$ cannot lie on the upper branch). On a rigorous level, 
minimizing the energy at fixed momentum has been used first in \cite{BOS} to 
construct a sequence of traveling waves with speeds tending to $ 0 $ 
in dimension $N\geq 3$. Minimizing the energy $E$ at fixed momentum 
$Q$ has the advantage of providing orbitally stable traveling waves, 
and this is intimately related  to the concavity of the curve $Q \mapsto E $. 
On the other hand, if $Q \mapsto E $ is convex, as it is the case 
for  the upper branch in figure \ref{diaggrospit} (b), one expects 
orbital instability.

In the case of the Gross-Pitaevskii equation, the curves describing the minimum of the energy at fixed 
momentum in dimension $2$ and $3$ have been obtained in \cite{BGS}, where 
the existence of minimizers of $E$  under the constraint $Q = q = constant$ 
is also proven for any $ q >0$ if $N=2$, 
respectively for any $ q \in (q_0, \infty)$ (with $q_0 >0$) if $N=3$. 
The proofs in \cite{BGS} depend on the special algebraic structure of 
the Gross-Pitaevskii nonlinearity and it seems difficult to extend them 
to other nonlinearities. The existence of minimizers in \cite{BGS} has been shown by 
considering the corresponding problem on tori $\left( \R / 2 n \pi \Z  \right)^N$, 
proving a priori bounds for minimizers on tori, then passing to the 
limit as $ n \lra \infty$. Although this method gives the existence of minimizers 
in $\R^N$, it does not imply the precompactness of all minimizing sequences, 
and therefore it  leaves the question of the orbital stability of 
minimizers completely open. 

\medskip

In  space dimension $N \geq 3$,
the existence of traveling waves for (\ref{1.1})  for any speed 
$c \in (0, v_s)$  and under general conditions 
on the nonlinearity  has been proven in \cite{M10} by minimizing the action 
$E - cQ$ under a Pohozaev constraint.
Although the traveling waves obtained in \cite{M10} minimize the action $ E - cQ $ 
among all traveling waves of speed $c$, the constraint used to prove their 
existence is not conserved by the flow of (\ref{1.1}) and consequently 
this method does not imply directly
their orbital stability (which is expected 
at least for sufficiently  small speeds  $c$). 

\medskip

In space dimension two, the existence of finite energy traveling waves is much more tricky. 
Of course, these solutions are still critical points of the functional $ E - cQ $ and they satisfy two Pohozaev identities 
(corresponding to scaling with respect to $x_1$ and $ x_2$, respectively). 
However, the geometry of level sets of this functional  is  more complicated. 
For instance, if $ \Po$  is  the set  of functions satisfying both Pohozaev identities, 
we are able to show that $E - cQ$ does not admit even local minimizers in  $ \Po$. 

One of the main goals of the present 
paper is to prove the existence of two-dimensional traveling-waves for  (\ref{1.1}) under general conditions 
on the nonlinearity $F$ (similar to the assumptions in \cite{M10}).
We use two approaches to show the existence of such solutions. 
 If the nonlinear potential 
$V$ is nonnegative, we consider the problem of minimizing the energy $E$  while
the momentum $Q $ is kept fixed.  If $F$ behaves nicely in a neighborood of $1$  we show that 
there exist minimizers for any $ q \in (0, \infty)$.  
The minimizers are traveling waves and their speeds are 
the Lagrange multipliers associated to the variational problem. 
These speeds  tend to zero as $ q \lra \infty$ and to $ v_s$ as $ q \lra 0$. 
For general nonlinearities we show that the energy-momentum diagram of these  minimizers is exactly as in figure  \ref{diaggrospit} (a).

If $V$ achieves negative values (this happens, for instance, in the case of  the cubic-quintic NLS) 
the above approach cannot be used because the infimum of the energy in the set of functions of 
constant momentum is always $ - \infty $. 
 In this case we minimize the functional 
$E - Q$ in the set of functions $ \psi $ satisfying 
$ \ii_ { \R^2} |\nabla \psi |^2 \, dx = k$.  
Under general assumptions we show that minimizers exist for all $k \in (0, k_{\infty})$ (with $ k_{\infty} = \infty $ if and only if $ V \geq 0$) 
and, after scaling, they give rise to traveling waves. 
The speeds of these traveling waves tend to $ v_s $ as $ k \lra 0 $, and to $0$ if $V$ is nonnegative and $ k \lra \infty$.

In space dimension two, even if $V$ takes negative values 
it is still possible to find local minimizers of  the energy 
under the constraint $ Q = q= constant $ for any $q$ in some 
interval $(0, q_{\infty})$.
The proof relies on the tools developed to show the existence of minimizers of $E- Q$ under the constraint 
$ \ii_ { \R^2} |\nabla \psi |^2 \, dx = k = constant$.
Clearly, these minimizers are traveling waves and their speeds  tend to $ v_s$ as $ q \lra 0 $.

\medskip

We work with general nonlinearities and we  consider only the following set of assumptions: 

\medskip

{\bf (A1) } The function $F$ is continuous on $[0, \infty)$, 
$C^1$ in a neighborhood of $1$, $F(1) = 0$ and $F'(1) < 0$.

\smallskip

{\bf (A2) } 
There exist $C > 0$ and $ p _0 < \frac{2}{N-2} $ (with $p_0 < \infty $ if $ N =2$) 
such that $|F(s) | \leq C(1 + s^{p_0}) $ for any $ s \geq 0$.

\smallskip

{\bf (A3) } There exist $C, \,  \al _0> 0$ and $r_* > 1 $  such that $F(s) \leq - C s^{\al _0} $ 
for any $ s \geq r_*$. 

\medskip

{\bf (A4) } $F$ is $C^2$ near $ 1 $ and
$$
F(s) = - ( s - 1) + \frac 12 F''( 1) ( s - 1) ^ 2 + \Oo (( s - 1)^3) 
\qquad \mbox{ for $s$ close to } 1. 
$$

For simplicity, we summarize in the next theorem our existence results in space dimension  two. 
Our methods work as well in higher dimension, as we will see later in Theorems \ref{T1.1} and \ref{T1.2}.

\setcounter{section}{0}

\begin{Theorem}
\label{T0.1}
Let $N=2$. Assume that (A1), (A2), (A4) are satisfied  and $ F''(1) \neq 3$. The following holds.

\medskip

(i) Assume in addition that $ V \geq 0 $ on $[0, \infty)$.  Then for any $ q \in (0, \infty)$ there is 
a traveling wave $ \psi $ of (\ref{1.1}) with  speed $ c = c(\psi) \in (0, v_s)$ such that $  Q ( \psi ) = q$. 
Moreover, $\psi $ minimizes the energy $E(\phi) $  among all functions $ \phi $ satisfying $ Q( \phi )= q$. 

\medskip

(ii) There is $ k_{\infty } > 0 $ (with $ k_{\infty} = \infty$ if $ V $ is nonnegative) such that 
for all $ k \in (0, k_{\infty})$  there is a traveling wave $ \tilde{\psi }$ of (\ref{1.1}) 
with speed $ c = c(\tilde \psi) \in (0, v_s)$ such that $ \ii_{\R^2} |\nabla \tilde{\psi} |^2 \, dx = k$. 
Furthermore, $\tilde{\psi} $ minimizes the quantity $ I_c( \phi) = - c Q( \phi) + \ii_{\R^2} V(|\phi |^2) \, dx $ (or, equivalently, $ E( \phi) - c Q( \phi)$)
  among all functions $ \phi $ satisfying $  \ii_{\R^2} |\nabla \phi |^2 \, dx = k$.    

\medskip

(iii) There is $ q_{\infty }^{\sharp} > 0 $ such that for any $ q^{\sharp} \in (0, q_{\infty}^{\sharp})$ there 
is a traveling wave $ \psi ^{\sharp} $ of (\ref{1.1}) with  speed $ c = c(\psi^{\sharp}) \in (0, v_s)$ 
such that $  Q ( \psi ^{\sharp}) = q^{\sharp}$ and $\psi ^{\sharp} $ is a local minimizer of $ E $ under 
the constraint $ Q(\phi)= q^{\sharp}$. 
\end{Theorem}

The behavior of $F$ near 1 (assumption (A4)) is important only for the existence of "small energy" traveling waves. 
 If (A4) holds but  $ F''(1)= 3 $
it follows from Proposition \ref{smallE} below  that statements (i), (ii) and (iii) in  Theorem \ref{T0.1} may hold only for 
$ q \in (q_0, \infty)$, $ k \in (k_0, k_{\infty} )$, and  $ q^{\sharp}  \in ( q_0^{\sharp}, q_{\infty}^{\sharp})$, respectively,
where $ q_0, \; k_0 $ and $ q_0^{\sharp} $ are strictly positive. 
Our proofs show that under assumptions  (A1) and (A2) alone  the three minimization problems  above  can still be solved, but only 
if  $q$, $k$ and $q^{\sharp}$ are sufficiently large.

\setcounter{section}{1} 

\medskip

Our results cover as well nonlinearities of Gross-Pitaevskii type, for which 
the potential $ V$ is nonnegative and both (i) and (ii) apply, 
and of cubic-quintic type, for which $V$ achieves negative values and (ii) and 
(iii) hold. To the best of our knowledge, all previous results in 
the literature about the existence of traveling waves for (\ref{1.1}) 
in space dimension two  are concerned only with the Gross-Pitaevskii equation 
and the proofs make use of the specific algebraic properties of this nonlinearity.

\medskip

It is possible to minimize the energy $E$ at fixed momentum $Q$ (provided that $V \geq 0$) 
or to minimize the functional $I(\phi ) = - Q( \phi) + \ii_{\R^N} V(|\phi|^2) \, dx $ at constant kinetic energy (i.e., $ \ii_{\R^N} |\nabla \phi |^2 \, dx = k$) 
in any space dimension $N\geq 2$. 
The minimizers give rise to traveling waves for (\ref{1.1}) (after scaling in the latter case). 
The existence of minimizers and the compactness of minimizing sequences are proven exactly as in dimension two. 
Last but not least, minimizing the energy at fixed momentum gives a set of solutions which is orbitally stable by the flow of (\ref{1.1}); 
this property is, in general, not true for minimizers under Pohozaev constraints. 
For these reasons and in view of subsequent work  we state and prove our results in any dimension $ N \geq 2$; see Theorems \ref{T1.1}, \ref{T1.2}, \ref{T1.3} below for precise statements.
  We  obtain the properties 
of the curve $ q \longmapsto E_{min}(q)$ representing  the minimum of the energy vs. momentum   for general nonlinearities $F$  such that $ V \geq 0$.
If $ N \geq 3 $ or  ($ N =2$ and $ F''(1) = 3$) there exists $ q_0 > 0 $ such that   $ E_{min}(q)   = v_s q$ for  $q \in (0, q_0)$,
and $E_{min}$ is represented by 
 the curve in figure \ref{diaggrospit} (b) below the line 
$ E_{min}(q) = v_s q $ on $ (q_0, \infty)$.   This is in full  agreement 
with the results in \cite{JR}, \cite{JPR} and \cite{BGS}.

\medskip

As already  mentioned, an important issue is the orbital  stability of traveling waves. 
We prove 
the precompactness of all minimizing sequences for all  variational  problems  presented above.
Since the energy $E$ and the momentum $Q$ are conserved  quantities for  (\ref{1.1}), 
a classical result (see \cite{CL}) implies that the set of minimizers given by statements  (i) and (iii) in  Theorem \ref{T0.1} (or by Theorems \ref{T1.1} and \ref{T1.3} below) are orbitally stable 
 by the flow of (\ref{1.1}).
In particular, our results apply to the two and three-dimensional Gross-Pitaevskii equation. 
It has been conjectured in \cite{JPR} that, in the case of  the Gross-Pitaevskii equation, 
two-dimensional traveling waves corresponding to the diagram in  figure \ref{diaggrospit}  (a) should be stable, 
as well as the three-dimensional solutions corresponding to the "lower branch" in  figure \ref{diaggrospit}  (b), while those on the upper branch should be unstable. 
As far as we know, the asymptotic stability of traveling waves in dimension $N \geq 2$ and the orbital stability or instability of traveling waves corresponding to the "upper branch" in figure  \ref{diaggrospit}  (b) are still open problems. 

\medskip

The minimization problems considered in this article are physically more relevant than minimization of $ E- cQ$ 
under a Pohozaev constraint, can be solved  in space dimension two  and minimization of the energy at 
fixed momentum gives directly the orbital stability of the set of solutions. Moreover, 
we get minimizers for any momentum in some interval $(q_0, \infty) $ or $(0 , q_{\infty}) $ 
and for any kinetic energy in some interval $(k_0, k_{\infty})$.
The price to pay is  that
 the speeds of the traveling waves  obtained in this way  are Lagrange multipliers, so we cannot 
guarantee that these speeds cover a whole interval 
(but in all cases  we get at least an uncountable set of speeds). 
We mention that  the two-dimensional traveling waves to 
(\ref{1.1}) have been studied numerically in \cite{CS}, in the case of  general nonlinearities (as those studied 
in dimension one in \cite{C1d}).
The numerical algorithms in \cite{CS} allow to perform the constrained 
minimization procedures used in the present paper. Numerical computations suggest  that for $N=2$, 
even if the nonlinearity  $F$ satisfies  the assumptions (A1), (A2) and (A4) with $ F''(1) \not= 3 $, 
it is not true in general that minimizing $E$ at 
fixed $Q$  
 provides  a single interval of speeds;   for instance, it may provide  the union 
of two disjoint intervals, even if we require $F$ to be decreasing.

One might ask whether there is a relationship between the families of 
traveling waves obtained by different  approaches. 
In  dimension $ N \geq 3$ we prove that all traveling waves  found in 
the present paper also minimize the action $ E - cQ$ under the Pohozaev 
constraint considered in \cite{M10}. The converse is, in general, not true. 
For instance, in the case of the Gross-Pitaevskii equation in dimension 
$ N \geq 3 $, it was shown in \cite{BGS, dL} that there are no traveling 
waves of small energy. 
In section \ref{sectionsmallenergy} we prove a sharp version of  that result, valid for general nonlinearities (see Proposition \ref{smallE} below).
This implies that  there is $ c_0 < v_s$ such that there are no traveling waves 
of speed $c \in (c_0, v_s)$ which minimize the energy at fixed momentum. 
However, if $N \geq 3$  the existence of traveling waves as minimizers of $E - cQ$ under 
a Pohozaev constraint has been proven   for any  $ c \in (0, v_s)$. This is in 
agreement with the energy-momentum diagram in  figure \ref{diaggrospit} (b), 
where the traveling waves with speed $ c $ close to the speed of sound $v_s$ 
are expected to lie on the upper branch. We also prove that all minimizers of 
the energy at fixed momentum are  (after scaling) minimizers of 
$ E - Q$ at fixed kinetic energy. It is an open question to find sufficient conditions on $F$ which guarantee that the converse is also true. 
Whenever the converse is true, the set of speeds of traveling waves that minimize 
the energy at fixed momentum is the interval $(0, v_s)$.

\medskip

{\bf Main results. } 
We work in the natural energy space $ \Eo $ associated to (\ref{1.1}).  
There are several equivalent definitions for $ \Eo$ 
(see the next paragraph). For the statement of the results, it suffices to know that
$$
\Eo = \{ \psi : \R^N \lra \C\; \big| \; 
\psi \mbox{ is measurable, } |\psi | - 1 \in L^2( \R^N), \nabla \psi \in L^2(\R^N) \} 
$$
and that $ \Eo $  is endowed with the semi-distance
$$
d_0 (\psi _1, \psi _2) = \| \nabla \psi _1 - \nabla \psi _2 \| _{L^2 (\R^N)} 
+  \| \, |\psi _1|  - |\psi _2 |\, \| _{L^2(\R^N)} .
$$
In Section \ref{mom} we give a rigorous definition of the momentum on the whole space $ \Eo$ 
and we study its properties. For now the formal definition 
$ Q  (\psi ) = \ii_{\R^N} \langle i \frac{ \p \psi }{\p x_1}, \psi \rangle \, dx $
is sufficient. 

Our most important results can be summarized as follows.

\medskip

\setcounter{Theorem}{0}

\begin{Theorem} 
\label{T1.1}
Assume that $N\geq 2$, (A1) and (A2) are satisfied and $V \geq 0$ on $[0, \infty)$. 
For $q \geq 0$, let 
$$
E_{min}(q) = \inf \{ E( \psi) \; \big| \; \psi \in \Eo, \; Q( \psi ) = q \}.
$$
Then: 

\medskip

(i) The function $E_{min}$ is concave, increasing on $[0, \infty)$, 
$ E_{min}(q) \leq v_s q $ for any $ q \geq 0$, the right derivative of 
$E_{min} $ at $0$ is $ v_s$,  $ \; E_{min}(q) \lra \infty$ and 
$ \frac{ E_{min}(q) }{q} \lra 0$ as $ q \lra \infty$.

\medskip

(ii) Let $ q_0 = \inf \{ q > 0 \; | \; E_{min}(q) < v_s q \}.$ 
For any $ q > q_0$, all sequences $(\psi_n)_{n \geq 1} \subset \Eo $ 
satisfying $Q( \psi _n) \lra q $ and $E(\psi _n) \lra E_{min}(q)$ are 
precompact for $d_0$ (modulo translations). 

The set $ \So _q = \{ \psi \in \Eo \; | \; Q( \psi ) = q,  \; E( \psi ) = E_{min}(q) \}$ 
is not empty and is orbitally stable (for the semi-distance 
$d_0$)  by the flow associated to (\ref{1.1}). 

\medskip

(iii) Any $ \psi_q \in \So _q $ is a traveling wave for (\ref{1.1}) of speed 
$ c( \psi _q) \in [ d^+ E_{min}(q), d^- E_{min}(q) ]$, where we denote by 
$d^- $ and $d^+$ the left and right derivatives. We have 
$  c( \psi _q) \lra 0 $ as $ q \lra \infty$. 

\medskip

(iv) If $ N \geq 3 $ we have always $ q_0 > 0 $. Moreover, if $ N =2$ and 
assumption (A4) is satisfied, we have $ q_0 = 0$ if and only if $F''(1) \not = 3 $, 
in which case $  c( \psi _q) \lra v_s $ as $ q \lra 0$. 

\end{Theorem}

If $V$ achieves negative values, the infimum of $E$ in the set 
$\{ \psi \in \Eo \; | \; Q( \psi ) = q \} $ is $- \infty$ for any $q$. 
In this case we prove the existence of traveling waves by minimizing the functional 
$I(\psi ) = - Q( \psi ) + \ii_{\R^N} V(|\psi |^2) \, dx $ 
(or, equivalently, the functional $E - Q$)
under the constraint $ \ii_{\R^N} |\nabla \psi |^2 \, dx = k.$ 
More precisely, we have the following results:

\begin{Theorem} 
\label{T1.2}
Assume that $N\geq 2$ and (A1), (A2) are satisfied. For $ k \geq 0 $, let
$$
I_{min}(k) = \inf \Big\{ I( \psi) \; \big| \; \psi \in \Eo, \; 
\ii_{\R^N} |\nabla \psi |^2 \, dx = k \Big\}.
$$
Then there is $k_{\infty} \in (0, \infty]$ such that the following holds:

\medskip

(i) For any $ k > k_\infty $, $I_{min}(k) = - \infty$. The function 
$I_{min}$ is concave, decreasing on $[0, k_\infty)$, 
$ I_{min}(k) \leq - k / v_s^2 $ for any $ k \geq 0$, the right derivative 
of $ I_{min} $ at $0$ is $ - 1/ v_s^2 $, and $ \frac{ I_{min}(k) }{k} \lra - \infty$ 
as $ k \lra \infty$.

\medskip

(ii) Let $ k_0 = \inf \{ k > 0 \; | \; I_{min}(k) <  - k / v_s^2 \} \in [ 0, k_\infty]$. 
For any $ k \in (k_0, k_{\infty})$, all sequences $(\psi_n)_{n \geq 1} \subset \Eo $ 
satisfying $\ii_{\R^N} |\nabla \psi_n |^2 \, dx \lra k $ and 
$ I(\psi _n) \lra I_{min}(k)$ are precompact for $d_0$ (modulo translations). 
If $ \psi_k \in \Eo $ is a minimizer for $I_{min}(k)$, there exists 
$ c = c( \psi _k) \in \left[ \sqrt{ - 1 / d^+ I_{min}(k)}, \sqrt{ - 1 / d^- I_{min}(k)} \right]$ 
such that $ \psi_k ( \frac{\cdot}{c} )$ is a 
traveling wave  of (\ref{1.1}) of speed $c(\psi_k)$.

\medskip

(iii) We have $ k_{\infty} < \infty$ if and only if ($N= 2$ and $ \inf V < 0$). 
If $ k_{\infty} = \infty$, the speeds of the traveling waves obtained 
from minimizers of $I_{min}(k)$ tend to $0$ as $ k \lra \infty$. 

\medskip

(iv) For $ N \geq  3$, we have $ k_0 > 0 $. If $N=2$ and assumption (A4) 
is satisfied we have $ k_0 = 0$ if and only if $F''(1) \not= 3 $, in which 
case the speeds of the traveling waves obtained from minimizers of $I_{min}(k)$ 
tend to $ v_s $ as $ k \lra 0 $. 

\end{Theorem}

Notice that statements (iii)  and (iv) in Theorem \ref{1.2} provide sufficient conditions to  have 
$ k_0 < k_\infty $. Actually, this is always the case if $ N \geq 3 $. 
In the case  $ N = 2$, we have $ k_0 < k_\infty $  if $\inf V \geq 0 $, 
or if ($\inf V < 0 $, $F$ verifies assumption (A4) and $F''(1) \not = 3$). 
The  main  physical example of nonlinearity satisfying  $ \inf V < 0 $ 
is the cubic-quintic nonlinearity, for which one has $  F''(1) \neq 3 $.

In space dimension two, the tools developed to prove Theorem \ref{T1.2} 
enable us to find minimizers of $E$ at fixed momentum on a subset of 
$ \Eo$ even if $V$ achieves negative values. We have: 

\begin{Theorem} 
\label{T1.3}
Assume that $N = 2$ and that (A1), (A2) are satisfied. Let
$$ E_{min}^{\sharp}(q) = \inf \Big\{ E( \psi) \; | \; \psi \in \Eo, \; Q( \psi ) = q \; 
\mbox{ and } \; \ii_{\R^2} V(|\psi|^2) dx \geq 0 \Big\}.$$
Then: 

\medskip

(i) The function $E_{min}^{\sharp}$ is concave, nondecreasing on $[0, \infty)$, 
$ E_{min}^{\sharp}(q) \leq v_s q $, $d^+E_{min}^{\sharp}(0) = v_s$  and
$ E_{min }^{\sharp}(q) \leq k_{\infty}$ for any $ q > 0$, where 
$  k_{\infty}$ is as in Theorem \ref{T1.2}.

\medskip

(ii) Let $ q_0^\sharp = \inf \{ q > 0 \; \big| \; E_{min}^\sharp(q) <  v_s q \} \in [ 0, \infty ) $ 
and $ q_\infty^\sharp = \sup \{ q > 0 \; \big| \; E_{min}^{\sharp}(q) <  k_\infty \} \in ( 0, \infty]$. 
Then $ q_0^\sharp \leq q_\infty^\sharp $ and for any $q \in ( q_0^\sharp , q_\infty^\sharp)$, 
all sequences $(\psi_n)_{n \geq 1} \subset \Eo $ satisfying 
$Q( \psi _n) \lra q $ and $E(\psi _n) \lra E_{min}^{\sharp}(q)$ are 
precompact for $d_0$ (modulo translations). 

The set $ \So _q^{\sharp} = \{ \psi \in \Eo \; \big| \; Q( \psi ) = q,  \; 
E( \psi ) = E_{min}^{\sharp}(q) \}$ is not empty and is orbitally stable 
by the flow of (\ref{1.1}) for the semi-distance $d_0$. 

\medskip

(iii) Any $ \psi_q \in \So _q ^{\sharp}$ verifies $ \ii_{\R^2} V(|\psi_q|^2) dx > 0 $, 
hence minimizes $E$ under  the constraint $Q=q$ in the open set 
$ \{ w\in \Eo \; \big| \; \ii_{\R^2} V(|w|^2) dx > 0 \} $. Therefore, it is a 
traveling wave for (\ref{1.1}) of speed 
$ c( \psi _q) \in [ d^+ E_{min}^\sharp(q), d^- E_{min}^\sharp(q) ]$. 

\medskip

(iv) If assumption (A4) is satisfied, we have $ q_0^\sharp = 0 $ if and 
only if $F''(1) \not = 3$, and in this case $ c( \psi _q) \lra v_s $ as $q \lra 0 $.

\end{Theorem}


The concavity of $E_{min}^{\sharp}$ is significantly more delicate than that of $E_{min}$ because there is an additional constraint. 
 In Theorem \ref{1.3} it might happen that 
$ q_0^\sharp = q_\infty^\sharp $, in which case (ii) never holds. 
Statement (iv) gives a sufficient condition (which is  satisfied by the cubic-quintic NLS) ensuring that  $ 0 = q_0^\sharp < q_\infty^\sharp $.

\medskip

 The stability results in Theorems \ref{T1.1} (ii) and \ref{T1.3} (ii) are  proven in section \ref{sectionorbistab}. 
We underline that these results
concern the set of traveling waves obtained as minimizers of the energy at fixed momentum. 
The uniqueness of these solutions (up to the invariances of the problem) is not known. 
In order to study the orbital stability of a single traveling wave $ \psi _{c_*}$ of speed $c_*$
one would need to prove first its nondegeneracy, 
which means that the linearized operator $L$ defined by 
$ L \phi = i c_* \frac{\p \phi}{\p x_1} + \Delta \phi + F(|\psi_{c_*}|^2) \phi + F'(|\psi_{c_*}|^2) \langle \psi_{c_*}, \phi \rangle \psi_{c_*}$
has its kernel spanned only by the derivatives of $  \psi _{c_*}$.
(The derivatives $ \frac{ \p \psi _{c_*}}{\p x_j}$ always belong to the kernel of $L$ because (\ref{1.3}) is translation invariant.)
This would also give more precise information on $  \psi _{c_*}$ (such as local uniqueness up to invariances of (\ref{1.3})
or the existence of a smooth curve of traveling waves $ c \longmapsto \psi _c$ for $ c $ near $ c_*$). 
Proving the nondegeneracy seems very challenging. 
To our knowledge such results were obtained in similar problems only in the radial case or in dimension one.  
Traveling waves of (\ref{1.1}) are clearly {\it not } radial. 

\medskip

 If (A1) and (A3) are satisfied, traveling waves are uniformly bounded
(cf. Proposition 2.2 p. 1078 in \cite{M8}) and 
 it is explained in the introduction of \cite{M10} 
how it is possible to modify $F$ in a neighborhood of infinity in such a 
way that the modified function $\tilde{F}$ satisfies (A1), (A2),  (A3) and, moreover,  (\ref{1.1}) has 
the same traveling waves as the equation obtained from it by replacing $F$ 
by $\tilde{F}$. If (A1) and (A2) hold, we get traveling waves as minimizers of 
some functionals under constraints. If (A1) and (A3) are verified 
but (A2) is not, the above argument still implies  the existence of 
such solutions, but they are minimizers only for some modified functionals.~We~get: 

\begin{Corollary} 

There exist finite energy traveling waves to (\ref{1.1}) under the same assumptions as in Theorems \ref{T1.1}, \ref{T1.2}, \ref{T1.3}, respectively, 
except that condition (A2) is replaced by (A3). 
\end{Corollary}

In section \ref{3fam} we investigate the relationship between the traveling waves given by 
Theorems \ref{T1.1}, \ref{T1.2}, \ref{T1.3} above and those found in \cite{M10}.
We show that minimizers of the energy at fixed momentum are (after scaling) minimizers of  $  E- Q$ at fixed kinetic energy,  and the
traveling waves of speed $c$  obtained by minimization at fixed kinetic energy are among the minimizers of the action $E - cQ$ under a Pohozaev constraint.

In section \ref{slow} we shall see that  the equation $\Delta \psi + F(|\psi |^2) \psi = 0 $  admits nontrivial  solutions in $ \Eo $ if and only if $V$ achieves negative values. In this case it admits solutions of minimum energy (also called {\it ground states}) and we show that
the traveling waves    which minimize $E - cQ$ under a Pohozaev constraint ($N \ge 3$) converge to these ground states  as  $ c \lra 0$.

\bigskip

We conclude with a result concerning the nonexistence of small energy solutions 
to ({\ref{1.3}). This is a sharp version of a result  proven in \cite{BGS} for the Gross-Pitaevskii 
nonlinearity in dimension $N=3$, then extended to $N \geq 4$ 
in \cite{dL}. 
The cases where $ q_0 > 0 $, $ k_0 > 0 $ or $ q_0^\sharp > 0 $ in the above 
theorems follow directly from this result.

\begin{Proposition}
\label{smallE}
Assume that $N \geq 2 $ and that $F$ verifies (A1) and ((A2) or (A3)). Suppose that either 

$\bullet$ $N \geq 3$, or

$\bullet$ $N = 2$, $F$ satisfies (A4) and $ F''(1) = 3 $.
\\
The following holds.

\medskip

(i) There is  $ k_* > 0$, depending only on $ N $ and $F$, such that 
if $ c \in [0, v_s] $ and if $ U \in \Eo $ is a solution 
to (\ref{1.3}) satisfying $ \ii_{\R^N} |\nabla U |^2 \, dx < k_*$, 
then $ U $ is constant. 

\medskip

(ii) Assume, moreover, that  $F$ satisfies (A2) with $ p_0 < \frac 2N$ or $F$ satisfies (A3).
There is  $ \ell_* > 0$, depending only on $ N $ and $F$, such that any solution $ U \in \Eo $ to (\ref{1.3}) 
with $ c \in [0, v_s]$ and $ \ii_{\R^N} \left( |U|^2 - 1 \right)^2 \, dx < \ell_* $ is constant.
\end{Proposition}

In the present paper we do not study the one-dimensional traveling waves of (\ref{1.1}). 
The existence of such solutions can be proved by using ODE techniques; 
we refer to \cite{C1d} for a thorough analysis of the 1D case 
and to \cite[Theorem 5.1 p. 1099]{M7} for nonexistence results in the supersonic case. 
It turns out that the energy-momentum diagrams depend strongly on the nonlinearity 
(this is also the case in space dimension two, see \cite{CS} for numerical results, 
and very probably in higher dimensions). Even for nice nonlinearities we may have a great variety 
of behaviors: multiplicity of solutions,  
branches of traveling waves that intersect each other or whose energy or momentum tend to infinity,  
nonexistence of traveling waves for some speed $ c_* \in (0, v_s)$, existence of a sonic traveling wave. 
It has been shown in \cite[Lemma 2]{BGS-survey} that minimizing the energy at fixed momentum is not 
possible in dimension 1. Here, by momentum, we mean the same quantity as in section 
\ref{mom} below. Actually, a renormalized momentum has been introduced in 
\cite{KivYan} to treat this problem and in the case of the Gross-Pitaevskii nonlinearity it 
has been proved in \cite{BGS-survey} that it is possible to minimize the energy at fixed renormalized 
momentum and the minimizers are precisely the traveling waves; 
see Theorem 2 in \cite{BGS-survey} which is the 1-D counterpart of our Theorem \ref{T1.1} for 
the Gross-Pitaevskii nonlinearity. The stability of 1-D traveling waves has been 
addressed in a series of papers. We refer to \cite{C1dStab} 
for a detailed study and for a review of the existing literature on this topic. 
Roughly speaking, being given a 1-D traveling wave $ \psi_{c_*} $ of speed $ c_* \in (0, v_s)$ one 
may construct a curve $ c \longmapsto \psi _c$ of traveling waves for $c$ in a neighborhood of $ c_*$.
Then $ \psi _{c _*}$ is orbitally stable if $ \frac{ dQ ( \psi _c)}{dc }_ {| c = c_*} < 0 $ and 
orbitally unstable if $ \frac{ dQ ( \psi _c)}{dc }_ {| c = c_*} > 0 , $ where $Q$ is the momentum. 
Equivalently, $ \psi _{c _*}$ is orbitally stable if  the mapping 
$ Q( \psi _c) \longmapsto E( \psi _c)$ is concave and orbitally unstable if this mapping is convex.

\medskip

{\bf Notation and function spaces.  }
Throughout the paper,  $ \Lo ^N$ is the Lebesgue measure on $ \R^N$ and  
$\mathcal{H}^{s}$ is the $s-$dimensional Hausdorff measure on $ \R^N$. 
For $ x = ( x_1, \dots, x_N) \in \R^N$, 
we denote  $ x' = ( x_2, \dots , x_N ) \in \R^{N-1}$. 
We write $\langle z_1, z_2 \rangle $ for the scalar product 
of two complex numbers $z_1, z_2$. 
Given a  function $f$ defined on $ \R^N$ and $ \la , \, \si  > 0$, 
we denote
\beq
\label{scale}
f_{\la, \si } (x) = f\left( \frac{x_1}{\la }, \frac{x'}{\si} \right).
\eeq
If $ 1 \leq p <N$, we write $ p^*$ for the Sobolev exponent 
associated to $p$, that is $ \frac{1}{p^*} = \frac 1p - \frac 1N$.

\medskip

If $F$ satisfies  (A1),    using Taylor's formula for $s$ in a neighborhood  of $ 1 $ we have
\beq
\label{1.4}
V(s) = \frac 12 V''(1) ( s - 1)^2 + ( s - 1)^2 \e (s - 1)
= \frac 12  ( s - 1)^2 + ( s - 1)^2 \e (s - 1) , 
\eeq
 where $ \e (t) \lra 0  \mbox{ as } t \lra 0.$ 
Hence for $|\psi |$ close to $ 1$, $V(|\psi |^2) $ can be approximated by 
the Ginzburg-Landau potential $\frac 12 ( |\psi |^2 - 1)^2$. 

\medskip

We fix an odd function $ \ph \in C^{\infty } (\R)$ such that $ \ph (s) = s $ 
for $ s \in [0, 2  ]$, $ 0 \leq \ph ' \leq 1 $ on $ \R$ and 
$ \ph (s) = 3   $ for $ s \geq 4  $. 
If assumptions (A1) and (A2) are satisfied,  it is not hard to see that 
there exist $C_1, \, C_2, \, C_3 > 0$  such that
\beq
\label{i1}
\begin{array}{l}
|V(s)| \leq C_1 ( s - 1)^2 \quad \mbox{ for any } s \leq 9 ; 
\\
\mbox{in particular, } 
| V(\ph ^2(\tau )) | \leq C_1 (\ph ^2(\tau ) - 1) ^2 \mbox{ for any  } \tau, 
\end{array}
\eeq
\vspace*{-11pt}
\beq
\label{i2}
|V(b) - V(a)| \leq C_2 |b-a| \max (a^{p_0}   , b^{p_0}   )
\qquad \mbox{ for any } a, \, b \geq 2 .
\eeq
Given  $ \psi \in H_{loc}^1 ( \R^N) $ and  an open set  $ \Om \subset \R^N$, 
the modified  Ginzburg-Landau  energy of $\psi $ in $ \Om $ is 
\beq
\label{1.5}
E_{GL }^{\Om } (\psi ) = \ds \int_{\Om } |\nabla \psi |^2 \, dx 
+ \frac 12 \int_{\Om } \left( \ph ^2(|\psi |) - 1 \right)^2 \, dx .
\eeq
We  simply write $ E_{GL}(\psi) $ instead of $E_{GL}^{\R^N} (\psi)$.
The modified Ginzburg-Landau energy will  play a central role in our analysis.

We denote 
$ \dot{H}^1(\R^N) = \{ \psi \in L_{loc}^1(\R^N) \; | \; \nabla \psi \in L^2(\R^N) \}$ and 
\beq
\label{E}
\begin{array}{rcl}
\Eo  & = & \{ \psi \in \dot{H}^1(\R^N) \; \big| \; \ph^2(|\psi |) - 1 \in L^2( \R^N) \} 
\\
 & = & \{ \psi \in \dot{H}^1(\R^N) \; \big| \; 
 E_{GL}(\psi) < \infty \}.
\end{array}
\eeq
Let $\DR$ be the completion of $ C_c^{\infty } (\R^N) $ 
for the norm $\|v \| = \|\nabla v \| _{L^2 (\R^N) }$ and let 
\beq
\label{X} 
\begin{array}{rcl}
\Xo  \! \! \! &  = &  \! \! \! \{ u \in \DR \; \big| \; \ph ^2(|1+ u|) - 1 \in L^2( \R^N) \} 
\\
 & = &  \! \! \! \{ u \in \dot{H}^1(\R^N) \; \big| \; 
u \in L^{2^*} (\R^N), \; E_{GL}(1+ u) < \infty \}  \quad \mbox{ if } N \geq 3, \mbox{  where }  2^*\!  = \! \frac{2N}{N-2}.
\end{array}
\eeq

If  $ N \geq 3 $ and $ \psi \in \Eo$, 
 there exists a constant $ z_0 \in \C$ such that  $ \psi - z_0 \in L^{2^*}(\R^N)$ 
 (see, for instance, Lemma 7 and Remark 4.2 pp. 774-775 in \cite{PG}).
It follows that $ \ph (|\psi | ) - \ph(|z_0|) \in  L^{2^*}(\R^N)$.
On the other hand, the fact that $E_{GL}(\psi ) < \infty $ implies $ \ph (|\psi | ) - 1 \in L^2(\R^N)$, 
thus necessarily $ \ph(|z_0|)  = 1$, that is $ |z_0 | = 1 $. 
Then it is easily  seen that there exist $ \al _0 \in [0, 2\pi)$ and $ u \in \Xo$, uniquely determined by $ \psi$, such that $ \psi = e^{i \al _0} ( 1+ u).$
In other words, if $ N \geq 3$ we have 
$\Eo = \{ e^{i \al _0 } (1+ u) \; | \; \al _0 \in [0, 2\pi), \; u \in \Xo \}.$
This simple description of $ \Eo $ is no longer true if $ N =2$.

It is not hard  to see that for $N \geq 2$ we have 
\beq
\label{L2}
\Eo = \{ \psi : \R^N \lra \C\; \big| \; 
\psi \mbox{ is measurable, } |\psi | - 1 \in L^2( \R^N), \nabla \psi \in L^2(\R^N) \}. 
\eeq
Indeed, we have $\big| \ph ^2 ( |\psi |) - 1\big| \leq 4  \big| \, | \psi | - 1 \big|$, hence
$ \ph ^2 ( |\psi |) - 1 \in L^2( \R^N)$ if $ |\psi | - 1 \in L^2( \R^N) $.
Conversely, let $ \psi \in \Eo$. 
If $ N =2$, it follows from Lemma \ref{L2.1} below that $ |\psi |^2 - 1 \in L^2 ( \R^2)$ 
and we have 
$\big| \, | \psi | - 1 \big| = \frac{1}{  |\psi | + 1}  \big| \, | \psi |^2 - 1 \big| \leq  \big| \, | \psi |^2 - 1 \big|.$
If $ N \geq 3$, we know that $ \ph (|\psi |  )- 1 \in L^2( \R^N)$ and 
$ 0 \leq |\psi | - \ph(|\psi |) \leq |\psi | \1_{\{ |\psi | \geq 2  \} } 
\leq 2 ( |\psi | - 1) \1_{\{ |\psi | \geq 2  \} }  
\leq 2 \big| \,  |\psi | - 1 \big| ^{ \frac{ 2^*}{2} } \1_{\{ |\psi | \geq 2  \} } $
and the last function belongs to $ L^2( \R^N)$ by the Sobolev embedding. 
Moreover, one may find bounds for $\| \, |\psi | - 1 \|_{L^2( \R^N)} $ in terms of $E_{GL}(\psi)$ (see Corollary \ref{C4.2} below).

Proceeding as in \cite{PG2}, section 1, one proves  that $ \Eo \subset L^2  + L^{\infty}(\R^N)$ 
and that $ \Eo$ endowed with the distance
\beq
\label{distance}
d_{\Eo}(\psi _1, \psi _2) = \| \psi _1 - \psi _2 \| _{L^2 + L^{\infty}(\R^N)} 
+  \| \nabla \psi _1 - \nabla \psi _2 \| _{L^2 (\R^N)} 
+  \| \, |\psi _1|  - |\psi _2 |\, \| _{L^2(\R^N)} 
\eeq
is a complete metric space. We recall that, given two Banach spaces $X$ and $Y$ of distributions 
on $ \R^N$,  the space $X + Y$  with norm defined by 
$\| w \|_{X + Y} = \inf \{ \| x \|_X + \| y \|_Y \; \big| \; w = x + y, x \in X , y \in Y \}$ is a Banach  space.

We will also consider the following semi-distance on $ \Eo$:
\beq
\label{semidistance}
d_0(\psi _1, \psi _2) =\| \nabla \psi _1 - \nabla \psi _2 \| _{L^2 (\R^N)} 
+  \| \, |\psi _1|  - |\psi _2 |\, \| _{L^2(\R^N)} .
\eeq
If $\psi _1, \psi _2 \in \Eo $ and $d_0(\psi _1, \psi _2) = 0$, 
then   $|\psi _1 | = |\psi _2 |$ a.e. on $\R^N$ 
and $ \psi _1 - \psi _2$ is a constant (of modulus not exceeding $ 2 $) a.e. on $\R^N$. 

In space dimension $ N = 2, 3, 4$, the Cauchy problem for the Gross-Pitaevskii 
equation has been studied by Patrick G\'erard (\cite{PG, PG2}) in the space 
naturally associated to that equation, namely
$$
\E = \{ \psi \in H_{loc}^1(\R^N) \; | \; \nabla \psi \in L^2( \R^N), |\psi |^2 - 1 \in L^2( \R^N) \}
$$
endowed with the distance
$$
d_{\E} (\psi _1, \psi _2) = \| \psi _1 - \psi _2 \| _{L^2 + L^{\infty}(\R^N)} 
+  \| \nabla \psi _1 - \nabla \psi _2 \| _{L^2 (\R^N)} 
+  \| \, |\psi _1| ^2 - |\psi _2 | ^2 \, \| _{L^2(\R^N)} .
$$
 If $ N = 2, 3 $ or $4$ 
it can be proved that $ \E = \Eo $ and the distances $ d_{\Eo} $ 
and $ d_{\E}$ are equivalent on $ \Eo$. 
Global well-posedness was shown in \cite{PG, PG2} (see section \ref{sectionorbistab}) if 
$  N \in \{ 2, 3 \}  $ or if $N= 4$ and the initial data is small.
In the case $ N=4$, global well-posedness for any initial data in $\E$ was recently proven in \cite{KOPV}.

\medskip

{\bf Some ideas in the proofs. } 
Theorem \ref{T1.1} is proven in Section \ref{minem}.
If $V$ is nonnegative, we show first that the energy $E$ can be estimated in terms of the Ginzburg-Landau energy $E_{GL}$, and conversely. 

If $ 0 < c < v_s = \sqrt{2}$, we may choose $ \e, \, \de > 0$ such that $ c < \sqrt{2}(1 - 2 \e) (1 - \de)$. 
Suppose that $ \psi \in \Eo $ satisfies $ 1 - \de \leq |\psi| \leq 1 + \de$. 
Then there is a lifting $ \psi = \rho e^{i \theta}$  and a simple computation shows that 
$$
\begin{array}{l}
|\nabla \psi |^2 = |\nabla \rho |^2 + \rho ^2 |\nabla \theta |^2 , \qquad
Q(\psi) = - \ds \int_{\R^N} ( \rho ^2 - 1) \frac{ \p \theta}{\p x_1} \, dx \qquad \mbox{ and } 
\\
\\
V(|\psi |^2) = V(\rho^2) = \frac 12 (\rho^2 -1) ^2 + o ((\rho^2 -1) ^2 ) \geq \frac{1-\e}{2} (\rho^2 -1) ^2
\end{array}
$$
provided that $ \de $ is sufficiently small.
Then we have 
\beq
\label{vanish0}
\begin{array}{l}
\ds |c Q(\psi)| \leq \sqrt{2} (1 - 2 \e) (1 - \de) \int_{\R^N} | \rho ^2 -1| \cdot \Big| \frac{ \p \theta}{\p x_1} \Big| \, dx 
\\
\\
\ds \leq (1 - 2 \e) \int_{\R^N} ( 1 - \de )^2 \Big| \frac{ \p \theta}{\p x_1} \Big| ^2 + \frac 12 ( \rho ^2 -1)^2 \, dx 
\\
\\
\ds \leq \int_{\R^N} ( 1 - 2 \e) \rho ^2 | \nabla \theta |^2 + V( \rho ^2) - \frac{\e}{2} ( \rho ^2 -1)^2  \, dx 
\leq E(\psi) - \e E_{GL}(\psi).
\end{array}
\eeq
Thus $E (\psi) \geq  | c Q ( \psi )| + \e E_{GL}(u) $ if $ |\psi|$ is sufficiently close to $1$ in the $L^{\infty}$ norm. 
Since $E_{GL}(\psi)$ measures, in some sense, the distance from $\psi $ to constant functions of modulus  $1$, 
we would like to establish a similar estimate for all functions with  small Ginzburg-Landau energy.  
However, $E_{GL}(\psi)$ does not control $\| \, |\psi| - 1 \|_{L^{\infty}}$. Moreover, there are functions with 
arbitrarily small Ginzburg-Landau energy which have  small-scale topological "defects" (e.g., dipoles). 
To overcome  these difficulties  we use a procedure of regularization by minimization, 
which is  studied  in Section \ref{reg}. 
Given $ \psi \in \Eo$, we minimize the functional 
$\zeta \longmapsto E_{GL}(\zeta) + \ds \frac{1}{h^2} \int_{\R^N} \ph(|\zeta-\psi | ^2) \, dx $ in the set 
$\{ \zeta \in \Eo \; | \;  \zeta - \psi \in H^1( \R^N) \}.$
It is shown that minimizers exist 
and any minimizer $ \zeta_h $ has nice  properties, for instance:

$\bullet $ $E_{GL}(\zeta_h) \leq E_{GL}(\psi)$, 

$\bullet $ $\| \zeta_h - \psi \|_{L^2} \lra 0 $ as $ h \lra 0$, and 

$\bullet $  $ \| \, |\zeta_h | - 1 \|_{L^{\infty}} $ can be estimated in terms of $h$ and $E_{GL}(\psi)$ 
and  is arbitrarily small if $E_{GL}(\psi)$  is sufficiently small. 

Using the regularization procedure described above, we show that an estimate of the form (\ref{vanish0}) is true for all functions in $ \Eo $ with sufficiently small Ginzurg-Landau energy. 
In particular, this implies that for all $ c \in (0, v_s)$ we have  $ E_{min}(q) \geq c|q|$ if $q$ is small enough. 
Using appropriate test functions we show that  $ E_{min}(q) \leq  v_s|q|$ for all $q$. 
The concavity of $E_{min}$ is proven  by using  test functions obtained from  "approximate minimizers"
 by  reflection  with respect to hyperplanes  in $\R^N$. 
It has been shown in \cite{CM} that in the limit $ c \lra v_s$, two and three dimensional traveling waves of (\ref{1.1})  have a lifting and their modulus and phase can be approximated (after scaling)  by the ground states of  the Kadomtsev-Petviashvili I (KP-I)  equation. 
If $N=2$, (A4) is satisfied and $ F''(1) \neq 3$, using test functions constructed from  two-dimensional ground states for KP-I we show that $E_{min}( q) < v_s q$ for all $ q> 0$  (see Theorem \ref{T4.13}).

We use the concentration-compactness principle and the regularization procedure in Section \ref{reg}  to show the existence of minimizers for $E_{min}(q)$. 
The hardest part is to show that minimizing sequences do not "vanish," that is their Ginzburg-Landau energy does not spread over $ \R^N$. 
Assume that $ (\psi _n)_{n \geq 1}$ is a vanishing minimizing sequence. 
We show that 
$
\ii_{\R^N} V(|\psi_n|^2) \, dx = \frac 12 \ii_{\R^N} \left( \ph ^2(|\psi_n |) -1 \right)^2 \, dx + o (1). 
$
Using Lemma \ref{vanishing} we construct a sequence $ h_n \lra 0 $ and for each $n$ we find a minimizer $ \zeta _n $ of the functional 
$ E_{GL}(\zeta) + \ds \frac{1}{h^2} \int_{\R^N} \ph(|\zeta-\psi _n | ^2) \, dx $ such that $ \| \, | \zeta _n | - 1 \|_{L^{\infty}} \lra 0 $. 
We show that $ Q( \zeta _n ) = Q( \psi _n ) + o(1) $, then  using (\ref{vanish0}) we get for all $ c \in (0, v_s)$, 
$$
E( \psi _n ) = E_{GL}( \psi _n) + o(1) \geq E_{GL}( \zeta _n)  + o(1) \geq c |Q( \zeta _n ) | + o(1) = c |Q( \psi _n ) | + o(1) .
$$
Passing to the limit as $ n \lra \infty$, then letting $ c \uparrow v_s$ we obtain $ E_{min}(q) \geq v_s q$ .
Hence minimizing sequences cannot vanish  if   $ E_{min}(q) < v_s q$ .

If "dichotomy" occurs, we must have $ E_{min}(q) = E_{min}(q_1) + E_{min}( q - q_1)$ for some $ q_1 \in (0, q)$. 
However, the concavity of $E_{min}$ implies that $ E_{min} ( q_1) \leq \frac{ q_1}{q} E_{min}(q) $ and 
 $ E_{min} ( q-q_1) \leq \frac{ q-q_1}{q} E_{min}(q) $, with equality if and only if $ E_{min}$ is linear on $[0, q]$. 
Taking into account the behavior of $E_{min}$ near the origin, that would imply $E_{min}(q) = v_s q$, a contradiction.

Since "vanishing" and "dichotomy" are ruled out,  we must  have "concentration."
Then Lemmas \ref{L4.10} and \ref{L4.11} are powerful tools that enable us to conclude 
that $(\psi _n)_{n \geq 1}$ has a convergent subsequence. 

\medskip

The starting point in the proof of Theorem \ref{T1.2} is a refinement of (\ref{vanish0}), namely: 
for $ c \in (0, v_s)$ and $ \e $ sufficiently small (depending on $c$) there holds
$$
E( \psi ) - \e E_{GL}(\psi) \geq |c Q( \psi)|
$$
for all functions $ \psi \in \Eo $ such that $ \| \nabla \psi \|_{L^2}$ is sufficiently small
 (see Lemma \ref{L5.1}).
This enables us to show that $ I_{min}(k) > - \infty$ if $ k $ is sufficiently small. 
If $ N \geq 3 $ or ($N=2$ and $V \geq 0$) we may use scaling to prove that $I_{min}(k)$ is finite for all $k$. 
If $ N = 2 $ and $ V$ achieves negative values, the value of $ k_{\infty}$ is given  in Lemma \ref{L5.3bis}.  
The other properties of the function $I_{min}$ are proven by using appropriate test functions and scalings. 

The boundedness of the Ginzburg-Landau energy for minimizing sequences of $I_{min}$ is not obvious. This is done in Lemma \ref{L5.4}. 
Then we use again  the concentration-compactness principle   and the analysis developed in Sections \ref{reg} and \ref{minem} to show the 
existence of minimizers.

\section{The momentum}
\label{mom}

The momentum (with respect to the $x_1$ direction) should be a functional 
defined on $ \Eo$ whose "G\^ateaux differential"\footnote{We did not introduce a 
manifold structure on $ \Eo$, although this can be done in a natural way, 
see \cite{PG, PG2}. However, it will be clear (see (\ref{gateaux})) what 
we mean here by "G\^ateaux differential." }  is $2 i \p_{x_1}$.
In dimension $ N \geq 3$, it has been shown in \cite{M10} how to define 
the momentum on $\Xo$ (and, consequently, on $\Eo$). The definition in \cite{M10} cannot be used directly in dimension $N=2$. 
In this section we will extend that definition in dimension two.

It is clear that on the affine space $ 1 + H^1(\R^N) \subset \Eo$, 
the momentum should be defined by 
$Q(1+ u) =  \ii_{\R^N} \langle i u _{x_1}, u \rangle \, dx .$
In order to define the momentum on the whole $ \Eo$, 
we introduce  the space $ \Yo = \{ \p _{x_1} \phi \; | \; \phi \in \dot{H}^1(\R^N) \}$.
It is easy to see that $\Yo $ endowed with the  norm 
$\| \p _{x_1} \phi \|_{\Yo}  = \| \nabla \phi \| _{L^2(\R^N )} $ is a Hilbert space.

In dimension $ N \geq 3$, it follows from Lemmas 2.1 and 2.2 in \cite{M10} 
that for any $ u \in \Xo $ we have $\langle i u_{x_1} , u \rangle \in L^1(\R^N) + \Yo.$ 
If $ N \geq 3 $ and $ \psi \in \Eo$, we have already seen there are 
$ u \in \Xo $ and $ \al _0 \in [0, 2 \pi) $ such that $\psi = e^{i \al _0 } (1+ u).$
An easy computation gives
$\langle i \psi_{x_1} , \psi  \rangle =  Im (u_{x_1} ) + \langle i u_{x_1} , u \rangle  $ and it is obvious that 
$  Im (u_{x_1} ) \in \Yo$, thus $ \langle i \psi_{x_1} , \psi  \rangle \in  L^1(\R^N) + \Yo.$
The next Lemma shows that a similar result holds if $N=2$.

\begin{Lemma} \label{L2.1}   
Let $ N \! = \!2. $ For any $ \psi \in \Eo $ we have $ |\psi |^2 - 1 \in L^2(\R^2)$ and~$\langle i \psi _{x_1}, \psi \rangle \in L^1(\R^2) + \Yo.$
\end{Lemma}

{\it {Proof.}}
The following facts, borrowed from \cite{brezis-lieb}, will be useful here and in the sequel: 
for any $ q \in [2, \infty )$ there is $C_q >0$ such that for all $ \phi \in L_{loc}^1(\R^2) $ satisfying
$\nabla \phi \in L^2(\R^2)$ and $\Lo ^2( supp (\phi ) ) < \infty $ we have 
\beq
\label{ineq1}
\| \phi \| _{L^q(\R^2)} \leq 
C_q \| \nabla \phi \|_{L^1(\R^2)} ^{\frac 2q} \| \nabla \phi \|_{L^2(\R^2)} ^{1- \frac 2q} 
\eeq
(see inequality (3.12) p. 108 in \cite{brezis-lieb}).
Since $ \nabla \phi = 0 $ a.e. on $ \{ \phi = 0 \}$, 
(\ref{ineq1}) and the Cauchy-Schwarz inequality give
\beq
\label{ineq2}
\| \phi \| _{L^q(\R^2)} \leq 
C_q \| \nabla \phi \|_{L^2(\R^2)} \left( \Lo ^2(\{ \phi (x) \neq 0  \} ) \right)^{\frac 1q}.
\eeq
Notice that (\ref{ineq2}), which is a variant of inequality (3.10) p. 107 in \cite{brezis-lieb},  holds for any $q \in [1, \infty)$. 

 Let $ \psi  \in \Eo$. It is clear that
\beq
\label{e1}
 \ds \int_{\{ | \psi | \leq 2  \}} (|\psi |^2 - 1)^2 \, dx 
= \int_{\{ | \psi | \leq 2  \}} (\ph ^2(|\psi | ) - 1 )^2 \, dx < \infty .
\eeq
Obviously, 
$\Lo ^2(\{ |\psi | \geq \frac 32   \}) < \infty$   (because $E_{GL}(\psi ) < \infty$) and
$ |\psi |^2 - 1 \leq C ( |\psi |- \frac 32  ) ^2 $ on $\{ |\psi  | \geq 2  \}$. 
Using (\ref{ineq2})  for $\phi =  ( |\psi |- \frac 32  ) _+$ (which satisfies 
$ |\nabla \phi | \leq | \nabla \psi | \1_{ \{ |\psi  | \geq \frac 32  \} }$ a.e.)~we~get
\beq
\label{e2}
\ds \int_{\{ | \psi | > 2  \}}  \left(|\psi |^2 - 1 \right)^2  dx 
\leq C \int  \left( |\psi | - \frac 32  \right)_+ ^4 \, dx 
\leq C \| \nabla \psi \|_{L^2(\R^2)} ^4 \Lo ^2( \{ |\psi | \geq \frac 32   \}) < \infty.
\eeq
Thus  $ |\psi |^2 - 1 \in L^2(\R^2)$.

It follows from Theorem 1.8 p. 134 in \cite{PG2} that there exist  $w \in H^1(\R^2)$ and
 a real-valued function $\phi$ on $\R^2$
such that $\phi \in L_{loc}^2(\R^2)$,   $\;  \p^{\al } \phi \in L^2(\R^2)$ 
for any $ \al \in \N^2$ with $|\al | \geq 1$  and
\beq
\label{2.0}
\psi =  e^{i \phi } + w.
\eeq
A simple computation gives
\beq
\label{2.1}
\langle  i \psi _{x_1}, \psi \rangle 
= -  \frac{\p \phi}{\p x_1}   
+  \frac{\p }{\p x_1} \left( \langle iw, e^{i \phi} \rangle \right)
 - 2  \langle \frac{ \p \phi}{\p x_1} e^{ i \phi } , w \rangle
+ \langle i w_{x_1}, w \rangle.
\eeq
The Cauchy-Schwarz inequality implies that 
$\langle \phi _{x_1} e^{i \phi } , w \rangle $ and $ \langle i w _{x_1}, w \rangle $ 
belong to $L^1(\R^2)$. It is obvious that $\frac{\p \phi}{\p x_1} \in \Yo$. 
We have $ \langle iw, e^{i \phi} \rangle \in L^2 (\R^2)$ and 
$$
\frac{\p }{\p x_j} \left( \langle iw, e^{i \phi} \rangle \right)
= \langle i\frac{\p w}{\p x_j}, e^{i \phi} \rangle + \langle w, \frac{\p \phi}{\p x_j} 
e^{i \phi} \rangle.
$$
The fact that $ w $ and $  \frac{ \p \phi}{\p x_j}$ belong to $H^1(\R^2)$ and the Sobolev 
embedding give
$w, \;  \frac{ \p \phi}{\p x_j} \in L^p(\R^2)$ for any $ p \in [2, \infty)$, hence 
$\langle w, \frac{\p \phi}{\p x_j} e^{i \phi} \rangle \in L^p(\R^2)$ for any 
$ p \in [1, \infty)$. 
Since $\langle i \frac{\p w }{\p x_j},  e^{i \phi} \rangle \in L^2(\R^2)$, we get 
$\frac{\p }{\p x_j} \left( \langle iw, e^{i \phi} \rangle \right) \in L^2(\R^2)$, hence 
$\langle iw, e^{i \phi} \rangle \in H^1(\R^2)$ and consequently 
$ \frac{\p }{\p x_1} \left( \langle iw, e^{i \phi} \rangle \right) \in \Yo.$
The proof of Lemma \ref{L2.1} is complete.
\hfill $\Box$

\medskip

For $ v \in L^1(\R^N)$ and $ w \in \Yo $, let $ L(v + w ) =  \ii _{\R^N} v (x) \, dx$. 
It follows from Lemma 2.3 in \cite{M10} that $L$ is well-defined and that it is a 
continuous linear functional on $  L^1(\R^N) + \Yo$. Taking into account Lemma 
\ref{L2.1} and the above considerations, for any $N \geq 2$ we  give the following

\begin{Definition}
\label{D2.2}
Given $ \psi \in \Eo$, the momentum of $\psi$ with respect to the
$x_1-$direction is 
$$
Q(\psi) = L( \langle i \psi_{x_1}, \psi \rangle  ).
$$
\end{Definition}

Notice that  the momentum (with respect to the
$x_1-$direction) has been defined  in \cite{M10} 
 for functions  $u \in  \Xo $  by $ \tilde{Q}(u) = L( \langle i \frac{ \p u}{\p x_1} , u \rangle)$.
If $\psi = e^{i \al _0} ( 1+ u)$, it is easy to see that 
$Q(\psi )= \tilde{Q}(u) $.

\medskip

If $ \psi \in \Eo $ is symmetric with respect to $ x_1$ (in particular, if $\psi $ is radial), then
$
Q(\psi) = Q( \psi (- x_1, x')) = - Q(\psi) , 
$ hence $ Q(\psi) = 0$.

\medskip

If $ \psi \in \Eo $ has a lifting $ \psi = \rho e^{ i \theta } $ with $ \rho ^2 - 1 \in L^2( \R^N)$ 
and $ \theta \in \dot{H}^1( \R^N)$ 
(note that if $ 2 \leq N \leq 4$ we have always $|\psi |^2 - 1 \in L^2 ( \R^N)$ by (\ref{L2}) and 
the Sobolev embedding), then 
\beq
\label{lift}
Q( \psi ) = L( - \rho ^2 \theta _{x_1}) = - \ds \int_{\R^N} ( \rho ^2 - 1) \theta _{x_1} \, dx.
\eeq

The next Lemma is an "integration by parts"  formula.

\begin{Lemma}
\label{L2.3}
For any $ \psi \in \Eo$ and $ v \in H^1 ( \R^N) $ we have 
$\langle i\psi _{x_1}, v \rangle \in L^1( \R^N)$, 
$\langle i\psi , v_{x_1}  \rangle  \in L^1( \R^N) + \Yo $ and 
\beq
\label{}
L(\langle i\psi_{x_1}, v \rangle +  \langle i\psi , v_{x_1}  \rangle) = 0.
\eeq
\end{Lemma}

{\it Proof.} If $ N \geq 3$ this follows immediately from Lemma 2.5 in \cite{M10}.
We give the proof in the case $N=2$. The Cauchy-Schwarz inequality implies 
$\langle i\psi_{x_1}, v \rangle  \in L^1(\R^2).$
Let $w \in H^1(\R^N)$ and $ \phi $ be  as in (\ref{2.0}), so that 
$ \psi =  e^{i \phi } + w$. Then 
\beq
\label{2.5}
 \langle i\psi , v_{x_1}  \rangle
=  \frac{\p }{\p x_1} \left( \langle i e ^{i \phi } , v \rangle \right) 
+  \langle  \phi _{x_1} e ^{i \phi } , v \rangle 
+ \langle i w, v _{x_1} \rangle .
\eeq
From the Cauchy-Schwarz inequality we have
$\langle  \phi _{x_1} e ^{i \phi } , v \rangle \in L^1(\R^2) $ and 
$ \langle i w, v _{x_1}\rangle \in L^1(\R^2) $.
As in the proof of Lemma \ref{L2.1} we  obtain 
$\langle i e ^{i \phi } , v \rangle  \in H^1(\R^2)$, hence 
$ \frac{\p }{\p x_1} \left( \langle i e ^{i \phi } , v \rangle \right)  \in \Yo$. 
We conclude that $\langle i\psi , v_{x_1}  \rangle  \in L^1( \R^N) + \Yo $. 
Using (\ref{2.0}), (\ref{2.5}) and the definition of $L$ we get 
$$
L(\langle i\psi_{x_1}, v \rangle +  \langle i\psi , v_{x_1}  \rangle)
= L(\langle i w_{x_1}, v \rangle +  \langle i w , v_{x_1}  \rangle) 
= \int_{\R^N} \langle i w_{x_1}, v \rangle +  \langle i w , v_{x_1}  \rangle \ dx 
$$
and the last quantity is zero by the standard integration by parts formula for 
functions in $H^1(\R^2)$ (see, e.g., \cite{brezis} p. 197). \hfill $\Box$

\begin{Corollary} 
\label{C2.4}
Let $ \psi_1, \; \psi_2  \in \Eo $ be such that $ \psi _1 -\psi _2 \in L^2( \R^N)$. 
Then 
\beq
\label{2.7}
|Q(\psi _1) - Q(\psi _2) | \leq \|\psi _1 - \psi _2\|_{L^2(\R^N)}
\left( \Big\| \frac{\p \psi _1}{\p x_1} \Big\| _{L^2(\R^N)} 
+ \Big\| \frac{\p \psi _2}{\p x_1} \Big\| _{L^2(\R^N)} \right)
\eeq
\end{Corollary}

{\it Proof. }  The same as the proof of Corollary 2.6 in \cite{M10}. 
\hfill 
$\Box$

\medskip

Let $ \psi \in \Eo$. 
It is easy to see that for any  function with compact support $ \phi \in H^1(\R^N)$ 
we have $ \psi + \phi \in \Eo$ and using Lemma \ref{L2.3} we get 
\beq
\label{gateaux}
\lim_{t \ra 0 } \frac 1t (Q(\psi + t \phi ) - Q(\psi ) )
= L(\langle i \psi_{x_1}, \phi \rangle + \langle i \phi_{x_1}, \psi \rangle )
= 2 \int_{\R^N} \langle i \psi _{x_1} , \phi \rangle \, dx.
\eeq

The momentum has a nice behavior with respect to dilations: 
for $ \psi \in \Eo$, $ \la, \; \si > 0$ we have
\beq
\label{2.9}
Q(\psi_{\la,  \si}) = \si ^{N-1} Q(\psi).
\eeq

\section{A regularization procedure} 
\label{reg}

The regularization procedure described below will be an important tool 
for our analysis. It  was first introduced in 
\cite{AB}, then developed in \cite{M10}, where it was   a key ingredient  in proofs. 
It enables us to 
get rid of the small-scale topological defects of functions and in 
the meantime to control the Ginzburg-Landau energy and the momentum  of the regularized functions.

In this section $\Om $ is an open set in $ \R^N$. 
We do not assume $ \Om $ bounded, nor connected. 
If $ \p \Om \neq \emptyset$, we assume that $ \p \Om $ is $ C^2$. 
Fix $ \psi \in \Eo $ and  $  h > 0$. 
We consider the functional
$$ 
G_{h, \Om }^{\psi} (\zeta) = \left\{ 
\begin{array}{ll}
\ds E_{GL}^{\Om } (\zeta) + \frac{1}{h^2} \int_{\Om} |\zeta - \psi |^2 \ dx  \qquad 
& \mbox{ if } N=2,
\\ \\
\ds E_{GL}^{\Om } (\zeta) + 
\frac{ 1}{h^2} \ds \int_{\Om } \ph \left( |\zeta - \psi |^2  \right) \, dx 
& \mbox{ if  } N \geq 3.
\end{array}
\right. $$
Note that $G_{h, \Om }^{\psi} (\zeta) $ may equal  $ \infty $ for some  $ \zeta \in \Eo $;
however, $G_{h, \Om }^{\psi} (\zeta) $ is finite whenever $ \zeta \in \Eo $ and 
$ \zeta - \psi \in L^2( \Om )$. 
We denote 
$ H_0^1(\Om ) = \{ u \in H^1( \R^N) \; | \; u = 0 \mbox{ on } \R^N \setminus \Om \}$ 
and 
$$
H_{\psi} ^1 ( \Om ) = \{ \zeta \in \Eo \; | \; \zeta - \psi \in H_0^1 ( \Om ) \}.
$$

Assume that $ N \geq 3$ and $ \psi = e^{i \al _0 } ( 1+ u ) \in \Eo$, 
where $ \al _ 0 \in [0, 2 \pi ) $ and  $ u \in \Xo$. 
Then 
$$
H_{\psi} ^1 ( \Om ) = \{ e^{i \al _0 } ( 1+ v) \; | \; v \in H_u^1(\Om) \}.
$$
Let 
$$
\tilde{G}_{h, \Om }^u (w) = E_{GL}^{\Om } (1+ w) + 
\frac{ 1}{h^2} \ds \int_{\Om } \ph \left( |w-u|^2  \right) \, dx. 
$$
It is obvious that $ \zeta = e^{i \al _0 } ( 1+ v)$ is a minimizer of 
$ G_{h, \Om }^{\psi} $ in $ H_{\psi}^1(\Om)$ if and only if $v$ is a minimizer 
of $\tilde{G}_{h, \Om }^u $ in $ H_u^1(\Om)$, hence  the results proved 
in \cite{M10} for minimizers of $\tilde{G}_{h, \Om }^u $ also hold 
for minimizers of  $ G_{h, \Om }^{\psi} $.

\medskip

The next three lemmas are analogous to Lemmas 3.1, 3.2 and 3.3 in \cite{M10}.
For the convenience of the reader we give  the full statements in any space 
dimension, but for the proofs in the case $ N \geq 3 $ we refer to \cite{M10}; 
we only indicate here what changes in proofs if $N=2$.

\begin{Lemma}
\label{L3.1}
(i) The functional $ G_{h, \Om }^{\psi} $ has a minimizer in $H_{\psi} ^1 ( \Om ) $. 

(ii) Let $ \zeta _h$ be a minimizer of $ G_{h, \Om }^{\psi}  $ in $H_{\psi} ^1 ( \Om ) $. 
There exist constants $C_i >0$, depending only on $N$,  such that:
\beq
\label{3.1}
E_{GL}^{\Om } ( \zeta_h )\leq E_{GL}^{\Om } ( \psi ); 
\eeq
\vspace{-10pt}
\beq
\label{3.2}
\| \zeta_h - \psi \|_{L^2( \Om ) } ^2  \leq
\left\{ 
\begin{array}{ll}
 h^2 E_{GL}^{\Om } ( \psi ) & \qquad \mbox{ if } N=2,
\\
\\
 h^2 E_{GL}^{\Om } ( \psi ) 
 + C _1 \left( E_{GL}^{\Om } ( \psi ) \right)^{1 + \frac 2N} h^{\frac 4N} 
& \qquad \mbox{ if } N \geq 3.
\end{array}
\right.
\eeq
\vspace{-10pt}
\beq
\label{3.3}
\ds \int_{\Om } \Big\vert \left( \ph ^2( |\zeta_h|) - 1 \right)^2 - 
  \left( \ph ^2( |\psi |) - 1 \right)^2 \Big\vert \, dx 
\leq C_2   h E_{GL}^{\Om } ( u ) ; 
\eeq
\vspace{-10pt}
\beq
\label{3.4} 
| Q(\zeta_h) - Q( \psi) | \leq 
\left\{ 
\begin{array}{ll}
 2h E_{GL}^{\Om } ( \psi ) & \qquad \mbox{ if } N=2,
\\
\\
C_3 \left(  h^2 +  \left( E_{GL}^{\Om } ( \psi ) \right)^{ \frac 2N} 
h^{\frac 4N} \right)^{\frac 12}
E_{GL}^{\Om } ( \psi ) & \qquad \mbox{ if } N \geq 3.
\end{array}
\right.
\eeq

(iii) For $ z \in \C$, denote 
$H(z) = \left( \ph ^2( |z |) - 1 \right) \ph ( |z |) \ph '( |z |) 
\frac{z }{ |z |}$  if $ z \neq 0$ and $H(0) = 0$. 
Then any minimizer $\zeta_h $ of $G_{h, \Om }^{\psi}  $ in $H_{\psi}^1( \Om )$ 
satisfies in $\Do '(\Om )$  the equation
\beq
\label{3.5}
\left\{ 
\begin{array}{ll}
\ds - \Delta \zeta_h +  H( \zeta_h ) 
+ \frac{ 1}{ h^2 } (\zeta _h - \psi ) = 0 & \qquad \mbox{ if } N=2,
\\
\\
\ds - \Delta \zeta_h +  H( \zeta_h ) 
+  \frac{ 1}{  h^2 } \ph ' \left( |\zeta _h - \psi|^2  \right) 
(\zeta _h - \psi ) = 0  & \qquad \mbox{ if } N\geq 3.
\end{array}
\right.
\eeq
Moreover, for any $ \omega \subset \subset \Om $ we have $ \zeta_h \in W^{2, p} ( \omega ) $ 
for   $ p \in [1, \infty )$; thus, in particular, $ \zeta_h \in C^{1, \al } (\omega )$
for $ \al \in [0, 1)$. 

(iv) For any $ h>0$, $ \de >0$ and $ R> 0$ there exists a constant 
$K = K ( N,  h, \de, R) > 0 $ such that  for any $ \psi \in \Eo $ with 
$E_{GL}^{\Om } ( \psi) \leq K$ and for any minimizer 
$\zeta_h$ of $G_{h, \Om }^{\psi} $ in $H_{\psi}^1( \Om )$~we~have 
\beq
\label{3.6}
1 - \de < |\zeta _h(x) | < 1 + \de 
\qquad \mbox{ whenever } x \in \Om \mbox{ and } dist(x, \p \Om ) > 4R.
\eeq
\end{Lemma}

{\it Proof. } Let $N =2$. 

(i)  The existence of a minimizer  is proven  exactly as in Lemma 3.1~in~\cite{M10}.

(ii) Let $\zeta _h$ be a minimizer.  We have $G_{h, \Om }^{\psi} (\zeta _h ) \leq 
G_{h, \Om }^{\psi} (\psi ) = E_{GL } (\psi ) $ and this gives (\ref{3.1}) and 
(\ref{3.2}). It is obvious that 
$$
\big\vert \left( \ph^2(|z_1|) - 1 \right)^2 -  \left( \ph^2(|z_2|) - 1 \right)^2 
\big\vert
\leq 6  \big\vert \ph (|z_1 |) -  \ph (|z_2 |)  \big\vert     \cdot 
\big\vert    \ph (|z_1|^2) + \ph (|z_2|^2) - 2  \big\vert
$$
and $ |\ph (|z_1 |) -  \ph (|z_2 |) | \leq |z_1 - z_2|$.  
Using the Cauchy-Schwarz inequality and (\ref{3.2}) we get 
$$
\begin{array}{l}
\ds \int_{\Om } \Big\vert \left( \ph ^2( |\zeta_h|) - 1 \right)^2 - 
  \left( \ph ^2( |\psi |) - 1 \right)^2 \Big\vert \, dx 
\\
\\
\leq 6  \| \zeta _h - \psi \|_{L^2(\Om )} 
\left( \ds \int_{\Om } 
\Big\vert  \ph ^2( | \zeta _h|) + \ph ^2( | \psi |) -2  \Big\vert ^2\, dx \right)^{\frac 12}
\\
\\
\leq 6   h \left( E_{GL }^{\Om } (\psi )  \right)^{\frac 12} \cdot 
\left( 2  \ds \int_{\Om }   \left( \ph ^2( |\zeta_h|) - 1 \right)^2 +
  \left( \ph ^2( |\psi |) - 1 \right)^2 \, dx \right)^{\frac12} 
\leq  12 \sqrt{2} h E_{GL }^{\Om } (\psi ) 
\end{array}
$$
and (\ref{3.3}) is proven. 
Finally, (\ref{3.4}) follows  from Corollary \ref{C2.4}, (\ref{3.1}) and (\ref{3.2}).

\medskip

(iii) For any 
$ \phi \in C_c^{\infty }(\Om )$ we have $ \zeta_h + \phi \in H_{\psi} ^1( \Om )$ 
and the function $ t \longmapsto G_{h , \Om }^{\psi} ( \zeta_h + t \phi )$ is 
differentiable and achieves its minimum at $ t =0$. 
Hence $\frac{d}{dt }_{\big\vert _{t =0}} \left(G_{h , \Om }^{\psi} ( \zeta_h  + t \phi )\right) =0$
for any $ \phi \in  C_c^{\infty }(\Om ) $ and this is precisely (\ref{3.5}).

For any $ z \in \C$ we have 
\beq
\label{3.7}
|H(z) | \leq 3  | \ph ^2 ( |z |) - 1| \leq 24  . 
\eeq
Since $ \zeta_h \in \Eo$, we have $ \ph ^2 (| \zeta _ h |) - 1 \in L^2 ( \R^2) $  
and the previous inequality  gives $ H( v_h ) \in L^{2}\cap L^{\infty } ( \R^2)$.
We have $\zeta_h, \psi \in H_{loc}^1(\R^2)$ and from the Sobolev embedding theorem 
we get $\zeta_h, \psi \in L_{loc}^p(\R^2)$ for any $ p \in [2, \infty)$. 
Using  (\ref{3.5}) we infer  that $ \Delta \zeta_h \in L_{loc}^p ( \Om )$ 
for any $ p \in [2, \infty)$. Then (iii) follows from standard elliptic 
estimates (see, e.g., Theorem 9.11 p. 235 in \cite{GT}).

\medskip

iv)   Using (\ref{3.7}) we get
$$
\|H(\zeta_h)\|_{L^2(\Om )} \leq 3    \| \ph ^2 ( |\zeta _h |) - 1  \|_{L^2(\Om )} 
\leq 3 \sqrt{2} \left( E_{GL}^{\Om } (\zeta _h ) \right)^{\frac 12} 
\leq 3 \sqrt{2} \left( E_{GL}^{\Om } (\psi ) \right)^{\frac 12}. 
$$
From (\ref{3.5}), (\ref{3.2}) and the above estimate we get
\beq
\label{3.8}
\| \Delta \zeta _h \| _{L^2(\Om ) } \leq \left( 3 \sqrt{2} + \frac 1h \right)  
\left( E_{GL}^{\Om } (\psi ) \right)^{\frac 12}. 
\eeq

For a measurable set  $ \omega \subset \R^N$ with $ \Lo ^N ( \omega ) < \infty$ and for 
 $ f \in L^1( \omega)$, we denote by 
$m( f, \omega) = \frac{1}{\Lo ^N ( \omega ) } \ds \int_{\omega} f (x) \, dx $ 
the mean value of $f$ on $ \omega$.  In particular, if $ f \in L^2(\omega)$ using 
the Cauchy-Schwarz inequality  we get 
$| m(f, \omega)| \leq \left( \Lo ^N (\omega) \right)^{- \frac 12 } \| f \| _{L^2(\omega)} $ 
and consequently
\beq
\label{3.9}
\| m(f, \omega) \| _{L^q (\omega) }
= \left( \Lo ^N (\omega) \right)^{\frac 1q } | m(f, \omega)| 
\leq  \left( \Lo ^N (\omega) \right)^{\frac 1q  - \frac 12} \| f \| _{L^2(\omega)} .
\eeq

Let $ x_0 $ be such that $B( x_0 , 4R) \subset \Om $. 
Using the Poincar\'e inequality and (\ref{3.1}) we have 
\beq
\label{3.10}
\| \zeta_h - m( \zeta_h, B( x_0, 4R)) \| _{L^2( B(x_0, 4R))} 
\leq C_P R \| \nabla \zeta_h \| _{L^2( B(x_0, 4R))}  
\leq C_P R \left( E_{GL}^{\Om } (\psi) \right)^{\frac 12}.
\eeq

It is well-known (see  Theorem 9.11 p. 235 in \cite{GT})
that for $ p \in (1, \infty )$ there exists $ C = C(N, r, p ) > 0$ such that 
for any $ w \in W^{2, p} ( B(a, 2r)) $  we have 
\beq
\label{3.11}
\| w \| _{W^{2, p}(B(a, r))} \leq C \left( 
\| w\| _{L^ p(B(a, 2 r))} + \| \Delta w\| _{L^ p(B(a, 2 r))} \right).
\eeq

From (\ref{3.8}), (\ref{3.10}) and (\ref{3.11}) we get 
\beq 
\label{3.12}
\|  \zeta_h - m( \zeta_h, B( x_0, 4R)) \| _{W^{2,2 } (B(x_0, 2R))}
 \leq C( h, R) \left( E_{GL}^{\Om } (\psi) \right)^{\frac 12}
\eeq
and in particular 
\beq
\label{3.13}
\forall 1 \leq i, j \leq 2, \quad \quad \quad 
\Big\| \frac{ \p ^2 \zeta _h}{\p x_i \p x_j } \Big\|_{L^2  (B(x_0, 2R))}
\leq C( h, R) \left( E_{GL}^{\Om } (\psi) \right)^{\frac 12}.
\eeq

We will use the following variant of the Gagliardo-Nirenberg inequality: 
\beq
\label{3.14}
\| w - m( w,B ( a, r) ) \|_ {L^p ( B(a, r))} 
\leq C(p, q, N, r) \|w \| _ {L^q ( B(a,2 r))} ^{\frac qp} 
\|\nabla w \| _ {L^N ( B(a,2 r))} ^{1- \frac qp} 
\eeq
for any $ w \in W^{1, N}(B(a, 2r))$, where $ 1 \leq q \leq p < \infty $ 
(see, e.g., \cite{kavian} p. 78). 
Using (\ref{3.14}) with $N=2$, $ p =4$, $ q = 2$, then (\ref{3.1}) and (\ref{3.13}) 
we find 
\beq 
\label{3.15}
\begin{array}{l}
\| \nabla \zeta _h - m( \nabla \zeta _h , B(x_0, R)) \|_{L^4(B(x_0, R)) }
\leq C   \| \nabla \zeta _h \| _{L^2 (B(x_0, 2R ) )} ^{\frac 12 } 
       \| \nabla ^2 \zeta _h \| _{L^2 (B(x_0, 2R ) )} ^{\frac 12 }
\\
\\
\leq C( h , R ) \left( E_{GL}^{\Om } (\psi) \right)^{\frac 12}.
\end{array}
\eeq
By (\ref{3.9}) and (\ref{3.1}) we have  
$\| m( \nabla \zeta _h , B(x_0, R)) \|_{L^4(B(x_0, R)) } \leq (\pi R ^2)^{- \frac 14} 
\left( E_{GL}^{\Om } (\psi) \right)^{\frac 12}.$
Together with (\ref{3.15}), this gives
\beq
\label{3.16}
\| \nabla \zeta _h   \|_{L^4(B(x_0, R) )}
\leq  C( h , R ) \left( E_{GL}^{\Om } (\psi) \right)^{\frac 12}.
\eeq

We will use the Morrey inequality which asserts that, for any 
$ w \in C^0 \cap W^{1, p}(B(x_0, r))$ with $ p >N$ we have 
\beq
\label{3.17}
| w(x) - w(y) | \leq C(p, N) | x-y|^{1 - \frac Np} \| \nabla w \|_{L^p(B(x_0, r))}
\qquad \mbox{ for any } x, y \in B(x_0, r)
\eeq
(see  the proof of Theorem IX.12 p. 166 in \cite{brezis}). 
The Morrey inequality and (\ref{3.16}) imply that 
\beq
\label{3.18}
|\zeta_h (x) - \zeta_h (y) | \leq  C_*( h , R ) 
\left( E_{GL}^{\Om } (\psi) \right)^{\frac 12} |x -y|^{\frac 12} 
\quad \mbox{ for any } x, y \in B(x_0, R).
\eeq

Fix $ \de > 0 $. 
Assume that there exists $ x_0 \in \Om $ such that $ dist( x_0, \p \Om ) > 4R $ 
and $\big| \; | \zeta _h (x_0)| - 1 \big| \geq \de $. 
Since 
$\Big| \; \big| \; | \zeta _h (x)| - 1 \big|  - \big| \; | \zeta _h (y)| - 1 \big|  \; \Big| 
\leq |\zeta_h (x) - \zeta _h (y) |$, 
using (\ref{3.18}) we infer that 
$$
\big| \; | \zeta _h (x)| - 1 \big| \geq \frac{\de }{2} 
\quad \mbox{ for any } x \in B(x_0, r_{\de }), 
$$
where $ r_{\de }
 =  \min \left( R, \frac{\de ^2}{4  C_* ^2 ( h , R )  E_{GL}^{\Om } (\psi)} \right).
$ 
Let 
\beq
\label{3.19}
\eta ( s) = \inf \{ ( \ph ^2 ( \tau ) - 1 )^2 \; | \; 
\tau \in ( - \infty, 1 - s ] \cup [1 + s, \infty) \}.
\eeq
It is clear that $ \eta $ is nondecreasing and positive on $(0, \infty )$. 
We have:
\beq
\label{3.20}
\begin{array}{l}
E_{GL}^{\Om }(\psi ) \geq  E_{GL}^{\Om }(\zeta _h) 
\geq \frac 12 \ds \int_{B(x_0,  r_{\de } )} \left( \ph ^2 (|\zeta _h |) - 1 \right)^2 \, dx 
\\
\\
\geq \frac 12   {\ds \int_{B(x_0,  r_{\de } )}  } \eta (\frac{\de }{2} ) \, dx = \frac{\pi}{2}  \eta (\frac{\de }{2} )   r_{\de}^2
= \frac{\pi}{2}   \eta (\frac{\de }{2} ) 
 \min \left( R, \frac{\de ^2}{4  C_* ^2 ( h , R )  E_{GL}^{\Om } (\psi)} \right)^2.
\end{array}
\eeq
It is clear that there exists a constant $ K= K(  h, R, \de ) $ such 
that (\ref{3.20}) cannot hold if $E_{GL}^{\Om }(\psi ) \leq K$. 
We infer that $ \big| \; | \zeta _h (x_0)| - 1 \big| < \de $ whenever 
$ x_0 \in \Om, $ $ dist( x_0, \p \Om ) > 4R$ and $E_{GL}^{\Om }(\psi ) \leq K$. 
\hfill $\Box $

\begin{Lemma}
\label{vanishing}
Let $ (\psi_n)_{n \geq 1} \subset \Eo $ be a sequence of functions satisfying: 

\smallskip

(a) $ (E_{GL}(\psi_n))_{n \geq 1} $ is bounded and

\smallskip

(b) $ \ds \lim_{n \ra \infty } \Big( \sup_{y \in \R^N} E_{GL}^{B(y, 1)} (\psi_n) \Big) =0.$

\smallskip

There exists a sequence $ h_n \lra 0 $ such that for any minimizer $ \zeta_n $ of 
$G_{h_n, \R^N}^{\psi _n}$ in $H_{\psi _n} ^1 (\R^N)$ we have 
$ \| \, | \zeta_n | - 1 \|_{L^{\infty }(\R^N)} \lra 0 $ as $ n \lra \infty $.
\end{Lemma}

{\it Proof. } Let  $N=2$.  We split the proof into several steps. 

\medskip

{\it Step 1. Choice of the sequence $(h_n)_{n \geq 1}$.} 
Let $M = \ds \sup_{n \geq 1} E_{GL}(\psi_n)$. 
For $ n \geq 1 $ and $ x \in \R^2$ we denote
$$
m_n (x) = m(\psi _n, B(x, 1)) = \frac{1}{ \pi } \int_{B(x, 1)} \psi _n (y) \, dy. 
$$
The Poincar\'e inequality implies that there exists $C_P > 0$ such that 
$$
\ds \int_{B(x, 1)} | \psi _n (y) - m_n (x) |^ 2 \, dy \leq C_P \int_{B(x,1)} | \nabla \psi _n |^2 \, dy.
$$
Using assumption (b) we find
\beq
\label{3.21}
\sup_{x \in \R^2} \| \psi _n - m_n (x) \| _{L^2(B(x,1))} \lra 0 \quad 
\mbox{ as } n \lra \infty.
\eeq
Proceeding exactly as in the proof of Lemma 3.2 in \cite{M10} (see the 
proof of (3.35) there) we get
\beq
\label{3.22}
\lim_{ n \ra \infty} \sup_{x \in \R^2} |H(m_n(x))| = 0.
\eeq
Let 
\beq
\label{3.23}
h_n = \max \left( \left( \sup_ {x \in \R^2}  \| \psi _n - m_n (x) \| _{L^2(B(x,1))} \right)^{\frac 13} , 
\sup _{x \in \R^2} |H(m_n(x))| \right).
\eeq
From (\ref{3.21}) and (\ref{3.22}) it follows that 
$ h_n \lra 0 $ as $ n \lra \infty $. 
Hence we may assume that $ 0 < h_n < 1$ for each $n$ 
(if $ h_n = 0$  then  $ \psi _n$ is constant  a.e. 
and any  minimizer $ \zeta_n $ of 
$G_{h_n, \R^2}^{\psi _n}$ equals $\psi _n$  a.e.). 

Let $ \zeta_n $ be a minimizer of $G_{h_n, \R^2}^{\psi _n}$,
as given   by Lemma \ref{L3.1} (i).
It follows from Lemma \ref{L3.1} (iii) that $\zeta_n $ satisfies (\ref{3.5})
 and $\zeta _n \in W_{loc}^{2,2}(\R^2)$.

\medskip

{\it Step 2.  We  prove that $\| \Delta \zeta _n \| _{L^2(B(x, \frac 12 ))}$ is 
bounded independently of $n$ and of $x$. }
There is no loss of generality to assume  that $ x = 0 $. 
Then we observe that (\ref{3.5}) can be written as
\beq
\label{3.24}
- \Delta \zeta _n  + \frac{1}{h_n^2} ( \zeta _n -  m_n(0)) = f_n \qquad \mbox{ in } 
\Do' (\R^2), 
\eeq
 where 
\beq
\label{3.25}
 f_n = \frac{1}{h_n^2} ( \psi _n - m_n(0)) -  (H(\zeta _n) - H( m_n(0))) -  H(m_n(0)).
\eeq
From  (\ref{3.2}) we have
$\| \zeta _n - \psi _n \| _{L^2(\R^2)} \leq h_n E_{GL}(\psi _n)^{\frac 12} 
\leq h_n M^{\frac 12}$
and from  (\ref{3.23}) we obtain  $ \| \psi _n - m_n(0)  \|_{L^2(B(0,1))} 
\leq h_n ^3 \leq h_n$,
hence
\beq
\label{3.26}
\| \zeta _n - m_n(0) \|_{L^2(B(0,1))} \leq (M^{\frac 12} +1) h_n.
\eeq
Since $H$ is Lipschitz, we get 
\beq
\label{3.27}
\| H(\zeta _n ) - H( m_n(0) )  \| _{L^2(B(0,1)} 
\leq C_1 \| \zeta _n - m_n(0) \|_{L^2(B(0,1))} 
\leq C_2 h_n.
\eeq
Using (\ref{3.25}), (\ref{3.23}) and (\ref{3.27}) we get
\beq
\label{3.28}
\begin{array}{l}
\| f_n \| _{L^2(B(0,1))} 
\\
\leq \frac{1 }{h_n^2} \| \psi _n - m_n (0) \| _{L^2(B(0,1))} 
+  \| H(\zeta _n ) - H( m_n(0) )  \| _{L^2(B(0,1))} 
+  \pi ^{\frac 12} | H(m_n(0) ) | 
\\
\leq C_3 h_n.
\end{array}
\eeq
It is obvious that for any bounded domain $ \Om \subset \R^2$, each term in 
(\ref{3.24}) belongs to $H^{-1}(\Om)$. 
Let $ \chi \in C_c^{\infty}(\R^2)$ be such that $ supp( \chi ) \subset B(0,1)$, 
$ 0 \leq \chi \leq 1$ and $ \chi = 1 $ on $B(0, \frac 12)$. 
Taking the duality product of (\ref{3.24}) by $ \chi ( \zeta _n - m_n(0))$ we find
\beq
\label{3.29}
\int_{\R^2} \chi | \nabla \zeta _n |^2 \, dx - \frac 12 \int_{\R^2} (\Delta \chi ) | 
\zeta _n - m_n(0) |^2 \, dx
+ \frac{1}{h_n ^2}  \int_{\R^2} \chi  |  \zeta _n - m_n(0) |^2 \, dx 
= \int_{\R^2} \langle f_n , \zeta_n - m_n(0) \rangle \chi  \, dx.
\eeq
Using (\ref{3.29}), the Cauchy-Schwarz inequality and (\ref{3.26}), (\ref{3.28}) we infer that 
\beq
\label{3.30}
\begin{array}{l}
\ds \frac{1}{h_n ^2}  \int_{B(0, \frac 12)}  |  \zeta _n - m_n(0) |^2 \, dx 
\\
\\
\ds \leq \| \Delta \chi \|_{L^{\infty} (\R^2) } \int_{B(0, 1)}  |  \zeta _n - m_n(0) |^2 \, dx
+ \| f_n \| _{L^2(B(0,1))} \| \zeta _n - m_n(0)  \|_{L^2(B(0,1))} 
\leq C_4  h_n ^2.
\end{array} 
\eeq
Now (\ref{3.24}), (\ref{3.28}) and (\ref{3.30}) imply that there is $C_5 >0 $ 
such that $\| \Delta \zeta _n \| _{L^2(B(0, \frac 12))}  \leq C_5$. Thus we 
have proved that for any $n$ and $x$,
\beq
\label{3.31}
\| \Delta \zeta _n \| _{L^2(B(x, \frac 12))}  \leq C_5,
\qquad \mbox{  where  } C_5 \mbox{  does not depend on $x$ and $n$.}
\eeq

\medskip

{\it Step 3. A H\"older estimate on $ \zeta _n$. } It follows from (\ref{3.11}) that 
\beq
\label{3.32}
\| \zeta _n  - m_n\| _{W^{2,2} (B(x, \frac 14) )} \leq 
C ( \| \Delta \zeta _n  \|_{L^2(B(x, \frac 12 ))} + 
\| \zeta _n - m_n \| _{L^2(B(x, \frac 12 ))}  ) \leq C _6.
\eeq
From (\ref{3.14}) and (\ref{3.32}) we find 
\beq
\label{3.33}
\| \nabla  \zeta _n  - m ( \nabla  \zeta _n , B(x, \frac 18 )) \| _{L^4(B(x, \frac 18))} 
\leq C \| \nabla  \zeta _n  \| _{L^2(B(x, \frac 14))} ^{\frac 12} 
 \| \nabla ^2 \zeta _n  \| _{L^2(B(x, \frac 14))} ^{\frac 12}  \leq C_7.
\eeq
It is clear that 
$|  m ( \nabla  \zeta _n , B(x, \frac 18 ))  | \leq \left( \Lo ^2(  B(x, \frac 18 )) \right)^{- \frac 12}
\|  \nabla  \zeta _n \| _{L^2 (  B(x, \frac 18 )) } \leq C_8.$
Then (\ref{3.33}) implies that $ \| \nabla  \zeta _n  \| _{L^4(B(x, \frac 18))}  $ 
is bounded independently of $n$ and of $x$. Using the Morrey inequality 
(\ref{3.17}) we infer  that there is $C_9 >0$ such that 
\beq
\label{3.34}
|\zeta _n (x) - \zeta _n (y) | \leq C_9 |x-y |^{\frac 12} \qquad 
\mbox{ for any $ n \in \N^*$ and any $x, y \in \R^2$ with } |x-y| < \frac 18.
\eeq

\medskip

{\it Step 4. Conclusion. } 
Let $ \de _n = \| \, |\zeta _n | - 1 \|_{L^{\infty}(\R^2)}$  if $\zeta _n $ 
is bounded, and $ \de _n = 1$ otherwise.
Choose $ x_0^n \in \R^2$  such that  $\big|\,  |\zeta _n ( x_0 ^n) | - 1 \big| \geq \frac{\de _n }{2} .$ 
From (\ref{3.34}) we infer that $\big|\,  |\zeta _n (x) | - 1 \big| \geq \frac{\de _n }{4}$ 
for any $ x \in B(x_0 ^n, r_n)$, 
where $ r_n = \min \left( \frac 18, \left(\frac{\de _n }{ 4 C_9 } \right)^2 \right).$
Let $ \eta $ be as in (\ref{3.19}).  Then we have 
\beq
\label{3.35}
\int_{B(x_0 ^n, r_n) } \left( \ph ^2(|\zeta _n |) - 1 \right)^2\, dx \geq 
\int_{B(x_0 ^n, r_n) }  \eta \left( \frac{\de _n}{4} \right) \, dx = 
\eta \left( \frac{\de _n}{4} \right) \pi r_n ^2.
\eeq

On the other hand, 
the function $ z\longmapsto \left( \ph ^2(| z|) - 1 \right)^2 $ is 
Lipschitz on $\C$. 
From this fact, the Cauchy-Schwarz inequality, (\ref{3.2}) and assumption 
(a) we get
$$
\begin{array}{l}
\ds \int_{B(x, 1)} \Big\vert 
\left( \ph ^2(|\zeta _n (y)|) - 1 \right)^2
- \left( \ph ^2(|\psi _n (y) |) - 1 \right)^2 \Big\vert  \, dy 
\\
\leq C \ds \int_{B(x, 1)} | \zeta_n (y) - \psi_n (y) | \, dy 
\leq C \pi ^{\frac 12} \|  \zeta_n - \psi_n   \| _{L^2 ( B(x, 1) )} 
\leq C  \pi ^{\frac 12} \|  \zeta_n - \psi_n   \|_{L^2 ( \R^2) } 
\leq C_{10} h_n  .
\end{array}
$$
Then using assumption (b) we infer that 
\beq
\label{3.36}
\ds \sup_{ x \in \R^2} \ds \int_{B(x, 1)}
\left( \ph ^2(|\zeta _n (y)|) - 1 \right)^2 \, dy 
\lra 0 \qquad \mbox{ as } n \lra \infty.
\eeq

From (\ref{3.35}) and (\ref{3.36}) we get 
$ {\ds \lim_{n \ra \infty }} \eta \left( \frac{\de_n}{4} \right) r_n ^2 = 0$
and this clearly implies $ \ds \lim_{n \ra \infty }  \de _n = 0$. 
This completes the proof of Lemma \ref{vanishing}.
\hfill
$\Box $

\medskip

The next result is based on Lemma \ref{L3.1} and will be very useful 
in the next sections to 
prove the "concentration" of minimizing sequences.
For $0 < R_1 < R_2 $ we denote 
$ \Om _{R_1, R_2} = B(0, R_2) \setminus \ov{B}(0, R_1)$. 

\medskip

\begin{Lemma}
\label{splitting}
Let $ A > A_3 > A_2 >1$. 
There exist  $ \e _0 >0$ and 
$C_i >0$, depending only on $  \; N, \; A, \; A_2, \; A_3$
(and $F$ for (vi)) 
 such that for any $ R \geq 1$, 
$ \e \in (0, \e _0)$ and $ \psi \in \Eo $  verifying
$E_{GL}^{\Om _{R, AR}} (\psi) \leq \e, $
there exist two functions $ \psi_1, \, \psi_2 \in \Eo $ and a constant $ \theta _0 \in [0, 2 \pi)$ 
satisfying the following properties:

\smallskip

(i) $\psi _1 = \psi $ on $B(0, R)$ and $ \psi _1 = e^{i \theta _0}$ on 
$\R^N  \setminus B(0, A_2 R) $,

\smallskip

(ii) $ \psi_2 = \psi  $ on $ \R^N \setminus B(0, AR)$ and 
$ \psi _2 =  e^{ i \theta _0 }= constant$
on $B(0, A_3R)$,  

\medskip

(iii) $ \ds \int_{\R^N} \Big\vert  \, 
 \Big\vert \frac{\p \psi}{\p x_j } \Big\vert ^2 - 
 \Big\vert \frac{\p \psi_1}{\p x_j } \Big\vert ^2 - 
 \Big\vert \frac{\p \psi_2}{\p x_j } \Big\vert ^2 \, \Big\vert \, dx 
 \leq C_1 \e  $ for $ j =1, \dots, N$, 
 
\medskip
 
(iv) $ \ds \int_{\R^N} \Big\vert  
\left( \ph ^2(| \psi |) - 1 \right)^2 
- \left( \ph ^2(|\psi_1|) - 1 \right)^2 
- \left( \ph ^2(|\psi_2|) - 1 \right)^2 
\Big\vert   \, dx \leq C_2 \e $, 

\medskip

(v) $ | Q(\psi) - Q(\psi_1) - Q( \psi_2) | \leq C_3 \e $, 

\medskip

(vi) If assumptions (A1) and  (A2) in the introduction hold, then
$$
 \ds \int_{\R^N}  \Big\vert V( |\psi  |^2) - V( |\psi _1|^2) - V( |\psi _2 |^2) \Big\vert \, dx 
\leq
\left\{
\begin{array}{l}
 C_4 \e + C_5 \sqrt{\e } \left( E_{GL}(\psi) \right)^{\frac{2^* -1}{2}} \quad \mbox{ if } N \geq 3, 
\\
\\
C_6 \e + C_7 \sqrt{\e } \left( E_{GL}(\psi ) \right) ^{ p_0 + 1} \quad \mbox{ if } N =2.
\end{array}
\right.
$$ 
Furthermore, the same estimate holds with $ V_+$ (respectively $V_-$) instead of $V$. 

\end{Lemma}

{\it Proof.  } 
If $ N \geq 3$, this is  Lemma 3.3 in \cite{M10}.

Let $ N =2$. 
Fix $ k > 0$,  $ A_1 $ and  $  A_4 $ such that 
$ 1 + 4k < A_1 < A_2 < A_3 < A_4 < A - 4k. $
Let $ h = 1 $ and $ \de = \frac{1}{2}$. 
Let   $K( N, h, \de, r )$  be  as in Lemma \ref{L3.1} (iv). 
We will prove that Lemma \ref{splitting} holds for 
$\e _0 = \min \left( K( 2, 1,  \frac{1}{2}, k ),   
 \frac{\pi}{8} \ln \left( \frac{A - 4k}{1 + 4k} \right)    \right)   .$

Fix $ \e < \e _0 $. Consider  $\psi \in \Eo $  such that 
$ E_{GL}^{\Om _{R, AR}}(\psi ) \leq \e$. 
Let $\zeta  $ be a minimizer of $G_{1, \Om _{R, AR}}^{\psi }$ in the space $H_{\psi }^1(\Om _{R, AR})$.
Such minimizers exist by Lemma \ref{L3.1} (but are perhaps not  unique). 
From Lemma \ref{L3.1} (iii) we have $\zeta  \in W_{loc}^{2,p} (\Om _{R, AR}) $ 
for any $ p \in [1, \infty)$, hence $ \zeta   \in C^1  (\Om _{R, AR}) $. 
Moreover, Lemma \ref{L3.1} (iv) implies that 
\beq
\label{3.37}
\frac{1}{2} \leq |\zeta  (x) | \leq  \frac{ 3 }{2} 
\qquad \mbox{ for any } x \mbox{ such that } R + 4k \leq |x| \leq AR - 4k.
\eeq
Therefore, the topological degree $ deg( \frac{ \zeta }{|\zeta |} , \p B(0, r))$ is well defined for any 
$ r \in [R+ 4k, AR - 4k]$ and does not depend on $r$.
It is well-known that $ \zeta $ admits a $C^1$ lifting $\theta$ (i.e.
$ \zeta  = |\zeta  |e^{i \theta} $) on $\Om_{R+ 4k, AR- 4k} $
if and only if
$deg( \zeta  , \p B(0, r))= 0$ for $ r \in (R + 4k, AR - 4k)$.
Denoting by $ \tau = (- \sin t, \cos t) $ the unit tangent vector at
$\p B(0, r)$ at a point $ r e^{i t} = (r \cos t, r \sin t ) \in \p B(0, r)$, we get
\beq
\label{3.38}
\begin{array}{l}
| deg( \zeta  , \p B(0, r)) |
= \bigg\vert \ds \frac{1}{2 i \pi} \int_0 ^{2 \pi}
\frac{\frac{\p }{\p t } ( \zeta  (r  e^{ i t}))}
{\zeta (r e^{i t })} \, d t
\bigg\vert
= \bigg\vert \ds \frac{r}{2 i \pi} \int_0 ^{2 \pi}
\frac{\frac{\p \zeta }{\p \tau} (r e^{ i t})}
{\zeta (r e^{i t })} \, d t
\bigg\vert
\\
\\
\leq \ds \frac{r}{2 \pi} \int_0^{2 \pi} 2  | \nabla \zeta  (r e^{i t})| \, d t
\leq \frac{r}{\pi } \sqrt{ 2 \pi} \left(
\int_0^{2 \pi} | \nabla \zeta  (re^{i t })| ^2 \, d t \right)^{\frac 12}.
\end{array}
\eeq
On the other hand,
$$
\ds \int_{ \Om_{R+ 4k, AR- 4k} } |\nabla \zeta  (x)|^2 \, dx =
\int_{R + 4k }^{AR - 4k  } r \int_0^{2 \pi} | \nabla \zeta  (re^{i t})| ^2 \, d t \, dr.
$$
We have
 ${\ds \int_{\Om_{R+ 4k, AR- 4k} }  } |\nabla \zeta  (x)|^2 \, dx \leq E_{GL}^{\Om_{R, AR}} (\zeta )
\leq   E_{GL}^{\Om_{R, AR}} (\psi)  < \e _0 \leq \frac{\pi }{8} \ln \left( \frac{AR - 4k}{R + 4k} \right)$
 and we infer that there exists
$r _* \in (R + 4k, AR - 4k) $ such that
$
\ds r _* \int_0^{2 \pi} | \nabla \zeta  (R_* e^{i t})| ^2 \, d t
< \frac{ \pi }{8} \frac {1}{r _*}.
$
From (\ref{3.38}) we get
$$
| deg( \zeta , \p B(0, r _* )) |  < \frac{r _*}{\pi } \sqrt{ 2 \pi}
\left(  \frac{ \pi }{8} \frac {1}{r _* ^2} \right)^{\frac 12}
= \frac 12.
$$
Since the topological degree is an integer, we have necessarily
$deg( \zeta , \p B(0, r _*))  =0$.
Consequently $deg( \zeta , \p B(0, r))  =0$ for any
$  r \in ( R+ 4k, AR- 4k) $ and $ \zeta  $ admits a $C^1$  lifting 
$\zeta = \rho e^{i  \theta }$. In fact, 
$\rho, \, \theta \in W_{loc}^{2, p}(\Om_{R+ 4k, AR- 4k} ) $ 
because $\zeta \in W_{loc}^{2, p}(\Om_{R+ 4k, AR- 4k}) $ 
(see Theorem 3 p. 38 in \cite{BBM}).

Consider  $ \eta _1, \eta _2 \in C ^{\infty }(\R)$ satisfying  
the following properties:
$$
\begin{array}{l}
\eta_1 = 1 \mbox{ on } (-\infty, A_1], \quad \eta_1 = 0 \mbox{ on } [A_2, \infty), 
\quad \eta_1 \mbox{ is nonincreasing, }
\\
\eta_2 = 0 \mbox{ on } (-\infty, A_3], \quad \eta_2 = 1 \mbox{ on } [A_4, \infty), 
\quad \eta_2 \mbox{ is nondecreasing. }
\end{array}
$$

Denote $ \theta _0 = m( \theta, \Om _{A_1 R, A_4 R })$.
We define $ \psi _1$ and $ \psi _2$ as follows: 
\beq
\label{3.39}
\psi _1 (x) = \left\{ \begin{array}{l} 
\psi (x)  \quad \mbox{ if } x \in \ov{B}(0, R), 
\\
\zeta  (x) \quad \mbox{ if } x \in {B}(0, A_1R) \setminus \ov{B}(0, R), 
\\
\left(1 + \eta_1 (\frac{|x|}{R} ) (\rho  (x)  - 1  )\right) 
e^{i \left(\theta _0 + \eta_1 (\frac{|x|}{R} ) (\theta (x) - \theta_0 )\right)}
\\
\qquad \qquad \qquad \qquad  \qquad  \quad \mbox{ if } x \in {B}(0, A_4R) \setminus B(0, A_1R),
\\
 e^{i \theta _0}  \quad \mbox{ if } x \in \R^2 \setminus {B}(0, A_4R ), 
\end{array}
\right.
\eeq
\beq
\label{3.40}
\psi _2 (x) = \left\{ \begin{array}{l} 
 e^{i \theta _0}  \quad \mbox{ if } x \in  \ov{B}(0, A_1R ), 
\\
\left(1 + \eta_2 (\frac{|x|}{R} ) (\rho (x) - 1 )\right) 
e^{i \left(\theta _0 + \eta_2 (\frac{|x|}{R} ) (\theta (x) - \theta_0 )\right)}
\\
\qquad \qquad \qquad \qquad \quad \qquad  \mbox{ if } x \in {B}(0, A_4R) \setminus \ov{B}(0, A_1R),
\\
\zeta  (x) \quad \mbox{ if } x \in {B}(0, AR) \setminus {B}(0, A_4 R),
\\
\psi (x)  \quad \mbox{ if } x \in \R^2 \setminus {B}(0, AR ). 
\end{array}
\right.
\eeq
Then $ \psi _1, \; \psi _2 \in \Eo $ and satisfy (i) and (ii). 
The proof of (iii), (iv) and (v) is exactly as in \cite{M10}. 
Next we prove (vi). 

Assume that (A1) and (A2) are satisfied and let $W(s) = V(s) - V(\ph^2(\sqrt{s}))$. 
Then $W(s) = 0 $ for $ s \in [0,4]$ and it is easy to see that $W$ satisfies
\beq
\label{i3}
|W(b^2) - W(a^2) | \leq C_3 |b-a| 
\left( a^{2p _0 +1} \1_{\{ a > 2  \} } + b^{2p _0 +1} \1_{\{ b > 2  \} }  \right)
\; \;  \mbox{ for any } a, b \geq 0.
\eeq
Using (\ref{i1}) and (\ref{i3}), then H\"older's inequality we obtain
\beq
\label{3.41}
\begin{array}{l}
\ds \int _{\R ^2  } 
\big\vert V( | \psi|^2) - V( | \zeta |^2) \big\vert \, dx
\\
\\
\leq 
\ds \int _{\Om _{ R, \, A R} }
\big\vert V( \ph ^2(| \psi |) ) \! - \! V( \ph ^2(|\zeta | )) \big\vert 
+ \big\vert W( |\psi |^2) \! - \! W( |\zeta |^2) \big\vert  dx 
\\
\\
\leq C \ds \int _{\Om _{ R, \, A R} }  
\left( \ph ^2(| \psi |) - 1 \right)^2 + \left( \ph ^2(| \zeta |) - 1 \right)^2 \, dx 
 \\
 \\
 \qquad + C \ds \int _{\Om _{ R, \, A R} } 
 \big\vert \, | \psi | - | \zeta  | \, \big\vert
 \left( | \psi| ^{2 p_0 + 1} \1 _{\{ | \psi | > 2  \} } 
 + | \zeta | ^{2 p_0 + 1} \1 _{\{ |\zeta | > 2  \} } \right)\, dx
\\
\\
\leq \! C' \e + 
\| \psi \! -  \! \zeta  \|_{L^2(\Om _{ R, \, A R} ) } \! \!
\left[ \! \!  {\left( \ds \int _{\Om _{ R, \, A R} }   \!  \! \! \! \! \! 
| \psi| ^{4 p_0 + 2} \1 _{\{ | \psi | > 2  \} }   dx \!   \right) \! \! } ^{\frac 12} 
\! \!+ \!
 {\left( \ds \int _{\Om _{ R, \, A R} }    \!  \! \! \! \! \! 
| \zeta | ^{4 p_0 + 2} \1 _{\{ | \zeta  | > 2  \} }   dx \! \right)  \! \!  }^{\frac 12}  \right]  \! \! .
\end{array}
\eeq

Using (\ref{ineq2})  we get 
\beq 
\label{3.42} 
\int_{\R^2 }  | \psi| ^{4 p_0 + 2} \1 _{\{ | \psi | > 2  \} }  \, dx  
\leq C \| \nabla \psi \| _{L^2(\R^2 )} ^{ 4 p_0 + 2} 
\Lo ^2\left( \{ x \in \R^2 \; | \; |\psi (x) | \geq 2  \} \right).
\eeq
On the other hand, 
\beq
\label{3.43}
9 
\Lo ^2\left( \{ x \in \R^2 \; \big| \; |\psi (x) | \geq 2   \} \right)
\leq \int_{\R^2} \left( \ph ^2 ( |\psi |) - 1 \right)^2 \, dx 
\leq 2 E_{GL}(\psi ) 
\eeq
and  a similar estimate holds for $ \zeta $. 
We insert (\ref{3.42}) and (\ref{3.43}) into (\ref{3.41}) to discover
\beq
\label{3.44}
 \int _{\R ^2  } 
\big\vert V( | \psi|^2) - V( | \zeta |^2) \big\vert \, dx \leq C '\e + C \sqrt{\e } \left( E_{GL}(\psi ) \right)^{p_0 +1 }.
\eeq
Proceeding exactly as in \cite {M10} (see the proof of (3.88) p. 144 there) we obtain
\beq
\label{3.45}
\int_{\R^2} \big\vert V(|\zeta |^2 ) - V( |\psi _1|^2) - V(|\psi _2|^2)  \big\vert \, dx \leq C \e.
\eeq
Then (vi) follows from (\ref{3.44}) and (\ref{3.45}).
\hfill
$\Box $

 \begin{Corollary} 
\label{C3.4}
For any $ \psi \in \Eo$, there is a sequence of functions
 $(\psi _n )_{n \geq 1} \subset \Eo $ satisfying:

(i) $\psi _n = \psi $ on $B(0, 2^n)$ and $ \psi _n =  e^{i \theta _n } = constant $ 
on $\R^N \setminus B(0, 2^{n+1})$, 

(ii) $\| \nabla \psi  _n - \nabla \psi \|_{L^2(\R^N)} \lra 0$ and 
$\| \ph ^2 (|\psi _n|) - \ph ^2(|\psi |) \| _{L^2(\R^N)} \lra 0$, 

(iii) $ Q(\psi _n ) \lra Q(\psi)$, $\ds \int_{\R^N} \big\vert V(|\psi _n|^2 ) - V(|\psi |^2) \big\vert \, dx \lra 0 $ 
and \\ 
$\ds \int_{\R^N} \big\vert \left( \ph ^2(|\psi _n |) - 1 \right)^2 
-  \left( \ph ^2(|\psi |) - 1 \right)^2  \big\vert \, dx \lra 0 $ as $ n \lra \infty.$

\end{Corollary}

{\it Proof. } Let $ \e _n = E_{GL}^{\R^N \setminus B(0, 2^n) }(\psi)$, so that $ \e _n \lra 0 $ as $ n \lra \infty$. 
Let $ A = 2$,  fix $ 1 < A_2 < A_3 < 2 $ and use Lemma \ref{splitting} with $ R = 2^n$
to obtain two functions $ \psi _1^n$, $\psi _2^n$ with properties (i)-(vi) in that Lemma. 
Let $ \psi _n = \psi _1^n$. It is then  straightforward to prove that $(\psi _n)_{n \geq 1}$ 
satisfies (i)-(iii) above. 
\hfill
$\Box$

\medskip

The next Lemma allows  to approximate functions in $ \Eo $ by functions 
with higher regularity. 

\begin{Lemma}
\label{L3.5}
(i) Assume that $ \Om = \R^N$ or that $ \p \Om $ is $ C^1$. 
Let $ \psi \in \Eo $. For each $ h > 0$,  let $ \zeta_h $ be a minimizer of 
$G_{h, \Om } ^{ \psi }$ in $H_{\psi }^1( \Om)$. Then 
$
\| \zeta _h - \psi  \| _{H^1( \Om )} \lra 0 $
as $ h \lra 0. $

\medskip

(ii) Let $ \psi \in \Eo $. For any $ \e > 0$ and any $ k \in \N$ there is $ \zeta \in \Eo $ such that 
$ \nabla \zeta \in H^k( \R^N )$, $E_{GL}(\zeta ) \leq E_{GL}(\psi) $ and 
$\| \zeta - \psi \|_{
H^1( \R^N)} < \e$. 

\end{Lemma}

{\it Proof. } (i) It suffices to prove that for any sequence $ h_n \lra 0 $ and any choice of a minimizer 
$ \zeta _{n}$ of  $G_{h _n, \Om } ^{ \psi }$ in $H_{\psi }^1( \Om)$, there is a subsequence 
$ (\zeta_{n_k})_{k \geq 1}$  such that 
$
\ds \lim_{k \ra \infty} \| \zeta _{h_{n_k}} - \psi  \| _{H^1( \Om )} = 0. $

 Let $ h_n \lra 0 $ and let $ \zeta _{n}$ be as above. 
By (\ref{3.2}) we have $ \zeta_{n} - \psi \lra 0 $ in $L^2( \Om )$ and it is clear that 
$ \zeta_{n} - \psi $ is bounded in $H_0^1( \Om)$. 
Then there are $v \in H_0^1( \Om )$ and a subsequence $ (\zeta_{n_k})_{k \geq 1}$  
such that 
$$
(\zeta_{n_k} - \psi) \rightharpoonup v \quad \mbox{ weakly in } H_0^1( \Om)
\quad \mbox{ and } \quad 
(\zeta_{n_k} - \psi) \lra v \quad \mbox{  a.e. on } \Om.
$$
Since $ \zeta_{n_k} - \psi \lra 0 $ in $L^2( \Om)$ we infer that $ v = 0 $ a.e., therefore 
$\nabla \zeta_{n_k}   \rightharpoonup  \nabla \psi $ 
weakly in $ L^2( \Om )$ and $ \zeta_{n_k} \lra \psi $ a.e on $ \Om$. 
By weak convergence we have 
$ 
\ii _{\Om } |\nabla \psi |^2 \, dx 
\leq {\ds \liminf _{k \ra \infty }}  \ii_{\Om } |\nabla  \zeta_{n _k} |^2 \, dx
$ 
and Fatou's Lemma gives
$
 \ii _{\Om }  \left( \ph^2(|\psi |) - 1 \right)^2 \, dx \leq 
{\ds \liminf _{k \ra \infty } } \ii _{\Om } \left( \ph^2(|\zeta _{n_k} |) - 1 \right)^2 \, dx .
$
Thus we get $ E_{GL}^{\Om} ( \psi ) \leq \ds \liminf _{k \ra \infty }  E_{GL}^{\Om} (\zeta_{n_k}).$
On the other hand we have $ E_{GL}^{\Om} (\zeta_{n_k}) \leq  E_{GL}^{\Om} ( \psi )  $ for all $k$. 
We infer that necessarily $ \ds \lim _{k \ra \infty }  E_{GL}^{\Om} (\zeta_{n_k}) = E_{GL}^{\Om} ( \psi ) $
and 
$ {\ds \lim _{k \ra \infty } }  \ii _{\Om } |\nabla  \zeta_{n _k} |^2 \, dx = \ii _{\Om } |\nabla \psi |^2 \, dx.$
Taking into account that $\nabla \zeta_{n_k}   \rightharpoonup  \nabla \psi $ 
weakly in $ L^2( \Om )$, we deduce that $\nabla \zeta_{n_k}   \lra  \nabla \psi $ 
strongly in $ L^2( \Om )$, thus $( \zeta_{n _k } - \psi ) \lra 0 $ in $H_0^1( \Om)$, as desired. 

\medskip

(ii) Let $ h > 0$ and let $ \zeta _h $ be a minimizer of $ G_{h, \R^N} ^{\psi }.$
Then $ \zeta _h$ satisfies (\ref{3.5}) in $ \Do '(\R^N)$, thus $ \Delta \zeta _h \in L^2( \R^N)$ 
and this implies $ \frac{ \p ^2 \zeta _h }{\p x_i \p x_j } \in  L^2( \R^N)$ for any $i, j$, 
hence $ \nabla \zeta _h \in H^1( \R^N)$. 
Moreover, if $ \nabla \psi \in H^{\ell} ( \R^N)$ for some $ \ell \in \N$, 
taking successively the derivatives of (\ref{3.5}) up to order $ \ell $ and repeating the above argument we get 
$ \nabla \zeta _h \in H^{\ell +1}( \R^N)$. 

Fix $ \psi \in \Eo$, $ k \in \N$ and $ \e > 0$. 
Using (i), there are $ h_1 > 0$ and a minimizer $ \zeta _1 $ of $G_{h_1, \R^N}^{\psi } $ such that 
$\| \zeta _1 - \psi \| _{H^1( \R^N ) } < \frac{\e}{2}$  and $ \nabla \zeta _1 \in H^1( \R^N)$. 
Then there are $ h_2 > 0$ and a minimizer $ \zeta _2 $ of $G_{h_2, \R^N}^ {\zeta _1}$ such that 
$\| \zeta _2 - \zeta _1 \|_{H^1( \R^N ) } < \frac{\e}{2^2}$  and $ \nabla \zeta _2 \in H^2( \R^N)$, 
and so on. 
After $k$ steps we find $ h_k $ and $ \zeta _k $ such that $ \zeta _k $ is a minimizer of $G_{h_k, \R^N}^{\zeta_{k-1}}$, 
$\| \zeta _k - \zeta _{k-1} \|_{H^1( \R^N ) } < \frac{\e}{2^k}$, and  
$ \nabla \zeta _k \in H^k( \R^N)$. Then 
$ \| \zeta _k - \psi \| _{H^1( \R^N ) } <  \| \zeta _k - \zeta_{k-1} \| _{H^1( \R^N ) } + 
\dots + \| \zeta _2 - \zeta _1 \| _{H^1( \R^N ) } + \| \zeta _1 - \psi \| _{H^1( \R^N ) } < \e$. 
Moreover, 
$ E_{GL}(\zeta_{k}) \leq  E_{GL}(\zeta_{k-1}) \leq \dots \leq E_{GL}(\psi).$ 
\hfill
$\Box$

\section{Minimizing the energy at fixed momentum}
\label{minem}

The aim of this section is to investigate the existence of minimizers of the 
energy $E$ under the constraint $Q = q >0$. If such minimizers exist, they are 
traveling waves to (\ref{1.1}) and their speed is precisely the Lagrange 
multiplier appearing in the variational problem.

We start with some useful properties of the functionals $E$, $E_{GL}$ and $Q$. 

\begin{Lemma}
\label{L4.1} If (A1) and (A2) in the Introduction hold, then 
 $V(|\psi|^2) \in L^1(\R^N)$ whenever $ \psi \in \Eo $.
Moreover, for any $ \de > 0$ there exist $C_1 (\de), \, C_2 ( \de ) >0$ such that 
for all $ \psi  \in \Eo $ we have
\beq
\label{4.1}
\begin{array}{l}
\ds \frac{1 - \de }{2}  \ds \int_{\R^N}  
\left( \ph^2(|\psi |) - 1 \right)^2 \, dx 
- C_1 ( \de ) \| \nabla \psi \|_{L^2(\R^N)}^{2^*} 
\leq 
\ds \int_{\R^N}  V(| \psi |^2 ) \, dx 
\\
\\
\ds \leq 
\frac{1 + \de }{2}  \ds \int_{\R^N}  
\left( \ph^2(|\psi |) - 1 \right)^2 \, dx 
+ C_2 ( \de ) \| \nabla \psi  \|_{L^2(\R^N)}^{2^*}  \
\qquad \mbox{ if } N \geq 3 ,
\end{array}
\eeq
respectively
\beq
\label{4.2}
\begin{array}{l}
\ds \left( \frac{1 - \de }{2} - C_1(\de) \| \nabla \psi \|_{L^2(\R^2)}  ^{2 p_0 + 2}  \right)
 \ds \int_{\R^2}   \left( \ph^2(|\psi |) - 1 \right)^2 \, dx 
\leq 
\ds \int_{\R^2}  V(| \psi |^2 ) \, dx 
\\
\\
\ds \leq 
\left( \frac{1 + \de }{2}  + C_2(\de) \| \nabla \psi \|_{L^2(\R^2)}  ^{2 p_0 + 2}  \right)
\ds \int_{\R^2}   \left( \ph^2(|\psi |) - 1 \right)^2 \, dx 
\qquad \mbox{ if } N =2.
\end{array}
\eeq
These estimates still hold if we replace the condition $ F \in C^0([0, \infty)) $ in (A1) by 
$F \in L_{loc}^1([0, \infty))$ and if we replace $V$ by $|V|$.
\end{Lemma}

{\it Proof. }
Inequality (\ref{4.1}) follows from Lemma 4.1 p. 144 in \cite{M10}. We only prove (\ref{4.2}).

Fix $ \de > 0$. There exists $ \beta = \beta (\de ) \in (0, 1]   $ such that 
\beq
\label{4.3}
\frac{ 1 - \de }{2}  ( s - 1) ^2 \leq V(s) \leq \frac{ 1 + \de }{2}  ( s - 1) ^2 
\qquad \mbox{ for any } s \in ( ( 1 - \beta )^2, (1 + \beta )^2).
\eeq

Let $ \psi \in \Eo$. It follows from (\ref{4.3}) that $V(|\psi |^2) \1 _{ \{ 1 - \beta \leq |\psi | \leq 1 + \beta \} } \in L^1(\R^2)$ and
\beq
\label{4.4}
\begin{array}{l}
\ds  \frac{1 - \de }{2}  \int_{\{ 1 - \beta \leq |\psi | \leq 1 + \beta \} } 
 \left( \ph^2(|\psi |) - 1 \right)^2 \, dx 
\leq  \int_{\{ 1 - \beta \leq |\psi | \leq 1 + \beta \} }  V(|\psi |^2) \, dx 
\\
\\
\ds \leq  \frac{1 + \de }{2}  \int_{\{ 1 - \beta \leq |\psi | \leq 1 + \beta \} } 
\left( \ph^2(|\psi |) - 1 \right)^2 \, dx .
\end{array}
\eeq
Using (A2) we infer that there exists $ C(\de ) > 0 $ such that 
\beq
\label{4.5}
\Big\vert V( s^2) - \frac{ 1 \pm \de }{2}  \left( \ph^2(s) - 1 \right)^2 \Big\vert
\leq C(\de ) \left( |s - 1 | - \frac 12 \beta  \right)^{ 2 p_0 + 2} 
\eeq
for any $ s \geq 0 \mbox{ satisfying } |s - 1 | \geq \beta.$
Let $K = \{ x \in \R^2 \; \big\vert \; \; \big| \, | \psi (x)| - 1 \big| \geq \frac{\beta }{2} \}. $
Let $ \eta $ be as in (\ref{3.19}). Then  $  \left( \ph^2(|\psi |) - 1 \right)^2 \geq \eta (\frac{ \beta }{2} )$ on $K$, hence
\beq
\label{4.6}
\Lo ^2 ( K) \leq \frac{ 1}{\eta (\frac{ \beta }{2} ) } \int_{\R^2}  \left( \ph^2(|\psi |) - 1 \right)^2 \, dx .
\eeq
Let $ \tilde{\psi } = \left( \big| \, | \psi | - 1 \big| - \frac{\beta }{2} \right)_+ \! .$
Then $\tilde{\psi} \in L_{loc}^1(\R^2)$, 
$|\nabla \tilde{\psi}| \leq |\nabla \psi |$ a.e. on $ \R^2$ and using (\ref{ineq2}) we get 
\beq
\label{4.7}
\int_{\R^2} |\tilde{\psi } |^{2 p_0 + 2} \, dx \leq C \| \nabla \tilde{\psi} \| _{L^2(\R^2)} ^{2 p_0 + 2} 
\Lo ^2(K).
\eeq
Using (\ref{4.5}),  (\ref{4.6})  and  (\ref{4.7})  we obtain 
\beq
\label{4.8}
\begin{array}{l}
\ds \int_{\R^2 \setminus \{ 1 - \beta \leq |\psi | \leq 1 + \beta \} } 
\Big\vert V( |\psi |^2) - \frac{ 1 \pm \de }{2} \left( \ph^2(|\psi |) - 1 \right)^2 \Big\vert \, dx 
\\
\\
\ds \leq C(\de) \int_{\R^2} |\tilde{\psi } |^{2 p_0 + 2} \, dx 
\leq C'(\de) \| \nabla \tilde{\psi} \| _{L^2(\R^2)} ^{2 p_0 + 2} 
\ds \int_{\R^2}   \left( \ph^2(|\psi |) - 1 \right)^2 \, dx .
\end{array}
\eeq
From (\ref{4.4}) and (\ref{4.8}) we infer that $V (|\psi |^2) \in L^1(\R^2)$ and (\ref{4.2}) holds.
\hfill
$\Box$

\medskip

The following result is a direct consequence of (\ref{4.2}).

\begin{Corollary}
\label{C4.1}
Assume that $ N =2$ and (A1) and (A2) hold. 
There is $ k_1 > 0 $ such that for any $ \psi \in \Eo $ satisfying 
$  \ii_{\R^2} |\nabla \psi |^2 \, dx \leq k_1 $ we have
$
 \ii_{\R^2} V( |\psi |^2 ) \, dx \geq 0.
$
\end{Corollary}

If $N \geq 3 $  and there exists $ s_0 \geq 0 $ satisfying $V(s_0 )<0$, Corollary \ref{C4.1} is not valid anymore.
Indeed, if $V$ achieves negative values it easy to see that there 
exists $\psi \in \Eo $ such that $ \ii_{\R^N} V( |\psi |^2 ) \, dx < 0.$
Then $ \ii_{\R^N} V( |\psi _{\si, \si }|^2 ) \, dx = \si ^N  \ii_{\R^N} V( |\psi |^2 ) \, dx < 0$
for any $ \si > 0 $ and 
$ \ii_{\R^N} |\nabla \psi _{\si, \si } |^2 \, dx 
=  \si ^{N-2} \ii_{\R^N} |\nabla \psi  |^2 \, dx \lra 0 $ as $ \si \lra 0 $.

\begin{Corollary}
\label{C4.2}
Let $ N \geq 2$. There is an increasing function $ m : \R_+ \lra \R_+$ such that 
$ \ds \lim_{\tau \ra 0 } m ( \tau ) = 0 $ and 
$$
\| \, | \psi | - 1 \|_{L^2(\R^N)} \leq m ( E_{GL}(\psi ) ) \qquad 
\mbox{ for any } \psi \in \Eo.
$$ 
\end{Corollary}

{\it Proof.} Let $ \tilde{F}(s) = \frac{ 1}{\sqrt{s}} - 1$. It is obvious that $ \tilde{F}$
satisfies the assumptions (A1) and (A2) in the introduction (except the continuity at $0$, but this plays no role here).
Let $ \tilde{V}(s) =  \ii_s ^{1} \tilde{F} (\tau ) \, d \tau $, so that 
$ \tilde{V}(s) = ( \sqrt{s} - 1) ^2 $ and $\ii_{\R^N} \tilde{V}(|\psi |^2) \, dx 
= \| \, |\psi | - 1 \|_{L^2(\R^N)}^2$. 
The conclusion follows by using the second inequalities in (\ref{4.1}) and (\ref{4.2})
with $ \tilde{F}$ and $\tilde{V}$ instead of $F$ and $V$. 
\hfill
$\Box $

\medskip

\begin{Lemma}
\label{L4.2}
(i) Let $ \de \in (0, 1)$ and let $ \psi \in \Eo $ be such that 
$ 1 - \de \leq |\psi | \leq 1 + \de $ a.e. on $\R^N$. Then 
$$
| Q(\psi) | \leq \frac{ 1}{\sqrt{2} (1 - \de )} E_{GL}(\psi).
$$

(ii) Assume that $ 0 \leq c < v_s $ and let $ \e \in (0, 1 - \frac{c}{v_s})$. 
There exists a constant $ K _1 = K_1( F, N, c, \e ) > 0$ 
such that for any $ \psi \in \Eo $ satisfying $E_{GL}(\psi) < K_1$ we have
$$
\ds \int_{\R^N} |\nabla \psi  |^2 \, dx +  \int_{\R^N}  V(| \psi |^2 ) \, dx 
- c | Q(\psi) | \geq \e E_{GL}(\psi). 
$$
\end{Lemma}

{\it Proof.} If $N\geq 3$, (i) is precisely Lemma 4.2 p. 145 and (ii) is Lemma 4.3 p. 146
in \cite{M10}. In the case $N=2$ the proof is similar and is left to the reader. 
\hfill
$\Box$

\medskip

For  any $ q \in \R$ we  denote 
$$
E_{min}(q) = \inf  \left\{ \ii_{\R^N} |\nabla \psi |^2\, dx + \ii_{\R^N} \big| V(|\psi |^2) \big| \, dx \; \; \; \Big| \; \psi \in \Eo, Q(\psi ) = q \right\}.
$$
Notice that if $ V \geq 0$, the above definition of $E_{min}$ is the same as the one given in Theorem \ref{T1.1}.
For later purpose we need this more general definition. To simplify the notation, we denote 
$$
\ov{E}(\psi) = \ii_{\R^N} |\nabla \psi |^2\, dx + \ii_{\R^N} \big| V(|\psi |^2) \big| \, dx \qquad \mbox{ for any } \psi \in \Eo.
$$

There are functions $ \psi \in \Eo $ such that $ Q(\psi ) \neq 0$
(see for instance Lemma 4.4 p. 147 in \cite{M10}). 
For any $ \psi \in \Eo$, the function $ \tilde{\psi }(x) = \psi ( - x_1, x')$ 
also belongs to $ \Eo$ and satisfies $\ov{E}(\tilde{\psi}) = \ov{E}({\psi})$,
$Q(\tilde{\psi}) = - Q(\psi)$.
Taking into account (\ref{2.9}), it is clear that for any $q$ the set 
$ \{ \psi \in \Eo \; | \; Q(\psi ) = q \}$ is not empty and 
$E_{min}(- q) = E_{min}(q).$
Thus it suffices to study $E_{min} (q)$ for $q \in [0, \infty)$.

\medskip

If there is $ s_0 $ such that $ V(s_0^2) <0$, then 
$ \inf  \{ E(\psi ) \; | \; \psi \in \Eo, Q(\psi ) = q \} = - \infty $ for all $ q \in \R$. 
(This is one reason why we use $\ov{E}$, not $E$, in the definition of $E_{min}$.)
Indeed, fix $ q \in \R$. From Corollary \ref{C3.4} and (\ref{2.9}) 
we see that there is $\psi _* \in \Eo$ such that 
$Q(\psi _* ) = q $ and $ \psi _* = 1$ outside a ball $B(0, R_*)$. 
It is easy to construct a radial, real-valued function $\psi _0$ such that 
$E(\psi _0) < 0 $ and $ \psi _0 = 1$ outside a ball $B(0, R_0)$ 
(for instance, take $ R_0$ sufficiently large, let $\psi _0 = s_0$ on $B(0, R_0 -1)$, 
$\psi _0 = 1$ on 
$ \R^N \setminus B(0, R_0) $ and $\psi _0$  affine in $|x|$ for $R_0 - 1 \leq |x|\leq R_0$). 
Then $ Q( \psi _0 )  = 0$.
Let $ e_1 = (1, 0 , \dots 0 )$. For $ n \geq 1$, we define $ \psi _n $ by 
$ \psi _n = \psi _* $ 
on $B(0, R_*)$, and $ \psi _n (x)  = \psi _0 (\frac xn - n ^2 (R_0 + R_*) e_1)$
on $ \R^N \setminus B(0, R_*)$.
 Then $Q( \psi _n ) = Q(\psi _* ) + n^{N-1} Q ( \psi_0) = q$ 
and $ E(\psi _n) = E(\psi _*) +  n^{N-2}  \ii_{\R^N} |\nabla \psi _0 |^2 \, dx 
+ n^N \ii _{\R^N} V(|\psi _0 |^2) \, dx \lra - \infty $ as $ n \lra \infty$.

\medskip


\medskip

The next Lemmas establish the properties of $E_{min}$. 

\begin{Lemma}
\label{L4.3}
Assume that $N \geq 2$. 
For any $ q > 0 $ we have $E_{min}(q) \leq v_s q. $
Moreover, there is a sequence $ (\psi _n)_{n \geq 1}$ such that 
$ \psi _n - 1 \in C_c^{\infty} (\R^N)$, $ V(\psi_n) \geq 0$, 
$Q(\psi _n ) = q $, $E(\psi_n ) \lra v_s q$, $E_{GL} (\psi_n ) \lra v_s q$
and $\ds \sup_{x \in \R^N} | \p ^{\al } \psi _n(x )| \lra 0 $ as $ n \lra \infty $ 
for any $ \al \in \N^N$, $ |\al | \geq 1$.
\end{Lemma}

{\it Proof.} 
Fix 
$ \chi \in C_c ^{\infty } (\R^N ),$ $ \chi \neq 0$. 
We will consider three parameters $ \e, \, \la , \, \si > 0 $ such that $ \e \lra 0 $, 
$ \la \lra \infty$,  $ \si \lra \infty$  and $\la \ll \si$.
We put 
$$
\rho_{\e, \la , \si } (x) = 1 + \frac{ \e }{\sqrt{2} \la } \frac{\p \chi }{\p x _1} \left(\frac{ x_1}{\la}, \frac{ x'}{\si } \right), 
\quad
\theta_{\la, \si } (x) = \chi \left( \frac{ x_1}{\la}, \frac{ x'}{\si } \right), 
\quad
\psi_{\e, \la, \si } (x) = \rho_{\e, \la , \si } (x) e^{-  i \e \theta_{\la, \si }  (x) }.
$$
It is clear that $V(\rho_{\e, \la , \si }^2) \geq 0$ if $\frac{\e}{\la}$ is small enough. A straightforward computation gives
$$
\int_{\R^N} \Big\vert \frac{ \p \rho_{\e, \la , \si } }{\p x_1 } \Big\vert ^2 \, dx
= \frac{ \e^2 \si^{N-1}}{2 \la^3 } \int_{\R^N} \Big\vert \frac{ \p ^2 \chi }{\p x_1 ^2}  \Big\vert ^2\, dx, 
$$
$$
\int_{\R^N} \Big\vert \frac{ \p \rho_{\e, \la , \si } }{\p x_j } \Big\vert ^2 \, dx
= \frac{ \e^2 \si^{N-3}}{2 \la  } \int_{\R^N} \Big\vert \frac{ \p ^2 \chi }{\p x_1  \p x_j}  \Big\vert ^2\, dx, 
\qquad j = 2, \dots, N, 
$$
$$
\int_{\R^N} \rho_{\e, \la , \si } ^2 \Big\vert \frac{ \p \theta_{ \la , \si } }{\p x_1 } \Big\vert ^2 \, dx
= \frac{\si^{N-1}}{\la } \int_{\R^N} \left( 1 +\frac{ \e}{ \sqrt{2} \la }  \frac{ \p \chi}{\p {x_1} } \right) ^2 
\Big\vert  \frac{ \p \chi}{\p {x_1} } \Big\vert^2 \, dx
\simeq
 \frac{\si^{N-1}}{\la } \int_{\R^N} \Big\vert  \frac{ \p \chi}{\p {x_1} } \Big\vert^2 \, dx, 
$$
$$
\int_{\R^N} \rho_{\e, \la , \si } ^2 \Big\vert \frac{ \p \theta_{ \la , \si } }{\p x_j } \Big\vert ^2 \, dx
= \si^{N-3} {\la } \int_{\R^N} \left( 1 +\frac{ \e}{ \sqrt{2} \la }  \frac{ \p \chi}{\p {x_1} } \right) ^2 
\Big\vert  \frac{ \p \chi}{\p {x_j} } \Big\vert^2 \, dx
\simeq
 {\si^{N-3}}{\la } \int_{\R^N} \Big\vert  \frac{ \p \chi}{\p {x_j} } \Big\vert^2 \, dx, 
$$
$$
\int_{\R^N} V( \rho_{\e, \la , \si } ^2 ) \, dx \simeq    \frac{\e ^2 \si^{N-1}}{\la } 
 \int_{\R^N} \Big\vert  \frac{ \p \chi}{\p {x_1} } \Big\vert^2 \, dx, 
$$
$$
\int_{\R^N} \left( \ph ^2( \rho_{\e, \la , \si })  - 1 \right)^2 \, dx 
\simeq   \frac{2 \e ^2\si^{N-1}}{ \la } 
 \int_{\R^N} \Big\vert  \frac{ \p \chi}{\p {x_1} } \Big\vert^2 \, dx, 
$$
$$
Q(\psi_{\e, \la, \si } )
 = \e  \int_{\R^N} ( \rho_{\e, \la , \si } ^2 - 1 ) \frac{ \p \theta_{ \la , \si } }{\p x_1} \, dx
\simeq   \frac{\sqrt{2}  \e ^2 \si^{N-1}}{ \la } 
 \int_{\R^N} \Big\vert  \frac{ \p \chi}{\p {x_1} } \Big\vert^2 \, dx. 
$$
Now fix $ q >0$. 
Then choose sequences of positive numbers $(\e _n)_{n \geq 1}$, 
$(\la _n)_{n \geq 1}$,  $(\si _n)_{n \geq 1}$ such that $ \e _n \lra 0 $, 
$ \la _n \lra \infty$, $ \si _n \lra \infty$, $\frac{ \la _n}{\si _n} \lra 0 $ 
and $Q( \psi_{\e_n, \la_n , \si _n } ) = q $ for each $n$. 
Such a choice is possible in view of the last estimate above. In particular, 
this gives 
$ \ds \frac{ \e _n ^2 \si _n ^{N-1}}{ \la _n } 
 \int_{\R^N} \Big\vert  \frac{ \p \chi}{\p {x_1} } \Big\vert^2 \, dx \lra \frac{q }{\sqrt{2}}.$
Let $ \psi _ n = \psi_{\e_n, \la_n , \si _n }$. 
It follows from the above estimates that 
$$
\ov{E}(\psi _n) = E(\psi _n) = \int_{\R^N} |\nabla \rho_{\e_n, \la _n , \si _n } | ^2
+ \e _n^2 \rho_{\e _n, \la _n , \si _n}   ^2 |\nabla \theta_{ \la _n , \si _n} |^2  
+ V(  \rho_{\e _n, \la _n, \si _n} ^2) \, dx 
\lra \sqrt{2} q = v_s q
$$
and similarly $
E_{GL} (\psi _n) \lra  v_s q  $ as $ n \lra \infty$.
The other statements are obvious.
Notice that a similar construction can be found in the proof of Lemma 3.3 p. 604 in   \cite{BGS}. 
\hfill
$\Box $

\begin{Lemma}
\label{L4.4}
Let $N\geq 2$. For each  $ \e > 0$ there is $ q_{\e}>0$ such that 
$$
E_{min}(q) > ( v_s - \e) q \qquad \mbox{ for any } q \in (0, q _{\e}).
$$
\end{Lemma}

{\it Proof.} 
Fix $ \e > 0$. 
It follows from Lemma \ref{L4.2} (ii) that  there is $ K_1( \e ) > 0 $ such 
that for any $ \psi \in \Eo $ satisfying $ E_{GL}(\psi ) < K_1 ( \e ) $ we have 
$$
\ov{E}(\psi )  \geq \left ( v_s - \frac{ \e }{2} \right) | Q(\psi )|.
$$

Using Lemma \ref{L4.1} we infer that there exists $ K_2 ( \e ) > 0 $ such 
that for any $ \psi \in \Eo $ 
satisfying $ \ov{E}(\psi ) < K_2 ( \e ) $ we have $ E_{GL}(\psi ) < K_1 ( \e ) $.

Take $ q _{\e}= \frac{ K_2 (\e)} { v_s + 1 }.$ 
Let $ q \in (0, q_{\e})$.
 There is $ \psi \in \Eo $ such that 
$ Q(\psi ) = q $ and $ \ov{E}(\psi ) < E_{min}(q) + q $. 
Since $E_{min}(q) \leq v_s q $ by Lemma \ref{L4.3},
for any such $ \psi $ we have $\ov{E}(\psi ) < ( v_s + 1) q_{\e} = K_2 ( \e ) $ and we infer that 
$E_{GL}(\psi ) < K_1 ( \e) $, thus 
$\ov{E}(\psi ) \geq  \left ( v_s - \frac{ \e }{2} \right) | Q(\psi )| = 
\left ( v_s - \frac{ \e }{2} \right) q .$
This clearly implies $E_{min}(q) \geq  \left ( v_s - \frac{ \e }{2} \right) q.$ 
\hfill
$\Box $

\begin{Lemma}
\label{L4.5}
Assume that  $N\geq 2$.

(i)  The function $E_{min} $ is subadditive:   for any $ q_1, \, q_2 \geq 0 $ we have
$E_{min}(q_1 + q_2) \leq E_{min}(q_1) + E_{min}(q_2).$

(ii) The function $E_{min} $ is nondecreasing on $[0, \infty)$,  concave, Lipschitz 
continuous and its best Lipschitz constant is $ v_s$. 
Moreover, for $ 0 < q_1 < q_2 $ we have 
$E_{min}(q_1) \leq \left( \frac{ q_1}{q_2} \right)^{\frac{N-2}{N-1}} E_{min}(q_2).$

(iii) For any $ q > 0 $ we have the following alternative: 

$\qquad \bullet $ either $E_{min}(\tau ) = v_s \tau $ for all $ \tau \in [0, q], $

$\qquad \bullet $ or $E_{min}(q) < E_{min}(\tau) + E_{min}(q - \tau ) \quad $ 
for all $ \tau \in (0, q)$.
\end{Lemma}

{\it Proof. } (i) Fix $ \e > 0$. From Corollary \ref{C3.4} and (\ref{2.9}) it 
follows that there exist 
$ \psi _1, \psi _2 \in \Eo$ such that $Q(\psi _i ) = q_i $, 
$\ov{E}(\psi _i) < E_{min}(q_i) + \frac{ \e}{2}$ 
and $ \psi _i = 1 $ outside  a  ball $B(0, R_i)$, $ i = 1,2$. 
 Let $ e \in \R^N$ be a vector of length $1$. Define 
$ \psi (x) = \left\{ \begin{array}{l} \psi _1 (x) \quad \mbox{ if } |x| \leq R_1, \\
\psi _2 (x- 4(R_1 + R_2) e)\quad \mbox{ otherwise. }
\end{array}
\right. $
Then $ \psi \in \Eo$, $Q(\psi ) = Q(\psi _1) + Q(\psi _2 ) = q_1 + q_2 $ and 
$E_{min}(q_1+ q_2) \leq \ov{E}(\psi ) =   \ov{E}(\psi _1) + \ov{E}(\psi _2 ) < 
E_{min}(q_1) + E_{min}(q_2) + \e.$ 
Letting $ \e \lra 0 $ we get $E_{min}(q_1 + q_2) \leq E_{min}(q_1) + E_{min}(q_2).$

\medskip

(ii) From  Lemma \ref{L4.3} 
we obtain 
$ 0 \leq E_{min }(q) \leq v_s q $ for any $ q \geq 0$. 

For $ \psi \in \Eo$ we have $ \psi _{\si, \si } = \psi \left( \frac{\cdot }{\si } \right) 
\in \Eo$, 
\beq
\label{4.9}
Q(\psi _{\si, \si } ) = \si ^{N-1} Q(\psi )\quad \mbox{ and } \quad 
\ov{E}(\psi _{\si, \si } ) = \si ^{N-2} \int_{\R^N} |\nabla \psi |^2 \, dx 
+ \si ^N \int_{\R^N}  \big| V(|\psi |^2) \big| \, dx.
\eeq

Assume that $0 < q_1 < q_2$. Let $ \si _0 =\left( 
\frac{ q_1}{q_2} \right) ^{\frac{1}{N-1}} < 1$. 
For any $ \psi \in \Eo$ satisfying $ Q(\psi ) = q_2$ we have $Q(\psi_{\si _0, \si _0}) = q_1 $ 
and from (\ref{4.9}) we see that $E_{min} (q_1) \leq \ov{E}(\psi _{\si _0, \si _0}) \leq 
\si _0 ^{N-2} \ov{E}(\psi )$. 
Passing to the infimum over all $\psi $ verifying $ Q(\psi ) = q_2$ we find 
$ E_{min}(q_1) \leq \left( \frac{ q_1}{q_2} \right)^{\frac{N-2}{N-1}} E_{min} (q_2).$
In particular,  $E_{min}$ is nondecreasing.
Using (i) and Lemma \ref{L4.3} we get 
$$
0 \leq  E_{min} (q_2) -  E_{min} (q_1) \leq  E_{min} (q_2- q_1) \leq v_s( q_2 - q_1).
$$
Hence $E_{min}$ is Lipschitz continuous and $ v_s$ is a Lipschitz constant for $E_{min}$.
 Lemma \ref{L4.4} implies that $ v_s$ is indeed  the best Lipschitz constant of  $E_{min}$.

Given a function $f$ defined on $ \R^N$ and $ t \in \R$, we denote 
by $S_t^+ f$ and $S_t^-f$, respectively, the functions 
\beq
\label{4.10}
S_t^+ f(x) = \left\{ 
\begin{array}{l}
f(x) \quad \mbox{ if } x_N \geq t, 
\\
f( x_1, \dots, x_{N-1}, 2t - x_N ) \mbox{ if } x_N < t, 
\end{array}
\right.
\eeq
\vspace*{-12pt}
\beq
\label{4.11}
S_t^- f(x) = \left\{ 
\begin{array}{l}
f( x_1, \dots, x_{N-1}, 2t - x_N )  \quad \mbox{ if } x_N \geq t, 
\\
f( x ) \mbox{ if } x_N < t. 
\end{array}
\right.
\eeq
It is easy to see  that for all $ \psi \in \Eo $ and $ t \in \R$ 
we have $S_t^+ \psi , S_t^- \psi  \in \Eo, \quad$ 
$\ov{E}(S_t^+\psi ) + \ov{E}(S_t^-\psi ) = 2 \ov{E}( \psi ) $ and 
$\langle i (S_t^{\pm} \psi )_{x_1}, S_t^{\pm} \psi  \rangle = 
S_t^{\pm} (\langle i \psi _{x_1}, \psi \rangle ).$ 
Moreover, if $ \phi \in \dot{H}^1(\R^N)$ then 
$S_t^{\pm} \phi \in \dot{H}^1(\R^N)$ and 
$ \p _{x_1} (S_t^{\pm} \phi) = S_t^{\pm}  ( \p _{x_1} \phi).$
If $ \psi \in \Eo$, there are $ \phi \in \dot{H}^1(\R^N)$ and $ g \in L^1 ( \R^N) $ 
such that $\langle i \psi _{x_1}, \psi \rangle = \p _{x_1}\phi + g$ (see Lemma \ref{L2.1} 
and the remarks preceding it). Then
$\langle i (S_t^{\pm} \psi )_{x_1}, S_t^{\pm} \psi  \rangle = S_t^{\pm} (\langle i \psi _{x_1}, \psi \rangle )
= \p_{x_1} ( S_t^{\pm}  \phi ) + S_t^{\pm}  g$  and Definition \ref{D2.2} gives
$Q(S_t^{\pm}  \psi ) = \ii _{\R^N} S_t^{\pm} g \, dx$. 
It follows that $Q(S_t^{+}\psi )  + Q(S_t^{-}\psi )  = 2 Q( \psi )$ 
and the mapping $ t \longmapsto Q(S_t^{+}\psi ) =  \ii _{\R^N} S_t^{+} g \, dx
= 2 \ii _{\{ x_n \geq t \} } g \, dx $ is continuous  on $ \R$, tends to $0 $ as $ t \lra  \infty $ 
and to $ 2  \ii _{\R} g \, dx = 2 Q( \psi ) $ as $t \lra - \infty$.

Fix $ 0 < q_1 < q_2 $ and $ \e > 0$.
Let $ \psi \in \Eo $ be such that $ Q( \psi) = \frac{ q_1 + q_2}{2} $ and 
$\ov{E} (\psi ) < E_{min}\left( \frac{q_1 + q_2}{2} \right) + \e$. 
The continuity of $  t \longmapsto Q(S_t^{+}\psi ) $ implies that 
there exists $ t_0 \in \R $ such that $Q( S_{t_0}^{+}\psi  ) = q_1$. 
Then necessarily $Q( S_t^{-}\psi  ) = q_2$ and we infer that 
$\ov{E}(  S_t^{+}\psi  ) \geq E_{min}( q_1) $, $\ov{E}(  S_t^{-}\psi  ) \geq E_{min}(q_2)$, 
and consequently
$$
 E_{min}\left( \frac{q_1 + q_2}{2} \right) + \e > \ov{E}(\psi ) 
= \frac 12 ( \ov{E}(  S_t^{+}\psi  ) + \ov{E}(  S_t^{-}\psi  ) )
\geq \frac 12 (  E_{min}( q_1) +  E_{min}( q_2)).
$$
Passing to the limit as $ \e \lra 0 $ in the above inequality we discover
\beq
\label{4.12}
 E_{min}\left( \frac{q_1 + q_2}{2} \right) \geq  \frac 12 (  E_{min}( q_1) +  E_{min}( q_2)).
\eeq
It is an easy exercise to prove that any continuous function satisfying (\ref{4.12}) 
is concave. 

\medskip

(iii)  Fix $ q > 0 $. 
By the concavity of $E_{min}$ we have $E_{min}(\tau ) \geq \frac{ \tau}{q} E_{min} (q) $ 
for any 
$ \tau \in (0, q) $ and equality may occur if and only if $E_{min} $ is linear on $[0, q]$. 
Therefore for any $ \tau \in (0, q ) $ we have  
$E_{min}(\tau ) + E_{min}(q - \tau ) \geq 
\frac{ \tau}{q} E_{min} (q) +  \frac{ q - \tau}{q} E_{min} (q) = E_{min} (q) $
and equality  occurs if and only if $ E_{min} $ is linear on $[0, q]$, that is 
$E_{min}(\tau ) = a \tau $ for $\tau \in [0, q]$ and some $ a \in \R$. 
Then  Lemma \ref{L4.3} gives $ a \leq v_s$ and Lemma \ref{L4.4} implies $ a \geq v_s - \e $ 
for any $ \e >0$, hence $ a = v_s$.
\hfill
$\Box $

\medskip

The function $ q \longmapsto \frac{E_{min}(q)}{q} $ is nonincreasing (because $E_{min}$ 
is concave), positive and by   Lemma 4.4 in \cite{M10}  there is 
a sequence $ q_n \lra \infty $ such that 
$ {\ds \lim_{n \ra \infty} }\frac{E_{min}(q_n)}{q_n }  = 0 $, hence 
$  {\ds \lim_{q \ra \infty} }\frac{E_{min}(q)}{q}  = 0 $. Let 
$$
q_0 = \inf \{ q > 0 \; | \; E_{min}(q) < v_s q \} ,
$$
so that $ q_0 \in [0, \infty)$, $E_{min}(q) = v_s q $ for $ q \in [0, q_0 ]$ 
and $E_{min}(q) < v_s q $ for any $q > q_0$.

\begin{Lemma} 
\label{L4.6}
Let $ N \geq 2 $. Assume that (A1), (A2) hold.
Then for any $ m , \, M > 0 $ there exist $ C_1(m), \, C_2(M) > 0 $ such that 
for all $ \psi \in \Eo $ satisfying $ m \leq \ov{E}( \psi ) \leq M$ we have
$$
C_1(m) \leq E_{GL}(\psi) \leq C_2(M). 
$$
\end{Lemma}

{\it Proof. } 
If $ N \geq 3$, Lemma \ref{L4.6} follows directly from (\ref{4.1}) with $|V|$ instead of $V$. 
If $ N=2$, the second inequality in (\ref{4.2}) implies that there is $C_1(m) > 0 $ 
such that $E_{GL}(\psi ) \geq C_1(m) $ if $ \ov{E}( \psi ) \geq m$. All we have 
to do is to prove that 
$ \ii _{\R^2} \left( \ph^2( |\psi |) - 1 \right)^2 \, dx  $ remains bounded 
if $ \ov{E}( \psi ) \leq M$. 
This would be trivial if $ \inf \{ V(s^2) \; | \; s \geq 0, \; |s - 1 | \geq \de \} > 0 $ 
for any $ \de > 0$; however, our assumptions do not prevent  $V$ to vanish somewhere 
on $[0, \infty)$ or to tend to zero at infinity. 
Since the proof is the same if $ N =2$ or if $ N \geq 2$, let us consider the general case.

Fix $ \de \in (0, 1 ] $ such that 
$V(s^2) \geq \frac 14  ( s^2 - 1) ^2 $ for $ s \in [1 - \de, 1 + \de]$. 
Consider $ \psi \in \Eo $ such that $ \ov{E}(\psi ) \leq M$. 
Clearly, $   \ii _{\{ | \, |\psi | - 1 | 
\leq \de \} } \left( \ph^2 (|\psi | ) - 1 \right)^2 \, dx 
\leq 4 \ii _{\{ | \, |\psi | - 1 | \leq \de \} } V( | \psi |^2) \, dx \leq 4M$ and 
we have to prove that 
$  \ii _{\{ | \, |\psi | - 1 | > \de \} } \left( \ph^2 (|\psi | ) - 1 \right)^2 \, dx $ 
is bounded. Since $ \ph $ is bounded, it suffices to prove that 
$ \Lo ^N ( \{ \big| \, |\psi | - 1 \big| > \de \} ) $ is bounded.

Let $ w = |\psi | - 1$. Then $|\nabla w | \leq |\nabla \psi | $ a.e., hence 
$ \nabla w \in L^2( \R^N)$, and $ \Lo ^N ( \{ | w| > \al \} ) $ is finite for all $ \al > 0$ 
(because $ \psi \in \Eo $). 
Let $ w_1 ( x) = \phi_1 (|x|) $ and $ w_2 ( x) = \phi_2 (|x|) $ be the symmetric 
decreasing rearrangements of $ w_+$ and $w_-$, respectively. 
Then $ \ph _1$ and $ \ph _2$ are finite,  nonincreasing on 
$ (0, \infty)$ and tend to zero at infinity. 
From Lemma 7.17 p. 174 in \cite{lieb-loss} it follows that
 $ \| \nabla w _1 \|_{L^2(\R^N)} \leq \| \nabla w _+ \|_{L^2(\R^N)} $ and 
 $ \| \nabla w _2 \|_{L^2(\R^N)} \leq \| \nabla w _- \|_{L^2(\R^N)} $. 
In particular, $ w_1, \, w_2 \in H^1 ( \Om _{R_1, R_2})$ for any $0< R_1 < R_2 < \infty$, 
where $ \Om _{R_1, R_2} = B(0, R_2) \setminus \ov{B(0, R_1)} $.
Using Theorem 2 p. 164 in \cite{EG} we infer that $ \phi _1 ,\,  \phi _2 \in H_{loc}^1((0, \infty)) $, 
hence are continuous on $(0, \infty)$.

Let $t_i = \inf \{ t \geq 0 \; | \; \phi_i(t) \leq \de \}$, $i = 1,2$,  so that 
$0 \leq \phi _i (t) \leq \de $ on $[t_i, \infty)$ and, if $ t_i >0$, then $ \phi_i (t_i) = \de$. 
It is clear that
\beq
\label{4.13} 
\begin{array}{l}
\Lo ^N ( \{ \big| \, | \psi | - 1 \big| > \de \}) 
= \Lo ^N ( \{ w_+ > \de \}) + \Lo ^N ( \{ w_ - > \de \}) 
\\
= \Lo ^N ( \{ w_1 > \de \}) + \Lo ^N ( \{ w_ 2 > \de \}) 
= ( t_ 1^N + t_2 ^N) \Lo ^N (B(0,1)).
\end{array}
\eeq

Define $h_1(s) = s^2 + 2  s$, $H_1(s) =  \frac 13 s^3 +  s^2$, 
$h_2(s) = - s^2 + 2  s$, $H_2(s) = - \frac 13 s^3 +  s^2$, so that $ H_1 ' = h_1$ and $ H_2' = h_2$.
If $ t_1 > 0$ we have:
\beq
\label{4.14}
\begin{array}{l}
\ov{E}(\psi ) \geq \ds \int_{\R^N} |\nabla \psi |^2 \, dx 
+ \frac{1}{4} \int_{\R^N}  \left( \ph^2( |\psi |) - 1 \right)^2 \1_{\{ 1 \leq |\psi |\leq 1 + \de \}}\, dx
\\
\\
\geq \ds \int_{\R^N} |\nabla w_+ |^2 \, dx 
+ \ds  \frac 14   \int_{\{ w_+ \leq \de \}} \left( (w_+ + 1 )^2 - 1 \right)^2 \, dx
\\
\\
\geq \ds \int_{\R^N} |\nabla w_1 |^2 \, dx 
+ \ds  \frac 14  \int_{\{ w_1 \leq \de \}} \left( (w_1 + 1 )^2 - 1 \right)^2 \, dx
\\
\\
\geq \ds \int_{\R^N \setminus B(0, t_1)} |\nabla w_1|^2 +  \frac 14  h_1^2 ( w_1) \, dx 
\\
\\
= | S^{N-1} | \ds \int_{t_1}^{\infty} 
\left( |\phi _1 '(s) | ^2 +  \frac 14 h_1^2 ( \phi _1 (s)) \right) s^{N-1} \, ds
\\
\\
\geq    t_1^{N-1}  |S^{N-1}| \ds \int_{t_1}^{\infty}  
 |\phi _1 '(s) | ^2 +  \frac 14 h_1^2 ( \phi _1 (s))  \, ds
 \\
 \\
\geq   t_1^{N-1} | S^{N-1} |\ds \int_{t_1}^{\infty}  
-  h_1( \phi_1 (s)) \phi _1'(s) \, ds 
\\
\\
=  t_1^{N-1} |S^{N-1} | \left[ -  H_1 ( \phi _1(s)) \right]_{s = t _0 }^{\infty} 
=   |S^{N-1} | H_1 ( \de )  t_1^{N-1} , 
\end{array}
\eeq
where $ |S^{N-1}| $ is the surface measure of the unit sphere in $ \R^N$. 
From (\ref{4.14}) we get $ t_1 ^{N-1} \leq C E(\psi)$, where $C$ depends only on $N$ and $V$.
It is clear that a similar estimate holds for $t_2$. Then using (\ref{4.13}) we obtain 
$$
\Lo ^N ( \{ \big| \, | \psi | - 1 \big| > \de \}) \leq C \left( \ov{E}(\psi) \right)^{\frac{N}{N-1}}, 
$$
where $C$ depends only on $N$ and $V$, and the proof of Lemma \ref{L4.6} is complete.
\hfill
$\Box $

We can now state the main result of this section, showing precompactness 
of minimizing sequences for $E_{min}(q)$ as soon as $q > q_0$. 

\begin{Theorem}
\label{T4.7}
Assume that $ q > q_0$, that is $E_{min}(q) < v_s q$. 
Let $ (\psi _n )_{n \geq 1 } $ be a sequence in $ \Eo $ satisfying
$$ 
Q( \psi _ n ) \lra q \quad \mbox{  and } \quad  \ov{E}( \psi _n ) \lra E_{min}(q).
$$ 
There exist a subsequence $(\psi_{n_k})_{k \geq 1}$, a sequence 
of points $(x_k)_{k \geq 1}  \subset \R^N$, 
and $ \psi \in \Eo $ such that $Q( \psi )\! = \!q$, $\ov{E}( \psi )\! =\! E_{min}(q) $, 
$\psi_{n _k} ( \cdot + x_k ) \lra \psi $ a.e. on $ \R^N$ and 
$ d_0( \psi_{n _k} ( \cdot + x_k ), \psi) \lra 0 $, that~is  
$$
\| \nabla \psi_{n _k} ( \cdot + x_k ) - \nabla \psi \|_{L^2(\R^N)} \lra 0 , 
\qquad 
\| \,  | \psi_{n _k} |( \cdot + x_k ) - | \psi  | \,  \|_{L^2(\R^N)} \lra 0 
\qquad \mbox{ as } k \lra \infty.
$$

\end{Theorem}

{\it Proof. }  
Since $\ov{E}(\psi _n) \lra E_{min}(q) > 0$, 
it follows from Lemma \ref{L4.6} that there are two positive constants $M_1, \, M_2$ 
such that $M_1 \leq E_{GL}(\psi _n ) \leq M_2$ for all sufficiently large $n$. 
Passing to a subsequence if necessary, we may assume that $E_{GL}(\psi _n) \lra \al _0> 0$.

We will use the concentration-compactness principle \cite{lions}.
We denote by $ \Lambda_n(t)$ the concentration function associated to 
$E_{GL}(\psi _n)$, that is 
\beq
\label{4.15}
\Lambda_n(t ) = \sup_{y \in \R^N} \int_{B(y, t)}
|\nabla \psi_n |^2 + \frac 12 \left( \ph ^2( |\psi _n |) - 1 \right)^2 \, dx. 
\eeq
Proceeding as in \cite{lions}, it is straightforward to prove that there 
exists a subsequence of $((\psi_n, \Lambda_n))_{n \geq 1}$, still denoted 
$((\psi_n, \Lambda_n))_{n \geq 1}$, there exists a nondecreasing 
function  $ \Lambda : [0, \infty ) \lra \R$ and there is $\al \in [0, \al _0] $ 
such that 
\beq
\label{4.16}
\Lambda_n (t) \lra \Lambda(t)\; \mbox{  a.e on }  [0, \infty ) \; \mbox{ as } n \lra \infty
\qquad \mbox{ and } \qquad
\Lambda(t) \lra \al \mbox{ as } t \lra \infty.
\eeq
As in the proof of Theorem 5.3 in \cite{M10}, we see that there is a nondecreasing sequence $t_n \lra \infty$ 
such that
\beq
\label{annulus}
 \lim_{n \ra \infty }\Lambda_n (t_n) =  \lim_{n \ra \infty }\Lambda_n \left(\frac{t_n}{2} \right) = \al.
\eeq

Our aim is to prove that $\al = \al _0$. 
The next lemma implies that $ \al > 0$. 

\begin{Lemma}
\label{L4.8}
Assume that $N \geq 2$ and assumptions (A1) and (A2) in the Introduction hold. 
Let $(\psi_n) _{n \geq 1} \subset \Eo $ be a sequence satisfying:

\medskip

(a) $   E_{GL}(\psi _n) \leq M $ for some positive constant $M$. 

\medskip

(b) $ \ds \liminf_{n \ra \infty} Q(\psi _n ) \geq q \in \R \cup \{ \infty \} $ 
as $ n \lra \infty.$

\medskip

c) $\ds \limsup_{n \ra \infty } E(\psi_n) < v_s q$. 
\\
Then there exists $ k > 0$ such that 
$ \ds \sup_{ y \in \R^N } E_{GL}^{B(y,1)} ( \psi _n ) \geq k 
$
for all sufficiently large $n$. 
\end{Lemma}

{\it Proof. } We argue by contradiction and we  suppose that the conclusion is false. 
Then there  is a subsequence (still denoted $(\psi_n) _{n \geq 1} $) such that 
\beq
\label{4.17}
\lim_ {n \ra \infty } \sup_{ y \in \R^N } \int_{B(y, 1)} 
|\nabla \psi_n |^2 + \frac 12 \left( \ph ^2( |\psi_n |) - 1 \right)^2 \, dx  =0.
\eeq

The first step is to prove that 
\beq
\label{4.18}
\lim_ {n \ra \infty } \int_{\R^N} 
\Big\vert V(|\psi _n |^2) - \frac 12 \left( \ph ^2(  |\psi _n |) - 1 \right)^2 \Big\vert \, dx = 0.
\eeq
If $ N \geq 3 $ this is done exactly as in the proof of Lemma 5.4 p. 156 in \cite{M10}.
We consider here only the case $ N =2$.

Fix $ \e > 0 $. By (A1) there is $ \de (\e) > 0 $ such that 
$$
\Big\vert V(s^2) - \frac 12 \left( \ph ^2( s) - 1 \right)^2 \Big\vert \leq \frac{\e }{2}  \left( \ph ^2( s) - 1 \right)^2
\qquad \mbox{ for any } s \in [ 1 - \de ( \e),  1 + \de ( \e)], 
$$
hence
\beq
\label{4.19}
\begin{array}{l}
\ds \int_{\{ 1 - \de ( \e) \leq |\psi | \leq   1 + \de ( \e) \} }
\Big\vert V(|\psi _n |^2) - \frac 12 \left( \ph ^2(  |\psi _n |) - 1 \right)^2 \Big\vert \, dx 
\\
\\
\ds \leq \frac{\e}{2}  \int_{\{1 - \de ( \e) \leq |\psi | \leq   1 + \de ( \e) \} }
\left( \ph ^2(  |\psi _n |) - 1 \right)^2 \, dx \leq \e M.
\end{array}
\eeq

Using (A2) we infer that there is $ C(\e) > 0 $ such that 
\beq
\label{4.20}
\Big\vert V(s^2) - \frac 12 \left( \ph ^2( s) - 1 \right)^2 \Big\vert \leq C(\e) ( |s| - 1 ) ^{2 p_0 + 2 } 
\quad \mbox{ for any } s 
\mbox{ satisfying } | s - 1 | \geq \de (\e).
\eeq
Let $ w_n = \big| \, |\psi _n | - 1 \big|$. 
Then $ w_n \in L_{loc}^1(\R^N)$ and  $|\nabla w_n | \leq | \nabla \psi _n|$ a.e., hence
$ \| \nabla w_n \|_{L^2( \R^2)} \leq \| \nabla \psi _n \|_{L^2( \R^2)}  \leq \sqrt{M}$. 
Using (\ref{ineq2}) for $ \left( w_n - \frac{ \de (\e )}{2} \right)_+ $ we obtain 
\beq
\label{4.21}
\begin{array}{l}
\ds \int_{\{ w_n > \de ( \e ) \} } w_n ^{2 p_0 + 2 } \, dx 
\leq 2^{ 2 p_0 + 2} \int_{\{ w_n > \de ( \e ) \} } \left( w_n - \frac{ \de (\e )}{2} \right)_+  ^{ 2 p_0 + 2} \, dx 
\\
\\
\leq C  \| \nabla w_n \|_{L^2( \R^2)} ^{ 2 p_0 + 2 } 
\Lo ^2 ( \{ w_n > \frac{\de ( \e )}{2}  \}  ) 
\leq C M^{ p_0 + 1} \Lo ^2 ( \{ w_n > \frac{\de ( \e )}{2} \}  ) .
\end{array}
\eeq
We claim that for any $ \de > 0 $ we have
\beq
\label{4.22}
\lim_{n \ra \infty }  \Lo ^2 ( \{ w_n > \de  \}  )   = 0.
\eeq
The proof of (\ref{4.22}) relies on Lieb's Lemma 
(see Lemma 6 p. 447 in \cite{lieb} or Lemma 2.2 p. 101 in \cite{brezis-lieb})
and is the same as the proof of (5.20) p. 157 in \cite{M10}, so we omit it. 

From (\ref{4.20}), (\ref{4.21}) and (\ref{4.22}) we get 
\beq
\label{4.23}
\int_{ \{ |\, |\psi | - 1 | > \de ( \e ) \} } 
\Big\vert V(s^2) - \frac 12 \left( \ph ^2( s) - 1 \right)^2 \Big\vert \, dx 
\leq
C(\e ) \int_{\{ w_n > \de ( \e ) \} } w_n ^{2 p_0 + 2 } \, dx \lra 0 
\eeq
as $ n \lra \infty.$ Then (\ref{4.18}) follows from (\ref{4.19}) and (\ref{4.23}).

From (\ref{4.17}) and Lemma \ref{vanishing} we infer that there exists a sequence $ h_n \lra 0 $ 
and for each $n $ there is a minimizer $ \zeta _n $ of $G_{h_n, \R^N} ^{\psi _n} $ in 
$H_{\psi _n }^1( \R^N) $ such that 
\beq
\label{4.24}
\de _n := \| \, |\zeta _n | - 1 \|_{L^{\infty} (\R^N)} \lra 0 \qquad \mbox{ as } n \lra \infty.
\eeq
Then Lemma \ref{L4.2} (i) implies 
\beq
\label{4.25}
E_{GL } ( \zeta_n ) \geq \sqrt{2} ( 1 - \de _n ) | Q( \zeta _n) |.
\eeq
From (\ref{3.4}) we obtain $ \ds \lim_{n \ra \infty} | Q( \zeta _n ) - Q( \psi _n) | = 0$, 
hence $\ds \liminf_{n \ra \infty}  Q( \zeta _n ) =  \liminf_{n \ra \infty}  Q( \psi _n)  \geq q.$
Using (\ref{4.18}), the fact that $E_{GL} ( \zeta _n) \leq E_{GL} ( \psi _n) $ and (\ref{4.25}) we get 
$$
\begin{array}{l}
E(\psi _n) = E_{GL} (\psi _n) + \ds  \int_{\R^N} V(|\psi_n |^2) - \frac 12 ( \ph ^2(|\psi _n|) - 1)^2\, dx 
\\
\geq 
E_{GL} (\zeta _n) + \ds  \int_{\R^N} V(|\psi_n |^2) - \frac 12 ( \ph ^2(|\psi _n|) - 1)^2\, dx 
\\
\geq 
\sqrt{2} ( 1 - \de _n ) |Q(\zeta _n )| + \ds  \int_{\R^N} V(|\psi_n |^2) - \frac 12 ( \ph ^2(|\psi _n|) - 1)^2\, dx .
\end{array}
$$
Passing to the limit as $ n \lra \infty $ in the above inequality we get 
$$
\liminf_{n \ra \infty} E(\psi _n) \geq \sqrt{2} q = v_s q, 
$$
which contradicts assumption c) in Lemma \ref{L4.8}. This ends the 
proof of Lemma \ref{L4.8}.
\hfill
$\Box $

\medskip

Next we prove that $ \al \not\in (0, \al _0)$. 
We argue again by contradiction and we assume that 
$ 0 < \al < \al _0$. 
Let $ t_n$ be as in (\ref{annulus}) and let $ R_n= \frac{t_n}{2}$. 
For each $n \geq 1$, fix $ y_n \in \R^N$ such that 
$E_{GL}^{B(y_n, R_n)} ( \psi_n) \geq \Lambda_n ( R_n) - \frac 1n$. 
Using (\ref{annulus}),  we have 
\beq
\label{4.27}
\e _n : = E_{GL}^{B(y_n, 2R_n) \setminus B(y_n, R_n)} ( \psi _n) 
\leq \Lambda_n ( 2R_n) - \left( \Lambda_n( R_n) - \frac 1n \right) \lra 0 
\mbox{ as } n \lra \infty . 
\eeq
After a translation, we may assume that $ y_n = 0$. 
Using Lemma \ref{splitting} with $ A = 2$, $ R = R_n$, $ \e = \e_n$, 
we infer that for all $n$  sufficiently large there exist 
two functions $ \psi_{n, 1}$, $ \psi_{n, 2}$ having the properties (i)-(vi) 
in Lemma \ref{splitting}. 
In particular, we have 
$E_{GL}(\psi _{n, 1} ) \geq E_{GL}^{B(0, R_n)}  ( \psi _n ) \geq \Lambda( R_n) - \frac 1n$, 
$E_{GL}(\psi _{n, 2} ) \geq E_{GL}^{\R^N \setminus B(0,2 R_n)} ( \psi _n ) \geq E_{GL}(\psi _{n} ) -  \Lambda(2 R_n) $
and
$ | E_{GL}(\psi _{n} ) - E_{GL}(\psi _{n, 1} ) - E_{GL}(\psi _{n, 2} )| \leq C \e _n \lra 0 $ as $ n \lra \infty.$
Taking into account (\ref{annulus}), we conclude that necessarily
\beq
\label{4.28}
E_{GL}(\psi _{n, 1} ) \lra \al 
\qquad 
\mbox{ and } 
\qquad
E_{GL}(\psi _{n, 2} ) \lra \al _0 - \al 
\qquad
\mbox{ as } n \lra \infty.
\eeq
From Lemma \ref{splitting} (iii)-(vi) we get 
\beq
\label{4.29}
| \ov{E}(\psi _n) -  \ov{E}(\psi _{n, 1}) - \ov{E}(\psi _{n, 2}) | \lra 0 \qquad \mbox{ and } 
\eeq
\beq
\label{4.30}
| Q(\psi _n) -  Q(\psi _{n, 1}) -  Q(\psi _{n, 2}) | \lra 0 \qquad \mbox{ as } n \lra \infty. 
\eeq
In particular,  $\ov{E}(\psi _{n, i}) $ is bounded, $ i =1,2$. 
Passing to a subsequence if necessary, we may assume that $ \ov{E}(\psi _{n, i}) \lra m_i \geq 0 $ 
as $ n \lra \infty$.
Since $\ds \lim_{n \ra \infty } E_{GL}(\psi _{n,i} ) > 0$, it follows from Lemma \ref{L4.1} that
 $ m_i > 0$, $i=1,2$.
Using (\ref{4.29}) we see that $ m_1 + m_2 = E_{min}(q)$, hence $ m_1, m_2 \in (0, E_{min}(q)).$

Assume that $\ds \liminf_{n \ra \infty } Q( \psi_{n,1}) \leq 0.$
Then (\ref{4.30}) implies  $\ds \limsup_{n \ra \infty } Q( \psi_{n,2}) \geq q.$
It is obvious that 
$$
\ov{E}(\psi_{n, 2} ) \geq E_{min}(Q(\psi _{n,2})).
$$
Passing to $\ds \limsup$  in the above inequality and using the continuity 
and the monotonicity of $E_{min}$ we get 
$ m_2 \geq E_{min}(q), $ a contradiction. 
Thus necessarily $\ds \liminf_{n \ra \infty } Q( \psi_{n,1})  > 0$ and similarly 
 $\ds \liminf_{n \ra \infty } Q( \psi_{n,2})  > 0$.
From (\ref{4.30}) we get 
$\ds \limsup_{n \ra \infty } Q( \psi_{n,i})< q$, $ i = 1,2$. 
Passing again to a subsequence, we may assume that 
$ Q( \psi_{n,i}) \lra q_i$ as $ n \lra \infty$, $ i = 1,2$, where $ q_1, q_2 \in (0, q)$.
Using  (\ref{4.30}) we infer that $ q_1 + q_2 = q$.
Since $\ov{E}(\psi_{n,i} ) \geq E_{min}(Q(\psi _{n,i}))$, 
passing to the limit we get $ m_i \geq E_{min}(q_i)$, $ i =1,2$ and consequently
$$
E_{min}(q) = m_1 + m_2 \geq  E_{min}(q_1) +  E_{min}(q_2).
$$
Since $ E_{min}(q) < v_s q$,  the above inequality  is in  contradiction 
with Lemma \ref{L4.5} (iii). Thus we cannot have $ \al \in (0, \al _0)$.

\medskip

So far we have proved that $ \al = \al _0$.
Then it is standard to prove that there  is a sequence $(x_n)_{n \geq 1} \subset \R^N$ such that 
for any $ \e > 0 $ there is  $ R_{\e} > 0 $  satisfying 
$E_{GL} ^{ \R^N\setminus B(x_n, R_{\e})} (\psi _n)  < \e$ for all sufficiently large $ n$. 
Denoting $ \tilde{\psi}{_n} = \psi _n ( \cdot + x_n)$, we see that 
for any $ \e > 0 $ there exist   $ R_{\e} > 0 $ and $ n_{\e } \in \N$ such that 
\beq
\label{4.31} 
E_{GL} ^{ \R^N\setminus B(0, R_{\e})}(\tilde{\psi}{_n})  < \e \qquad \mbox{ for all } n \geq n_{\e}.
\eeq
Obviously,  $ (\nabla \tilde{\psi}_{n})_{n \geq 1} $ is bounded in $L^2(\R^N)$ and it is easy to see that
$( \tilde{\psi}{_n})_{n \geq 1} $ is bounded in $L^2 (B(0,R))$ for any $R >0$ (use (\ref{e1}) and 
(\ref{e2}) if $ N=2$, respectively (\ref{e1}) and the Sobolev embedding if $ N \geq 3$).
By a standard argument,  there exist a function $ \psi \in H_{loc}^1(\R^N)$
such that $ \nabla \psi \in L^2(\R^N)$ and a subsequence $ (\tilde{\psi}_{n _k})_{k \geq 1} $
 satisfying 
 \beq
 \label{4.32}
 \begin{array}{l}
 \nabla \tilde{\psi}_{n_k } \rightharpoonup \nabla \psi \qquad \mbox{ weakly in } L^2( \R^N), 
 \\
 \tilde{\psi}_{n_k } \rightharpoonup  \psi \qquad \mbox{ weakly in } H^1(B(0,R)) \mbox{ for all } R > 0, 
 \\
 \tilde{\psi}_{n_k } \lra  \psi \; \;   \mbox{ strongly in } L^p(B(0,R)) \mbox{ for  } R > 0 \mbox{ and }
 p \in [1, 2^* ) \; ( p \in [1, \infty) \mbox{ if } N =2), 
 \\
  \tilde{\psi}_{n_k } \lra  \psi \; \;   \mbox{  almost everywhere on } \R^N.
 \end{array}
 \eeq

By weak convergence we have 
\beq
\label{4.33}
\int_{\R^N} |\nabla \psi |^2 \, dx 
\leq \liminf_{k \ra \infty } 
\int_{\R^N} |\nabla \tilde{\psi }_{n_k} |^2 \, dx .
\eeq
The a.e. convergence and Fatou's Lemma imply
\beq
\label{4.34}
\int_{\R^N}
\left( \ph^2 (|\psi |) - 1 \right)^2 \, dx 
\leq \liminf_{k \ra \infty } 
\int_{\R^N} \left( \ph^2 (|\tilde{\psi }_{n_k}|) - 1 \right)^2 \, dx 
\quad \mbox{ and } 
\eeq
\beq
\label{4.35}
\quad
\int_{\R^N} \big| V(|\psi |^2) \big| \, dx \leq \liminf_{k \ra \infty } \int_{\R^N} \big| V(|\tilde{\psi} _{n_k}|^2) \big| \, dx.
\eeq
From (\ref{4.33}), (\ref{4.34}) and (\ref{4.35}) we obtain 
\beq
\label{4.36}
E_{GL} ( \psi ) \leq \liminf_{k \ra \infty }  E_{GL} ( \tilde{\psi } _{n_k}) = \al _0 
\quad
\mbox{ and } 
\quad
\ov{E} ( \psi ) \leq \liminf_{k \ra \infty }  \ov{E} ( \tilde{\psi } _{n_k})  = E_{min}(q).
\eeq
Similarly, for any $ \e > 0 $ we get 
\beq
\label{4.37}
E_{GL}^{\R^N \setminus B(0, R_{\e})} ( \psi ) 
\leq \liminf_{k \ra \infty }  E_{GL}^{\R^N \setminus B(0, R_{\e})}  ( \tilde{\psi } _{n_k})
\leq \limsup_{k \ra \infty }  E_{GL}^{\R^N \setminus B(0, R_{\e})}  ( \tilde{\psi } _{n_k}) \leq \e.
\eeq

The following holds. 

\begin{Lemma}
\label{L4.10}
Assume that $N \geq 2$ and assumptions (A1) and (A2)  are  verified. 
Let $ ( \g _n ) _{n \geq 1} \subset \Eo $ be a sequence  satisfying:

\smallskip

(a) $ (E_{GL}( \g _n)) _{n \geq 1}$ is bounded and for any $ \e > 0 $ 
there are $ R_{\e } > 0$ and $n_\e \in \N$ such that 
$E_{GL}^{\R^N \setminus B(0, R_{\e })} ( \g _n ) < \e$ for $n \geq n_\e$. 

\smallskip

(b) There exists $ \g \in \Eo $ such that 
$ \g _n \lra \g $ strongly in $L^2( B(0, R))$ for any $ R > 0$, 
and $ \g_n \lra \g $ a.e. on $ \R^N$ as $ n \lra \infty$. 

\smallskip

Then $ \| \,  |\g_n | -  | \g | \, \| _{L^2( \R^N)} \lra 0 $ and
$ \| V( |\g _n |^2) - V( |\g |^2) \| _{L^1( \R^N)} \lra 0 $ as $ n \lra \infty$. 

\end{Lemma}

{\it Proof of Lemma \ref{L4.10}.}
Fix $ \e > 0 $. Let $ R_{\e}$ and $n_\e \in \N$ be as in assumption (a). 
Then 
\beq
\label{4.38}
\| \ph ( |\g _n |) - 1 \|_{L^2(\R^N \setminus B(0, R_{\e }))} ^2
\leq
 \int_{\R^N \setminus B(0, R_{\e })}  \left( \ph ^2 ( |\g _n |)  - 1 \right)^2 dx 
\leq 2 \e
\eeq
for $ n \geq n_\e $. It is clear that a similar estimate holds for $ \g$. Let 
$$
\begin{array}{ll}
\tilde{\g }_n = |\g _n | - \ph (|\g _n|), & \tilde{\g } = |\g  | - \ph (|\g |), 
\\
\\
A_n = \{ x \in \R^N \; \big| \;  \; | \g _n (x) | \geq 2  \}, & A = \{ x \in \R^N \; \big| \; \; | \g (x) | \geq 2  \}, 
\\
\\
A_n ^{\e }= \{ x \in \R^N \setminus B(0, R_{\e}) \; \big| \; | \g _n (x) | \geq 2  \}, 
& A ^{\e} = \{ x \in \R^N \setminus B(0, R_{\e}) \; \big| \; | \g (x) | \geq 2  \}.
\end{array}
$$
We have 
$$
9  \Lo ^N ( A_n^{\e}) 
\leq  \int_{\R^N \setminus B(0, R_{\e })}  \left( \ph ^2 ( |\g _n |)  - 1 \right)^2 dx 
\leq  2 E_{GL} ^{\R^N \setminus B(0, R_{\e })}  ( \g _n )  
\leq 2 \e
$$
and similarly $ 9 \Lo ^N ( A^{\e}) \leq 2 \e.$
In the same way $  \Lo ^N ( A_n ) \leq \frac{2}{9 } E_{GL}(\g_n)$ and 
$ \Lo ^N ( A ) \leq \frac{2}{9 } E_{GL}(\g)$.

Since $ 0 \leq \ph' \leq 1 $, it is easy to see that $|\nabla \tilde{\g} _n | \leq | \nabla \g _n|$ a.e. and 
$|\nabla \tilde{\g}  | \leq | \nabla \g |$ a.e., hence 
$(|\nabla  \tilde{\g} _n | )_{n \geq 1}$ and $ \nabla \tilde{\g}$ are bounded in $L^2( \R^N)$. 
If $ N \geq 3$, the Sobolev embedding implies that $(\tilde{\g} _n)_{n \geq 1}$ 
is bounded in $ L^{2^*}(\R^N)$. 
Then using the fact that $ \tilde{\g}_n = 0 $ on $ \R^N \setminus A_n$ and H\"older's inequality we infer that 
 $ \tilde{\g}_n$ is bounded in $ L^p(\R^N)$ for $1 \leq p \leq 2^*$. 
If $ N = 2$, by (\ref{ineq2}) we get 
$$
\| \tilde{\g}_n \|_{L^p(\R^2)} ^p \leq C_p^p \| \nabla \tilde{\g }_n \|_{L^2(\R^2)} ^p \Lo ^ 2( A_n), 
$$
hence $(\tilde{\g} _n)_{n \geq 1}$  is bounded in $L^p(\R^N)$ for any $ 2 \leq p < \infty . $
Let $ p = 2^* $ if $ N \geq 3$ and let $p > 2  p_0 + 2 $ if $ N = 2$. 
Using H\"older's inequality ($p> 2  p_0 + 2 > 2$) we have
\beq
\label{4.39}
\| \tilde{\g }_n \|_{L^2(\R^N \setminus B(0, R_{\e }))} ^2 
=  \int_{A_n^{\e} } |\tilde{\g}_n |^2 \, dx 
\leq 
\| \tilde{\g }_n \|_{L^p (\R^N)} ^2 \Lo ^N ( A_n ^{\e } )^{1 - \frac 2p} \leq C _1\e ^{ 1 - \frac 2p}, 
\eeq
where $C_1 $ does not depend on $n$. 
It is clear that a similar estimate holds for $ \tilde{\g}$. 

In the same way, using (A2) and H\"older's inequality ($p> 2  p_0 + 2 $) we get
\beq
\label{4.39bis}
\begin{array}{l}
\ds \int_{(\R^N \setminus B(0, R_{\e })) \cap \{ |\g _n | \geq 4  \} } | V ( |{\g } _n | ^2 ) | \, dx 
\leq C' \int_{(\R^N \setminus B(0, R_{\e })) \cap \{ |\g _n | \geq 4  \} } |\g _n |^{ 2 p_0 + 2} \, dx 
\\
\\
\leq C''
\ds \int_{A_n^{\e} } |\tilde{\g}_n |^{2 p_0 + 2}  \, dx 
\leq C '' \| \tilde{\g }_n \|_{L^p (\R^N)} ^{2 p_ 0 + 2}  \Lo ^N ( A_n ^{\e } )^{1 - \frac{ 2p_0 + 2}{p}} 
\leq C_2 \e ^{ 1 - \frac{2p_0 + 2}{p} }, 
\end{array}
\eeq
and (A1) implies 
\beq
\label{4.39ter}
\ds \int_{(\R^N \setminus B(0, R_{\e })) \cap \{ |\g _n | \leq 4  \} } | V ( |{\g } _n | ^2 ) | \, dx 
\leq C''' \int_{\R^N \setminus B(0, R_{\e }) } 
 \left( \ph ^2 ( |\g _n |)  - 1 \right)^2 dx 
\leq C_3 \e, 
\eeq
where the constants $C_2, C_3$ do not depend on $n$. 
The same estimates are obviously valid for $ \g$. 

From (\ref{4.38}) and (\ref{4.39}) we get 
\beq
\label{4.39qua}
\begin{array}{l}
\| \, | \g _n | - |\g | \, \|_{L^2( \R^N \setminus B(0, R_{\e}\! ))}
\! \leq \!
\| \ph (|\g _n |) \! - \! 1 \| _{L^2( \R^N \setminus B(0, R_{\e}))}
\! + \! \| \ph (|\g  |) \! - \! 1  \| _{L^2( \R^N \setminus B(0, R_{\e}\! ))}  
\\
\quad + \| \tilde{\g} _n \| _{L^2( \R^N \setminus B(0, R_{\e}))}   
+ \| \tilde{\g}  \| _{L^2( \R^N \setminus B(0, \R_{\e}))}   
\leq 2 \sqrt{ 2 } \e + 2 C_1 \e ^{1 - \frac 2p }.
\end{array}
\eeq
Using (\ref{4.39bis}) and (\ref{4.39ter}) we obtain
\beq
\label{4.39cinq}
\ds \int_{\R^N \setminus B(0, R_{\e})} | V(|\g _n |^2) | \, dx 
\leq  C_2 \e ^{ 1 - \frac{2p_0 + 2}{p} } + C_3 \e.
\eeq
It is obvious that $ \g $ also satisfies (\ref{4.39cinq}).

 Since $ |\g _n | = \ph (|\g _n |) + \tilde{\g }_n $ is bounded in $L^p( B(0,R))$ for any $ p \in [2, 2^*]$ if $ N \geq 3$, 
 respectively $ p \in [2, \infty )$ if $ N = 2$, and 
 $  \g_n \lra \g$ in $L^2( B(0, R))$ by assumption (b), using interpolation we infer that  
 $  \g_n \lra \g$ in $L^p( B(0, R))$ for any $p \in [1, 2^*)$ (with $ 2^* = \infty$ if $N=2$).
This implies that $V(|\g _n |^2) \lra V(|\g |^2) $ in $L^1( B(0, R))$ 
(see, for instance, Theorem A2 p. 133 in \cite{willem}).
Thus we have $ \| \, |\g _n | - |\g | \, \|_{L^2( B(0, R_{\e}))} \leq \e $ 
and $ \| V( |\g _n |^2)  - V(|\g |^2) \|_{L^1( B(0, R_{\e}))} \leq \e $  for all sufficiently large $n$. 
Together with  inequalities (\ref{4.39qua}) and (\ref{4.39cinq}), this implies
$\| \, |\g _n | - |\g | \, \|_{L^2(\R^N)} \leq 2 \sqrt{2} \e + 2 C_1 \e ^{1 - \frac 2p } + \e$
and
$ \| V( |\g _n |^2)  - V(|\g |^2) \|_{L^1( \R^N)}
\leq 2 C_2 \e ^{ 1 - \frac{2p_0 + 2}{p} } + ( 2C_3 +1) \e $
for all  sufficiently large $n$. 
Since $ \e $ is arbitrary, Lemma \ref{L4.10} follows.
\hfill
$\Box $

\medskip

We come back to the proof of Theorem \ref{T4.7}.
From (\ref{4.31}), ({\ref{4.32}) and Lemma \ref{L4.10} we obtain 
$\| \, |\tilde{\psi}_{n _k} | - |\psi | \, \| _{L^2( \R^N) } \lra 0 $ as $ k \lra \infty.$ 
Clearly, this implies 
$\| \ph^2(| \tilde{\psi}_{n_k } |) - \ph^2(| \psi  |) \|_{L^2(\R^N)} \lra 0$.

\medskip

We will use the following result: 

\begin{Lemma}
\label{L4.11}
Let $N \geq 2$ and assume that $ ( \g _n ) _{n \geq 1} \subset \Eo $ is 
a sequence satisfying:

\smallskip

(a) $ (E_{GL}( \g _n) ) _{n \geq 1} $ is bounded and for any $ \e > 0 $ 
there are $ R_{\e } > 0$ and $n_\e \in \N $ such that 
$E_{GL}^{\R^N \setminus B(0, R_{\e })} ( \g _n ) < \e$ for $n \geq n_\e$. 

\smallskip

(b) There exists $ \g \in \Eo $ such that 
$ \nabla \g _n \rightharpoonup \nabla \g $ weakly in $ L^2(\R^N)$ and  
$ \g _n \lra \g $ strongly in $L^2( B(0, R))$ for any $ R > 0$ 
as $ n \lra \infty$. 

\smallskip

Then $ Q( \g_n ) \lra Q( \g )$ as $ n \lra \infty$. 

\end{Lemma}

We postpone the proof of  Lemma \ref{L4.11} and we complete the proof of Theorem \ref{T4.7}. 
From (\ref{4.31}), (\ref{4.32}) and Lemma \ref{L4.11} it follows that 
$Q(\psi ) = \ds \lim_{k \ra \infty} Q(\tilde{\psi}_{n_k } ) = q.$
Then necessarily 
$ \ov{E}(\psi ) \geq E_{min}(q) = \ds \lim_{k\ra \infty} \ov{E}(\tilde{\psi}_{n_k } ).$
From (\ref{4.36}) we get $ \ov{E}(\psi ) = E_{min}(q)$, hence $ \psi $ is a minimizer of $\ov{E}$ 
under the constraint $ Q (\psi) = q$. Taking into account (\ref{4.33}), (\ref{4.35}) 
and the fact that $ \ov{E}(\tilde{\psi }_{n _k} ) \lra \ov{E}( \psi ) $, 
we infer that 
$\ii _{\R^N}  |\nabla \tilde{\psi }_{n_k}| ^2 \, dx \lra \ii _{\R^N}  |\nabla \psi | ^2 \, dx .$
Together with the weak convergence 
$\nabla \tilde{\psi}_{n_k } \rightharpoonup \nabla \psi $ in $ L^2( \R^N)$, 
this gives the strong convergence 
$\| \nabla \tilde{\psi}_{n_k } - \nabla \psi \|_{L^2(\R^N)} \lra 0 $ as $ k \lra \infty$
and Theorem \ref{T4.7} is proven. 
\hfill
$\Box$

\medskip

{\it Proof of Lemma \ref{L4.11}.}
It follows from Lemma \ref{L4.1} and Lemma \ref{L4.2} (ii) that 
there are $ \e _ 0 > 0 $ and  $C_0>0$ such that for 
any $ \phi \in \Eo $ satisfying $E_{GL}( \phi ) \leq \e_0 $ we have 
\beq
\label{4.40}
|Q(\phi )| \leq C _0 E_{GL}(\phi). 
\eeq

Fix $ \e \in (0, \frac{ \e _0}{2})$. 
Let $R_{\e}$ and $n_\e$ be as in assumption (a). 
We will use the conformal transform. Let
\beq
\label{4.41}
v_k (x) = \left\{
\begin{array}{ll}
\g _k (x) & \quad \mbox{ if } |x| \geq R_{\e} , 
\\
\g _k  \left( \frac{ R_{\e}^2}{ |x|^2 } x \right) & \quad \mbox{ if } |x| < R_{\e},
\end{array}
\right.
\qquad
v (x) = \left\{
\begin{array}{ll}
{\g } (x) & \quad \mbox{ if } |x| \geq R_{\e} , 
\\
{\g} \left( \frac{ R_{\e}^2}{ |x|^2 } x \right) & \quad \mbox{ if } |x| < R_{\e}.
\end{array}
\right.
\eeq

A straightforward  computation gives
\beq
\label{4.43}
\int_{B(0, R_{\e})} |\nabla v_k |^2 \, dx = \!
\int_{\R^N \setminus B(0, R_{\e})}  \! |\nabla \g_k (y) |^2 
\left( \frac{ R_{\e} ^2}{|y |^2} \right)^{N-2} \!  dy 
\leq \!
\int_{\R^N \setminus B(0, R_{\e})} |\nabla \g_k  (y) |^2 \, dy , 
\eeq
\beq
\label{4.44}
\begin{array}{l}
\ds \int_{B(0, R_{\e})} \left( \ph ^2 (| v_k | ) - 1 \right)^2 \, dx = 
\int_{\R^N \setminus B(0, R_{\e})} \left( \ph ^2 (|\g _k (y) | ) - 1 \right)^2 
\left( \frac{ R_{\e} ^2}{|y |^2} \right)^{N}   dy 
\\
\\
\leq
\ds \int_{\R^N \setminus B(0, R_{\e})} \left( \ph ^2 (|\g _k   | ) - 1 \right)^2  \, dy, 
\end{array}
\eeq
so that $ v_k \in \Eo $ and $E_{GL}(v_k) < 2\e < \e _0 $.
Similarly  $ v \in \Eo $ and $E_{GL}(v) < 2\e.$
From (\ref{4.40}) we get
\beq
\label{4.45}
| Q( v_k) | \leq 2 C _0 \e \qquad \mbox{ and } \qquad | Q( v) | \leq 2 C _0 \e.
\eeq

Since $ \nabla \g _k \rightharpoonup \nabla \g $ weakly in $L^2( \R^N)$,  
a simple change of variables shows that for any fixed $ \de \in (0, R_{\e}) $ we have 
$ \nabla v _k \rightharpoonup \nabla v $ weakly in $L^2( B(0, R_{\e}) \setminus B(0, \de ))$. 
On the other hand, 
$$
\ds \int_{B(0, \de)} |\nabla v_k |^2 \, dx = \!
\int_{\R^N \setminus B(0, \frac{R_{\e}^2}{\de})}  \! |\nabla \g _k  (y) |^2 
\left( \frac{ R_{\e} ^2}{|y |^2} \right)^{N-2} \!  dy 
\leq \int_{\R^N \setminus B(0, \frac{R_{\e}^2}{\de})}  \! |\nabla \g  _k  (y) |^2 \, dy 
$$
and $ \ds \sup_{k \geq 1} 
\int_{\R^N \setminus B(0, \frac{R_{\e}^2}{\de})}  \! |\nabla \g _k  (y) |^2 \, dy 
\lra 0 $ 
as $ \de \lra 0 $ by assumption (a). We conclude that 
\beq
\label{4.46}
 \nabla v_k \rightharpoonup \nabla v \qquad \mbox{ weakly in } L^2( B(0, R_{\e})).
\eeq

Since $ \g _k \lra \g $ in $L^2( B(0,R)) $ for any  $R>0$, 
we have  for any fixed $ \de \in (0, R_{\e})$,
$$
\int_{B(0, R_{\e}) \setminus B(0, \de)} |v_k - v|^2 \, dx 
= \int_{B(0, \frac{R_{\e}^2}{\de})  \setminus  B(0, R_{\e})} |\g _k  (y)- \g (y) |^2 
\left( \frac{ R_{\e} ^2}{|y |^2} \right)^{N}   dy 
\lra 0 \mbox{ as } k \lra \infty. 
$$
It is easy to see that there is $ p > 2$ such that 
$\left( \left(|v_k | - 2  \right)_+ \right)_{k \geq 1}$ 
is bounded in $L^p(\R^N)$. 
(If $ N \geq 3$ this follows for $ p = 2^*$  from the Sobolev embedding 
because  $\|\nabla v_k \|_{L^2(\R^N)}^2 \leq E_{GL} (v_k)\leq 2 \e$. 
If $ N =2$, the fact that $E_{GL} (v_k)\leq 2 \e$ implies  that
 $\Lo ^2( \{ | v_k | \geq 2  \})$ and $\|\nabla v_k \|_{L^2(\R^2)} $
are bounded  and the conclusion follows from (\ref{ineq2}).)
Using H\"older's inequality we obtain 
$$
\int_{B(0, \de )} 
\left(|v_k | - 2  \right)_+  ^2 \, dx \leq \| \left(|v_k | - 2  \right)_+ \| _{L^p(\R^N)} ^{2}
\left( \Lo ^N ( B(0, \de )) \right)^{1 - \frac 2p} 
$$
and the last quantity tends to zero as $ \de \lra 0 $ uniformly with respect to $k$. 
This implies 
$$
\int_{B(0, \de )}  | v_k |^2 \, dx \lra 0 \qquad \mbox{ as } \de \lra 0 
\mbox{ uniformly with respect to } k
$$
and we conclude that 
\beq
\label{4.47}
v_k \lra v \qquad \mbox{ in  } L^2(B(0, R_{\e})).
\eeq
Let 
\beq
\label{4.42}
w_k = \g _k  - v_k, \qquad w = \g - v.
\eeq
It is obvious that $ w_k, \, w \in H_0^1 (B(0, R_{\e}))$, $ \g _k = v_k + w_k $ and $ \g = v + w$. 
From  assumption (b), (\ref{4.46})  and (\ref{4.47}) it follows that
\beq
\label{4.48}
w_k \lra w \quad \mbox{ strongly $ \quad $ and } \quad 
\nabla w_k \rightharpoonup \nabla w \quad \mbox{  weakly in } L^2(B(0, R_{\e})).
\eeq

Using Definition \ref{D2.2} we have 
\beq
\label{4.49}
\begin{array}{l}
| Q( \g _k ) - Q( \g ) | \leq |Q(v_k ) - Q(v) |
+ \big\vert L ( \langle i \frac{ \p v_k}{\p x_1 } , w_k \rangle - \langle i \frac{ \p v}{\p x_1 } , w \rangle  ) \big\vert
\\
\\
+ \big\vert L ( \langle i \frac{ \p w_k}{\p x_1 } , v_k \rangle - \langle i \frac{ \p w}{\p x_1 } , v \rangle )  \big\vert
+ \big\vert L ( \langle i \frac{ \p w_k}{\p x_1 } , w_k \rangle - \langle i \frac{ \p w}{\p x_1 } , w \rangle ) \big\vert.
\end{array}
\eeq
From (\ref{4.45}) we get  $ |Q(v_k ) - Q(v)  | \leq 4 C _0 \e$. 
Since $ w_k = 0 $ and $ w = 0 $ outside $\ov{B}(0, R_{\e})$, using  the definition of $L$ we obtain 
$$
\begin{array}{l}
 L ( \langle i \frac{ \p v_k}{\p x_1 } , w_k \rangle - \langle i \frac{ \p v}{\p x_1 } , w \rangle  ) 
= {\ds \int _{B(0, R_{\e})} }
\langle i \frac{ \p v_k}{\p x_1  } - i \frac{ \p v}{\p x_1  }, w \rangle 
+ \langle i \frac{ \p v_k}{\p x_1  }, w_k - w \rangle \, dx 
\lra 0 \quad \mbox{ as } k \lra \infty 
\end{array}
$$
because $\frac{ \p v_k}{\p x_1} - \frac{ \p v}{\p x_1} \rightharpoonup 0 $ 
weakly and $ w_k - w \lra 0 $ strongly  in $L^2(B(0, R_{\e}))$.
Similarly the last two terms in (\ref{4.49}) tend to zero as $ k \lra \infty$. 
Finally we get 
$
| Q( \g _k ) - Q( \g ) | \leq  ( 4 C _0 + 1) \e 
$ for all sufficiently large $ k$.
Since $ \e \in (0, \frac{\e_0}{2})$ is arbitrary, the conclusion of Lemma \ref{L4.11} follows.
\hfill
$\Box $

\begin{Corollary}
\label{C4.13}
Assume that $ N \geq 2 $ and (A1), (A2) are satisfied. 
If $ (\g _n ) _{n \geq 1} \subset \Eo $,  $ \g \in \Eo $ are such that 
$ d_0 ( \g _n , \g ) \lra 0 $, then 
$\ds \lim_{n \ra \infty}  Q( \g _n ) = Q( \g)$ and 
$ \ds \lim_{n \ra \infty} \| V(|\g _n |^2) - V(|\g |^2) \|_{L^1( \R^N) } = 0 .$

In particular,  $Q$ and $E$ are continuous functionals on $ \Eo $ endowed with 
the semi-distance $d_0$.

\end{Corollary}

{\it Proof.} 
We have $ \nabla \g _n \lra \nabla \g $ and $(|\g _n| - |\g |) \lra 0 $ 
in $L^2( \R^N)$ as $ n \lra \infty$, hence 
$
|\nabla \g _n |^2 + \frac 12  \left( \ph^2( |\g _n |) - 1 \right)^2 
\lra 
|\nabla \g  |^2 + \frac 12 \left( \ph^2( |\g |) - 1 \right)^2 
$
in $L^1( \R^N)$, and consequently $ (\g _n) _{n \geq 1}$ satisfies 
assumption (a) in Lemma \ref{L4.11}.

Consider a subsequence $(\g_{n_{\ell}})_{\ell \geq 1}$ of $(\g_n)_{n \geq 1}$.
Then there exist a subsequence 
$(\g_{n_{\ell _k}})_{k \geq 1} $ and $ \g _0 \in \Eo $  that  satisfy (\ref{4.32}).
Since $ \nabla \g _{n _{\ell _k}} \rightharpoonup \nabla \g _0 $ weakly 
in $L^2( \R^N)$ and $ \nabla \g _{n _{\ell _k}} \lra \nabla \g  $  in $L^2( \R^N)$ 
we see that $\nabla \g _0 = \nabla \g$ a.e. 
on $ \R^N$, hence there is a constant $ \beta \in \C$ such that 
$ \g_0 = \g + \beta $ a.e. on $ \R^N$. The convergence 
$|\g_{n _{\ell _k}}|  \lra |\g _0|$ in $L_{loc}^2( \R^N) $ gives 
$ |\g _0 | = |\g | $ a.e. on $ \R^N$. By the definition of $Q$ it follows that
$Q(\g _0) = Q( \g + \beta ) = Q(\g)$. 
Using Lemma \ref{L4.11} we get $Q( \g_{n _{\ell_k}}) \lra Q( \g _0 ) = Q(\g)$ as $ k \lra \infty$
and Lemma \ref{L4.10} implies that 
$V(| \g_{n_{\ell _k}} |^2) \lra V(|\g _0|^2) = V(|\g |^2) $ in $ L^1( \R^N)$ as $ k \lra \infty$.
Hence any subsequence $(\g_{n_{\ell}})_{\ell \geq 1}$ of $(\g_n)_{n \geq 1}$ contains a subsequence 
$(\g_{n_{\ell _k}})_{k \geq 1} $ such that $Q( \g_{n _{\ell_k}}) \lra Q( \g  )$
and $ \|  V(| \g_{n_{\ell _k}} |^2) - V(|\g |^2) \|_{L^1( \R^N)} \lra 0 $, 
and this
 clearly implies the desired conclusion.
%
\hfill
$\Box$

\medskip

Assume that for some $ q >0$ there is $ \psi \in \Eo $ such that $ Q( \psi ) = q \, $ 
and $ \ov{E}(\psi) = E_{min}(q)$. 
Using Corollary \ref{C4.13}, for any sequence $(\psi_n )_{n \geq 1} \subset \Eo $ 
such that $ d_0 ( \psi _n, \psi ) \lra 0 $ and 
for any sequence of points $(x_n )_{n \geq 1} \subset \R^N$ we have 
$Q(\psi _n(\cdot + x_n)) \lra q $ and $\ov{E}(\psi _n(\cdot + x_n)) \lra E_{min}(q). $
Hence the convergence result provided by Theorem \ref{T4.7} for minimizing 
(sub)sequences of $\ov{E}$ under the constraint $Q= q$ is optimal.

\medskip

Next we show that if $ V \geq 0 $ on $[0, \infty)$, the minimizers of $ \ov{E} = E$ at fixed momentum are traveling waves to (\ref{1.1}).
We denote by $ d^- E_{min}(q) $ and $ d^+ E_{min}(q) $ the left and right derivatives 
of $E_{min}$ at $ q>0 $ (which exist and are finite for any $ q > 0 $ because 
$E_{min}$ is concave). We have: 

\begin{Proposition}
\label{P4.12}
Let $ N \geq 2 $ and $ q > 0$. 
Assume that $V(s) \geq 0 $ for any $ s \geq 0$ and $\psi $ is a minimizer of $E$ in the set 
$ \{ \phi \in \Eo\; | \; Q( \phi ) = q \}$. Then: 

\smallskip

(i) There is $ c \in [d^+ E_{min}(q),  d^- E_{min}(q) ] $ such that $ \psi $ satisfies
\beq
\label{4.50}
i c \psi _{x_1} + \Delta \psi + F(|\psi |^2) \psi = 0 \qquad \mbox{ in } \Do' ( \R^N).
\eeq

\smallskip

(ii) Any solution $\psi \in \Eo $ of (\ref{4.50}) satisfies 
$ \psi \in W_{loc}^{2, p } ( \R^N)$ and $\nabla  \psi \in W^{1, p } ( \R^N)$ 
for any $ p \in [2, \infty)$, $\psi $ and $ \nabla \psi $ are bounded 
and $ \psi \in C^{1, \al } (\R^N)$ for any $ \al \in [0, 1)$.

\smallskip

(iii) After a translation, $ \psi $ is axially symmetric with respect 
to the $ x_1-$axis if $ N \geq 3$. 
The same conclusion holds for $N=2$ if we assume in addition that $F$ is $C^1$. 

(iv) For any $ q > q_0 $ there are $ \psi ^+, \psi ^- \in \Eo $ such that 
$Q( \psi ^+ ) = Q( \psi ^- )=p$, $E( \psi ^+ ) = E( \psi ^- )= E_{min}(p)$ and 
$ \psi ^+, \psi ^-$ satisfy (\ref{4.50}) with speeds $c^+ = d^+E_{min}(p)$ and 
$c^- = d^-E_{min}(p)$, respectively.

\end{Proposition} 

{\it Proof. } 
(i) It is easy to see that 
$ \Delta \psi + F( |\psi |^2 ) \psi \in H^{-1}(\R^N)$, $ i \psi _{ x_1} \in L^2( \R^N)$ 
and for any $ \phi \in C_c^{\infty }( \R^N)$ we have $ \psi + \phi \in \Eo $, 
$ \ds \lim_{t \ra 0 } \frac 1t ( Q( \psi + t \phi ) - Q( \psi)) 
= 2 \langle i \psi _{x _1}, \phi \rangle_{L^2( \R^N)}$ and 
\beq
\label{gateaux2}
\begin{array}{rcl}
\ds \lim_{t \ra 0 } \frac 1t ( E( \psi + t \phi ) - E( \psi)) & = &
2 \ds \int_{\R^N} \langle \nabla \psi , \nabla \phi \rangle 
-  F( |\psi |^2 ) \langle \psi, \phi \rangle \, dx 
\\
& = &  -2 \langle \Delta \psi + F( |\psi |^2 ) \psi , \phi \rangle_{H^{-1}(\R^N), H^1 (\R^N)}.
\end{array}
\eeq
Denote $ E' ( \psi) . \phi = -2 \langle \Delta \psi + F( |\psi |^2 ) \psi , \phi \rangle_{H^{-1}(\R^N), H^1 (\R^N)}$ and 
$ Q' ( \psi) . \phi = 2 \langle i \psi _{x _1}, \phi \rangle_{L^2( \R^N)}$.
We have 
$ E( \psi + t \phi ) \geq E_{min}(Q( \psi + t \phi))$, hence for all $ t > 0$
\beq
\label{4.52}
\frac 1t ( E( \psi + t \phi ) - E ( \psi ) ) \geq \frac 1t ( E_{min}( Q(\psi + t \phi )) 
- E_{min} ( q ) ).
\eeq
If $  Q' ( \psi) . \phi >0$, we have $ Q( \psi + t \phi ) > Q ( \psi ) = q $ 
for $ t > 0$ and $t$ close to $0$, 
then passing to the limit as $ t \downarrow 0 $ in (\ref{4.52}) we get 
$ E' ( \psi) . \phi \geq d^+ E_{min} (q)  Q' ( \psi) . \phi $. 
If $  Q' ( \psi) . \phi < 0$, we have $ Q( \psi + t \phi ) < Q ( \psi ) = q $ 
for $t$ close to $0$ and $ t > 0$, then
passing to the limit as $ t \downarrow 0 $ in (\ref{4.52}) we get 
$ E' ( \psi) . \phi \geq d^- E_{min} (q)  Q' ( \psi) . \phi $. 
Putting $ - \phi $ instead of $ \phi$ in the above, we discover 
\beq
\label{4.53}
\begin{array}{l}
d^+ E_{min} (q)  Q' ( \psi) . \phi  \leq E' ( \psi) . \phi 
\leq d^- E_{min} (q)  Q' ( \psi) . \phi
\qquad \mbox{ if }  Q' ( \psi) . \phi > 0, \mbox{ and } 
\\
d^- E_{min} (q)  Q' ( \psi) . \phi  \leq E' ( \psi) . \phi 
\leq d^+ E_{min} (q)  Q' ( \psi) . \phi
\qquad \mbox{ if }  Q' ( \psi) . \phi < 0.
\end{array}
\eeq

Let $ \phi _0 \in C_c^{\infty}(\R^N)$ be such that $ Q' ( \psi) . \phi _0 = 0$. We claim that $ E' ( \psi) .  \phi _0 = 0$. 
To see this, consider $ \phi \in C_c^{\infty}(\R^N)$  such that $ Q' ( \psi) . \phi \neq 0$. 
(Such a $\phi$ exists for otherwise, we would have 
$ 0 = Q' ( \psi) . \phi = 2 \langle i \psi _{x _1}, \phi \rangle_{L^2( \R^N)} $ 
for any $ \phi \in C_c^{\infty}(\R^N)$, yielding $ \psi _{x _1} = 0 $, 
hence $ Q ( \psi) = 0 \not = q $.) 
Then for any $n\in \N$ we have $ Q' ( \psi) . ( \phi + n \phi _0 ) = Q' ( \psi) .  \phi  $. 
From (\ref{4.53}) it follows that 
$E' ( \psi) . (\phi + n \phi _0 )=  E' ( \psi) .  \phi + n  E' ( \psi) .  \phi _0 $ 
is bounded, thus necessarily $ E' ( \psi) .  \phi _0 = 0$. 

Take $ \phi _1  \in C_c^{\infty}(\R^N)$ such that $ Q' ( \psi) . \phi _1 = 1$.  
Let $ c  =  E' ( \psi) .  \phi _1.$ 
Using (\ref{4.53}) we obtain $ c \in [d^+ E_{min}(q),  d^- E_{min}(q) ] $. 
For any $ \phi \in C_c^{\infty} (\R^N)$ we have 
$Q' ( \psi) .( \phi - (Q' ( \psi) .\phi) \phi _1 ) = 0$, hence
$E' ( \psi) .( \phi - (Q' ( \psi) .\phi) \phi _1 ) = 0$, that is 
$E' ( \psi) . \phi = c Q' ( \psi) .\phi$ 
and this is precisely (\ref{4.50}).

\medskip

(ii) If $ N \geq 3$ this is Lemma 5.5 in \cite{M10}. 
If $ N =2$ the proof is very similar and we omit it. 

\medskip

(iii) If $ N \geq 3$, the axial symmetry  follows from the fact that the 
minimizers are $C^1$ and from  Theorem 2' p. 329 in \cite{M7}.
We use an argument due to  O. Lopes \cite{lop1} to give a proof which requires 
$F$ to be  $C^1$, but works also for  $N =2$. 
Let $S_t^+$ and $S_t ^-$ be as in (\ref{4.10}) and (\ref{4.11}), respectively.
Proceeding as in the proof of Lemma \ref{L4.5} (ii), 
we find $ t \in \R$ such that $ Q( S_t^+ \psi) =  Q( S_t^- \psi) = q.$
This implies $E(S_t^+ \psi ) \geq E_{min}(q) $ and $E(S_t^- \psi ) \geq E_{min}(q) $. 
On the other hand $E(S_t^+ \psi ) + E(S_t^- \psi ) = 2 E( \psi )= 2 E_{min}(q)$, 
thus necessarily $E(S_t^+ \psi ) = E(S_t^- \psi ) = E_{min}(q) $ 
and $S_t^+ \psi $ and $S_t^- \psi $ are also minimizers.
Then  $S_t^+ \psi $ and $S_t^- \psi $  satisfy (\ref{4.50}) 
(with some coefficients $c_+$ and $c_-$ instead of $c$)
and have the regularity properties given by (ii). 
Since $S_t^+ \psi = \psi $ on $\{ x_N > t  \}$ and $S_t^- \psi = \psi $ 
on $\{ x_N < t  \}$, we infer that necessarily $ c_+ = c_- = c$. 
Let $ \phi _ 0 (x) = e^{ \frac{ i c x_1}{2}} \psi (x)$, 
$ \phi _ 1 (x) = e^{ \frac{ i c x_1}{2}}S_t^+  \psi (x)$, 
$ \phi _ 2 (x) = e^{ \frac{ i c x_1}{2}} S_t^- \psi  (x)$. 
Then $ \phi _0, \; \phi _1 $ and $ \phi _2$ are bounded, belong 
to $W_{loc}^{2, p}(\R^N)$ for any $ p\in [2,\infty)$  and solve the equation 
$$
\Delta \phi + \left( \frac{ c^2}{4} + F(|\phi |^2) \right) \phi = 0 \qquad \mbox{ in } \R^N. 
$$
Since $F$ is $ C^1$ and $ \phi _0$, $\phi _1$ are bounded, the function 
$ w = \phi _1 - \phi _0$ satisfies an equation 
$$
\Delta w + A(x) w = 0  \qquad \mbox{ in } \R^N, 
$$
where $A(x)$ is a $2 \times 2$ matrix and $ A \in L^{\infty}(\R^N)$. 
Since $w \in H_{loc}^2( \R^N)$ and $w = 0 $ in $ \{ x_N > t \}$, the 
Unique Continuation Theorem (see, for instance, the appendix of \cite{lop1}) 
implies that $ w = 0 $ on $ \R^N$, that is 
$S_t^+ \psi = \psi $ on $ \R^N$. We have thus proved that $ \psi $ is 
symmetric with respect to the hyperplane  $\{ x_N = t \}$. 
Similarly we prove that for any $ e \in S^{N-1}$ orthogonal to $ e_1 = (1, 0 , \dots, 0)$ 
there is $ t_e \in \R$ such that $ \psi $ is symmetric with respect to the hyperplane 
$ \{ x \in \R^N \; | \; x.e = t  _e \}$.
Then it is easy to see that  after a translation $ \psi $  is symmetric 
with respect to $O x_1$. 

\medskip

iv) Consider a sequence $ q_n \uparrow q $.  We may assume $ q_n > q_0 $ 
for each $n$. By Theorem \ref{T4.7} there is $ \psi _n \in \Eo $ such that 
$ Q(\psi _n ) = q_n \lra q $ and $E(\psi _n) = E_{min}(q_n ) \lra E_{min}(q)$ 
by continuity of $E_{min}$. Since $ q > q_0 $ we have $E_{min}(q) < v_s q$ and using Theorem \ref{T4.7} again
 we infer that there are a subsequence $(\psi_{n_k})_{k \geq 1}$, 
a sequence $(x_k)_{k \geq 1} \subset \R^N$ and $ \psi ^- \in \Eo$ such that 
$Q( \psi ^- ) = q$, $E(\psi^-) = E_{min}(q)$ and, denoting 
$\tilde{\psi}_{n_k} = \psi_{n_k} ( \cdot + x_n)$, we have
$\tilde{\psi}_{n_k} \lra \psi^- $ a.e. on $ \R^N$ and 
$ d_0 (\tilde{\psi}_{n_k} , \psi ^-) \lra 0 $ as $ k \lra \infty$.

By (i) we know that each $ \tilde{\psi} _{n_k}$ satisfies (\ref{4.50}) for some 
$c_{n_k} \in [ d^+E_{min}(q_{n_k}), d^-E_{min}(q_{n_k}) ]$. 
Since $E_{min}$ is concave, we have $c_{n_k} \lra d^-E_{min}(q)$ as $ k \lra \infty$. 
It is easily seen that $ \tilde{\psi} _{n_k} \lra \psi ^-$ and 
$F(| \tilde{\psi} _{n_k} |^2) \tilde{\psi} _{n_k} \lra F(|\psi ^-|^2) \psi ^-$ 
in $ \Do '(\R^N)$. Writing (\ref{4.50}) for each $\tilde{\psi} _{n_k}$ and 
passing to the limit as $ k\lra \infty$ we infer that 
$\psi ^- $ satisfies (\ref{4.50}) in $\Do'( \R^N)$ with $ c = d^- E_{min}(q)$. 

The same argument for a sequence $ q_n \downarrow q$ gives the existence 
of $ \psi ^+$.  
\hfill
$\Box $

\medskip

If $F$ satisfies assumption (A4) in the introduction and $ F''(1) \neq 3$, 
we prove that in space  dimension $N=2$ we have $ q_0 = 0$. This implies that we can minimize $E$ under 
the constraint $Q= q$ for any $ q >0$. The traveling waves obtained 
in this way have small energy and  speed tending to $ v_s $ as $ q \lra 0 $. 
For the two-dimensional Gross-Pitaevskii equation, the numerical and formal 
study in \cite{JR} suggests that these traveling waves are rarefaction 
pulses asymptotically described by the ground states of the Kadomtsev-Petviashvili I  (KP-I) equation. 
The rigorous convergence, up to rescaling and renormalization, of the 
traveling waves of (\ref{1.1}) in the transonic limit 
to the ground states of the (KP-I) equation 
has been proven  in \cite{BGS1} in the case of  the two-dimensional Gross-Pitaevskii equation. 
That result has been extended in \cite{CM} to a general 
nonlinearity satisfying (A1), (A2) and (A4) with $ F''(1) \neq 3 $.

A  result similar to Theorem \ref{T4.13} below  is {\it not } true in higher dimensions: 
in view of Proposition \ref{smallE}   we have $ q_0 > 0 $ for any $N \geq 3$. 
If $ N \geq 3$, the existence of traveling waves with speed close to $ v_s$ is guaranteed by 
Theorem 1.1 and Corollary 1.2 p. 113 in \cite{M10}. 
In space dimension three, the convergence of the traveling waves constructed in \cite{M10} 
 to the ground states of the  three-dimensional (KP-I) equation as $ c \to v_s$ has   been rigorously justified 
 under the same assumptions as in dimension two (see Theorem 6 in \cite{CM}). 
It was also shown  in \cite{CM}  that these solutions have high energy and momentum 
(of order $1/\sqrt{v_s ^2 - c^2 } $ as $ c \to v_s$) and thus lie on the upper branch in figure  \ref{diaggrospit} (b).

\begin{Theorem}
\label{T4.13} 
Suppose that $ N =2$, the assumption (A4) in the introduction holds
and  $F''(1) \neq 3 $. Then $E_{min}(q) < v_s q $ for any $ q > 0$. 
In other words, $q_0 = 0$.
\end{Theorem}

\begin{remark} \label{R4.16} \rm
If $N = 2$, $ V \geq 0$  and (A1), (A2) and (A4) hold with $ F''(1) \neq 3 $, 
it follows from 
Theorems \ref{T4.13} and \ref{T4.7} that for any $ q > 0 $ there is 
$ \psi _ q \in \Eo $ such that $ Q( \psi _q ) = q $ and $ E( \psi _p) = E_{min}(q)$. 
Proposition \ref{P4.12} (i) implies that $ \psi_q $ is a traveling wave of (\ref{1.1}) 
of speed $ c ( \psi _q ) \in [ d^+E_{min}(q), d^-E_{min}(q) ].$
Using Lemmas \ref{L4.3} and \ref{L4.4} we infer that $ c( \psi _q) \lra  v_s $ 
as $ q \lra 0 $. In particular, we see that there are traveling waves of 
arbitrarily small energy whose speeds are arbitrarily close to $ v_s$.

\medskip

In view of the formal asymptotics given in \cite{JR}, it is natural to try to 
prove Theorem \ref{T4.13}  by using test functions constructed from 
an ansatz related to the (KP-I) equation.


\end{remark} 

{\it Proof of Theorem \ref{T4.13}.}
Fix $ \g > 0 $ (to be chosen later). 
We consider the (KP-I) equation 
\beq
\label{KP-I}
u_t - \g u u _x + \frac {1}{ v_s ^2} u_{xxx} - \p_x^{-1} u _{yy} = 0 , 
\qquad t \in \R, ( x, y ) \in \R^2 ,
\eeq
where $u$ is real-valued. Let $ Y$ be the completion of 
$ \{ \p_x \phi \; | \; \phi \in C_c^{\infty}(\R^2, \R) \}$ for the norm 
$\| \p _ x \phi \|_Y ^2  =  \| \p _ x \phi \|_{L^2(\R^2)} ^2 + 
v_s ^2 \| \p _ y \phi \|_{L^2(\R^2)} ^2 + \| \p _ {xx} \phi \|_{L^2(\R^2)} ^2. $
A traveling wave for (\ref{KP-I}) moving with velocity $ \frac{1}{v_s ^2}$ is a 
solution of the form 
$u( t, x, y ) = v ( x - \frac{t }{ v_s ^2}, y)$, where $ v \in Y$. 
The traveling wave profile $v$  solves the equation 
$$
\frac{1}{v_s ^2} v_x + \g v v_x - \frac{1}{v_s ^2} v_{xxx} + \p_x^{-1} v _{yy} = 0 
\qquad \mbox{ in } \R^2, 
$$
or equivalently, after integrating in $x$, 
\beq
\label{4.55}
\frac{1}{v_s ^2} v + \frac{\g}{2} v ^2  - \frac{1}{v_s ^2} v_{xx} + \p_x^{-2} v _{yy} = 0 
\qquad \mbox{ in } \R^2. 
\eeq
It is a critical point of the functional (called the {\it action})
$$
\Sr (v) = \int_{\R^2 } \frac{1}{v_s ^2} |v|^2 + \frac{1}{v_s ^2} |v _x|^2 + | \p _x ^{-1} v _y |^2 \, dx \, dy 
+ \frac{ \g }{3} \int_{\R^2} v ^3 \, dx \, dy
= \frac{1}{v_s ^2} \| v \|_Y ^2 + \frac{ \g }{3} \int_{\R^2} v ^3 \, dx \, dy.
$$
Equation (\ref{4.55}) is indeed nonlinear if $\g \not = 0 $. The existence of a 
nontrivial traveling wave solution $w$ for (KP-I) follows 
from Theorem 3.1 p. 217 in \cite{dBSIHP}. The solution found in \cite{dBSIHP} 
minimizes $ \| \cdot \|_Y$ in the set 
$ \{ v \in Y \; | \; \ds \int_{\R^2} v ^3 \, dx \, dy = \int_{\R^2} w ^3 \, dx \, dy\}$.
It was also proved (see Theorem 4.1 p. 227 in \cite{dBSIHP}) that
$ w \in H^{\infty}(\R^2) := \ds \cap_{m \in \N} H^m( \R^2)$,  $ \p _x ^{-1}w _y \in H^{\infty}(\R^2) $
and $w$ minimizes the action $ \Sr $ among all nontrivial solutions of (\ref{4.55})
(that is, $w$ is a {\it ground state}). 
Moreover, $w$ satisfies the following integral identities: 
\beq
\label{identites}
\left\{\begin{array}{ll}
\ds \int_{\R^2} \frac{1}{v_s^2}  w^2 + \frac {\g}{2} w^3 + \frac{1}{v_s^2} | \p _x w |^2 
+ | \p _x^{-1} w _y |^2 \, dx \, dy = 0, 
\\  
\\ 
\ds \int_{\R^2} \frac{1}{v_s^2}  w^2 + \frac {\g}{3} w^3 - \frac{1}{v_s^2} | \p _x w |^2 
+ 3 | \p _x^{-1} w _y |^2 \, dx \, dy = 0,
\\
\\
\ds \int_{\R^2} \frac{1}{v_s^2}  w^2 + \frac {\g}{3} w^3 + \frac{1}{v_s^2} | \p _x w |^2 
- | \p _x^{-1} w _y |^2 \, dx \, dy = 0.
\end{array}\right.
\eeq
The first identity is obtained by multiplying (\ref{4.55}) by $w$ and integrating, 
while the two other are Pohozaev identities associated to the scalings in $x$, 
respectively in $y$. They are formally obtained by multiplying (\ref{4.55}) by 
$ xw$, respectively by $ y \p _x^{-1} w _y$ and integrating by parts; see the 
proof of Theorem 1.1 p. 214 in \cite{dBSIHP} for a rigorous justification. 

Comparing $ \Sr(w) $ to the last equality in (\ref{identites}) we get 
\beq
\label{4.57}
\int_{\R^2}  | \p _x^{-1} w _y |^2 \, dx \, dy = \frac 12 \Sr (w).
\eeq
In particular, $ \Sr (w) > 0$. 
Then from the three identities (\ref{identites}) we obtain
\beq
\label{4.58}
\frac{1}{v_s ^2} \! \int_{\R^2} \! |w|^2 \, dx \, dy = \frac 32 \Sr (w), \quad
\frac{1}{v_s ^2} \! \int_{\R^2} \! |w_x |^2 \, dx \, dy =  \Sr (w), \quad
\frac{ \g }{6} \! \int_{\R^2} \! w^3 \, dx \, dy = - \Sr (w).
\eeq

Let $ w$ be as above and let $ \phi = v_s \p _x ^{-1} w$, so that $ \p _x \phi = v_s w$.
For $ \e > 0 $ small we define 
$$
\rho _{\e } (x, y ) =  1 + \e ^2 w ( \e x , \e ^2 y), \qquad
\theta_{\e } (x,y) = \e \phi ( \e x, \e ^2 y ) , \qquad
U_{\e} = \rho _{\e } e^{-i \theta _{\e}}. 
$$
Then $ U_{\e} \in \Eo $ (because $w \in H^{\infty}(\R^2)$). 
For $\e$ sufficiently small we have $V(|U_{\e}|^2) = V(\rho_{\e}^2) \geq 0$, hence $\ov{E}(U_{\e}) = {E}(U_{\e})$.
A straightforward computation and (\ref{4.57}), (\ref{4.58}) give
$$
\int_{\R^2} \Big\vert \frac{\p \rho_{\e}}{\p x } \Big\vert ^2 dx \, dy 
= \e ^3  \int_{\R^2} \Big\vert \frac{\p w}{\p x } \Big\vert ^2 dx \, dy  
= \e ^3  v_s ^2 \Sr (w) = 2  \e ^3   \Sr (w), 
$$

$$
\int_{\R^2} \Big\vert \frac{\p \rho_{\e}}{\p y } \Big\vert ^2 dx \, dy  
= \e ^5  \int_{\R^2} \Big\vert \frac{\p w}{\p y } \Big\vert ^2 dx \, dy ,
$$

$$
\begin{array}{l}
\ds \int_{\R^2} \rho _{\e } ^2  \Big\vert \frac{\p \theta_{\e}}{\p x } \Big\vert ^2 dx \, dy 
= \e  \int_{\R^2} ( 1 + \e ^2 w ) ^2 \vert \phi _x \vert ^2 dx \, dy 
\\
= \ds \e   v_s ^2 \int_{\R^2} ( 1 + \e ^2 w ) ^2  w^2 dx \, dy
= \frac 32 v_s ^4 \Sr (w) \e - \frac{12}{\g } v_s ^2 \Sr (w) \e ^3
+ v_s ^2 \e ^5   \int_{\R^2} w^4  dx \, dy,
\end{array}
$$

$$
\begin{array}{l}
\ds \int_{\R^2} \rho _{\e } ^2  \Big\vert \frac{\p \theta_{\e}}{\p y } \Big\vert ^2 dx \, dy 
= \e ^3  \int_{\R^2} ( 1 + \e ^2 w ) ^2 \vert \phi _y \vert ^2 dx \, dy 
\\
= \ds \e^3   v_s ^2 \int_{\R^2} ( 1 + \e ^2 w ) ^2  | \p _x ^{-1}w_y |^2 dx \, dy
\\
= \ds \frac 12  v_s ^2 \Sr (w) \e ^3 
+ 2 \e ^5  v_s ^2 \int_{\R^2} w | \p _x ^{-1}w_y |^2 dx \, dy 
+ \e^ 7   v_s ^2 \int_{\R^2} w ^2| \p _x ^{-1}w_y |^2 dx \, dy .
\end{array}
$$
Using (\ref{lift}) we get 
\beq
\label{4.59}
\begin{array}{l}
Q( U_{\e}) = \ds \int_{\R^2}  ( \rho _{\e }^2 - 1 ) \frac{ \p \theta _{\e}}{\p x } dx \, dy
= \e  \int_{\R^2} ( 2 w + \e^2 w^2) \phi _x \, dx \, dy
\\
= \e   v_s \ds \int_{\R^2} ( 2 w + \e^2 w^2) w \, dx \, dy
= 3  v_s ^3 \Sr (w) \e - \frac { 6}{\g }  v_s \Sr (w) \e ^3. 
\end{array}
\eeq

If (A4) holds  we have the expansion 
\beq
\label{4.60}
V(s) = \frac 12 ( s - 1)^2 - \frac 16 F''(1) ( s - 1)^3 + H(s), 
\eeq
where $|H(s) | \leq C(s - 1 )^4 $ for $ s $ close to $ 1$.
Using (\ref{4.60}) and the fact that $ w \in L^p( \R^2)$ for any $ p \in [2, \infty]$, 
for small $ \e $  we may expand $V(\rho _{\e})$ and integrate to get 
$$
\begin{array}{l} 
\ds \int_{\R^2} V( \rho _{\e }^2) \, dx \, dy 
= 2 \e  \int_{\R^2} w^2 dx \, dy 
+ \e ^3 \left( 2 - \frac 43 F''( 1)  \right) \int_{\R^2} w^3 dx \, dy  + \Oo ( \e ^5)
\\
\\
= \frac 32  v_s ^4 \Sr (w) \e - \frac{6}{\g } \left(  v_s ^2 - \frac 43 F''( 1)  \right)\Sr (w) \e^3
+ \Oo ( \e ^5) .
\end{array}
$$
From the previous computations we find 
\beq
\label{4.61}
E( U_{\e}) - v_s Q(U_{\e}) =   v_s ^2 \Sr (w) \left( \frac 32 - \frac{12 - 4 F''(1)}{\g } \right) \e ^3 + \Oo ( \e^5). 
\eeq
If $ F''(1) \neq 3$, choose $ \g \in \R$ such that $\frac 32 - \frac{12 - 4 F''(1)}{\g } < 0$ 
(take, for instance, $ \g = 6 - 2 F''(1)$). 
Let $w$ be a  ground state of (\ref{4.55}) for this choice of $ \g$. 
It follows from (\ref{4.61}) that there is $ \e _0 > 0 $ such that 
$ E( U_{\e }) - v_s Q(U_{\e}) < 0 $ for any $ \e \in (0, \e _0)$ (since $\Sr (w) > 0 $). 
On the other hand, using (\ref{4.59}) we infer that there is $ \e _1 < \e _0 $ 
such that the mapping $ \e \longmapsto Q(U_{\e})$ is a homeomorphism from $(0, \e _1)$ 
to an interval $(0, q_1)$. Since $E_{min}(Q(U_{\e})) \leq \ov{E}(U_{\e}) = E(U_{\e}) < v_s Q (U_{\e})$, 
we see that $ E_{min}(q) < v_s q $ for any $ q \in (0, q_1)$. 
Then the concavity of $E_{min}$ implies $ E_{min}(q) < v_s q $ for any $ q >0$. 
\hfill
$\Box$

\medskip

We pursue with some qualitative properties of $E_{min}$ for large $q$.
Theorem \ref{T4.16} (a) below implies that the speeds of traveling waves 
obtained from Theorem \ref{T4.7} tend to $0$ as $ q \lra \infty$.

\begin{Theorem}
\label{T4.16}
$\; $ 

(a) If (A1) holds and $N\geq 2$, there is $C > 0$ such that 
$ E_{min}(q) \leq C q^{ \frac{N-2}{N-1}} \ln q \; $ for large $q$. 

(b) If $ N\geq 2$ and (A1) and (A2) hold we have 
$ \ds \lim_{q \ra \infty} E_{min}(q) = \infty$. Moreover, if 
$N \geq 3$ there is $ C> 0$ such that $ E_{min}(q) \geq C q^{ \frac{N-2}{N-1}}$.

\end{Theorem}

{\it Proof. } 
(a) Using Lemma 4.4 p. 147 in \cite{M10}  we see that there is 
a continuous mapping $R \longmapsto v_R$ from $[2, \infty )$ to 
$H^1( \R^N)$ and constants $C_i > 0$, $i=1$, $2$, $3$, such that 
\beq
\label{4.62}
\int_{\R^N} |\nabla v_R| ^2 dx \leq C_1 R^{N-2} \ln R, 
\quad
 \int_{\R^N} \big|  V(|1+ v_R|^2) \big| \, dx  \leq C_2 R^{N-2},  
\eeq
\beq
\label{4.63}
C_3 (R-2) ^{N-1} \leq Q( 1+ v_R) \leq C_3 R^{N-1}.
\eeq
Let $ q_R =  Q( 1+ v_R)$. 
The set $ \{ q_R \; | \; R \geq 2 \}$ is an interval of the form $[q_*, \infty)$. 
By (\ref{4.63}) we have 
$ C_3 ^{ - \frac{1}{N-1}} q_R^{\frac{1}{N-1}} 
\leq	R \leq 
2 + C_3 ^{ - \frac{1}{N-1}} q_R^{\frac{1}{N-1}} .$
Then using (\ref{4.62}) we get for $R$ sufficiently large
$$
E_{min}(q_R) \leq \ov{E}( 1+ v_R) \leq C_1 R^{N-2} \ln R + C_2 R^{N-2} 
\leq C q_R^{\frac{N-2}{N-1}} \ln q_R.
$$

(b) As in the proof of  Lemma \ref{L4.5} (ii), using (\ref{4.9}) we get 
 $E_{min}(q_2) \geq \left( \frac{ q_2}{q_1} \right) ^{\frac{N-2}{N-1}} E_{min}(q_1)$
for any $ q_2 > q_1 > 0$. This is the second statement of (b), and it  
implies that $ \ds \lim_{q \ra \infty} E_{min}(q) = \infty$ if $N \geq 3$.

Let $ N =2$. 
We argue by contradiction and we assume that $ \ds \lim_{q \ra \infty} E_{min}(q) $ 
is finite. Using Theorem \ref{T4.7} for $q$ sufficiently large, we may choose 
$ \psi _q \in \Eo $ such that $ Q( \psi _q ) = q $ and $ \ov{E} ( \psi _q) = E_{min}(q)$. 
Consider a sequence $ q_n \lra \infty$. From Lemma \ref{L4.6} it follows that 
$E_{GL} ( \psi _{q_n} )$ is bounded and stays away from $0$. 
Passing to a subsequence we may assume that $E_{GL} ( \psi _{q_n} )\lra \al _0 > 0$. 
Let $ \Lambda_n(t)$ be the concentration function associated to $E_{GL}(\psi _{q_n})$ 
(as in (\ref{4.15})). Arguing as in the proof of Theorem \ref{T4.7} and passing 
to a subsequence (still denoted $(q_n)_{n \geq 1}$), we see that there exist 
a nondecreasing function $ \Lambda : [0, \infty ) \lra \R$, $ \al \in [0, \al _0 ] $ 
and a sequence $ t_n \lra \infty $ satisfying (\ref{4.16}) and (\ref{annulus}). 
Then we use Lemma \ref{L4.8}  to infer that $ \al > 0$. 

If $ \al \in (0, \al _0)$, proceeding as in the proof of Theorem \ref{T4.7} and using 
Lemma \ref{splitting}  for $ \psi _{q_n}$ we see that there exist functions 
$ \psi _{n,1}, \psi _{n,2} \in \Eo $ such that (\ref{4.28})$-$(\ref{4.30}) hold. 
Passing to a subsequence if necessary, we may assume that 
$ \ov{E}(\psi _{n, i}) \lra m_i \geq 0 $ as $ n \lra \infty$. 
Since $\ds \lim_{n \ra \infty } E_{GL}(\psi _{n,i} ) > 0$, it follows from Lemma 
\ref{L4.1} that $ m_i > 0$, $i=1,2$. Using (\ref{4.29}) we see that 
$ m_1 + m_2 = \ds \lim_{q \ra \infty} E_{min}(q)$, hence 
$ 0 < m_i  <  \ds \lim_{q \ra \infty} E_{min}(q).$ 
Since $ Q( \psi _{q_n} ) = q_n \lra \infty$, from (\ref{4.30}) it follows that at 
least one of the sequences $(Q(\psi_{n, i}))_{n \geq 1}$ contains a subsequence 
$(Q(\psi_{n_k, i}))_{k \geq 1}$ that tends to $\infty$. Then 
$ \ov{E}(\psi_{n_k, i}) \geq E_{min}(Q(\psi_{n_k, i})) $ and passing to the limit 
as $ k \lra \infty $ we find $ m_i \geq \ds \lim_{q \ra \infty} E_{min}(q)$, 
a contradiction. Thus we cannot have $ \al \in (0, \al _0)$.

We conclude that necessarily $ \al = \al _0 $. 
Proceeding again as in the proof of Theorem \ref{T4.7} we infer that 
there is a sequence $(x_n)_{n \geq 1} \subset \R^N$  such that 
$\tilde{ \psi}_n = \psi _{q_n} ( \cdot + x_n)$ satisfies (\ref{4.31}). 
Then there exist a subsequence $(\tilde{ \psi}_{n _k})_{k \geq 1}$ and 
$ \psi \in \Eo $ such that (\ref{4.32}) holds. 
Using Lemma \ref{L4.11} we infer that $ Q( \tilde{ \psi}_{n _k}) \lra Q( \psi ) \in \R$ 
and this is in contradiction with $ Q( \tilde{ \psi}_{n _k})= q_{n _k}  \lra \infty.$
Thus necessarily $ E_{min}(q) \lra \infty $ as $ q \lra \infty.$
\hfill
$\Box $

\medskip

An alternative proof of the fact that $  E_{min}(q) \lra \infty $ as $ q \lra \infty$ 
is to show that for $ \psi \in \Eo $ we may write 
$ \langle i \psi _{x_1}, \psi \rangle = f + g$, where $ g \in \Yo $ and $ f $ 
is bounded in $ L^1( \R^N)$ if $E_{GL}(\psi )$ is bounded, then to use Lemma 
\ref{L4.6} to infer that $Q( \psi )$ remains bounded if $\ov{E}(\psi )$ is bounded.

\medskip

From Theorem \ref{T4.16} and Lemma \ref{L4.6} we obtain the following:

\begin{Corollary}
\label{C4.17}
For all $ M >0$, the functional $Q$ is bounded  on the set
$ \{ \psi \in \Eo \; | \; E _{GL} ( \psi ) \leq M \} .$ \\
If (A1) and  (A2) hold,  $Q$ is also bounded on the set
$\{ \psi \in \Eo \; | \; \ov{E} ( \psi ) \leq M \} .$
\end{Corollary}

\section{Minimizing the action at fixed kinetic energy}
\label{another}

Although in many important physical applications the nonlinear potential $V$ 
may achieve negative values (this happens, for instance, for the cubic-quintic NLS), 
there are no results in the literature that imply the existence of 
finite energy traveling waves for (\ref{1.1}) in space dimension two for this 
kind of nonlinearity. We develop here a method that  works if $N \geq 2$ and 
$V$ takes negative values. The method used in \cite{M10} 
(minimization of $E_c$ under a Pohozaev constraint) does not require any  assumption on sign of the potential $V$,
 hence can be applied for
the cubic-quintic NLS if $ N \geq 3$, but  does not work in space dimension two.
 Throughout this section we assume that (A1) and 
(A2) are satisfied. 

We begin with a refinement of Lemma \ref{L4.2}.

\begin{Lemma} 
\label{L5.1}
Assume that $|c| < v_s $ and let $ \e \in (0, 1 - \frac{|c|}{v_s})$. There is $k > 0$ such that for any 
$ \psi \in \Eo $ satisfying $  \ii _{\R^N} |\nabla u |^2 \, dx \leq k$ we have 
$$ E(\psi) - \e E_{GL}(\psi) \geq |c Q(\psi)|. $$
\end{Lemma} 

{\it Proof. } 
Fix $ \e _1 > 0 $ such that $ \e + \e _1 < 1 - \frac{|c|} {v_s}$. 
It follows from Lemma \ref{L4.1} that there is $ k_1 > 0 $ such that 
\beq
\label{5.1}
( 1 - \e _1 ) E_{GL} (\psi ) \leq E( \psi ) \qquad 
\mbox{ for any } \psi \in \Eo \mbox{ satisfying } \ds \int_{\R^N} |\nabla \psi |^2 \,dx \leq k_1.
\eeq

Let 
$ \tilde{F}(s) =  ( 1 - \ph ^2 ( \sqrt{s}) ) \ph ( \sqrt{ s}) \ph '( \sqrt{ s})  \frac{1}{\sqrt{s}}$.
Then $ \tilde{F}(s) =  1 - s $ for $ s \in [0, 4 ]$ and $\tilde{F}$ satisfies 
(A1) and (A2). 
Let $ \tilde{V} (s) =  \ii _s ^{ 1} \tilde{F} (\tau ) \, d \tau 
= \frac 12 \left( \ph ^2( \sqrt{s}) - 1 \right) ^2. $
Using Lemma \ref{L4.2} (ii) with $\tilde{F}$ and $\tilde{V}$ instead of $F$ and $V$ 
we infer that there is $ k \in (0, \frac{k_1}{2}) $ such that for any $ \psi \in \Eo $ with 
$ E_{GL}(\psi ) \leq 2 k$ we have
\beq
\label{5.2}
( 1 - \e - \e _1) E_{GL}(\psi ) \geq | c Q (\psi )|.
\eeq

Let $ \psi \in \Eo $ be such that $  \ii _{\R^N} |\nabla \psi |^2 \, dx \leq k .$

If $ \frac 12  \ii _{\R^N} \left( \ph ^2 (|\psi |) ^2 - 1 \right)^2 \, dx \leq k$ we have 
$E_{GL}(\psi ) \leq 2 k$ and (\ref{5.2}) holds.
Then using (\ref{5.1}) we obtain 
$ E(\psi) - \e E_{GL}(\psi) \geq ( 1 - \e - \e _1 ) E_{GL}(\psi ) \geq | c Q ( \psi ) | . $

If $ \frac 12  \ii_{\R^N} \left( \ph ^2 (|\psi |) ^2 - 1 \right)^2 \, dx > k$, let
 $ \si = \left(  \ii _{\R^N} |\nabla \psi |^2 \, dx \right)^{\frac 12}
\left( \frac 12 \! \ii _{\R^N} \! \left( \ph ^2 (|\psi |) ^2 - 1 \right)^2 \, dx \! \right)^{ - \frac 12} \!$.
Then $ \si \in (0, 1)$ and 
$$\frac 12 \ds \int_{\R^N} \left( \ph ^2 (|\psi _{\si, \si } |) ^2 - 1 \right)^2 \, dx
= \int_{\R^N} |\nabla \psi _{\si, \si } |^2 \, dx 
= \frac 12 E_{GL}(\psi_{\si, \si})
= \si ^{N-2} \int_{\R^N} |\nabla \psi |^2 \, dx < k.$$
 Using (\ref{5.1})  and (\ref{5.2}) we get $E(\psi ) \geq ( 1 - \e_1) E_{GL}(\psi )$ and 
$ ( 1 - \e - \e _1) E_{GL}(\psi _{\si, \si }) \geq | c Q (\psi _{\si, \si })|.$
Then we have
$$
\begin{array}{l}
 E(\psi) - \e E_{GL}(\psi) - |c Q(\psi)| \geq ( 1 - \e - \e _1 ) E_{GL}(\psi) - |c Q(\psi)| 
 \\
\\
\geq \ds ( 1 - \e - \e _1 ) \left( \frac{1}{\si ^{N-2}} \! \int_{\R^N} \! |\nabla \psi _{\si, \si } |^2 \, dx
+ \frac{1}{2 \si ^N} \! \int_{\R^N} \! \left( \ph ^2 (|\psi _{\si, \si } |) ^2 - 1  \right)^2 \, dx \right)
- \frac{1}{\si ^{N-1}} |c Q(\psi _{\si, \si})|
\\
\\
\geq \ds \frac{1 - \e - \e _1}{2} \left( \frac{1}{\si ^{N-2}}  
+ \frac{1}{\si ^N} \right) E_{GL}(\psi _{\si, \si}) - \frac{1- \e - \e _1}{\si^{N-1} } E_{GL}(\psi _{\si, \si})
\geq 0.
\end{array}
$$
$\; $

\vspace*{-18pt}
\hfill
$\Box$

\medskip

Let $ I( \psi ) = - Q( \psi ) +  \ii_{\R^N} V(|\psi |^2) \, dx = E( \psi) - Q( \psi ) - \ii _{\R^N} |\nabla \psi |^2 \, dx.$
\smallskip
\\
We will minimize $I(\psi)$ under the constraint $\| \nabla \psi \|_{L^2( \R^N) } \!= \! constant.$
For any $ k > 0$~we~define
$$
I_{min}(k) = \inf \Big\{ I(\psi ) \; \big| \;  \; \psi \in \Eo, 
\int_{\R^N} |\nabla \psi |^2 \, dx = k \Big\}.
$$

The next Lemmas establish the basic  properties of the function $I_{min}$. 

\begin{Lemma}
\label{L5.2}
(i) For any $ k > 0 $ we have $I_{min}(k) \leq - \frac{ 1}{v_s ^2} k $.

(ii) For any $ \de > 0 $ there is $ k( \de ) > 0 $ such that 
$I_{min}(k) \geq - \frac{ 1 + \de }{v_s ^2} k $ for any $ k \in (0,  k ( \de )).$
\end{Lemma}

{\it Proof.} 
i) Let $ N \geq 3$. Let $ q = 2 k v_s ^{ N-3}$. In the proof of Lemma \ref{L4.3} 
we have constructed a sequence $(\psi _n )_{n \geq 1} \subset \Eo $ such that
$$
 Q( \psi _n ) = q, \qquad \int_{\R^N} |\nabla \psi _n |^2 \, dx \lra \frac 12 v_s q = k v_s^{N-2} ,
\qquad  \int_{\R^N} V(|\psi _n |^2) \, dx \lra \frac 12 v_s q 
$$
and $ \psi _n $ is  constant outside a large ball. 
Let $ \si _n = k^{\frac{1}{N-2}} \left( \ii _{\R^N} |\nabla \psi _n |^2 \, dx \right)^{- \frac{1}{N-2}} .$
Then $ \si _n \lra \frac{1}{v_s}$ as $ n \lra \infty$. We get
$$
\int_{\R^N} |\nabla ( (\psi _n )_{\si _n, \si _n }) |^2 \, dx = \si _n ^{N-2} \int_{\R^N} |\nabla \psi _n |^2 \, dx  = k, 
$$
$$
Q((\psi_n)_{\si_n, \si _n}) = \si _n ^{N-1} Q( \psi _n ) \lra 
\frac{q}{v_s ^{N-1}} = \frac{ 2 k }{v_s ^2}, 
$$
$$
\int_{\R^N}  V(|(\psi _n)_{\si _n , \si _n }|^2) \, dx 
= \si _n ^N \int_{\R^N}  V(|\psi _n |^2) \, dx 
\lra \frac{1}{v_s ^N} \cdot \frac{v_s q }{2} = \frac{k}{v_s^2}.
$$
We have $ I_{min}(k) \leq I( (\psi _n )_{\si _n, \si _n })$ for each $n$ 
and passing to the limit as $ n \lra \infty $ we obtain $ I_{min}(k) \leq  - \frac{1}{v_s ^2} k .$

If $ N =2$, let $ q = \frac{2k}{v_s}$, choose $ \psi _n $ as in the proof 
of Lemma \ref{L4.3} such that 
$$
\int_{\R^2} |\nabla \psi _n |^2 \, dx  = k, \qquad
Q( \psi _n ) \lra q = \frac{ 2k}{v_s} \qquad \mbox{ and } \qquad 
\int_{\R^2}  V(|\psi _n |^2) \, dx \lra k .
$$
Let $ \si = \frac{1}{v_s}.$
Then $ \ii _{\R^2} |\nabla ((\psi _n)_{\si, \si}) |^2 \, dx  = k$, 
$ Q((\psi _n)_{\si, \si}) = \si Q( \psi _n) \lra  \frac{ 2k}{v_s ^2}$ and  \\
${ \ii _{\R^2} } V(|(\psi _n )_{\si, \si} |^2) \, dx 
=\si ^2 { \ii _{\R^2} } V(|\psi _n |^2) \, dx \lra \frac{k}{v_s ^2}$, 
hence  $I( (\psi _n)_{\si, \si}) \lra - \frac{ k}{v_s ^2}$.

\medskip

(ii) Fix $ \de > 0 $ and let $ c = \frac{v_s}{\sqrt{ 1 + \de}}.$
Lemma \ref{L5.1} implies that there is $ k > 0 $ such that for any $ \psi \in \Eo $ 
with $ \ii _{\R^N} |\nabla \psi |^2 \, dx \leq k $ we have
\beq
\label{5.3}
\ds \int_{\R^N} |\nabla \psi |^2 \, dx - c Q( \psi) + \int_{\R^2}  V(|\psi  |^2) \, dx \geq 0.
\eeq
Let $ \psi \in \Eo $ be such that ${ \ii _{\R^N} } |\nabla \psi |^2 \, dx \leq \frac{k}{c^{N-2}}.$
Then 
$ \ii _{\R^N} |\nabla \psi _{c,c }|^2 \, dx = c^{N-2}  \ii _{\R^N} |\nabla \psi |^2 \, dx \leq k$, 
hence $ \psi_{c, c}$ satisfies (\ref{5.3}), that is 
$ c^{N-2} \ii _{\R^N} |\nabla \psi |^2 \, dx + c^N I(\psi ) \geq 0$ or equivalently
$$
I(\psi) \geq - \frac{1}{c^2} \ds \int_{\R^N} |\nabla \psi |^2 \, dx = - \frac{1 + \de}{v_s ^2} \ds \int_{\R^N} |\nabla \psi |^2 \, dx.
$$
Hence (ii) holds with $ k( \de) = \frac{k}{c^{N-2}}.$
\hfill
$\Box $

\medskip

We give now global properties of $I_{min}$.

\begin{Lemma}
\label{L5.3}
The function $I_{min}$ has the following properties:

\medskip

(i) $I_{min}$ is concave, decreasing on $[0, \infty)$ and 
${\ds \lim_{k \ra \infty}} \frac{ I_{min}(k) }{k} = - \infty. $

\medskip

(ii) $I_{min}$ is subadditive, that is $I_{min}(k_1 + k_2) \leq I_{min} (k_1)+ I_{min} (k_2) $
for any $ k_1, k_2 \geq 0$. 

\medskip

(iii) If either $N\geq 3$ or ($ N=2$ and $V\geq 0 $ on $[0, \infty)$), 
we have $I_{min}(k) > - \infty$ for any $ k > 0$. 

\medskip

(iv) If ($ N = 2$ and $ \inf V < 0$), then $ I_{min}(k) = - \infty$ for 
all sufficiently large $k$. 

\medskip

(v) Assume that $ N =2$,  (A4) holds and $ F''(1) \not= 3 $. 
Then $I_{min}(k) < - \frac{1}{v_s ^2} k $ for any $ k > 0$. 
\end{Lemma}

{\it Proof. } 
i) We prove that for any $ k > 0$, 
\beq
\label{5.4}
I_{min}(k) \geq \limsup_{h \downarrow k } I_{min}(h).
\eeq
Fix $ \psi \in \Eo $ such that $\ii _{\R^N} |\nabla \psi |^2\, dx = k$. 
At least one of the mappings $ t \longmapsto  \ii_{\R^N} |\nabla \psi_{t,t} |^2\, dx$, 
$ t \longmapsto  \ii_{\R^N} |\nabla \psi_{1,t} |^2\, dx$ or
$ t \longmapsto  \ii_{\R^N} |\nabla \psi_{t,1} |^2\, dx$
is (strictly) increasing on $[1, \infty)$. 
Let $ \psi ^t $ be either $\psi_{t,t}$ or $\psi_{1,t}$ or $\psi_{t,1}$, in such a way that 
$ t \longmapsto  \ii_{\R^N} |\nabla \psi ^{t} |^2\, dx$
is continuous and increasing on $[0, \infty)$. 
It is easy to see that $ I(\psi ^t) \lra I(\psi )$ as $ t \lra 1$. 
Let $(k_n)_{n\geq 1}$ be a sequence satisfying $ k_n \downarrow k$. 
There is a sequence $ t_n \downarrow 1 $ such that $ \ii _{\R^N} |\nabla \psi ^{t_n} |^2 \, dx = k_n$. 
For each $n$ we have $ I_{min}(k_n) \leq I(\psi^{t_n}) $ 
and passing to the limit as $ n \lra \infty $ we find $\ds \limsup_{n \ra \infty}  I_{min}(k_n) \leq I(\psi)$. 
Since this is true for any sequence $ k_n \downarrow k$ and any $ \psi \in \Eo $ satisfying 
$ \ii _{\R^N} |\nabla \psi |^2\, dx = k$, (\ref{5.4}) follows.

Proceeding exactly as in the proof of Lemma \ref{L4.5} (see the proof of (\ref{4.12}) there) we find
\beq
\label{5.5}
I_{min}\left(\frac{k_1 + k_2}{2} \right) \geq \frac 12 (I_{min}(k_1) + I_{min}(k_2)) 
\qquad \mbox{ for any } k_1, \, k_2 >0.
\eeq
Let $ 0 \leq k_1 < k_2$. Using (\ref{5.5}) and a straightforward induction we find 
\beq
\label{5.6}
I_{min} ( \al  k_1 + (1 - \al ) k_2) \geq \al I_{min}(k_1) + (1 - \al ) I_{min}(k_2) 
\qquad \mbox{ for any } \al \in [0,1] \cap \Q.
\eeq
Let $ \al \in (0, 1). $ Consider a sequence $ (\al _n)_{n \geq 1} \subset [0,1] \cap \Q$ 
such that $ \al _n \uparrow \al $. Using (\ref{5.4}) and (\ref{5.6}) we get 
$$
\begin{array}{l}
I_{min} ( \al  k_1 + (1 - \al ) k_2) 
\geq  \ds \limsup_{n \ra \infty} I_{min} ( \al _n k_1 + (1 - \al _n ) k_2) 
\\
\geq \ds \limsup_{n \ra \infty}  \left(\al_n  I_{min}(k_1) + (1 - \al _n ) I_{min}(k_2) \right)
= \al I_{min}(k_1) + (1 - \al ) I_{min}(k_2).
\end{array}
$$
Thus $ I_{min} $ is concave on $[0, \infty)$. Since $I_{min}(0) = 0$, 
by Lemma \ref{L5.2} $I_{min}$ is continuous at $0$ and negative on an 
interval $(0, \de)$ and we infer that $I_{min}$ is negative and decreasing 
on $(0, \infty)$.

The concavity of $I_{min}$ implies that the function 
$k \longmapsto \frac{I_{min}(k)}{k}$ is nonincreasing on $(0, \infty)$. 
Using Lemma 4.4 in \cite{M10}  we find a sequence 
$(\psi_n)_{n \geq 3}\subset \Eo $ such that 
$$
 k_n := \ds \int_{\R^N } |\nabla \psi _n |^2 \, dx \leq C_1 n^{N-2} \ln n, \quad
 \Big\vert \int_{\R^N } V(|\psi _n |^2) \, dx \Big\vert \leq C_2 n^{N-2} \;  \mbox{ and } \; \; 
 Q( \psi _n ) \geq C_3 n^{N-1}, 
$$
where $C_1,\,  C_2, \, C_3 >0$ do not depend on $n$. Then 
${\ds \lim_{k \ra \infty}} \frac{ I_{min}(k) }{k} 
\leq {\ds \lim_{n \ra \infty}} \frac{ I(\psi _n)}{k_n } = - \infty.$

\medskip

(ii) By concavity we have $I_{min}(k_i ) \geq \frac{k_i}{k_1 + k_2}I_{min}(k_1 + k_2), \; i = 1,2$, 
and the subadditivity follows. 

\medskip

(iii) Consider first the case $ N \geq 3$. Fix $ k >0$.
Argue by contradiction and assume that there is a sequence $(\psi_n)_{n\geq 1}\subset \Eo $ such that 
$ \ii _{\R^N} |\nabla \psi _n|^2 \, dx = k $ and 
\beq
\label{5.7}
I(\psi _n) = - Q( \psi _n) + \int_{\R^N} V(|\psi _n|^2) \, dx  \lra - \infty \qquad \mbox{ as } n \lra \infty . 
\eeq
Let $ c = \frac{v_s}{2}.$
By Lemma \ref{L5.1}  there exists $ k _2 > 0 $ such that $ \frac{k_2}{k} < \left( \frac{v_s}{2} \right)^{N-2}$ 
and  (\ref{5.3}) holds for any $ \psi \in \Eo $ with 
$ \ii _{\R^N} |\nabla \psi |^2 \, dx \leq k_2 $.  
Let $ \si = k_2^{\frac{1}{N-2}} k^{-\frac{1}{N-2}} <  \frac{v_s}{2}$.
Then $\ii _{\R^N} |\nabla((\psi_n)_{\si, \si }|^2 \, dx = k_2$, 
hence $(\psi _n)_{\si, \si} $ satisfies (\ref{5.3}), that is
\beq
\label{5.8}
\int_{\R^N} |\nabla \psi _n|^2 \, dx  - \si \frac{v_s}{2} Q( \psi _n) 
+ \si ^2 \int_{\R^N} V(|\psi _n|^2) \, dx \geq 0.
\eeq
From (\ref{5.7}) and (\ref{5.8}) we get 
$$
- \int_{\R^N} |\nabla \psi _n|^2 \, dx 
+ \left(\si \frac{v_s}{2} - \si ^2 \right)\int_{\R^N} V(|\psi _n|^2) \, dx 
\lra - \infty, 
$$
which implies  $  \ii_{\R^N} V(|\psi _n|^2) \, dx \lra - \infty $ as $ n \lra \infty$. 
Since $  \ii _{\R^N} |\nabla \psi _n|^2 \, dx  = k$, this contradicts the first inequality in (\ref{4.1}).

Next assume that $ N =2$ and $V\geq 0 $ on $[0, \infty)$. Fix $ k >0$. 
By Corollary \ref{C4.17} there is $q_k > 0$ such that $ |Q(\psi ) | \leq q_k $ 
for any $ \psi \in \Eo $ satisfying $E(\psi) \leq k+1$. Let $ \psi \in \Eo $ 
be such that $ \ii _{\R^2} |\nabla \psi |^2 \, dx = k$.
If $ \ii _{\R^2} V(|\psi |^2) \, dx = 0$ we infer that $ |\psi | = 1 $ a.e. 
on $ \R^2$ and then (\ref{lift}) implies $Q(\psi ) = 0$, hence $I(\psi)  = 0$. 
If $ \ii _{\R^2} V(|\psi |^2) \, dx > 0$ 
let $ \si = \left(  \ii _{\R^2} V(|\psi |^2) \, dx \right)^{- \frac 12}$ and 
$ \tilde{\psi} = \psi_{\si, \si}$, so that 
$\ii _{\R^2} V(|\tilde{\psi }|^2) \, dx = 1 $ and 
$ \ii _{\R^2} |\nabla \tilde{\psi} |^2 \, dx = k$.
We infer that $ |Q(\tilde{\psi})| \leq q_k $. 
Since $ \psi = \tilde{\psi}_{\frac{1}{\si}, \frac{1}{\si}}$ we have 
by scaling $ I(\psi) = \si^{-2} \ii _{\R^2} V(|\tilde{\psi} |^2) \, dx 
- \si^{-1} Q (\tilde{\psi}) \geq \si^{-2} -  \si^{-1} q_k \geq - \frac{q_k^2}{4}. $ 
We conclude that $I_{min}(k ) \geq - \frac{q_k^2}{4} > - \infty.$

\medskip

iv)  If $V$ achieves negative values, it is easy to see that there exists 
$ \psi _1 \in \Eo $ such that $ \ii _{\R^2} V(|\psi _1 |^2) \, dx < 0$. 
Let $k_1 =  \ii _{\R^2} |\nabla \psi _1|^2 \, dx .$ Then, for any $t > 0 $, 
$ \ii _{\R^2} |\nabla (\psi _1)_{t, t}|^2 \, dx = \ii _{\R^2} |\nabla \psi _1|^2 \, dx 
= k_1 $ because  $N=2$, thus
$$
I_{min}(k_1) \leq I((\psi_1)_{t, t}) = 
- t Q( \psi _1) + t^2 \ds \int_{\R^2} V(|\psi _1 |^2) \, dx \lra - \infty
$$ 
as $t \to  \infty $. By concavity we have $I_{min}(k) = - \infty $ for any $ k \geq k_1$.

\medskip

v) The proof relies on the comparison maps constructed in the proof 
of Theorem \ref{T4.13} from the (KP-I) ground state. Notice first 
that if $ \psi \in \Eo $ is such that $ \ii_{\R^2} |\nabla \psi |^2 \, dx = k$, $\ii_{\R^2} V(|\psi  |^2) \, dx > 0$ 
and $ Q( \psi ) > 0 $, then the function $ t \longmapsto I(\psi_{t,t}) = 
t^2 \ii_{\R^2} V(|\psi  |^2) \, dx - t Q( \psi ) $ achieves its minimum at 
$ t_0 = \frac 12 Q( \psi ) \left( \ii _{\R^2} V(|\psi  |^2) \, dx  \right)^{-1}$. 
Since $ \int_{\R^2} |\nabla \psi |^2 \, dx = \int_{\R^2} |\nabla \psi_{t,t} |^2 \, dx 
= k $ in dimension $N=2$, it follows   that 
\beq
\label{rajout}
I_{min} (k) \leq \inf_{t > 0 } I(\psi_{t,t}) 
= I(\psi_{t_0, t_0}) = - \frac { Q^2( \psi)}{ 4 \int_{\R^2} V(|\psi  |^2) \, dx  } .
\eeq
Fix $ \gamma \not = 0 $ (to be chosen later), and let  $w$ be a ground state 
for (\ref{4.55}). Then, for $\e$ small enough, we have seen in the proof 
of Theorem \ref{T4.13} how to construct from $w$ a comparison map $U_\e \in \Eo $ 
satisfying
$$
\begin{array}{l}
\ds Q( U_{\e}) = 3  v_s ^3 \Sr (w) \e - \frac { 6}{\g }  v_s \Sr (w) \e ^3 ,
\\ \\
\ds \ii_{\R^2} V(|U_\e|^2) \ dx = \frac32  v_s^4 \Sr (w) \e 
- \frac{6}{\gamma}\left(  v_s^2 - \frac43 F''(1)  \right) \Sr (w) \e^3 
+ \mathcal{O}(\e^5) ,
\\ \\ 
\ds \ii_{\R^2} | \nabla U_\e|^2 \ dx = \frac32  v_s^4 \Sr (w) \e 
+  v_s^2 \Sr (w) \left( \frac32 - \frac{12}{\g} \right) \e^3 
+ \mathcal{O}(\e^5) .
\end{array}
$$
Let $k_\e = \ii_{\R^2} | \nabla U_\e|^2 \ dx $. Since $ Q(U_\e) > 0 $ 
and $ \ii_{\R^2} V(|U_\e|^2) \ dx > 0 $ for $\e$ small, we infer from 
(\ref{rajout})  that 
$$
\begin{array}{rl}
\ds I_{min} (k_\e) \leq & \ds - \frac { Q^2( U_\e)}{ 4 \int_{\R^2} V(|U_\e |^2) \, dx  } 
= - \frac32  v_s^2 \Sr (w) \e 
+ \frac{ 4 F''(1) }{\gamma} \Sr (w) \e^3 
+ \mathcal{O}(\e^5) 
\\ \\
= & \ds - \frac{1}{v_s^2} \left[ \frac32  v_s^4 \Sr (w) \e 
- v_s^2  \frac{4F''(1) }{\gamma} \Sr (w) \e^3 
+ \mathcal{O}(\e^5) \right] .
\end{array}$$
Therefore, we have
$$ I_{min} (k_\e) < - \frac{k_\e}{v_s^2} $$
for all $\e $ sufficiently small
provided that $- \frac{4}{\g} F''(1) \geq \frac 32 - \frac{12}{\g}$, that is $ \frac{4( 3 - F''(1))}{\g } > \frac 32$
(take, for instance, $\g = 3 - F''(1)$). 
\hfill
$\Box$

\medskip

Let 
\beq
\label{k}
k_0 = \inf \Big\{ k \geq 0 \; \big| \; I_{min}  (k) < - \frac{1}{v_s ^2} k \Big\} 
\qquad \mbox{ and } \qquad
k_{\infty} = \inf \{ k > 0 \; \big| \; I_{min}  (k) = - \infty \}.
\eeq
By Lemmas \ref{5.2} and \ref{5.3} (i) we have $0 \leq k_0 < \infty$ and 
$ 0 < k_{\infty} \leq \infty$.
It is clear that $ k_0 \leq k_{\infty}$. If either $ N \geq 3$ or $N =2$ 
and $ V \geq 0 $ on $[0, \infty)$ we have $ k_{\infty} = \infty$, while if 
$ N=2$ and (A4) holds with $ F''(1) \not= 3 $, we have $ k_0 = 0$; obviously, in all these cases we have $k_0 < k_{\infty}$.
 The next Lemma gives further 
information in the case when $N = 2$ and $V$ achieves negative values. 
It brings into light  the relationship   between $k_\infty $ and the Dirichlet energy of 
the stationary solutions  of (\ref{1.1}) with minimal energy, the so-called ground 
states or bubbles.

\begin{Lemma}
\label{L5.3bis}
Assume that $N = 2$, (A1), (A2) are satisfied and $\inf V < 0$. Let 
$$
T = \inf \left\{ \int_{\R^2}|  \nabla \psi |^2 \, dx \; \Big\vert  \; \psi \in \Eo, \; 
|\psi | \mbox{ is not constant and } \int_{\R^2 } V(|\psi |^2) \, dx \leq 0 \right\}.$$
Then:

\medskip

(i) We have  $T > 0$ and the infimum is achieved for some $ \psi _0 \in \Eo $. 
Moreover, any such $ \psi _0 $ satisfies the equation $ \Delta \psi _0 + \si ^2 F(|\psi _0|^2) \psi _0 = 0 $ in $ \Do '(\R^2)$
for some $ \si > 0$, $ \int_{\R^2} V(|\psi _0|^2) \, dx = 0$, 
$\psi _0 $ belongs to $C^{1, \al }(\R^2)$ for any $ \al \in (0, 1)$ and, after a translation, 
$\psi _0$ is radially symmetric.

\medskip

(ii) For any $k < T $ and any $M >0$, $E_{GL}$ is bounded on the set 
$$
\Eo_{k, M} := \left \{ \psi \in \Eo \; \Big\vert \; \int_{\R^2} |\nabla \psi |^2 \, dx \leq k, \; 
\int_{\R^2}  V(|\psi |^2) \, dx \leq M \right\}.
$$

(iii) We have $ k_{\infty } = T.$

\end{Lemma}

{\it Proof.} 
(i) It follows from Corollary \ref{C4.1} that $ T>0$. 
The proof of the existence and regularity of  minimizers is rather classical and is similar to the proof of 
Theorem 3.1 p. 106 in \cite{brezis-lieb}, so we omit it. 
This is also an immediate consequence of Theorem 1.2 in \cite{MN}.
Notice that any minimizer of the considered problem is also a minimizer of 
$\int_{\R^2} V(|\psi |^2) \, dx $
under the constraint $\int_{\R^2} |\nabla \psi |^2 \, dx = T$ and then the radial symmetry follows from Theorem 
2 p. 314 in \cite{M7}.

\medskip

(ii)  
Fix $ \beta \in (0, 1]$ such that 
\beq
\label{5.13}
V(s^2) \geq \frac{1}{4} ( s^2 - 1)^2 \qquad \mbox{ for any } s \in ( ( 1 - \beta)^2, ( 1 + \beta )^2). 
\eeq

It suffices to prove that for any sequence $(\psi _n ) _{n \geq 1 } \subset \Eo _{k, M}$, 
$E_{GL}(\psi _n)$ is bounded. Let $(\psi _n ) _{n \geq 1 } \subset \Eo _{k, M}$.
Let $ K_n = \{ x \in \R^2 \; \big | \; \; \big| \, |\psi _n (x)| - 1 \big| \geq \frac{ \beta}{2} \}.$
We claim that it suffices to prove that $ \Lo ^2( K_n)$ is bounded. 

Indeed, assume $ \Lo ^2( K_n)$  bounded. 
Let $ \tilde{\psi} _n = \left( \big| \, |\psi _n | - 1 \big| - \frac{ \beta}{2} \right)_+.$
Then $ \tilde{\psi }_n \in L_{loc}^1 (\R^2),$  $|\nabla \tilde{\psi }_n  | \leq | \nabla \psi _n|$ 
a.e. on $ \R^2$ and by (\ref{ineq2}) we have
$$
\int_{\R^2} |\tilde{\psi}_n|^{2 p_0 + 2} \, dx 
\leq C_{2 p_0 + 2} \| \nabla \tilde{\psi } _n \|_{L^2(\R^2)}^{2 p_0 + 2} \Lo ^2 ( K_n) 
\leq C_{2 p_0 + 2} \| \nabla {\psi } _n \|_{L^2(\R^2) }^{2 p_0 + 2} \Lo ^2 ( K_n) .
$$
By (A1) and (A2) there is $C_0 >0$ such that $|V(s^2) | \leq C_0 \left( | s - 1 | - \frac{\beta}{2} \right)^{2 p_0 + 2}$ 
for any $ s $ satisfying $|s - 1 | \geq \beta$. Hence
$$
\int_{\R^2 \setminus \{ 1 - \beta \leq |\psi _n | \leq 1 + \beta \} }
| V(|\psi _n |^2) | \, dx  
\leq C_0 \int_{\R^2} |\tilde{\psi}_n|^{2 p_0 + 2} \, dx  \leq C_0 C_{2p_0 + 2} 
 \| \nabla {\psi } _n \|_{L^2 (\R^2) }^{2 p_0 + 2} \Lo ^2 ( K_n) 
$$
and the last quantity is bounded. Since $ \int_{\R^2} V(|\psi _n|^2) \, dx $ is bounded, 
we infer that \\
$  \int_{\{ 1 - \beta \leq |\psi _n | \leq 1 + \beta \} } V(|\psi _n|^2) \, dx $ is bounded, 
and by (\ref{5.13}), 
$  \int_{\{ 1 - \beta \leq |\psi _n | \leq 1 + \beta \} } \left( \ph ^2( |\psi_n |) - 1 \right)^2 \, dx $ is bounded. 
On the other hand, 
$  \int_{\R^2 \setminus \{ 1 - \beta \leq |\psi _n | \leq 1 + \beta \} } 
\left( \ph ^2( |\psi_n |) - 1 \right)^2 \, dx \leq \int_{ K_n}\left( \ph ^2( |\psi_n |) - 1 \right)^2 \, dx \leq 
 64 \Lo ^2 ( K_n)$ 
and the conclusion follows.

\medskip

It remains to prove the boundedness of $ \Lo ^2 ( K_n)$. Let 
$$
\psi _n ^+ = \left\{ 
\begin{array}{ll} 
|\psi _n | & \quad \mbox{ if } |\psi _n | \geq 1 
\\
1 & \quad \mbox{ otherwise, }
\end{array}
\right. 
\qquad \qquad \mbox{and} 
\qquad \qquad
\psi _n ^- = \left\{ 
\begin{array}{ll} 
|\psi _n |  & \quad \mbox{ if } |\psi _n | \leq 1
\\
1  & \quad \mbox{ if } |\psi _n | \geq 1.
\end{array}
\right.
$$
It is clear that $ \psi_n^+, \psi _n ^- \in \Eo$, 
$ \int_{\R^2} |\nabla \psi _n ^+|^2 + |\nabla \psi _n ^-| ^2 \, dx 
=  \int_{\R^2} |\nabla | \psi _n | \,|^2 \, dx \leq k $ and 
$  \int_{\R^2} V(|\psi _n ^+|^2) + V(|\psi _n ^-|^2) \, dx = \int_{\R^2} V(|\psi _n |^2)\, dx.$
If $  \int_{\R^2} V(|\psi _n ^+|^2)  \, dx <0$, 
by (i) we have $\int_{\R^2} |\nabla \psi _n^+|^2 \, dx \geq T > k$, 
 a contradiction. 
Thus necessarily  $  \int_{\R^2} V(|\psi _n ^+|^2)  \, dx \geq 0$
and similarly $  \int_{\R^2} V(|\psi _n ^-|^2)  \, dx \geq 0$, hence
$  \int_{\R^2} V(|\psi _n ^{\pm}|^2)  \, dx \in [0, M].$

Let $K_n ^+ = \{ x \in \R^2 \; \big| \; \; |\psi _n (x) | \geq 1 + \frac{ \beta}{2} \}, $
$K_n ^- = \{ x \in \R^2 \; \big| \; \; |\psi _n (x) | \leq 1 - \frac{ \beta}{2} \}. $
Let $w_n ^+ = \phi_n^+(|x|) $ and $w_n ^- = \phi_n^-(|x|) $ be the symmetric decreasing rearrangements of 
$(|\psi_n| - 1 )_+ = \psi_n ^+ - 1$ and of $(|\psi_n| - 1 )_ -= 1 - \psi _n ^-  $, respectively.
As in the proof of Lemma \ref{L4.6} we have $ \phi_n^{\pm} \in H_{loc}^1( (0, \infty)).$ Let 
$$
\begin{array}{c}
t_n  = \inf \{ t \geq 0 \; | \; \phi_n ^+(t)< \frac{ \beta}{2} \} \qquad \mbox{ and } \qquad
s_n = \inf \{ t \geq 0 \; | \; \phi_n ^- (t)< \frac{ \beta}{2} \} .
\end{array}
$$
Then $ \Lo ^2 (K_n^+) = \Lo ^2 ( \{ (|\psi_n| - 1 )_+  \geq \frac{\beta}{2} \} )
= \Lo ^2 ( \{ w_n ^+ \geq \frac{\beta}{2} \} ) = \Lo ^2 ( \ov{B} (0, t_n )) = \pi t_n ^2$ 
and similarly $  \Lo ^2 (K_n^-) =  \pi s_n ^2$, so that $ \Lo ^2( K_n) = \pi( t_n^2 + s_n ^2)$.

Assume that there is a subsequence $t_{n_j} \lra \infty$. 
Let $ \tilde{w}_j = (w_{n_j}^+)_{\frac{1}{t_{n_j}}, \frac{1}{t_{n_j}} } 
= \phi_{n_j} ^+ ( t_{n_j} |\cdot|), $
so that $\tilde{w}_j \geq \frac{\beta}{2}$ on $ \ov{B} (0, 1)$ and 
$ 0 \leq \tilde{w}_j  < \frac{\beta}{2}$ on $ \R^2 \setminus \ov{B} (0, 1)$.
Then $ \int_{\R^2} |\nabla \tilde{w}_j  |^2 \, dx 
= \int_{\R^2} |\nabla {w}_{n_j}  |^2 \, dx \leq k $ 
and using (\ref{ineq2}) we see that $(\tilde{w}_j  -  \frac{\beta}{2} )_+$ 
is uniformly bounded in $L^p( B(0,1)) $ for any $ p < \infty$, and consequently 
 $(\tilde{w}_j)_{j \geq 1} $ is bounded in $L^p(B(0, R))$ 
 for any $ p < \infty $ and any $ R \in(0, \infty)$. 
Then there is a subsequence of $(\tilde{w}_j )_{j \geq 1}$, still denoted $(\tilde{w}_j )_{j \geq 1}$,
and there is $ \tilde{w} \in H_{loc}^1(\R^2)$ such that $ \nabla \tilde{w} \in L^2 (\R^2)$ and 
$(\tilde{w}_j )_{j \geq 1}$, $\tilde{w}$ satisfy (\ref{4.32}).
It is easy to see that $ 1 + \tilde{w} \in \Eo $ and $ \tilde{w} \geq \frac{\beta}{2}$ on $\ov{B}(0,1).$
By weak convergence we have 
$$
\tilde{k} := \ds \int_{\R^2} |\nabla \tilde{w} |^2\, dx \leq \liminf_{ j \ra \infty} 
 \int_{\R^2} |\nabla \tilde{w} _j |^2\, dx \leq k.
$$
Using (A2),  the convergence $ \tilde{w} _j \lra \tilde{w}$ in $L^{2 p_0 + 2}(B(0, 1))$ 
 and Theorem A2 p. 133 in \cite{willem} we get 
 $ \int_{B(0,1)} V((1 + \tilde{w}_j )^2) \, dx \lra \int_{B(0,1)} V((1 + \tilde{w} )^2) \, dx.$
Since $ \tilde{w}_j \in [0, \frac{\beta}{2}]$ on $\R^2 \setminus {B}(0,1) $ and 
$ V( s^2 ) \geq 0$ for $ s \in [1 , 1 + \frac{\beta}{2} ]$, 
using Fatou's Lemma we obtain 
$  \int_{\R^2 \setminus B(0,1)} V((1 + \tilde{w} )^2) \, dx
\leq {\ds \liminf_{j \ra \infty}} \int_{\R^2 \setminus B(0,1)} V((1 + \tilde{w}_j )^2) \, dx.$
Therefore
$$
\begin{array}{l}
\ds \int_{\R^2} V((1 + \tilde{w} )^2) \, dx \leq 
\liminf_{j \ra \infty} \int_{\R^2 } V( (1 + \tilde{w}_j )^2) \, dx
= \liminf_{j \ra \infty} \frac{1}{t_{n_j} ^2} \int_{\R^2 } V((1 + {w}_{n_j}^+ )^2) \, dx
\\
= \ds \liminf_{j \ra \infty} \frac{1}{t_{n_j} ^2} \int_{\R^2 } V((\psi_{n_j}^+ )^2) \, dx \leq 0 
\end{array}
$$
because  $\int_{\R^2 } V((\psi_{n_j}^+ )^2) \, dx \leq M$ and $t_{n_j} \to +\infty$ 
by our assumption. Since $1 + \tilde{w} \geq 1 + \frac{ \beta}{2}$ on $B(0,1)$ 
we infer that $ \int_{\R^2} |\nabla \tilde{w} |^2 \, dx \geq T > k$, a contradiction. 

So far we have proved that $(t_n)_{n \geq 1}$ is bounded. Similarly $(s_n)_{n \geq 1}$ is bounded, 
thus $(\Lo ^2( K_n))_{n \geq 1}$ is bounded and the proof of (ii) is  complete.

\medskip

(iii) Consider a radial function $ \psi _0 \in \Eo $ such that 
$|\psi _0|$ is not constant, $\int_{\R^2} V(|\psi _0|^2) \, dx = 0 $ and 
$\int_{\R^2} |\nabla \psi _0 |^2 \, dx = T$. 
Since $F(|\psi _0 |^2  ) \psi _0$ does not vanish a.e. on $ \R^2$, 
 there exists a radial function $ \phi \in C_c^{\infty}(\R^2)$ such that 
$\int_{\R^2} \langle F(|\psi _0 |^2 ) \psi _0 , \phi \rangle \, dx > 0$. 
It follows that $\frac{d}{dt} _{| t =0} \int_{\R^2} V(|\psi _0 + t \phi |^2) \, dx 
= - 2 \int_{\R^2} \langle F(|\psi _0 |^2 \psi _0 , \phi \rangle \, dx  <0$, 
consequently there is $ \e > 0$ such that 
$ \int_{\R^2} V(|\psi _0 + t \phi |^2) \, dx  <  \int_{\R^2} V(|\psi _0  |^2) \, dx = 0$ for any $ t \in (0, \e)$.
Denote $ k(t) = \int_{\R^2} |\nabla ( \psi _0 + t \phi ) |^2 \, dx.$
It follows from  the proof of Lemma \ref{L5.3} (iv) that
$I_{min}(k(t) ) = - \infty $ for any $ t \in (0, \e)$, thus 
$k_\infty \leq k(t) $ for any $ t \in (0, \e)$. Since $ k(t) \lra T $ 
as $ t \lra 0$, we infer that $ k_{\infty} \leq T$.

Let $ k < T$. 
Consider $ \psi \in \Eo $ such that $\int_{\R^2} |\nabla \psi |^2 \, dx = k$. 
If $|\psi | = 1 $ a.e. we have $V(|\psi |^2) = 0 $ a.e. and 
$ Q( \psi ) = 0 $ by (\ref{lift}), hence $I(\psi) = 0$. 
If  $|\psi |$ is not constant,  then   we have necessarily $\int_{\R^2} V(|\psi|^2) \, dx > 0$. 
If $Q( \psi ) \leq 0$, it is obvious that $I(\psi ) > 0$. 
If $ Q(\psi ) > 0$ we have ${\ds \inf_{t > 0}} I(\psi_{t,t}) 
= - \frac 14 Q^2(\psi) \left( \int_{\R^2} V(|\psi|^2) \, dx \right)^{-1} $ and the infimum is achieved for 
$t_{min} = \frac 12 Q(\psi) \left( \int_{\R^2} V(|\psi|^2) \, dx \right)^{-1} $.
There exists $ t_1 > 0$ such that $\int_{\R^2} V(|\psi _{t_1, t_1}|^2) \, dx = 1$. Then 
$\int_{\R^2} |\nabla \psi _{t_1, t_1 } |^2 \, dx = k$ and
$$
I(\psi) \geq \inf_{t > 0} I(\psi_{t,t})  
 = - \frac {Q^2(\psi) }{4 \int_{\R^2} V(|\psi|^2) \, dx }
 =  - \frac { Q^2(\psi_{t_1, t_1}) }{ 4\int_{\R^2} V(|\psi_{t_1, t_1}|^2) \, dx }
 = - \frac 14 Q^2(\psi_{t_1, t_1}).
$$
This implies $I(\psi) \geq \inf \{ - \frac 14 Q^2( \phi ) \; | \; \phi \in \Eo _{k, 1} \}.$
By (ii) we know that $E_{GL}$ is bounded on $\Eo _{k, 1}$ and Corollary \ref{C4.17} implies that $Q$ 
is also bounded on $\Eo _{k, 1}$. We conclude that $ I_{min}(k) > - \infty$, hence $ k < k_{\infty }$. 
Since this is true for any $ k < T$, we infer that $ k_{\infty } \geq T$. 
Thus $ k_{\infty} = T$. 
\hfill
$\Box$

\begin{Lemma}
\label{L5.4}
Assume that $ 0 < k < k_{\infty }$ and $(\psi_n)_{n \geq 1} \subset \Eo $ is a sequence such that 
$ \int_{\R^N} |\nabla \psi _n |^2 \, dx \leq k $ for all $n$.
Suppose that  $(I(\psi_n))_{n \geq 1}$ is bounded in the case $N \geq 3$, 
respectively that $I(\psi _n ) < 0$ for all $n$ in the case $N = 2$.

Then $(Q(\psi_n))_{n \geq 1}$,  $ \left( \int_{\R^N} V(|\psi _n |^2) \, dx  \right)_{n \geq 1}$ and 
$(E_{GL}(\psi_n))_{n \geq 1}$  are bounded.
\end{Lemma}

{\it Proof.} 
Consider first the  case $ N \geq 3$. 
Let us show  that $  \int_{\R^N} V(|\psi _n |^2) \, dx $ is bounded from above.
We argue by contradiction and assume that this is false. 
Then there is a subsequence, still denoted $(\psi_n)_{n \geq 1}$, such that 
$ s_n : =  \int_{\R^N} V(|\psi _n |^2) \, dx \lra \infty $ as $ n \lra \infty$.

Let $ \si _n = s_n ^{- \frac 1N}$. Since 
$  \int_{\R^N}|\nabla ((\psi _n )_{\si _n, \si _n}) |^2 \, dx 
= \si_n ^{N-2} \int_{\R^N}|\nabla \psi _n  |^2 \, dx\lra 0 $ as $ n \lra \infty$,
 Lemma \ref{L5.1} implies that $(\psi_n)_{\si_n, \si _n}$ satisfies (\ref{5.3}) with $ c = \frac{v_s}{2}$ 
 for all sufficiently large $n$,  that is 
$$
\int_{\R^N}|\nabla \psi _n  |^2 \, dx - \frac{v_s}{2} s_n ^{- \frac 1N} ( s_n - I(\psi _n)) + s_n^{ 1 - \frac{2}{N}} \geq 0.
$$
Since $ \int_{\R^N}|\nabla \psi _n  |^2 \, dx $ and $ I(\psi _n)$ are bounded and $ s_n \lra \infty$, 
the left-hand side of the above inequality tends to $- \infty $ as $ n \lra \infty$, a contradiction. 
We conclude that there is $M > 0 $ such that $ \int_{\R^N} V(|\psi _n |^2) \, dx  \leq M$ for all $M$. 
Then (\ref{4.1}) implies that $  \ii_{\R^N} \left( \ph ^2( |\psi _n|) - 1 \right)^2 \, dx $ is bounded.
By (\ref{4.1}), $\ii _{\R^N} V(|\psi|^2) \, dx $ is bounded from below. 
 Using Corollary \ref{C4.17} we infer that $(Q(\psi_n))_{n \geq 1}$ is bounded.

\medskip

Consider next the case $ N =2$. 
Since $ \ii _{\R^2}|\nabla \psi _n  |^2 \, dx \leq k < k_{\infty}$, 
using Lemma \ref{L5.3bis} (i) and (iii) 
we see that either   $ \ii _{\R^N} V(|\psi _n |^2) \, dx > 0$ 
or $ \ii _{\R^N} V(|\psi _n |^2) \, dx = 0 $ and $ |\psi _n | = 1 $ a.e. on $\R^2$. 
In the latter case (\ref{lift}) implies $ Q( \psi _n) = 0$, hence $ I( \psi _n) = 0$, 
contrary to the assumption that $  I( \psi _n ) < 0$.
Thus necessarily $ 0<  \ii _{\R^N} V(|\psi _n |^2) \, dx < Q( \psi _n) $ for all $n$ because 
$ I( \psi _n) < 0$. 

Since $\ii _{\R^2} |\nabla (\psi _n) _{\si, \si}  |^2 \, dx = \ii _{\R^2} |\nabla \psi _n|^2 \, dx $
for any $ \si >0$, as in the proof of Lemma \ref{L5.3bis} (iii) we have
$$
 - \frac {Q^2(\psi _n ) }{4 \int_{\R^2} V(|\psi _n |^2) \, dx }
= \inf_{\si > 0} I((\psi_n)_{\si,\si})  \geq I_{min} 
\left( \int_{\R^2} |\nabla \psi _n |^2 \, dx  \right) \geq I_{min}(k)
$$
and this implies
$$
 Q^2 ( \psi _n)  \leq -4 I_{min}(k) \int_{\R^N} V(|\psi _n |^2) \, dx .
$$
Combining this with the inequality $ 0<  \int_{\R^2} V(|\psi _n |^2) \, dx < Q( \psi _n) $,  we get 
\beq
\label{5.12}
0< \ds \int_{\R^2} V(|\psi _n |^2) \, dx < Q( \psi _n) \leq -4 I_{min}(k).
\eeq
We have thus proved that
$(Q(\psi_n))_{n \geq 1}$ and  $ \left( \int_{\R^2} V(|\psi _n |^2) \, dx  
\right)_{n \geq 1}$ are bounded.
The boundedness of $ \int_{\R^2} \left( \ph ^2( |\psi_n |) - 1 \right)^2 \, dx $ 
follows from Lemma \ref{L4.6}  if  $V \geq 0 $ on $[0, \infty)$, 
respectively from Lemma \ref{L5.3bis} (ii) if $V$ achieves negative values. 
\hfill
$\Box$

\medskip

We now state the main result of this section, which shows precompactness 
of minimizing sequences for $I_{min}(k)$ as soon as $ k_0 < k < k_\infty $.

\begin{Theorem}
\label{T5.5}
Assume that $N \geq 2$ and (A1), (A2) hold. Let $ k \in ( k_0 , k_{\infty})$ and let 
$(\psi_n )_{n \geq 1} \subset \Eo $ be a sequence such that 
$$
\int_{\R^N} |\nabla \psi _n |^2 \, dx \lra k \qquad \mbox{ and } \qquad I(\psi _n) \lra I_{min}(k).
$$

There exist a subsequence $(\psi_{n_k})_{k \geq 1}$, a sequence of points $(x_k)_{k \geq 1}  \subset \R^N$, 
and $ \psi \in \Eo $ such that $ \int_{\R^N} |\nabla \psi  |^2 \, dx = k$, 
$I( \psi ) = I_{min}(k) $, $\psi_{n_k} ( x_k + \cdot ) \lra \psi $ a.e. on $\R^N$  and 
$$
\| \nabla \psi_{n _k} ( \cdot + x_k ) - \nabla \psi \|_{L^2(\R^N)} \lra 0 , 
\qquad 
\| \,  | \psi_{n _k} |( \cdot + x_k ) - | \psi  | \,  \|_{L^2(\R^N)} \lra 0 
\qquad \mbox{ as } k \lra \infty.
$$
\end{Theorem}

{\it Proof. }  
Since $I_{min}(k) < 0$, we have $I(\psi_n) < 0$ for all sufficiently large $n$. 
By Lemma \ref{L5.4} the sequences 
$(Q(\psi_n))_{n \geq 1}$,  $ \left( \int_{\R^N} V(|\psi _n |^2) \, dx  \right)_{n \geq 1}$ and 
$(E_{GL}(\psi_n))_{n \geq 1}$  are bounded. Passing to a subsequence 
if necessary, we may assume that $E_{GL}(\psi _n) \lra \al _0 \geq k > 0$ 
and  $Q( \psi _n) \lra q $ as $ n \lra \infty$.  

We use the Concentration-Compactness Principle  (\cite{lions}) 
and we argue  as in the proof of Theorem \ref{T4.7}.
Let $\Lambda_n (t)$ be the concentration function associated to $E_{GL}(\psi _n) $, as in (\ref{4.15}).
It is standard to prove  that there exist a subsequence of 
$((\psi_n, \Lambda_n))_{n \geq 1}$, still denoted $((\psi_n, \Lambda_n))_{n \geq 1}$, 
a nondecreasing function  $ \Lambda : [0, \infty ) \lra \R$, $ \al \in [0, \al _0 ]$, 
and a nondecreasing sequence $t_n \lra \infty$ such that (\ref{4.16}) and (\ref{annulus}) hold.  
The next result implies that $ \al >0$.

\begin{Lemma}
\label{L5.6}
Let $(\psi_n) _{n \geq 1} \subset \Eo $ be a sequence satisfying:

\medskip

(a) $   E_{GL}(\psi _n) \leq M $ for some positive constant $M$. 

\medskip

(b) $ \int_{\R^N} |\nabla \psi _n |^2\, dx \lra k $ and $ Q(\psi _n ) \lra q $ 
as $ n \lra \infty.$

\medskip

(c) ${\ds \limsup_{n \ra \infty }}\,  I (\psi_n) < - \frac{1}{v_s ^2} k$. 
\\
Then there exists $ \ell > 0$ such that 
$ \ds \sup_{ y \in \R^N } E_{GL}^{B(y,1)} ( \psi _n ) \geq \ell $ 
for all sufficiently large $n$. 
\end{Lemma}

{\it Proof. } 
It is obvious that the sequence $(\psi_n) _{n \geq 1}$ satisfies the 
conclusion of Lemma \ref{L5.6} if and only if 
$((\psi_n)_{v_s, v_s}) _{n \geq 1}$ satisfies the same conclusion.

By (a) we have $E_{GL}((\psi_n)_{v_s, v_s}) \leq \max( v_s^{N-2}, v_s^N) M = 2^{\frac N2} M. $
Assumption (b) implies 
$$\ds \int_{\R^N} |\nabla (\psi _n)_{v_s, v_s} |^2\, dx \lra v_s^{N-2}k \qquad \mbox{  and } 
\qquad  Q((\psi _n) _{v_s, v_s} ) \lra  v_s^{N-1} q = \tilde{q} .
$$
Using (c) we find
$$
\begin{array}{l}
\ds \limsup_{n \ra \infty}  E((\psi _n)_{v_s, v_s}) - v_s \tilde{q} 
= \ds \limsup_{n \ra \infty}  \left( \! v_s^{N-2} \! \! \int_{\R^N} |\nabla \psi _n |^2\, dx 
+ v_s^N \! \! \int_{\R^N} V(|\psi _n|^2) \, dx - v_ s ^N Q( \psi _n) \! \right) 
\\
= v_s ^N \left( \frac{1}{v_s ^2} k  + \ds \limsup_{n \ra \infty}  I( \psi _n )  \right) < 0.
\end{array}
$$
Then the result follows directly from Lemma \ref{L4.8}.
\hfill 
$\Box $

\medskip

Next we prove that $ \al \not\in (0, \al _0)$. 
We argue again by contradiction and we assume that 
$ 0 < \al < \al _0$. 
Arguing as in the proof of Theorem \ref{T4.7} and using Lemma \ref{splitting} 
for each $n $ sufficiently large we construct two functions 
$\psi_{n,1}, \psi_{n,2} \in \Eo $ such that 
\beq
\label{5.14}
E_{GL}(\psi _{n, 1} ) \lra \al 
\qquad 
\mbox{ and } 
\qquad
E_{GL}(\psi _{n, 2} ) \lra \al _0 - \al , 
\eeq
\vspace*{-15pt}
\beq
\label{5.15}
\int_{\R^N} \big\vert |\nabla \psi_{n} |^2 - |\nabla \psi_{n, 1} |^2 - |\nabla \psi_{n, 2} |^2 \big\vert \, dx \lra 0, 
\eeq
\vspace*{-15pt}
\beq
\label{5.16}
\int_{\R^N} \big\vert V (| \psi_{n} |^2) - V( | \psi_{n, 1} |^2 ) - V( | \psi_{n, 2} |^2 ) \big\vert \, dx \lra 0, 
\eeq
\vspace*{-15pt}
\beq
\label{5.17}
| Q(\psi _n) -  Q(\psi _{n, 1}) -  Q(\psi _{n, 2}) | \lra 0 \qquad \mbox{ as } n \lra \infty. 
\eeq
Passing to a subsequence if necessary, we may assume that 
$ \ii _{\R^N} |\nabla \psi_{n, i} |^2 \, dx \lra k_i \geq 0 $ 
as $ n \lra \infty$, $ i=1,2$. By (\ref{5.15}) we have $ k_1 + k_2 = k$. We claim that $ k_1 > 0 $ and $ k_2 >0$. 

To prove the claim assume, for instance, that $ k_1 = 0$. From (\ref{5.14}) it follows that \\
$ \frac 12  \ii _{\R^N} \left( \ph^2 (|\psi_{n,1}|) - 1 \right)^2 \, dx \lra \al.$
Using Lemma \ref{L4.1}   we find $  \ii _{\R^N}  V( | \psi_{n, 1} |^2 ) \, dx \lra \al $. 
From Lemma \ref{L4.2} (ii)  we infer that there is $ \kappa > 0$ 
such that 
$
E(\psi) \geq \frac{ v_s}{2} | Q( \psi )| 
$
 for any $ \psi \in \Eo $ satisfying $ E_{GL}(\psi ) \leq \kappa.$
It is clear that there are $ n_0 \in \N$ and  $ \si _0 > 0$ such that 
$E_{GL}( (\psi_{n,1})_{\si , \si }) \leq \kappa $ for any $ n \geq n_0 $ and any $ \si \in (0, \si _0]$. 
Then $E((\psi_{n,1})_{\si , \si }  ) \geq 
\frac{ v_s}{2} | Q( (\psi_{n,1})_{\si , \si }) |$, that is
$$
\frac{ v_s}{2} | Q( \psi_{n,1}) | \leq \frac{1}{\si } \int_{\R^N} |\nabla \psi_{n, 1} |^2 \, dx 
   +   \si \int_{\R^N}  V (|\psi_{n,1}| ^2)  \, dx 
$$
 for any $n \geq n_0 $  and  $ \si \in (0, \si_0].$
Passing to the limit as $ n \lra \infty $ in the above inequality we discover
$\frac{v_s}{2} \ds \limsup_{n \ra \infty} | Q( \psi_{n,1}) | \leq \si \al $ for any $ \si \in (0, \si _0]$, 
that is $\ds \lim_{n \ra \infty} | Q( \psi_{n,1}) | =0$. 
As a consequence we find $\ds \lim_{n \ra \infty}  I( \psi_{n,1}) = \al$. 
Since $|I(\psi _n ) -  I(\psi _{n, 1}) -  I(\psi _{n, 2}) | \lra 0 $ by (\ref{5.16}) and (\ref{5.17}), 
we infer that $I(\psi _{n, 2}) \lra I_{min}(k) - \al$ as $ n \lra \infty$. 
Since $ \ii _{\R^N}  |\nabla \psi_{n, 2} |^2 \, dx \lra k_2 = k$, 
this contradicts the definition of $I_{min}$ and the fact that $I_{min}$ is continuous at $k$. 
Thus necessarily $ k_1 >0$. Similarly we have $ k_2 >0$, that is $ k_1, \, k_2 \in (0, k)$. 

\medskip

We have $I(\psi_{n,i}) \geq I_{min} ( \int_{\R^N} |\nabla \psi_{n, i} |^2 \, dx )$ and passing to the limit we get 
$\ds \liminf_{n \ra \infty} I(\psi_{n,i}) \geq I_{min} (k_i )$, $ i =1,2$.
Using (\ref{5.16}), (\ref{5.17}) and the fact that $I(\psi _n) \lra I_{min} (k)$ 
we infer that $I_{min} (k) \geq I_{min} (k_1) + I_{min} (k_2)$.
On the other hand, the concavity of $I_{min} $ implies $I_{min}(k_i) \geq \frac{k_i}{k}  I_{min} (k)$, 
hence $I_{min}(k_1) + I_{min}(k_2) \geq I_{min}(k) $ and equality may occur if and only if $I_{min}$
is linear on $[0, k]$. Thus there is $ A \in \R$ such that $I_{min}(s) = As $ for any $ s \in [0, k]$. 
By Lemma \ref{L5.2} we have $ A = - \frac{1}{v_s ^2}$, hence $I_{min}(k) = - \frac{k}{v_s ^2}$, 
contradicting the fact that $ k > k_0$. 
Thus we cannot have $ \al \in (0, \al _0)$, and then necessarily $ \al = \al _0$.

\medskip

As in the proof of Theorem \ref{T4.7}, 
 there  is a sequence $(x_n)_{n \geq 1} \subset \R^N$ such that 
for any $ \e > 0 $ there is  $ R_{\e} > 0 $  satisfying 
$E_{GL} ^{ \R^N\setminus B(x_n, R_{\e})} (\psi _n)  < \e$ for all  $ n$ sufficiently large. 
Let $ \tilde{\psi}{_n} = \psi _n ( \cdot + x_n)$. Then 
for any $ \e > 0 $ there exist   $ R_{\e} > 0 $ and $ n_{\e } \in \N$ such that 
$(\tilde{\psi}_n)_{n \geq 1}$ satisfies (\ref{4.31}).
It is  standard to prove that there is  a function $ \psi \in H_{loc}^1(\R^N)$
such that $ \nabla \psi \in L^2(\R^N)$ and a subsequence $ (\tilde{\psi}_{n _j})_{j \geq 1} $
 satisfying (\ref{4.32})-(\ref{4.34})  and (\ref{4.37}).

\medskip

 Lemmas \ref{L4.10} and \ref{L4.11} imply that 
$ \| \, |\tilde{\psi} _{n _j} | - |\psi | \, \|_{L^2( \R^N) } \lra 0 $, 
$Q (\tilde{\psi} _{n_j})  \lra Q( \psi) $  and
$  \ii _{\R^N} V(|\tilde{\psi }_{n_j} |^2) \, dx  \lra  \ii _{\R^N} V(|\psi  |^2) \, dx $ 
as $ j \lra \infty. $ 
Therefore  $I (\tilde{\psi} _{n_j}) \lra I( \psi) $, and consequently $I(\psi ) = I_{min}(k).$
On the other hand, by (\ref{4.33}) we have $\ii _{\R^N} |\nabla \psi |^2 \, dx \leq k$. 
Since $I_{min}$ is strictly decreasing, we infer that necessarily  
$\ii _{\R^N} |\nabla \psi |^2 \, dx =  k 
= {\ds \lim_{j \ra \infty}} \ii _{\R^N} |\nabla \tilde{\psi}_{n_j} |^2 \, dx.$
Combined with the weak convergence $ \nabla \tilde{\psi}_{n_j} \rightharpoonup \nabla \psi$ in 
$L^2( \R^N)$, this gives the strong convergence 
$ \nabla \tilde{\psi}_{n_j} \lra \nabla \psi$ in $L^2( \R^N)$
and the proof of Theorem \ref{T5.5} is complete. 
\hfill
$\Box $

\medskip

Denote by $ d^- I_{min}(k) $ and $ d^+ I_{min}(k) $ the left and right 
derivatives of $I_{min}$ at $ k>0 $ (which exist and are finite for 
any $ k > 0 $ because $I_{min}$ is concave). We have: 

\begin{Proposition}
\label{P5.7}
(i) Let $ c > 0 $. Then the function $ \psi  $ is a minimizer of $I$ 
in the set $ \{ \phi \in \Eo \; | \;  \int_{\R^N} |\nabla \phi |^2 \, dx = k \}$ 
if and only if $ \psi_{c,c}$ minimizes the functional 
$$ I_c (\phi) = - cQ ( \phi ) +  \int_{\R^N} V(|\phi|^2) \, dx $$ 
in the set $ \{ \phi \in \Eo \; | \;  \int_{\R^N} |\nabla \phi |^2 \, dx = c^{N-2} k \}$.

\medskip

(ii) If $ \psi \in \Eo $ satisfies $  \int_{\R^N} |\nabla \psi |^2 \, dx = k $ 
and $ I(\psi) = I_{min}(k)$, there is $ \vartheta \in [d^+ I_{min}(k),  d^- I_{min}(k) ] $ 
such that 
\beq
\label{5.18}
i  \psi _{x_1} - \vartheta \Delta \psi + F(|\psi |^2) \psi = 0 \qquad 
\mbox{ in } \Do' ( \R^N).
\eeq
Then for $ c = \frac{1}{\sqrt{- \vartheta}}$ the function $ \psi _{c,c}$ 
satisfies (\ref{4.50}) and minimizes $E_c = E - cQ $ in the set 
$ \{ \phi \in \Eo \; | \;  \int_{\R^N} |\nabla \phi |^2 \, dx = c^{N-2} k \}$.
Moreover,  $ \psi \in W_{loc}^{2, p } ( \R^N)$ and 
$\nabla  \psi \in W^{1, p } ( \R^N)$ for any $ p \in [2, \infty)$. 

\medskip

(iii) After a translation, $ \psi $ is axially symmetric with respect 
to the $ x_1-$axis if $ N \geq 3$. The same conclusion is true if 
$N=2$ and we assume in addition that $F$ is $C^1$. 

\medskip

(iv) For any $ k \in (k_0, k_{\infty}) $ there are $ \psi ^+, \psi ^- \in \Eo $ 
such that 
$\int_{\R^N} |\nabla \psi ^+|^2 \, dx = \int_{\R^N} |\nabla \psi ^-|^2 \, dx = k$, 
$I(\psi ^+) = I(\psi ^-) = I_{min}(k)$ and $\psi ^+$, $\psi ^-$ satisfy (\ref{5.18}) 
with $ \vartheta^+ = d^+ I_{min}(k)$ and $ \vartheta^- = d^- I_{min}(k)$, respectively.

\end{Proposition} 

{\it Proof. } 
For any $ \phi \in \Eo $ we have $I_c( \phi_{c,c}) = c^N I(\phi)$, 
$ \ii _{\R^N} |\nabla \phi_{c,c} |^2 \, dx = c^{N-2} \ii _{\R^N} |\nabla \phi |^2 \, dx$
and (i) follows. 
The proofs of (ii), (iii) and (iv) are very similar to the proof of Proposition \ref{P4.12} 
and we omit them.
\hfill
$\Box $ 

\medskip

We will establish later (see Proposition \ref{P5.9} below) a  relationship between the traveling waves constructed in 
section \ref{minem} and those given by 
Theorem \ref{T5.5} and Proposition \ref{P5.7} above.
The next remark shows that, in some sense, 
there is equivalence between the inequalities $E_{min}(q) < v_s q$ and $I_{min}(k) < - \frac{k}{v_s^2}$.

\begin{remark} 
\label{R5.8}
 (i) \rm 
Let $ \psi \in \Eo $ be such that $ E(\psi ) < v_s Q(\psi)$ and let 
$ k = \ii _{\R^N} |\nabla \psi |^2 \, dx .$  
Then $I_{min}(\frac{k}{v_s^{N-2}}) < - \frac{k}{v_s^ N}.$
Indeed, we have  
$ { \ii _{\R^N} } |\nabla \psi _{\frac{1}{v_s}, \frac{1}{v_s}} |^2 \, dx = \frac{1}{v_s ^{N-2}} k $ and
$$
I_{min}\left(\frac{k}{v_s^{N-2}} \right) + \frac{k}{v_s^ N} \leq I \left( \psi_{\frac{1}{v_s}, \frac{1}{v_s}} \right)
+ \frac{1}{v_s ^2} \int_{\R^N}  |\nabla \psi _{\frac{1}{v_s}, \frac{1}{v_s}} |^2 \, dx 
= \frac{1}{ v_ s ^N} ( E(\psi) - v_s Q(\psi ) ) < 0.
$$

{\it (ii)} Conversely, let $\psi \in \Eo $ be such that 
$I(\psi ) < - \frac{1}{v_s ^2}   \ii _{\R^N}  |\nabla \psi |^2 \, dx$ 
and denote $ q = v_s^{N-1} Q(\psi)$.  Then $ E_{min}(p) < v_s q.$
Indeed, we have $ Q(\psi_{v_s, v_s}) = v_s^{N-1} Q(\psi) = q $ and
$$
E_{min}(q) - v_s q \leq E(\psi_{v_s, v_s}) - Q (\psi_{v_s, v_s}) 
= v_s ^N \left( \frac{1}{v_s ^2} \int_{\R^N}  |\nabla \psi |^2 \, dx + I(\psi) \right) 
< 0.
$$
\end{remark}

\section{Local minimizers of the energy at fixed momentum ($N=2$)}
\label{local}

We will use the  results in the previous section 
to find traveling waves to (\ref{1.1}) in space dimension $N=2$  which 
are {\it local} minimizers of the energy at fixed momentum
even  when $V$ achieves negative values. 

If $ N =2$ and $q \geq 0$, define
\beq
\label{5.40}
E_{min}^{\sharp }(q) = \inf \Big\{ E( \psi ) \; \big| \; \psi \in \Eo, \; Q( \psi) = q 
\mbox{ and } \int_{\R^2} V(|\psi |^2) \, dx \geq 0 \Big\}.
\eeq
This definition agrees with the one given in section \ref{minem}  in the case $ V \geq 0.$

\begin{Lemma}
\label{L5.10}
Assume that $N=2$ and (A1), (A2) are satisfied. 
The function $E_{min}^{\sharp }$ has the following properties:

\medskip

(i) $E_{min}^\sharp(q) \leq v_s q $ for any $ q \geq 0. $

\medskip

(ii) For any $ \e > 0 $ there is $ q_{\e } > 0 $ such that 
$E_{min}^\sharp (q) > (v_s - \e ) q $ for any $ q \in (0, q_{\e})$. 

\medskip

(iii) $E_{min}^\sharp$ is subadditive on $[0, \infty)$, nondecreasing, 
Lipschitz continuous and its best Lipschitz constant is $v_s$. 

\medskip

(iv) If $ \inf V < 0 $, then for any $ q > 0 $ we have 
$E_{min}^\sharp (q) \leq k_{\infty}$, where $k_{\infty}$ is as in (\ref{k}) or in Lemma \ref{L5.3bis} (iii).

\medskip

(v) $E_{min}^\sharp$ is concave  on $[0, \infty)$.

\end{Lemma}

{\it Proof. } If $ V \geq 0$ on $[0, \infty)$, the statements of Lemma 
\ref{L5.10} have already been proven in section \ref{minem}. 
We only consider here the case when $V$ achieves negative values.
The estimate (i) follows from Lemma \ref{L4.3}. 
For (ii) proceed as in the proof of Lemma \ref{L4.4} and use Lemma \ref{L5.1} instead of Lemma \ref{L4.2}. 
The proof of (iii) is the same as that of Lemma \ref{L4.5}  (i).

\medskip

(iv) Let $ q  > 0$. Fix $ \e > 0$, $ \e $ small. 
By (ii) there is $ \psi \in \Eo $ such that $ \ii _{\R^2} |\nabla \psi |^2 \, dx \leq \frac{\e }{4} $ and 
$ Q( \psi ) \geq \frac{  \e}{ 8 v_s }$. 
It is obvious that
$ \ii _{\R^2} |\nabla (\psi _{\si, \si}) |^2 \, dx = \ii _{\R^2} |\nabla \psi |^2 \, dx \leq \frac{\e }{4} $
for any $ \si > 0 $ and there is $ \si _0 > 0 $ such that $ Q( \psi _{\si _0, \si _0 }) > q$. 
Using Corollary \ref{C3.4} and (\ref{2.9}), we see that there is $ \psi _1 \in \Eo $ such that 
$ Q( \psi _1 ) = q$, $ \ii _{\R^2} |\nabla \psi _1 |^2 \, dx \leq \frac{\e }{2} $ and
$ \psi _1 = 1 $ outside a large ball $B(0, R_1)$. 
Let $M_1 = \ii _{\R^2} V( |\psi _1 |^2 ) \, dx. $

Let $ \psi _0 $ be  as in Lemma \ref{L5.3bis} (i). 
Proceeding as in the proof of Lemma \ref{L5.3bis} (iii) we see  that there exists a radial function 
$ \phi \in C_c^{\infty} ( \R^2)$ and there is $ \e _1 > 0$ such that 
$\ii _{\R^2} V( |\psi _0 + t \phi |^2) \, dx < 0 $ for any $ t \in (0, \e_1)$. 
Taking $ t \in (0, \e _1) $ sufficiently small and using a radial cut-off  and scaling 
it is not hard to construct a radial function $ \psi _2 \in \Eo $ such that 
$ \ii _{\R^2} |\nabla \psi _2 |^2 \, dx \leq k_{\infty} + \frac{\e}{4}$,
$ \ii _{\R^2} V( |\psi _2 |^2) \, dx = - M_2  < 0 $ 
and $ \psi _2 = 1 $ outside a large ball $B(0, R_2)$.
Since $ \psi _2$ is radial, we have $ Q( \psi _2 ) = 0$. 

Let $ t = \left( \frac{ M_1 - \frac{ \e}{4}}{M_2} \right)^{\frac 12}$. 
Choose $ x_0 \in \R^2$ such that $|x_0| > 2( R_1 + t R_2)$ and define
$$ 
\psi _* (x) = \left\{ 
\begin{array}{ll} 
\psi _1 (x) & \quad \mbox{ if } |x| \leq R_1, \\
\psi _2 \left( \frac{  x - x_0}{t } \right) & \quad \mbox{ if } |x| > R_1. 
\end{array}
\right. 
$$
Then $ \psi _* \in \Eo$, 
$ Q( \psi _*) = Q( \psi_1) + t Q(\psi_2) = q $, 
$ \ii _{\R^2} |\nabla \psi _* |^2 \, dx = \ii _{\R^2} |\nabla \psi _1 |^2 \, dx 
+ \ii _{\R^2} |\nabla \psi _2 |^2 \, dx \leq k_{\infty} + \frac{ 3 \e}{4}$, 
and 
$\ii _ {\R^2} V(|\psi _*|^2) \, dx = \ii _ {\R^2} V(|\psi _1|^2) \, dx + t^2 \ii _ {\R^2} V(|\psi _2|^2) \, dx
= M_1 - t^2 M_2 =\frac{ \e}{4 } > 0$. 
Thus $E_{min}^\sharp(q) \leq E(\psi _*) \leq k_{\infty} + \e.$ 
Since $ \e $ is arbitrary, the conclusion follows. 

\medskip

(v) The idea is basically the same as in the proof of Lemma \ref{L4.5} (ii) but we have to be more careful 
because the  functions $ \psi \in \Eo $ that satisfy $ \ii_{\R^2} V(|\psi |^2) \, dx \geq 0 $ do not necessarily satisfy
 $ \ii_{\R^2} V(|S_t^{\pm}\psi |^2) \, dx \geq 0 $ for all $t$, where $S_t^{\pm}$ are as in (\ref{4.10}) - (\ref{4.11}).

Let $E^\sharp =  \ds \sup_{q \geq 0 } E_{min}^\sharp(q).$ 
By (iv) we have $ E^\sharp \leq k_{\infty}$. Denote 
\beq
\label{qsharp}
 q^\sharp = \sup \{ q > 0 \; | \; E_{min}^\sharp(q) < E^\sharp \}.
\eeq

Define $E_{min}^{\sharp,-1} (k) = \sup \{ q \geq 0 \; | \; E_{min}^\sharp(q) \leq k \}.$
Then $E_{min}^{\sharp,-1}$ is finite, increasing, right continuous on $[0, E^\sharp )$ and
$E_{min}^{\sharp} (E_{min}^{\sharp,-1} (k)) = k $ for all $ k \in [0, E^{\sharp})$.
By convention, put $E_{min}^{\sharp,-1} (k) = 0$ if $ k <0$.
For any $ \phi \in \Eo $ with $\ii _{\R^2} V(|\phi |^2) \, dx \geq 0$ 
 we have $ E_{min}^{\sharp }( Q(\phi) ) \leq E(\phi)$, thus
\beq
\label{Einverse}
Q( \phi ) \leq E_{min}^{\sharp,-1} (E(\phi)).
\eeq

We will prove that for any fixed $ q \in (0, q^{\sharp})$ there are $ q_1 < q$ and $ q_2 > q$ such that $E_{min}^{\sharp}$
is concave on $[q_1, q_2]$. 

Let $ q \in (0, q^\sharp )$. Fix an arbitrary $ \e > 0 $ such that 
$ E_{min}^\sharp (q) + 4 \e < E^\sharp .$ 
Choose $ \psi \in \Eo $ such that $ Q( \psi ) = q $, 
$ \ii _{\R^2} V(|\psi |^2) \, dx > 0 $ and 
$E(\psi) < E_{min}^\sharp(q) + \e$.

We may assume that $ \psi $ is symmetric with respect to $ x_2$. 
Indeed, let $S_t^+$ and $S_t^-$ be as in (\ref{4.10})-(\ref{4.11}).
Arguing as in the proof of Lemma \ref{L4.5} (ii), there is $ t_0 \in \R$ such that 
$\ii _{\R^2} |\nabla (S_{t_0 }^+(\psi) ) |^2 \, dx = \ii _{\R^2} |\nabla (S_{t_0 }^-(\psi) ) |^2 \, dx 
= \ii _{\R^2} |\nabla \psi |^2 \, dx < k_{\infty}$.
After a translation, we may assume that $ t_0 = 0$. 
Let $ \psi _1 =  S_{0 }^-(\psi)$, $ \psi _2 = S_{0 }^+(\psi)$, denote $q_i = Q( \psi _i)$ and 
$ v_i = \ii _{\R^2} V(|\psi _i |^2) \, dx $, $ i = 1,2$ and $ v= \ii _{\R^2} V(|\psi  |^2)$, 
so that $ q_1 + q_2 = 2 Q(\psi) = 2q $ and $ v_1 + v_2 = 2 v$. Since 
$\ii _{\R^2} |\nabla \psi_i ) |^2 \, dx < k_\infty = T $, by Lemma 
{\ref{L5.3bis}  we have $ v_1 \geq 0$ and $ v_2 \geq 0$ and consequently $v_1, v_2 \in [0, 2v]$.  
If $ q_ 1 \leq 0 $ we have $ q_2 \geq 2 q$ and then for 
$\si _2 = \frac{ q}{q_2} \leq \frac 12$, we get
$Q((\psi_2)_{\si _2, \si_2}) = q$ and $E((\psi_2)_{\si _2, \si_2}) \leq E(\psi) 
< E_{min}^\sharp ( q ) + \e $, 
hence we may choose $(\psi_2)_{\si _2, \si_2} $ instead of $ \psi$, and $(\psi_2)_{\si _2, \si_2} $   is 
symmetric with respect to $x_2$. A similar argument works if $ q_2 \leq 0$. 
If $ q_1 > 0 $ and $ q_2 > 0$, let $\si _1 = \frac{ q}{q_1}$ and 
$\si _2 = \frac{ q}{q_2}$, so that $ \frac{1}{\si _1} + \frac{1}{\si_2} = 2$. 
We claim that  there is $ i \in \{ 1, 2 \}$ such that $ \si _i ^2 v_i \leq v$, 
and then we may choose $(\psi_i)_{\si_i, \si _i}$, 
which is symmetric with respect to $x_2$,  instead of $ \psi$ .
 Indeed, if the claim is 
false we have $ v_i > \frac{1}{\si _i ^2} v$ and taking the sum we get 
$ 2 > \frac{1}{\si_1^2} + \frac{1}{\si_2^2}$, 
which is impossible because  $\frac{1}{\si _1} + \frac{1}{\si_2} = 2$.

Since $ \psi $ is symmetric with respect to $ x_2$, we have 
$Q(S_0^{\pm } \psi ) = q $ and 
$E(S_0^{\pm } \psi ) = E( \psi) < k_{\infty} - 3 \e .$ As in Lemma 
\ref{L4.5} (ii), the mapping $ t \longmapsto E(S_t^- \psi) $ is continuous 
and tends to $ 2 E(\psi)$ as $ t \lra \infty$.
Let 
$$
 t_{\infty} = \inf \{ t \geq 0  \; | \; E(S_t^- \psi)  \geq k_{\infty} \} 
\qquad (\mbox{with possibly }  t_{\infty} = \infty  \mbox{ if } E(\psi) \leq \frac 12 k_{\infty}).
$$ 
For any $ t \in [0, t_{\infty})$ we have $E(S_t^- \psi)  < k_{\infty} .$
If there is $ t \in [0, t_{\infty})$  such that $\ii _{\R^2} V( |S_t^-\psi |^2) \, dx = 0$, 
we have necessarily $ \ii _{\R^2} |\nabla (S_t^-\psi )|^2 \, dx \geq k_{\infty}$, 
thus $E(S_t^- \psi)  \geq k_{\infty} $, a contradiction. 
We infer that the function $t \longmapsto \ii _{\R^2} V( |S_t^-\psi |^2) \, dx  $ 
is continuous, positive at $ t =0$ and cannot vanish on $[0, t_{\infty})$, 
hence $ \ii_{\R^2} V(|S_t^-\psi |^2) \, dx  >0$ for all $ t \in [0, t_{\infty})$. Consequently we have
\beq
\label{5.41}
E(S_t^- \psi) \geq E_{min}^\sharp (Q(S_t^- \psi)) \qquad 
\mbox{ for any } t \in [0, t_{\infty}).
\eeq
For any $ t \geq 0$ we have 
$\ii _{\R^2} |\nabla (S_t^+ \psi ) |^2 \, dx = 2 \ii _{\{x_2 \geq t \} } |\nabla \psi |^2 \, dx 
\leq 2 \ii _{\{ x_2 \geq 0 \}} |\nabla \psi |^2 \, dx \leq E(\psi ) < k_{\infty}$, 
hence $\ii _{\R^2} V( |S_t^+\psi |^2) \, dx \geq 0$ (by Lemma \ref{L5.3bis}) 
and therefore
\beq
\label{5.42}
E(S_t^+ \psi) \geq E_{min}^\sharp (Q(S_t^+ \psi)) \qquad \mbox{ for any } t \geq 0.
\eeq
The mapping $ t \longmapsto Q(S_t^+ \psi)$ is continuous, tends to $0$ as 
$ t \lra \infty$ and $Q(S_0^+ \psi) =q$. If $ t_{\infty}= \infty$, 
for any $ q_1 \in (0, q)$ there is $ t_{q_1} > 0$ such that 
$Q(S_{t_{q_1}}^+ \psi) = q_1.$ Then $Q(S_{t_{q_1}}^- \psi) = 2q -q_1$ 
and using (\ref{5.40}), (\ref{5.41}) we get
$$
E_{min}^\sharp(q) + \e > E(\psi) = \frac 12 \left( E(S_{t_{q_1}}^+ \psi) + E(S_{t_{q_1}}^- \psi) \right) 
\geq \frac 12 \left( E_{min}^\sharp( q_1) + E_{min}^\sharp(2q - q_1) \right) .
$$

In the case  $ t_{\infty } < \infty$ we have $E(S_{t_{\infty}}^ - \psi) = k_{\infty}$, 
hence $E(S_{t_{\infty}}^ + \psi) = 2 E(\psi) - E(S_{t_{\infty}}^ - \psi) < 2 E_{min}^{\sharp}(q) + 2 \e - k_{\infty} < E_{min}^\sharp(q) $ 
and by (\ref{Einverse}) it follows that 
$$
Q(S_{t_{\infty}}^ + \psi)  \leq E_{min}^{\sharp,-1} ( 2 E_{min}^\sharp(q) 
+ 2 \e - k_{\infty}) < q.
$$
For any $ q_1 \in [Q(S_{t_{\infty}}^ + \psi) , q]$ there is 
$ t_{q_1 } \in [0, t_{\infty}]$ such that $Q(S_{t_{q_1}}^+ \psi) = q_1.$ 
As above, we obtain
$$
E_{min}^\sharp(q) + \e > \frac 12 \left( E_{min}^\sharp( q_1) 
+ E_{min}^\sharp(2q - q_1) \right) 
\qquad \mbox{ for any } q_1 \in [E_{min}^{\sharp,-1} ( 2 E_{min}^\sharp(q) 
+ 2 \e - k_{\infty}), \; q ].
$$
Since $ \e \in (0, \; \frac 14( E^\sharp - E_{min}^\sharp(q) )) $ 
is arbitrary and $E_{min}^{\sharp,-1} $ is right continuous we infer 
that for any $ q \in (0, q^\sharp)$ there holds 
\beq
\label{5.43}
E_{min}^\sharp(q) \geq \frac 12 \left( E_{min}^\sharp( q_1) + E_{min}^\sharp(2q - q_1) \right) 
\quad \mbox{ for all } q_1 \in  (E_{min}^{\sharp,-1} 
( 2 E_{min}^\sharp(q)  - k_{\infty}), \; q ].
\eeq

The function $ q \longmapsto E_{min}^{\sharp,-1} ( 2 E_{min}^\sharp(q)  - k_{\infty}) $ is nondecreasing and right continuous on $(0, q^\sharp)$. Fix 
$ q_* \in (0, q^\sharp)$. We have 
$$
\lim_{q \downarrow q_* } E_{min}^{\sharp,-1} ( 2 E_{min}^\sharp(q)  - k_{\infty}) 
= E_{min}^{\sharp,-1} ( 2 E_{min}^\sharp(q _* )  - k_{\infty}) < q_* 
$$
because $ 2 E_{min}^{\sharp} ( q_*) - k_{\infty} < E_{min}^{\sharp} ( q_*)$.
It is then easy to see that  there are $ q_* ' < q_* $ and $ q_* '' \in (q_* , q^\sharp)$ 
such that for any $ q \in [q_* ', q_* ''],$
\beq
\label{5.44}
E_{min}^{\sharp,-1} ( 2 E_{min}^\sharp(q)  - k_{\infty}) < q_* ' .
\eeq
Using (\ref{5.43}) we see that for any $ q_1, q_2 \in [ q_* ', q_* ''] $ 
we have
$$
E_{min}^\sharp \left( \frac{ q_1 + q_2}{2} \right) \geq 
\frac 12 (E_{min}^\sharp( q_1) + E_{min}^\sharp( q_2) ).
$$
Since $E_{min}^\sharp$ is continuous, we infer that $E_{min}^\sharp$ 
is concave on $ [q_* ', q_* ''] $. Thus any point $ q_* \in (0, q^\sharp)$ 
has a neighborhood where $E_{min}^\sharp$ is concave 
and then it is not hard to see that $E_{min}^\sharp$ is concave on $[0, q^\sharp)$. 
If $ q^\sharp < \infty$ we have $E_{min}^\sharp = E^\sharp $ on 
$[q^\sharp, \infty)$, hence $E_{min}^\sharp$ is concave on $[0, \infty)$.
\hfill
$\Box$

\medskip

Let 
\beq
\label{5.45}
q_0^\sharp = \inf \{ q > 0 \; | \; E_{min}^{\sharp} (q) < v_s q \} 
\qquad \mbox{ and } \qquad 
q_{\infty}^{\sharp } = \sup \{ q > 0 \; | \; E_{min}^\sharp(q) < k_{\infty} \}.
\eeq
It is obvious that $ q_0^\sharp \leq q_{\infty}^\sharp $ and 
$ q_{\infty} > 0$ because $E_{min}^{\sharp }(q) \to 0 < k_\infty$ as $q \lra 0$. 
If $F$ satisfies assumption (A4) and $ F''(1) \not= 3 $, it follows from Theorem \ref{T4.13} 
that $ q_0^\sharp = 0$ (notice that 
the test functions 
$U_{\e}$ constructed in the proof of Theorem \ref{T4.13} satisfy $ V(|U_{\e}|^2)  \geq 0$ in $ \R^2$). 

\medskip

Our next result shows the 
precompactness of minimizing sequences for  $E_{min}^\sharp (q)$.

\begin{Theorem}
\label{T5.11}
Assume that $N=2$,  (A1), (A2) are satisfied, and $ \inf V < 0$. 
Let $ q \in (q_0^\sharp, q_{\infty}^\sharp)$ and assume that 
$ (\psi _n )_{n \geq 1 } \subset \Eo $ is a sequence satisfying
$$ 
\ii_{\R^2} V(|\psi_n|^2) \, dx \geq 0, \quad Q( \psi _ n ) \lra q 
\quad \mbox{  and } \quad  E( \psi _n ) \lra E_{min}^\sharp(q).
$$ 
There exist a subsequence $(\psi_{n_k})_{k \geq 1}$, a sequence of 
points $(x_k)_{k \geq 1}  \subset \R^N$, 
and $ \psi \in \Eo $ such that $Q( \psi ) = q$, $E( \psi ) = E_{min}^\sharp(q) $, 
$\psi_{n_k} ( x_k + \cdot ) \lra \psi $ a.e. on $ \R ^2$ and 
$ \ds \lim_{k \ra \infty} d_0(\psi_{n_k} ( x_k + \cdot ) , \;  \psi ) = 0 .$ 
Furthermore, $ \ii_{\R^2} V(|\psi|^2) \, dx > 0 $, hence 
$\psi \in \Eo $ is a local minimizer in the sense that
$$ E(\psi) = E_{min}^\sharp (q) = \inf \Big\{ E(w)\; \big| \; 
w \in \Eo, \qquad Q(w) = q , \qquad \ii_{\R^2} V(|w|^2) \, dx > 0 \Big\} .$$

Moreover, the conclusions of Proposition \ref{P4.12} hold  true with 
$E_{min}$ replaced by $E_{min}^\sharp$. 

\end{Theorem}

{\it Proof. } Fix $ k_1, k_2$ such that  $ 0< k_1 < E_{min}^\sharp(q) < k_2 < k_{\infty}$. 
We may assume that $ k_1 < E( \psi _n) < k_2$ for all $n$. 
By Lemma \ref{L4.1} there is $ C_1( k_1) > 0 $ such that 
$E_{GL}( \psi _n) \geq C_1( k_1).$ Since 
$\ii_{\R^2} V(|\psi_n|^2) \, dx \geq 0, $
we have $ \psi_n \in \Eo_{k_2, k_2}$ and using Lemma \ref{L5.3bis} we 
infer that $E_{GL}( \psi _n)$ is bounded. Passing to a subsequence 
if necessary, we may assume that $E_{GL}(\psi _n) \lra \al _0 > 0$. 
Then we proceed as in the proof of Theorem \ref{T4.7} and we use the Concentration-Compactness Principle 
for the sequence of functions 
$ f_n = |\nabla \psi _n |^2 + \frac 12 \left( \ph^2( |\psi _n|)- 1 \right)^2$. 

We rule out vanishing thanks to Lemma \ref{L4.8}.

If dichotomy occurs for a subsequence (still denoted $(\psi_n)_{n \geq 1}$), 
using Lemma \ref{splitting} for all $n$ sufficiently large we construct 
two functions $ \psi_{n, 1}, \; \psi_{n , 2} \in \Eo $ 
such that 
$\big\vert \ii _{\R^2} |\nabla \psi_n |^2 \,dx - \ii _{\R^2} |\nabla \psi_{n , 1} |^2 \,dx
- \ii _{\R^2} |\nabla \psi_{n , 2} |^2 \,dx \big\vert \lra 0$, and 
(\ref{4.28}), (\ref{4.29}), (\ref{4.30}) hold for some $ \al \in (0, \al _0)$. 
In particular, we have $\ii _{\R^2} |\nabla \psi_{n , i} |^2 \,dx < k_2 < k_{\infty }$, $ i =1,2$ for all $n$ sufficiently large and this implies 
$\ii _{\R^2} V(|\psi_{n, i}|^2) \, dx \geq 0$, 
so that $E(\psi_{n, i}) \geq E_{min}^\sharp (Q (\psi_{n, i}))$. 
Since $ q \in (q_0^\sharp, q_{\infty}^\sharp)$, using 
the concavity of $E_{min}^{\sharp }$ and Lemma \ref{L5.10} (i) and (ii)
we infer that $E_{min}^\sharp(q) < E_{min}^\sharp(q') + E_{min}^\sharp(q- q')$ 
for any $ q' \in (0, q)$. 
Then arguing  as in the proof of Theorem \ref{T4.7} we rule out dichotomy and we conclude that concentration occurs.
 
Hence there  is a sequence  $(x_n)_{n \geq 1}  \subset \R^N$
such that, denoting $ \tilde{\psi}_{n} = \psi_{n} ( x_{n} + \cdot ) $, 
(\ref{4.31}) holds.
Consequently there are  a subsequence $(\tilde{\psi}_{n_k})_{k \geq 1}$ and 
$ \psi \in \Eo $ that satisfy (\ref{4.32}) and (\ref{4.33}). 
Using Lemmas \ref{L4.10} and \ref{L4.11} we get 
$ {\ds \lim_{k \ra \infty} } \| \, | \tilde{\psi}_{n_k} | - | \psi | \, \| _{L^2( \R^N)} = 0,  $
\beq
\label{5.46}
\lim_{k \ra \infty} \ii _{\R^2} V( | \tilde{\psi}_{n_k} |^2) \, dx = \ii _{\R^2} V( | \psi |^2) \, dx 
\quad \mbox{ and } \quad \lim_{k \ra \infty} Q(  \tilde{\psi}_{n_k}) = Q( \psi).
\eeq 
In particular, we have $\ii _{\R^2} V( | \psi |^2) \, dx \geq 0$, $Q( \psi) = q$ 
and this implies $E(\psi) \geq E_{min}^\sharp(q)$. 
Combining this information with (\ref{4.33}) and (\ref{5.46}) we see that necessarily 
$\ii _{\R^2} |\nabla \tilde{\psi}_{n_k} |^2\, dx \lra \ii _{\R^2} |\nabla \psi |^2\, dx.$
Together with the weak convergence $\nabla \tilde{\psi}_{n_k} \rightharpoonup \nabla \psi$ in $L^2( \R^2)$, 
this implies the strong convergence 
$\| \nabla \tilde{\psi}_{n_k} - \nabla \psi \| _{L^2( \R^2)} \lra 0 $.  
Hence $ d_0( \tilde{\psi}_{n_k} , \psi) \lra 0$ as $ k \lra \infty$. 
The fact that $ \ii _{\R^2} V( | \psi |^2) \, dx > 0 $ comes from 
the fact that 
$|\psi |$ is not constant (because $Q(\psi ) = q >0$) and 
$ \ii_{\R^2} |\nabla \psi|^2 \, dx < k_\infty $.

The last part is proved in the same way as Proposition \ref{P4.12}.
\hfill
$\Box$

\medskip

If $ q_{\infty}^\sharp < \infty$ we have $E_{min}^\sharp(q) = k_{\infty} $ 
for all $ q \geq q_{\infty}^\sharp$. 
The conclusion of Theorem \ref{T5.11} is not valid for $ q \geq q_{\infty}^\sharp $. 
Indeed, for such $q$ the argument used in the proof of Lemma \ref{L5.10} (iv) 
leads to the construction of a minimizing sequence $(\psi_n)_{n \geq 1} \subset \Eo$ satisfying the assumptions of Theorem \ref{T5.11}, 
but $E_{GL}(\psi_n) \lra \infty$. Furthermore, if 
$ \ii_{\R^2} |\nabla \psi|^2 \, dx \geq k_\infty $, Lemma {\ref{L5.3bis} 
does not guarantee that the potential energy $ \ii _{\R^2} V(|\psi |^2) \, dx $ is positive.

\section{Orbital stability} 
\label{sectionorbistab}

It is beyond the scope of the present paper to study the Cauchy problem associated 
to (\ref{1.1}). Instead, we will content ourselves to assume in the sequel that 
the nonlinearity $F$ satisfies (A1), (A2) and is such that the following holds:

\medskip

{\bf (P1)} (local well-posedness) For any $ M > 0$ there is $T(M) > 0$ 
such that for any $ \psi _0 \in \Eo $ with $E_{GL}(\psi _0) \leq M$ there exist
$T_{\psi_0} \geq T(M)  $ and a unique solution 
$ t \longmapsto \psi (t) \in C([0, T_{\psi_0}), (\Eo, d))$ such that $ \psi(0) = \psi_0$. 
Moreover, $\psi(\cdot)$ depends continuously on the initial data in the following sense: 
if $ d( \psi _0 ^n, \psi _0) \lra 0$ and $ t \longmapsto \psi _n(t)  $ is the solution 
of (\ref{1.1}) 
with initial data $\psi_0^n$, 
then for any $T <  T_{\psi_0} $ we have $T < T_{\psi_0^n}$ for all sufficiently large $n$ 
and $d( \psi _n(t), \psi (t)) \lra 0$ uniformly on $[0, T]$ as $ n \lra \infty$. 

\medskip

{\bf (P2)} (conservation of phase at infinity) 
We have $\psi(\cdot) - \psi _0 \in C([0, T_{\psi_0}), H^1( \R^N)).$

\medskip

{\bf (P3)} (conservation of energy) We have $E(\psi(t) ) = E( \psi _0) $ for any 
$t \in [0, T_{\psi_0})$.

\medskip

{\bf (P4)} (regularity) If $\Delta \psi _0 \in L^2( \R^N)$, then 
$ \Delta \psi(\cdot) \in C([0, T_{\psi_0}), L^2( \R^N) ) $.

\medskip

In space dimension $N = 2,3,4$, the Cauchy problem for the Gross-Pitaevskii equation 
(that is (\ref{1.1}) with $F(s) = 1- s$) 
has been studied in \cite{PG, PG2} and it was proved that the flow has the properties 
(P1)-(P4) above. Moreover, the solutions found in \cite{PG, PG2} are global in 
time if $ N=2, 3$ or if $N=4$ and the initial data has sufficiently small energy. 
This comes from the conservation of energy and from the fact that the Gross-Pitaevskii 
equation is subcritical if $N=2, \; 3$ and it is critical  if $N=4$. 
It seems that the proofs in \cite{PG, PG2} can be easily adapted to more general 
subcritical nonlinearities provided that the associated nonlinear potential $V$ 
is nonnegative on $[0, \infty)$. Notice that any nonlinearity satisfying (A2) 
is subcritical.

Recently it has been proved in \cite{KOPV} that  the Gross-Pitaevskii equation is globally well-posed on the whole energy space $ \Eo$ in space dimension $N=4$
and that  the cubic-quintic NLS is globally 
well-posed on $ \Eo$ if $N=3$, 
 despite the fact that both problems are critical.

\medskip

Assume that (P1) and (P3) hold. 
If $ V \geq 0$, using the conservation of energy and Lemma \ref{L4.6} it is 
easy to prove that all solutions are global.

If $N=2$ and $ \inf V < 0 $, any solution $ t \longmapsto \psi (t)$ 
with initial data $ \psi _0 $ satisfying 
$ \ii _{ \R^2} |\nabla \psi _0 | ^2 \, dx < k_{\infty}$ and 
$E( \psi _0) < k_{\infty}$ is global. 
Indeed, the mapping $ t \longmapsto \ii_{\R^2} V(|\psi (t) |^2) dx $ is continuous; 
if it changes sign at some 
$ t_0 \in (0, T_{\psi _0})$, there are two possibilities: 
either $ \psi(t_0)$ is constant (and then $E(\psi (t_0)) = 0$, 
hence $E(\psi (t)) = 0 $ for all $t$ and $ \psi (t)$ is constant)
or  Lemma \ref{L5.3bis} (i) implies that 
$\ii_{\R^2} V(|\psi (t_0) |^2) dx = 0  $ and 
$ \ii _{ \R^2} |\nabla \psi _0 | ^2 \, dx \geq k_{\infty}$, thus $E(\psi(t_0) ) \geq k_{\infty}$, 
contradicting the fact that, by conservation of the energy, 
$E(\psi(t_0) ) = E( \psi _0 ) < k_{\infty}$.
Consequently $ 0 \leq \ii_{\R^2} V(|\psi (t) |^2) dx \leq E( \psi _0)$ and
$ 0 \leq \ii_{\R^2} |\nabla \psi (t) |^2 dx \leq E( \psi _0)$ 
as long as the solution exists. Then Lemma \ref{L5.3bis} (ii) implies that $E_{GL}( \psi (t))$ 
remains bounded and using (P1) we see that the solution is global.

\medskip

In the case of more general nonlinearities, the Cauchy problem for (\ref{1.1})
has been considered by C. Gallo in \cite{gallo}. 
In space dimension $N=1,\, 2,\, 3, \, 4 $ and under suitable assumptions on $F$, he proved the following
(see Theorems 1.1  and 1.2 pp. 731-732  in \cite{gallo}):

\medskip

{\bf (P1')}  For any $ \psi _0 \in \Eo$ and any $ u_0 \in H^1 ( \R^N)$, there exists a unique global 
solution $ \psi _0 + u(t)$, where $ u(\cdot) \in C([0, \infty), H^1( \R^N))$ 
and $ u(0) = u_0$. 
The solution depends continuously on the initial data $ u_0 \in H^1( \R^N)$.

\medskip

Notice that the solutions in \cite{gallo}  satisfy (P2) by construction
and they also satisfy  (P3) and (P4). 
Moreover, it is proved (see Theorem 1.5 p. 733 in \cite{gallo})  that 
any solution $ \psi \in C([0, T], \Eo )$ 
automatically satisfies (P2).

\begin{Lemma}
\label{L6.1} (conservation of the momentum)
Assume that $F$ is such that (A1), (A2), ((P1) or (P1')) and (P2)$-$(P4) hold. 
Let $ \psi _0 \in \Eo $ and let $ \psi $ be the solution of (\ref{1.1}) 
with initial data $ \psi _0$, as given by (P1) or (P1')). Then
$$
Q( \psi (t)) = Q( \psi _0) \qquad \mbox{ for any } t \in [0, T_{\psi_ 0}).
$$
\end{Lemma}

{\it Proof. } 
Assume  that $\psi _0 \in \Eo $ is such that $\Delta \psi _0 \in L^2( \R^N)$. 
Let $\psi(\cdot)$ be the solution of (\ref{1.1}) with initial data $ \psi _0$. 
By (P1) and (P4) we have $  \psi_{x_j}(\cdot ) \in C([0, T_{\psi_0}), H^1( \R^N ))$, $ j = 1, \dots, N$. 
Let $t, t+s \in [0, T_{\psi _0} )$. 
Since $ \psi( t+s) - \psi (t) \in H^1( \R^N)$ by (P2), 
the Cauchy-Schwarz inequality implies 
 $\langle i \psi_{x_1}( t+s) + i \psi _{x_1}(t) ,  \psi( t+s) - \psi (t) \rangle \in L^1( \R^N)$.
Using the definition of the momentum and Lemma \ref{L2.3} we get 
$$
\begin{array}{l}
\frac 1s \left( Q (\psi( t+s)) - Q( \psi (t)) \right) 
= \frac 1s L ( \langle i \psi_{x_1}( t+s) + i \psi _{x_1}(t) , \;  \psi( t+s) - \psi (t) \rangle )
\\
\\
= {\ds \ii _{\R^N} } \langle i \psi_{x_1}( t+s) + i \psi _{x_1}(t) , \;  \frac 1s ( \psi( t+s) - \psi (t) ) \rangle\, dx.
\end{array}
$$
Letting $ s \lra 0 $ in the above equality and using (\ref{1.1}) we get 
\beq
\label{6.1}
\frac{d}{dt} ( Q( \psi (t))) = 2 \ii_{\R^N } \langle \frac{ \p \psi (t)}{\p x_1 } , \; 
\Delta \psi (t) + F( |\psi |^2) \psi (t) \rangle \, dx.
\eeq
Since $ \frac{ \p \psi (t)}{\p x_j } \in H^1( \R^N)$, using the integration by parts formula for $H^1$ functions
(see, e.g., \cite{brezis} p. 197) we have
\beq
\label{6.2}
\ii_{\R^N } \langle \frac{ \p \psi (t)}{\p x_1 } , \; \Delta \psi (t) \rangle \, dx 
= - \ii_{\R^N} \sum_{j = 1}^N 
\langle \frac{ \p ^2 \psi (t) }{ \p x_1 \p x_j } , \; \frac{ \p \psi (t) }{ \p x _j} \rangle \, dx 
= - \frac 12 \ii_{\R^N} \frac{ \p } {\p x_1} \left( |\nabla \psi(t) |^2 \right) \, dx. 
\eeq
We have $|\nabla \psi(t) |^2  \in L^1( \R^N)$ and 
$  \frac{ \p } {\p x_k} \left( |\nabla \psi(t) |^2 \right)= 2 \sum_{j = 1}^N  
\langle \frac{ \p ^2 \psi (t)}{\p x_k \p x_j },\;  \frac{ \p  \psi (t)}{ \p x_j } \rangle \in L^1( \R^N)$, 
hence $ |\nabla \psi(t) |^2  \in W^{1,1}( \R^N)$.
It is well-known that for any  $ f \in W^{1,1}(\R^N)$ we have $ \ii _{\R^N} \frac{\p f}{\p x_j} (x) \, dx = 0$ 
and using (\ref{6.2}) we get 
$\ii_{\R^N } \langle \frac{ \p \psi (t)}{\p x_1 } , \; \Delta \psi (t) \rangle \, dx = 0.$

On the other hand, $ 2 \langle  \psi _{x_1} (t) ,\; F( |\psi |^2) \psi (t) \rangle 
= - \frac{ \p }{\p x_1} \left( V(|\psi(t )|^2) \right).$
We have $ V(|\psi(t )|^2)  \in L^1( \R^N)$ by Lemma \ref{L4.1}.  Using the fact that 
$ \psi _{x_j} ( t) \in H^1( \R^N)$, (A1), (A2) and the Sobolev embedding it is easy to see that 
$  \frac{ \p }{\p x_j} \left( V(|\psi(t )|^2) \right)  = -2 \langle  \psi _{x_j} (t) ,\; F( |\psi |^2) \psi (t) \rangle \in L^1( \R^N)$ 
for all $j$, hence $ V(|\psi(t )|^2)  \in W^{1,1}( \R^N)$ and therefore 
$\ii _{\R^N} \frac{ \p }{\p x_1} \left( V(|\psi(t )|^2) \right) \, dx = 0.$ 
Then using (\ref{6.1}) we obtain $\frac{d}{dt} ( Q( \psi (t))) = 0$ for any $t$, 
consequently $  Q( \psi ( \cdot))$ is constant on $[0, T_{\psi _0}).$ 

\medskip

Let $\psi _0 \in \Eo$ be arbitrary.
By Lemma \ref{L3.5}, there is a sequence $ (\psi _0 ^n)_{n \geq 1} \subset \Eo $ 
such that $ \nabla \psi _0 ^n \in H^2 ( \R^N)$
and $\| \psi _0 ^n - \psi _0 \|_{H^1( \R^N)} \lra 0 $ as $ n \lra \infty$
(thus, in particular, $d( \psi _0 ^n, \psi _0) \lra 0 $). 
Fix $ T \in (0, T_{\psi _0}).$ It follows from (P1) or (P1')
that for all sufficiently large $n$, the solution $\psi _n (\cdot ) $ of (\ref{1.1})  
with initial data $ \psi _0 ^n$ exists at least on $[0, T]$ and $ d ( \psi _n (t) , \psi (t) ) \lra 0 $ 
uniformly on $[0, T]$. 
Using Corollary \ref{C4.13} we infer that for any fixed $ t \in [0, T]$ we have 
$Q( \psi _n (t) ) \lra Q(\psi (t))$. 
From the first part of the proof and Corollary \ref{C2.4} we get 
$Q( \psi _n(t)) = Q (\psi _0^n) \lra Q( \psi _0) $ as $ n\lra \infty$. 
Hence $Q( \psi (t)) = Q( \psi _0)$. 
\hfill
$\Box $

\medskip

We now state our orbital stability result, which is  based on the argument 
in \cite{CL}.

\begin{Theorem}
\label{T6.2} 
Assume that (A1), (A2), ((P1) or (P1')) and (P2)$-$(P4) hold. 

$\bullet$ We assume $N \geq 2$ and  $V\geq 0$ on $[0, \infty )$.  Let 
$ q > q_0 $, and define $ \So_q = \{ \psi \in \Eo \; | \; Q( \psi ) = q, 
\mbox{ and } E( \psi ) = E_{min}(q) \}. $ 

\smallskip

Then, $\So _q$ is not empty 
and is orbitally stable by the flow of (\ref{1.1}) for the semi-distance $d_0$ 
in the following sense: for any $ \e > 0$ there is $ \de _{\e } >0$ 
such that any solution of (\ref{1.1}) with initial data $\psi _0$ 
such that $d_0( \psi_0, \So _q)< \de _{\e} $ is global and 
satisfies $d_0 (\psi(t), \So _q) < \e$ for any $ t > 0 $. 

\smallskip

$\bullet$ Assume that $N=2$ and $ \inf V < 0 $.  Let 
$ q \in (q_0^\sharp, q_{\infty}^\sharp)$, where $q_0^\sharp, q_{\infty}^\sharp $ are as in (\ref{5.45}),  and define 
$ \So_q^\sharp = \{ \psi \in \Eo \; | \; Q( \psi ) = q, \; 
\ii_{\R^2 } V(|\psi|^2) dx \geq 0  \mbox{ and } E( \psi ) = E_{min}^\sharp(q) \}. $ 

\smallskip

Then $\So _q^\sharp$ is orbitally stable by the flow of (\ref{1.1}) 
for the semi-distance $d_0$.

\end{Theorem}

{\it Proof. } We argue by contradiction and we assume that the 
statement is false. Then there is some $ \e _0> 0$ such that 
for any $ n \geq 1$ there is $ \psi _0^n \in \Eo$ satisfying 
$ d_0( \psi_0^n, \So _q) < \frac 1n$ (resp. $ d_0( \psi_0^n, \So_q^\sharp) < \frac 1n$)
 and there is $ t_n >0$ such that 
$d_0 ( \psi_n( t_n) , \So _q ) \geq \e _0 $ (resp. 
$d_0 ( \psi_n( t_n) , \So _q^\sharp ) \geq \e _0 $), where $\psi_n$ is 
the solution of the Cauchy problem associated to (\ref{1.1}) with 
initial data $ \psi _0^n$. 

\medskip

We claim that $Q( \psi _0^n ) \lra q $ and $E(\psi _0^n) \lra E_{min}(q) $ 
(resp. $E(\psi _0^n) \lra E_{min}^\sharp(q) $). 
Indeed,  for each $n$ there is $ \phi _n \in \So _q$ (resp. $\in \So _q^\sharp$) 
such that $ d_0 ( \psi _0^n , \phi _n ) < \frac 2n$. 
If $N=2$ and $V$ achieves negative values, we have 
$$
\limsup_{n \ra \infty}  \ii _{\R^2} |\nabla \psi _0^n | ^2 \, dx 
= \limsup_{n \ra \infty}  \ii _{\R^2} |\nabla \phi ^n | ^2 \, dx 
\leq \limsup_{n \ra \infty} E(\phi_n) = E_{min}^{\sharp}(q) < k_{\infty } , 
$$
hence $ \ii _{\R^2} V(|\psi _0^n |^2) \, dx \geq 0$ for all sufficiently large $n$. 
Consider an arbitrary subsequence $(\psi_0^{n_{\ell}})_{\ell \geq 1} $ 
of $(\psi_0^n)_{n \geq 1}$. Using either Theorem \ref{T4.7} or 
Theorem \ref{T5.11} we infer that there exist a subsequence 
$(\phi_{n_{\ell_k}})_{ k \geq 1}$ of $ (\phi_n)_{n \geq 1}$, 
a sequence $ (x_k)_{k \geq 1} \in \R^N$ and $ \phi \in \So_q$ (resp. $\in \So_q^\sharp$) such that 
$ d_0( \phi_{n_{\ell_k}} ( \cdot + x_k ), \, \phi ) \lra 0$ as $ k \lra \infty$. 
Then $ d_0( \psi_0^{n_{\ell_k}} ( \cdot + x_k ), \, \phi ) \leq 
 d_0( \phi_{n_{\ell_k}} ( \cdot + x_k ), \, \phi ) + \frac{2}{n_{\ell_k} } \lra 0$ 
and using Corollary \ref{C4.13} we get 
$ Q(  \psi_0^{n_{\ell_k}}  ) = Q(\psi_0^{n_{\ell_k}} ( \cdot + x_k ) ) 
\lra Q( \phi) = q $ and
$ E(  \psi_0^{n_{\ell_k}}  ) = E(\psi_0^{n_{\ell_k}} ( \cdot + x_k ) ) \lra E( \phi) = E_{min}(q) $ 
(resp. $ E(  \psi_0^{n_{\ell_k}}  )  \lra E( \phi) = E_{min}^\sharp (q)$). Since any subsequence 
of $(\psi_0^n)_{n \geq 1}$ contains a subsequence as above, the claim follows. 

\medskip

By (P3) and Lemma \ref{L6.1} we have 
$E(\psi_n(t_n) ) = E(\psi _0^n) \lra E_{min}(q)$ (resp. $E(\psi_n(t_n) ) \lra E_{min}^\sharp (q)$) 
and $Q(\psi_n(t_n) ) = Q(\psi _0^n) \lra q.$ 
Moreover, if $N=2$ and $ \inf V < 0 $, we have already seen that 
$ \ii _{\R^2} V(|\psi _n (t) |^2) \, dx$ cannot change sign, hence
$ \ii _{\R^2} V(|\psi _n (t_n) |^2) \, dx \geq 0$. 
Using again either Theorem \ref{T4.7} or Theorem \ref{T5.11} we see that 
there are a subsequence $ (n_k)_{k \geq 1}$, $ y _k \in \R^N$ and 
$ \zeta \in \So _q $ (resp. $\in \So_q^\sharp$) such that  
$ d_0 (\phi_{n_{k}} (t_{n_k}) , \, \zeta ( \cdot - y_k) ) \lra 0  $ 
as $ k \lra \infty$, and this contradicts the assumption 
$d_0 ( \psi_n( t_n) , \So _q ) \geq \e _0 $ (resp. 
$d_0 ( \psi_n( t_n) , \So _q^\sharp ) \geq \e _0 $) for all $n$. 
The proof of Theorem \ref{T6.2} is thus complete.
\hfill
$\Box $

\section{Three families of traveling waves} 
\label{3fam}

If the assumptions (A1), (A2) are satisfied and $V\geq 0$ on $[0, \infty)$, Theorem \ref{T4.7} and Proposition \ref{P4.12}
provide finite energy traveling waves to (\ref{1.1}) with any momentum $ q > q_0$; denote by  $ \mathscr{M}$  the family of these traveling waves. 
Theorem \ref{T5.5} and Proposition  \ref{P5.7} provide traveling waves that minimize the action $E - cQ$ at constant kinetic energy; 
let $ \mathscr{K}$ be the family of those solutions. 
If $N=2$, we have also a family $ \mathscr{M}^{\sharp}$ of traveling waves given by Theorem \ref{T5.11}.
Finally, let $\mathscr{P} $ be the family of traveling waves found in \cite{M10}; 
the elements of $\mathscr{P} $ are minimizers of the action 
$E - cQ$ under a Pohozaev constraint (see Theorem \ref{T4.19} below for a precise statement).
Our next goal is to establish  relationships between these families of solutions. 
We will prove that $ \mathscr{M} \subset \mathscr{K}$ and $ \mathscr{K} \subset \mathscr{P}$ if $ N \geq 3$, and that 
 $ \mathscr{M} \subset \mathscr{K}$ and $ \mathscr{M}^{\sharp} \subset \mathscr{K}$ if $ N =2$.
Besides, we find interesting characterizations of the minima of the associated functionals. 

\medskip

Let 
\beq
\label{functionals}
A(\psi ) \!=\! \int_{\R^N} \sum_{j = 2}^N \Big\vert \frac{ \p \psi}{\p x_j }  \Big\vert ^2  dx, 
\quad \;
E_ c ( \psi ) = E(\psi ) - c Q ( \psi), 
\quad \;
P_c(\psi) = E_c( \psi) - \frac{2}{N-1} A(\psi).
\eeq
It follows from Proposition 4.1 p. 1091 in \cite{M7} that any finite-energy 
traveling wave $\psi $ of speed $c$  of (\ref{1.1})  satisfies the Pohozaev 
identity $P_c(\psi) = 0$. Denote
\beq
\label{Tc}
\Co _c= \{ \psi \in \Eo \; | \; \psi \mbox{ is not constant and } P_c(\psi ) = 0 \} 
\quad \mbox{ and } \quad 
T_c = \inf \{ E_c( \psi ) \; | \; \psi \in \Co _c\}.
\eeq 
We summarize below the main results in \cite{M10}. 

\begin{Theorem}\label{T4.19}(\cite{M10}) Assume that $ N \geq 3 $ and (A1) and (A2) hold. Then: 

(i) For any $ c \in (0, v_s)$ the set $ \Co _c $ is not empty and $ T_c > 0 $. 

\medskip

(ii) Let $(\psi_n)_{n \geq 1} \subset \Eo $ be a sequence such that
$$
P_c(\psi_n ) \lra 0 \qquad \mbox{ and } \qquad E_c (\psi _n ) \lra T_c \qquad \mbox{ as } n \lra \infty.
$$
If $ N=3$ we assume in addition that there is a positive constant $ d $ such that 
$$
D(\psi_n) \lra d \quad \mbox{ as } n \lra \infty, \qquad \mbox{ where }
D(\phi ) =  \ii_{\R^N} \Big\vert \frac{\p \phi }{\p x_1} \Big\vert  ^2  + \frac 12  \left(\ph^2(|\phi|) - 1 \right)^2  \, dx .
$$

Then there exist a subsequence $(\psi_{n_k})_{k \geq 1}$, a sequence 
 $(x_k)_{k \geq 1}  \subset \R^N$, 
and $ \psi \in \Co_c  $ such that $E_c ( \psi ) = T_c$, that is, $\psi $ is a minimizer of $E_c$ in $\Co_c  $, 
$\psi_{n _k} ( \cdot + x_k ) \lra \psi $ in $L_{loc}^p(\R^N) $ for $1 \leq p < \infty$ and a.e. on $ \R^N$ and 
$$
\| \nabla \psi_{n _k} ( \cdot + x_k ) - \nabla \psi \|_{L^2(\R^N)} \lra 0 , 
\qquad 
\| \,  | \psi_{n _k} |( \cdot + x_k ) - | \psi  | \,  \|_{L^2(\R^N)} \lra 0 
\qquad \mbox{ as } k \lra \infty.
$$

\medskip

(iii) Let $ \psi $ be a minimizer of $E_c$ in $\Co_c  $. 
Then $\psi $ satisfies (\ref{1.3}) if $N \geq 4$, respectively there exists $ \si > 0 $ such that $\psi_{1, \si }$ satisfies 
(\ref{1.3}) if $N =3$. Moreover,   $ \psi $ (respectively $\psi_{1, \si }$)
is a minimum action solution of (\ref{1.3}), that is it minimizes the action $E_c$ among all finite energy solutions.
Conversely, any minimum action  solution to (\ref{1.3}) is a minimizer of $E_c$ in $ \Co _c$. 

\end{Theorem}

Part (i) is Lemma 4.7 in \cite{M10}, 
part (ii) follows from Theorems 5.3 and 6.2 there and part (iii) follows from Propositions 5.6 and 6.5 
in the same paper and from the fact that any solution $\psi$ satisfies the Pohozaev identity $P_c(\psi) = 0$. 

\begin{remark}\label{R4.20} \rm
As already mentioned in \cite{M10} p. 119, all the conclusions of Theorem \ref{T4.19} above are valid if $ c = 0 $ provided that the set
$ \Co_0 = \{ \psi \in \Eo \; | \; \psi \mbox{ is not constant and } P_0 (\psi ) = 0 \}$ 
is not empty. 
We will see later in section \ref{slow} that $ \Co _0 \neq \emptyset $ if and only if $V$ achieves negative values. 
\end{remark}


\begin{Proposition}
\label{P4.18}
Assume that $N\geq 3$, (A1) and (A2) hold and $V\geq 0$ on $[0, \infty)$. Then:

(i) $ T_c \geq E_{min}(q) - cq $ for any $ q > 0 $ and $ c \in (0, v_s)$. 

(ii) $ T_c \lra \infty $ as $ c \lra 0$. 

(iii) Let $ \psi \in \Eo $ be a minimizer of $E$ under the constraint 
$ Q = q_* > 0$. Assume that $ \psi $ satisfies an Euler-Lagrange equation 
$ E'(\psi ) = c Q'( \psi)$ for some $c \in (0,v_s)$. Then $\psi $ is a 
minimizer of $ E_c$ in $ \Co _c$.
\end{Proposition}

{\it Proof. } For $ \psi \in \Eo$ denote 
\beq
\label{Bc}
B_c( \psi ) =  \ii_{\R^N} \Big| \frac{ \p \psi}{\p x_1 }  \Big| ^2 \, dx - c Q ( \psi ) + \ii_{\R^N} V( |\psi |^2 ) \, dx.
\eeq 
Then $E_c( \psi) = A( \psi ) + B_c( \psi) = \frac{2}{N-1} A(\psi) + P_c( \psi)$ and $P_c( \psi) = \frac{N-3}{N-1} A( \psi) + B_c(\psi)$.

\medskip

i) Consider first the case $ N \geq 4$. Fix $ \psi \in \Co _c$. 
It is clear 
that $A(\psi) > 0 $, hence 
$ B_c( \psi ) = P_c( \psi) - \frac{N-3}{N-1} A( \psi) 
= - \frac{N-3}{N-1} A( \psi) <0$. Since $ V \geq 0 $ by hypothesis, 
it follows that $ c Q( \psi ) = \ii_{\R^N} V( |\psi |^2 ) \, dx 
+ \ii_{\R^N} \big\vert \frac{ \p \psi}{\p x_1 }  \big\vert ^2 \, dx 
- B_c(\psi) > 0$, hence $ Q( \psi ) > 0 $ because $c>0$. It is easy to see that the function 
$ \si \longmapsto E_c( \psi_{1, \si}) = \si^{ N-3} A( \psi) + \si^{N-1} B_c( \psi)$ 
achieves its maximum at $ \si = 1$. Fix $ q > 0$. Since 
$ Q( \psi_{1, \si} ) = \si ^{N-1} Q( \psi)$, there is $ \si _q > 0$ such that 
$ Q(\psi_{1, \si_q }) = q$. We have obviously $ E( \psi _{1, \si _q }) \geq E_{min}(q) $ 
and
$$
E_{min}(q) - cq \leq E( \psi _{1, \si _q }) - c Q(  \psi _{1, \si _q }) 
= E_c( \psi _{1, \si _q }) \leq E_c( \psi_{1,1}) = E_c(\psi). 
$$
Taking the infimum as $ \psi \in \Co _c$, then the supremum as $ q >0$ in the 
above inequality we get $\ds \sup_{q >0} ( E_{min}(q) - cq) \leq T_c.$

Now consider the case $N=3$. Let $ \psi \in \Co _c$. Then 
$P_c(\psi) = B_c( \psi) = 0 $, $Q( \psi ) > 0$  and 
$E_c ( \psi_{1, \si}) = A ( \psi) + \si ^2 B_c( \psi) = A( \psi)$ for any 
$ \si >0$. Fix $ q >0$. Since $Q( \psi_{1, \si }) = \si ^2 Q( \psi)$, there 
is $ \si _q >0$ such that $ Q(\psi_{1, \si_q }) = q$ and this implies 
$ E( \psi _{1, \si _q }) \geq E_{min}(q) $. We have 
$$
E_{min}(q) - cq \leq E( \psi _{1, \si _q }) - c Q(  \psi _{1, \si _q }) 
= E_c( \psi _{1, \si _q }) = A( \psi) = E_c( \psi_{1,1}) = E_c(\psi). 
$$
Since this is true for any $ \psi \in \Co _c$ and any $ q>0$, we conclude 
  again that  $\ds \sup_{q >0} ( E_{min}(q) - cq) \leq T_c.$

\smallskip

(ii) Fix $ q > \frac{1}{v_s}$. We have $E_{min}(q) - cq > E_{min}(q) - 1 $ for 
any $ c \in (0, \frac 1q)$. Using (i) we get
$$
T_c \geq E_{min}(q) - cq > E_{min} (q) - 1 \qquad \mbox{ for any } c \in (0, \frac 1q). 
$$
Since $E_{min}(q) \lra \infty $ as $ q \lra \infty $ by Theorem \ref{T4.16} (b), 
the conclusion follows. 

\smallskip

(iii) We know that $\psi $ is a traveling wave of speed $c$ and by 
Proposition 4.1 p. 1091 in \cite{M7} we have $P_c(\psi ) = 0$, 
that is $ \psi \in \Co _c$. Using (i) we obtain
\beq
\label{4.68}
E_c(\psi) \geq T_c \geq \sup_{q >0} ( E_{min}(q) - cq) .
\eeq
On the other hand, we have
$$
E_c( \psi)= E( \psi) - c Q( \psi) = E_{min}(q_*) - c q_* .
$$
Therefore all inequalities in (\ref{4.68}) have to be equalities. 
We infer that $\psi$ minimizes $ E_c$ in $ \Co _c$,  $ T_c = E_{min}(q_*) - c q_*$ 
and the function $ q \longmapsto E_{min}(q) - cq $ achieves its maximum 
at $ q_*$. 
\hfill
$\Box$

\medskip

The next result shows that the  minimizers of $E_{min}$ or $E_{min}^\sharp$ 
are also minimizers for $I_{min}$ (after scaling).

\begin{Proposition}
\label{P5.9}
Let $ N \geq 2$. Assume that (A1), (A2) hold and either

(a)   $V\geq 0$ on $[0, \infty )$ and  $ q > q_0$, or 

(b)  $N=2$, $ \inf V < 0 $ and  $ q \in (q_0^\sharp, q_{\infty}^\sharp).$

Consider  $ \psi \in \Eo$ such that  $Q(\psi ) = q$ and  $E(\psi) = E_{min}(q) $ in case (a), respectively $E(\psi) = E_{min}^{\sharp }(q) $ in case (b),
and $ \psi $ satisfies (\ref{4.50}) for some $ c \in ( 0, v_s )$
(the  existence of $\psi$ follows from 
Theorem \ref{T4.7} in case (a) and from Theorem \ref{T5.11} in case (b)).
Let $ k = \int_{\R^N} |\nabla \psi |^2 \, dx $.  

\smallskip

Then $\frac{ k }{c^{N-2}} > k_0$ 
and $ \psi_{\frac 1c, \frac 1c}$ is a minimizer of $I$  in the set 
$ \{ \phi \in \Eo \; | \; \int_{\R^N} | \nabla \phi |^2 \, dx =  \frac{k}{c^{N-2}} \}$, that is 
$ I( \psi_{\frac 1c, \frac 1c}) = I_{min} \left( \frac{k}{c^{N-2}} \right) . $ 

Equivalently, $\psi $ is a minimizer of $I_c$ (and of $E_c $) in the set 
$ \{ \phi \in \Eo \; | \; \int_{\R^N} | \nabla \phi |^2 \, dx =k \}$.

Moreover, 
the function $I_{min}$ is differentiable  at $\frac{k}{c^{N-2}}$ and a function $ \zeta \in \Eo$ is a minimizer for $I_{min}\left(\frac{k}{c^{N-2}}\right) $
if and only if  $ \zeta_{c, c}$ is a minimizer for $E_{min} (q)$ and  a traveling wave of speed $c$.
\end{Proposition}

{\it Proof. } By Remark \ref{R5.8} (i) we have $I_{min}(\frac{k}{v_s^{N-2}}) < - \frac{k}{v_s^{N}}$ and 
Proposition \ref{P4.12} (i) implies $ c \in (0, v_s)$, hence 
$ \frac{k}{c^{N-2}} > \frac{k}{v_s^{N-2}} > k_0$.
Using Theorem \ref{T5.5} we infer that there is a minimizer 
$\tilde{\psi} \in \Eo$ of the functional $I$ 
under the constraint $ \int_{\R^N} |\nabla \tilde{\psi} | ^2 \, dx = \frac{k}{c^{N-2}}$. 
By Proposition \ref{P5.7} (ii) there is $ c_1 \in (0, v_s)$ such that $\tilde{\psi}_{c_1, c_1}$ satisfies 
(\ref{4.50}) with $c_1$ instead of $c$.

Let $ \psi _1 = \tilde{\psi}_{c,c}$, so that $\int_{\R^N} |\nabla \psi _1 |^2 \, dx = 
c^{N-2} \int_{\R^N} |\nabla \tilde{ \psi } |^2 \, dx = k= \int_{\R^N} |\nabla  \psi  |^2 \, dx$. 
Denote $ q_1 = Q( \psi _1) = c^{N-1} Q( \tilde{\psi} )$. 

It follows from Proposition 4.1 p. 1091-1092 in \cite{M8} that $\psi $ and $\tilde{\psi}_{c_1, c_1}$ satisfy
the following Pohozaev identities:
\beq
\label{5.19} 
-(N-2)  \int_{\R^N} |\nabla  \psi  |^2 \, dx + c (N-1) Q( \psi) = N  \int_{\R^N} V(|\psi |^2 ) \, dx , 
\eeq
respectively 
$ 
-  (N-2)  \ii _{\R^N} |\nabla \tilde{ \psi }_{c_1, c_1} |^2 \, dx 
+ c _1 (N-1) Q( \tilde{\psi}_{c_1, c_1}) = N  \ii _{\R^N} V(|\tilde{\psi }_{c_1, c_1}|^2 ) \, dx .
$
Since $ \tilde{\psi}_{c_1, c_1} = (\psi _1)_{\frac{c_1}{c}, \frac{c_1}{c}}$, 
the latter equality is equivalent to 
\beq
\label{5.20}
- (N-2) \frac{ c_1^{N-2}}{c^{N-2}} \int_{\R^N} |\nabla  \psi _1 |^2 \, dx 
+  (N-1) \frac{ c_1^N}{c^{N-1}} Q( \psi _1) = N  \frac{ c_1^N}{c^N} \int_{\R^N} V(|\psi _1|^2 ) \, dx .
\eeq

Since $\ii  _{\R^N} |\nabla \psi_{\frac 1c, \frac 1c} |^2 \, dx = \frac{k}{c^{N-2}}
= \ii  _{\R^N} |\nabla \tilde{\psi }|^2 \, dx $
we have $I(\tilde{\psi}) \leq I(\psi_{\frac 1c, \frac 1c})$, that is 
\beq
\label{5.21}
- \frac{1}{c^{N-1}} Q( \psi _1) + \frac{1}{c^N} \int_{\R^N} V(|\psi _1 |^2) \, dx 
\leq - \frac{1}{c^{N-1}} Q( \psi ) + \frac{1}{c^N} \int_{\R^N} V(|\psi  |^2) \, dx .
\eeq
Replacing $\ii _{\R^N} V(|\psi  |^2) \, dx  $ and $\ii _{\R^N} V(|\psi _1 |^2) \, dx  $
from (\ref{5.19}) and (\ref{5.20})  into (\ref{5.21}) we get
\beq
\label{5.22}
cq + (N-2 ) k \leq c q_1 + (N-2) \frac{c^2}{c_1^2} k.
\eeq

Let $\si = \left( \frac{q}{q_1} \right)^{\frac{1}{N-1}}$. 
Then $ Q ( (\psi_1 )_{\si, \si } ) = q$  and consequently
$E(\psi ) \leq E((\psi_1 )_{\si, \si } )$, that is
\beq
\label{5.23}
k + \int_{\R^N} V(|\psi  |^2) \, dx  \leq \si^{N-2} k + \si ^N \int_{\R^N} V(|\psi  _1 |^2).
\eeq
We plug (\ref{5.19}) and (\ref{5.20}) into (\ref{5.23}) to obtain 
\beq
\label{5.24}
c q_1 + (N-2) \frac{c^2}{c_1^2} k \leq N c q_1 - \frac{N-1}{\si ^N} cq 
+ \left( \frac{N}{\si ^2} - \frac{2}{\si ^N} \right) k.
\eeq
Combining (\ref{5.24}) with (\ref{5.22}) we infer that 
$cq + (N-2) k \leq N c q_1 - \frac{N-1}{\si ^N} c q 
+ \left( \frac{N}{\si ^2} - \frac{2}{\si ^N} \right) k.$
Since $ q = \si^{N-1} q_1$, the last inequality can also be written as
\beq
\label{5.25}
\frac{ c q_1}{\si } ( \si ^N - N \si + N -1) 
+ \frac{ k}{\si ^N} ( ( N-2) \si ^N - N \si^{N-2} + 2) \leq 0.
\eeq

If $ N=2$, (\ref{5.25}) is equivalent to $\frac{ cq_1}{\si} ( \si - 1 )^2 \leq 0$ 
and it implies that $ \si = 1$, thus $ q = q_1$. 

If $ N \geq 3 $ we have $ \si ^N - N \si + N -1 = ( \si - 1)^2 \ds \sum_{j = 0}^{N-2} (N-1-j) \si ^j $ and
$$
( N-2) \si ^N - N \si^{N-2} + 2 = (\si - 1)^2 \left[ (N-2) \si^{N-2} + 2 \ds \sum_{j=0}^{N-3} (j+1) \si ^j \right].
$$
Inserting these identities into (\ref{5.25}) and using the fact that 
$ \si, c, q_1, k$ are positive we infer that 
$ \si = 1$, hence $ q = q_1$.  Then using (\ref{5.22}) we obtain $ c_1 ^2 \leq c^2$. 
On the other hand, from (\ref{5.24}) and the fact that $ q = q_1$, $\si = 1$ 
we obtain $ c^2 \leq c_1^2$. 
Since $ c$ and $ c_1$ are positive, we have necessarily $ c = c_1$. 

Since $ q = q_1$  (and  $ c = c_1$ in  the case $ N \geq 3$),  
 using (\ref{5.19}) and (\ref{5.20}) it is easy to see that $I(\psi _{\frac 1c, \frac 1c} ) = I(\tilde{\psi})$, 
hence $ I(\psi _{\frac 1c, \frac 1c} ) = I_{min}(\frac{k}{c^{N-2}})$, as desired. 

Moreover, in the case $ N \geq 3$ we have  proved that {\it any} minimizer $\tilde{\psi }$ of $I$ under 
the constraint $\ii _{\R^N} |\nabla \tilde{\psi}|^2 \, dx = \frac{k}{c^{N-2}}$
satisfies (\ref{5.18}) with $ \vartheta = - \frac{1}{c^2}$. 
It follows from Proposition \ref{P5.7} (iv) that 
$d^+ I_{min}\left( \frac{k}{c^{N-2}} \right) = d^- I_{min}\left( \frac{k}{c^{N-2}} \right) $, hence
$I_{min} $ is  differentiable at $\frac{k}{c^{N-2}}$ and 
$I_{min}' (\frac{k}{c^{N-2}}) = - \frac{1}{c^2}. $

It remains to show that $I_{min}$ is differentiable at $k$ in the case $N=2$. 
 We already  know that $ \psi _{\frac 1c, \frac 1c} $ is a minimizer for $I_{min}(k)$. 
Let $ \phi $ be any other minimizer for $I_{min}(k)$. 
Let $ \phi _ 1 = \phi _{\si, \si}$, where $ \si = \frac{ q}{Q( \phi)}$, so that $ Q( \phi _1 ) = q = Q ( \psi ). $
By (\ref{rajout}) we have 
$$
- I_{min}( k) = \frac{ Q^2 (\psi) }{4 \int_{\R^2} V( |\psi| ^2 )\, dx }  
= \frac{ Q^2 (\phi) }{4 \int_{\R^2} V( |\phi| ^2 )\, dx } 
= \frac{ Q^2 (\phi_1) }{4 \int_{\R^2} V( |\phi_1| ^2 )\, dx } .
$$
We infer that $  \ii_{\R^2} V( |\psi| ^2 )\, dx = \ii_{\R^2} V( |\phi_1| ^2 )\, dx$. 
Since $  \ii_{\R^2} |\nabla \phi _1 |^2 \, dx = k =  \ii_{\R^2} |\nabla \psi |^2 \, dx $ we have $ E( \phi _1 ) = E( \psi ) = E_{min}(k)$, 
hence $ \phi _1 $ is a minimizer for $E_{min}(k)$. 
It follows that there exists $ c_2 \in (0, v_s) $ such that $ \phi _1$ is a traveling wave of speed $ c_2. $
Writing the Pohozaev identities (\ref{5.19}) for $ \psi $ and $ \phi _1$ we find 
$ cq = 2 \ds \int_{\R^2} V( |\psi| ^2 )\, dx $ and $ c_ 2 q = 2 \ds \int_{\R^2} V( |\phi_1| ^2 )\, dx$, respectively. 
Hence $ c = c_2$. 

Proposition \ref{P5.7} (ii)  implies that there is $ c _ 3 $ such that $ \phi _{{c_3}, {c_3}} $ is 
a traveling wave of speed $ c_3$. Using Lemma \ref{L4} below it follows that $ c_3 = c$. 
Since this is true for any minimizer for $I_{min}(k)$, using Proposition \ref{P5.7} (iv) we 
get the desired conclusion. 

The last statement follows easily: if $ \zeta $ is {\it any } minimizer for $I_{min}\left(\frac{k}{c^{N-2}}\right) $ we already know that $ \zeta_{c, c } $ 
is a traveling wave of speed $c$, hence satisfies (\ref{5.19}). Furthermore, $ \ii_{\R^N} |\nabla \zeta_{c, c} |^2 \, dx = k = \ii_{\R^N} |\nabla \psi|^2 \, dx$ and $ I(\zeta) = I_{min}\left(\frac{k}{c^{N-2}}\right) = I(\psi_{\frac 1c, \frac 1c})$, and these equalities clearly imply $ Q( \zeta_{c,c }) = q= Q( \psi)$ and 
$\ii_{\R^N} V(|\zeta_{c, c} |^2) \, dx = \ii_{\R^N} V(|\psi |^2) \, dx$.
\hfill
$\Box$

\medskip

Notice that Proposition \ref{P5.9} does {\it not } imply directly  the differentiability of $I_{min}$ throughout on $(k_0, k_{\infty}).$
For instance, it is possible that for some $q $ there exist two minimizers $ \psi _1, \psi _2 $ for $E_{min}(q)$ with 
$ k_1 = \ii_{\R^N} |\nabla \psi _1 |^2 \, dx <  \ii_{\R^N} |\nabla \psi _2 |^2 \, dx = k_2$ 
and there are no minimizers with kinetic energy between $ k_1$ and  $ k_2$. 
Then $ \psi _1, \psi _2$ are traveling waves with speeds $ c_1, c_2$, respectively. 
By (\ref{5.19}) we have $ c_1 > c_2$ and  Proposition \ref{P5.9} implies that $I_{min} $ is differentiable at $ \frac{k_i}{ c_i ^{N-2}}$, $i= 1,2$, 
but gives no information about the differentiability of  $I_{min} $  on $(  \frac{k_1}{ c_1 ^{N-2}},  \frac{k_2}{ c_2 ^{N-2}}).$
If $N=2$, (A1), (A2), (A4) held with $ F''(1) \neq 3$, $V \geq 0$  and $I_{min}$ were differentiable on $(0, \infty)$, Theorem \ref{T1.2} would give the existence of finite energy traveling waves for any speed  $ c \in (0, v_s)$. 

\begin{Lemma}
\label{L4}
Let $ N \geq 2$. 
 If $ \zeta \in \Eo $ is a traveling wave of speed $c_1$ for (\ref{1.1}), $ Q( \zeta) \neq 0 $ and there is $ \tau > 0 $ 
such that $ \zeta_{\tau, \tau }$ is a traveling wave of speed $ c_2$, then necessarily $ \tau = 1$ and $ c_1 = c_2$. 
\end{Lemma}

{\it Proof. }
Indeed, $ \zeta $ satisfies the equations 
\beq
\label{s}
i c_1 \frac{ \p \zeta }{\p x_1} + \Delta \zeta + F(|\zeta |^2) \zeta = 0 \qquad
\mbox{ and } \qquad
i \frac{c_2}{\tau}  \frac{ \p \zeta }{\p x_1} + \frac{1}{\tau ^2} \Delta \zeta + F(|\zeta |^2) \zeta = 0
\qquad \mbox{ in } \R^N. 
\eeq
If $ \tau = 1$ we get $ ( c_ 2 - c_1) \frac{ \p \zeta}{\p x_1} = 0 $. Since $ \frac{ \p \zeta}{\p x_1} \not\equiv 0 $  
(because $Q (\zeta) \neq 0 $) we have $ c_1 = c_2$ and the Lemma is proven. 
We argue by contradiction and assume that $ \tau \neq 1$. 
Writing the Pohozaev identities corresponding to the $ x_1-$ direction for the two equations in (\ref{s}) (see Proposition 4.1  in \cite{M8})
and using the notation (\ref{functionals})  we find
$$
\int_{\R^N } \Big| \frac{ \p \zeta}{\p x_1} \Big|^2 \, dx - A(\zeta) 
 = \int_{\R^N} V(|\zeta |^2) \, dx 
= \frac{1}{\tau ^2} \left( \int_{\R^N } \Big| \frac{ \p \zeta}{\p x_1} \Big|^2 \, dx - A(\zeta)
\right) .
$$
If $ \tau \neq 1$ we infer  that necessarily 
$ \ii_{\R^N } \Big| \frac{ \p \zeta}{\p x_1} \Big|^2 \, dx = A(\zeta)  
$ 
and  $ \ii_{\R^N} V(|\zeta |^2) \, dx  = 0$. 
Then writing the Pohozaev identities with respect to $(x_2, \dots, x_N)$ we get 
$$
(N-1) \int_{\R^N } \Big| \frac{ \p \zeta}{\p x_1} \Big|^2 \, dx  + (N-3) A(\zeta) 
 - (N-1) c_ 1 Q( \zeta ) = 0, 
$$
 respectively 
 $$ 
\frac{1}{\tau ^2} \left((N-1) \int_{\R^N } \Big| \frac{ \p \zeta}{\p x_1} \Big|^2 \, dx  + (N-3) A(\zeta) 
\right) - (N-1) \frac{ c_2}{\tau } Q( \zeta ) = 0. 
$$
It follows that $ c_2 \tau = c_1$. 
Subtracting the two equations in (\ref{s}) we get $ i c_1 \frac {\p \zeta }{\p x_1} + \Delta \zeta = 0 $ and
 $ F( |\zeta |^2) \zeta = 0 $.
Since $ F(1) = 0 $ and $ F'(1) \neq 0$, there is $ \e > 0 $ such that $ F(s) \neq 0 $ for all 
$ s \in [( 1 - \e )^2, ( 1 + \e )^2 ] \setminus \{ 1 \}.$
Hence $ | \zeta ( x) | \not\in [1-\e,  1+ \e]\setminus \{ 1 \} $ for all $ x \in \R^N$. 
Since $ \zeta $ is continuous and tends to $1$ at infinity, this implies that necessarily 
$ |\zeta | = 1$ in $ \R^N$, and then $ Q( \zeta ) = 0 , $ a contradiction. Lemma  \ref{L4} is thus proven. 
\hfill
$\Box$

\medskip

The next  result establishes the relationship, if  $N\geq 3$, between the traveling 
waves obtained from minimizers of $I_{min}$ and the traveling wave solutions 
given by Theorem \ref{T4.19}.

\begin{Proposition}
\label{P5.13} Assume that  $N \geq 3$ and (A1), (A2) hold.
Let $ \Co _c$ and $ T_c$ be as in (\ref{Tc}). Then: 

(i) $ T_c \geq k + c^N I_{min}\left( \frac{ k}{c^{N-2}} \right) $ for any 
$ k >0 $ and any $ c \in ( 0, v_s)$. 

\medskip

(ii) Let $ \psi $ be a minimizer of $I$ under the constraint 
$ \ii_ {\R^N} |\nabla \psi |^2 \, dx = k$ 
and let $ c \in ( 0, v_s)$ be such that $\psi_{c,c}$ satisfies  (\ref{4.50}).
Then $\psi _{c,c}$  minimizes $E_c = E - cQ$ in $ \Co _c$. 
\end{Proposition}

{\it Proof. } 
We keep the same notation as in the proof of Proposition \ref{P4.18}. 

\medskip

(i) Consider  the case $N \geq 4$. Fix $ \psi \in \Co _c$ and $ k >0$. 
Since $A(\psi) > 0 $, the function 
$ \si \longmapsto \ii _ {\R^N} |\nabla \psi _{1, \si } |^2 \, dx 
= \si ^{N-3}A ( \psi ) + \si^{N-1} \ii _{\R^N} | \frac{ \p \psi}{\p x_1} |^2 \, dx $ is one-to-one from 
$(0, \infty )$  to $(0, \infty )$, 
so there is $ \si _k$ such that $\ii _ {\R^N} |\nabla \psi _{1, \si _k} |^2 \, dx  = k, $
that is $\ii _ {\R^N} |\nabla \psi _{\frac 1c, \frac { \si _k}{c}} |^2 \, dx  = \frac{ k}{c^{N-2}} $. 
This implies $I\left( \psi _{\frac 1c, \frac { \si _k}{c}} \right) \geq I_{min} \left (\frac{ k}{c^{N-2}} \right)$. 
We have $ 0 = P_c( \psi) = A(\psi) + B_c( \psi)$, thus $ A(\psi) > 0 > B_c(\psi ) $ and  the function 
$ \si \longmapsto E_c ( \psi _{1, \si} ) = \si^{N-3} A( \psi ) + \si ^{N-1} B_c( \psi ) $
achieves its maximum at $ \si = 1$. Then we have
$$
\begin{array}{l}
E_c( \psi) = E_c( \psi_{1, 1}) \geq E_c( \psi_{1, \si_k}) 
= \ii _ {\R^N} |\nabla \psi _{1, \si _k} |^2 \, dx  + I_c( \psi_{1, \si _k}) 
\\
\\
= k + c^N I ( \psi _{\frac 1c, \frac { \si _k}{c}})  \geq k + c^N I_{min}\left( \frac{ k}{c^{N-2}} \right).
\end{array}
$$
The above inequality is valid for any $\psi \in \Co _c$ and $ k >0$, hence 
$T_c \geq {\ds \sup_{k >0} } \left( k + c^N I_{min}\left( \frac{ k}{c^{N-2}} \right) \right).$

Next consider the case $N =3$. 
Let $ \psi \in \Co _c $ and let $ k > 0$. 
Then $P_c( \psi) = B_c( \psi) = 0 $ and for any $ \si > 0 $ we have $ E_c(\psi_{1, \si}) = E_c ( \psi ) = A( \psi )$ and 
$\ii _ {\R^3} |\nabla \psi _{1, \si } |^2 \, dx  = A( \psi ) + 
\si ^2  \ii _{\R^3} | \frac{ \p \psi}{\p x_1} |^2 \, dx .$
If $ A( \psi ) \geq k$ we have, taking into account that $I_{min}$ is negative on $(0, \infty)$, 
$$
E_c ( \psi) = A( \psi) \geq k > k + c^3 I_{min}\left( \frac{ k}{c} \right).
$$
If $A(\psi ) < k$, there is $ \si _k > 0 $ such that 
$\ii _ {\R^3} |\nabla \psi _{1, \si _k} |^2 \, dx  = k, $ which means 
$\ii _ {\R^3} |\nabla \psi _{\frac 1c, \frac { \si _k}{c}} |^2 \, dx  = \frac{ k}{c} $.
This implies
 $I_c (\psi _{1, \si _k } ) = c^3  I\left( \psi _{\frac 1c, \frac { \si _k}{c}} \right) 
 \geq c^3 I_{min} \left (\frac{ k}{c} \right)$. 
Thus we get 
$$
E_c( \psi ) = E_c (\psi _{1, \si _k} ) = 
\ii _ {\R^3} |\nabla \psi _{1, \si _k} |^2 \, dx  + I_c (\psi _{1, \si _k} )
\geq k + c^3 I_{min} \left (\frac{ k}{c} \right) .
$$
Hence $ E_c ( \psi ) \geq k + c^3 I_{min} \left (\frac{ k}{c} \right)  $ for any $ \psi \in \Co _c$ and
$ k>0$, and the conclusion follows.

\medskip

(ii) Since $ \psi_{c, c}$ satisfies (\ref{4.50}), 
by Proposition 4.1 p. 1091 in \cite{M7} we have $ \psi_{c,c} \in \Co _c$. Then 
\beq
\label{5.50}
E_c( \psi _{c,c}) \geq T_c \geq \sup_{\kappa >0} 
\left( \kappa + c^N I_{min}\left( \frac{ \kappa}{c^{N-2}} \right) \right).
\eeq
On the other hand, 
$$
E_c( \psi _{c,c}) =  c^{N-2} \ii_{\R^N} |\nabla \psi |^2 \, dx + c^N I( \psi ) 
= c^{N-2} k + c^N I_{min}(k) \leq \sup_{\kappa > 0 } 
\left( \kappa + c^N I_{min}\left( \frac{ \kappa}{c^{N-2}} \right) \right).
$$
Therefore all inequalities in (\ref{5.50}) are equalities, 
$ \psi_{c,c}$ minimizes $ E_c $ in $ \Co _c$, 
$T_c = c^{N-2} k + c^N I_{min}(k)  $ and the function 
$ \kappa \longmapsto \kappa + c^N I_{min}\left( \frac{ \kappa}{c^{N-2}} \right) $ 
achieves its maximum at $ \kappa = c^{N-2} k$. 
\hfill
$\Box$

\section{Small speed traveling waves} 
\label{slow}

Theorem \ref{T4.16} implies that  $ \frac{E_{min}(q)}{q} \lra 0$ as $ q \lra \infty$.  
Since $E_{min}$ is concave and positive, necessarily
$ d^+ E_{min}(q) \lra 0 $ and $  d^- E_{min}(q) \lra 0$  as $ q \lra \infty$ and we infer that 
the traveling waves provided by Theorem \ref{T4.7} and Proposition \ref{P4.12} have speeds close to zero as $ q \ra \infty$. 
Similarly, using Lemma \ref{L5.3} (i) and (iii) we find that $I_{min}$ is finite for all $ k > 0 $ and 
$ d^+ I_{min}(k) \lra - \infty $,  $  d^- I_{min}(k) \lra -\infty$  as $ k \lra \infty$
if either $ N \geq 3$ or ($N=2$ and $V\geq 0$). 
Hence the traveling waves given by Theorem \ref{T5.5} and Proposition \ref{P5.7} 
 have speeds that tend to zero as $ k \lra \infty$. 
This section is a first step in   understanding  the behavior of traveling waves in the limit $ c \lra 0 $.
As one would expect, this
is related to the existence of finite energy solutions to the stationary version of (\ref{1.1}), namely to the equation 
\beq
\label{5.650} 
\Delta \psi + F(|\psi |^2) \psi = 0  \qquad \mbox{ in } \R^N .
\eeq
Clearly, the solutions of (\ref{5.650}) are precisely the critical points of $E$. 
We call {\it ground state } of (\ref{5.650}) a solution that minimizes the energy $E$ among all nontrivial solutions. 

Assume that $ N \geq 2 $ and the assumptions (A1) and (A2) are satisfied. 
Then (\ref{5.650}) admits nontrivial solutions $ \psi \in \Eo $ if and only if the nonlinear potential 
$V$ achieves negative values. 
The existence follows from Theorem 2.1 p. 100 and Theorem 2.2 p. 103 in \cite{brezis-lieb} if $ N \geq 3$, 
respectively from Theorem 3.1 p. 106 in  \cite{brezis-lieb} if $ N =2$.
Moreover, the solutions found in \cite{brezis-lieb} are ground states. 

On the other hand, any solution $ \psi \in \Eo $ of (\ref{5.650}) has the regularity provided by Proposition \ref{P4.12} (ii) 
and this is enough to prove that $ \psi $ satisfies the Pohozaev identity 
\beq
\label{5.66}
(N-2) \ii_{\R^N} |\nabla \psi |^2 \, dx + N \ii_{\R^N} V(|\psi |^2) \, dx = 0 
\eeq
(see Lemma 2.4 p. 104 in \cite{brezis-lieb}). In particular, (\ref{5.66}) implies that (\ref{5.650}) cannot have finite energy solutions if $ V \geq 0$.

We will prove in the sequel that if $N\geq 3$ and $V$ achieves negative values, the traveling waves constructed in this paper 
tend to the ground states of (\ref{5.650}) as their speed goes to zero. 
If $ N \geq 3$,  we have shown in section \ref{3fam} that all traveling waves found here 
also belong to the family of traveling waves given by Theorem \ref{T4.19}, hence it suffices to establish the result for 
the solutions provided by Theorem \ref{T4.19}. 

If $ N=2$ and $V $ takes negative values, we were not able to prove that $d^{\pm} I_{min}(k) \lra -\infty$ as $ k \lra k_{\infty}$. 
Numerical computations  in \cite{CS} indicate that this is indeed the case, at least for some model nonlinearities (including the cubic-quintic one).
If $\ds \lim_{k \uparrow k_{\infty}} d^{\pm} I_{min}(k) = -\infty $, the speeds of the traveling waves given by Theorem \ref{T5.5} and Proposition \ref{P5.7} 
tend to zero as $  k \lra k_{\infty}$ and then we are able to prove a result similar to Proposition \ref{P5.14} below 
(although the proof is very different because minimization under Pohozaev constraints is no longer possible). 

If $V\geq 0$ on $[0, \infty)$, equation (\ref{5.650}) does not have finite energy solutions. 
Then the traveling waves of (\ref{1.1}) have large energy (see Proposition \ref{P4.18} (ii)) and 
are expected to develop vortex structures in the limit $ c \lra 0 $.
This is the case for the traveling waves to the Gross-Pitaevskii equation: in dimension two the solutions found in 
\cite{BS} have two vortices of opposite sign located at a distance of order $\frac 2c$, 
and in dimension three the traveling waves found in \cite{BOS} and \cite{chiron} have vortex rings. 
If $ V\geq 0$, a rigorous description of the behavior of traveling waves in the limit $ c \lra 0 $ 
is still missing. 

\begin{Proposition}
\label{P5.14} 
Let $ N \geq 3$. Suppose that (A1) and (A2) are satisfied and there exists $ s_0 \geq 0 $ such that $V( s_0 ) < 0$. 
Let $ (c_n)_{n \geq 1} $ be any sequence of numbers in $(0, v_s)$ such that $ c_n \lra 0 $. 
For each $ n $, let $ \psi _n \in \Eo$ be any minimizer of $ E_{c_n} = E - c_n Q$ in $ \Co _{c_n}$ such that 
$\psi _n$ is a traveling wave of (\ref{1.1}) with speed $ c_n $.   Then: 

(i) There are a subsequence $(c_{n_k})_{k \geq 1}$, a sequence $ (x_k ) _{k \geq 1} \subset \R^N$ 
and a ground state $ \psi $ of 
(\ref{5.650}) such that $\psi_{n _k} ( \cdot + x_k ) \lra \psi $ in $ L_{loc}^p( \R^N) $ for $ 1 \leq p < \infty $ and a.e. on $ \R^N$ and 
$$
\| \nabla \psi_{n _k} ( \cdot + x_k ) - \nabla \psi \|_{L^2(\R^N)} \lra 0 , 
\qquad 
\| \,  | \psi_{n _k} |( \cdot + x_k ) - | \psi  | \,  \|_{L^2(\R^N)} \lra 0 
\qquad \mbox{ as } k \lra \infty.
$$

(ii) There is a sequence $(a_k)_{k \geq 1}$ of complex numbers of modulus $1$ such that $a_k \lra 1$ as $ k \lra \infty$ and
$$
\| a_k \psi_{n_k} ( \cdot + x_k) - \psi \|_{W^{2, p }(\R^N)} \lra 0 \qquad \mbox{ as } k \lra \infty \quad \mbox{ for any } p \in [2^*, \infty).
$$
In particular, 
$\| a_k \psi_{n_k} ( \cdot + x_k) - \psi \|_{C^{1, \al }(\R^N)} \lra 0 $ as $ k \lra \infty $  for any $ \al \in [0,1).$
\end{Proposition}


If $F $ is $C^k$ it can be proved that the convergence in (ii) holds in $W^{k+2, p }(\R^N)$, $ 2^* \leq p < \infty.$

\medskip

{\it Proof.  }
(i) Let $ \psi _0 $ be any ground state of (\ref{5.650}). By (\ref{5.66}) we have 
$ \ii_{\R^N} V(|\psi _0 |^2) \, dx = - \frac{N-2}{N} \ii_{\R^N} |\nabla \psi |^2 \, dx <0$. 
It is shown in  \cite{brezis-lieb} that $ \psi _0$ is a minimizer of the functional $ J(\phi ) = \ii_{\R^N} |\nabla \phi|^2 \, dx $ 
subject to the constraint $ \ii_{\R^N} V(|\phi  |^2) \, dx = \ii_{\R^N} V(|\psi _0 |^2) \, dx$; 
conversely, any minimizer of this problem is a ground state to of (\ref{5.650}), and Proposition \ref{P4.12} (ii) 
implies that any minimizer is $C^1$ on $ \R^N$. 
It follows from Theorem 2 p. 314 in \cite{M7} that any ground state of  (\ref{5.650}) is, after translation, radially symmetric. 
In particular, the radial symmetry implies that  $Q(\psi _0 ) = 0$.

Let $A$, $E_c = E - cQ$, $P_c$ be as in (\ref{functionals}) and $ \Co _c $ and $ T_c$ as in ({\ref{Tc}).
Since $ \psi _0$ is a solution of (\ref{5.650}), it satisfies the Pohozaev identity $P_0 (\psi _0) = 0 $ and then we get 
$ P_c( \psi _0) = P_0 (\psi _0) - c Q( \psi _0) = 0$ for any $c$, that is $ \psi _0 \in \Co _c $ for any $c$. 
Therefore
\beq
\label{5.67}
A( \psi _n) = \frac{N-1}{2} \left( E_{c_n}(\psi _n) - P_{c_n}(\psi_n) \right) 
= \frac{N-1}{2} E_{c_n}(\psi _n) =  \frac{N-1}{2} T_{c_n} \leq \frac{N-1}{2} E_{c_n}(\psi _0) 
= A( \psi _0).
\eeq
On the other hand, by Proposition 10 (ii) in \cite{CM} the function $ c \longmapsto T_c$ is decreasing on $(0, v_s)$.
Fix $ c_* \in (0, v_s)$. For all sufficiently large $n$ we have $ c_n < c^*$, hence
\beq
\label{5.68}
A(\psi _n) =  \frac{N-1}{2} T_{c_n} \geq  \frac{N-1}{2} T_{c_*} > 0.
\eeq

Consider first the case $N\geq 4$. We claim that $E_{GL}(\psi_n)$ is bounded. 
To see this we argue by contradiction and we assume that there is a subsequence, still denoted $(\psi_n)_{n\geq1}$, such that 
$E_{GL}(\psi_n)\lra \infty$. 
By (\ref{5.67}) we have
\beq
\label{5.69}
D(\psi_n) = \ii_{\R^N} \Big| \frac{\p \psi_n}{\p x_1} \Big| ^2  + \frac 12  \left( \ph ^2 (|\psi_n|)-1 \right)^2 \, dx \lra \infty \qquad \mbox{ as } n \lra \infty.
\eeq
Using Lemma \ref{L4.2} (ii) we see that there are two positive constants $ k_0, \ell _0$ such that 
for any $ \psi \in \Eo$ satisfying $ E_{GL}(\psi) = k_0 $ and for any $c \in (0, c_*) $ 
(where $c_*$ is as in (\ref{5.68})) there holds 
\beq
\label{5.70}
E_c(\psi) \geq E(\psi) - c |Q(\psi)| \geq \ell_0.
\eeq
It is easy to see that for each $n$ there is $ \si _n >0$ such that 
\beq
\label{5.71}
E_{GL}((\psi_n)_{\si_n, \si_n}) = 
\si_n ^{N-3} A( \psi_n) + \si_n ^{N-1} D(\psi_n)
= k_0.
\eeq

In particular, $(\psi_n)_{\si_n, \si_n} $ satisfies (\ref{5.70}).

We recall that the functional $B_c$ was defined in (\ref{Bc}). We have 
$ B_{c_n}(\psi_n) = P_{c_n}( \psi_n) - \frac{N-3}{N-1}A(\psi_n)$.  
Then the fact that $P_{c_n}(\psi _n) = 0 $ and  (\ref{5.67}) imply that $ B_{c_n}(\psi_n)$ is bounded.
From (\ref{5.69}) and (\ref{5.71}) it follows that $ \si _n \lra 0 $ as $ n \lra \infty$, 
hence
$$
E_{c_n} ( (\psi _n)_{\si_n, \si_n}) = \si_n ^{N-3} A( \psi_n) + \si_n ^{N-1} B_{c_n} (\psi_n) \lra 0 \qquad \mbox{ as } n \lra \infty.
$$
This contradicts the fact that $ E_{c_n} ( (\psi _n)_{\si_n, \si_n}) \geq \ell_0$ for all $n$ and the claim is proven.

Using Corollary \ref{C4.17} we infer that $Q(\psi_n)$ is bounded. 
Since $ c_n \lra 0 $,  using (\ref{5.67}) we find 
\beq
\label{5.72}
P_0(\psi_n) = P_{c_n}(\psi_n) + c_n Q(\psi_n) \lra 0 \qquad \mbox{ and } 
\eeq
\beq
\label{5.73}
\begin{array}{l}
E(\psi_n) = E_{c_n}(\psi_n) + c_n Q(\psi_n)= \frac{2}{N-1}A(\psi_n) + P_{c_n}(\psi_n) + c_n Q(\psi_n)
\\
\\
\leq \frac{2}{N-1}A(\psi_0) + c_n Q(\psi_n)= E( \psi _0 ) +  c_n Q(\psi_n) \lra E(\psi _0)
\qquad \mbox{ as } n \lra \infty.
\end{array}
\eeq
Then the  conclusion follows from Theorem \ref{T4.19} (with $c=0$)  and Remark \ref{R4.20}.

\medskip

Next consider the case $ N = 3$. For all $n$ and all $ \si > 0$ we have 
$$
P_{c_n} ((\psi_n)_{1, \si}) = \si^2 P_{c_n} ( \psi_n ) = 0 
\; \;  \mbox{ and } \; \; 
E_{c_n} ((\psi_n)_{1, \si}) = A ((\psi_n)_{1, \si}) + P_{c_n} ((\psi_n)_{1, \si})  = A(\psi_n) = T_{c_n}, 
$$
hence $(\psi_n)_{1, \si} $ is also a minimizer of $ E_{c_n}$ in $ \Co_{c_n}$.
For each $n$ there is $ \si _n > 0 $ such that $D((\psi_n)_{1, \si_n}) = \si _n ^2 D(\psi_n) = 1$. 
We denote $\tilde{\psi}_n= (\psi_n)_{1, \si_n}$.
Then $\tilde{\psi}_n $ is a minimizer of $E_{c_n}$ in $ \Co_{c_n}$, 
$E_{GL} (\tilde{\psi}_n) = A(\tilde{\psi}_n) +1 = A(\psi_n) + 1$ is bounded  by (\ref{5.67}) 
and then Corollary \ref{C4.17} implies that $Q(\tilde{\psi}_n) $ is bounded. 
As in the case $N\geq 4$ we find that $(\tilde{\psi}_n)_{n\geq 1}$ satisfies (\ref{5.72}) and (\ref{5.73}).
From Theorem \ref{T4.19}  and Remark \ref{R4.20} it follows that there exist a subsequence 
$(\tilde{\psi}_{n_k})_{k\geq 1}$, a sequence $ (x_k)_{k \geq 1} \subset \R^3$ and a minimizer $\tilde{\psi}$ of $ E$ in $ \Co _0$ that 
satisfy the conclusion of Theorem \ref{T4.19} (ii). 
Moreover, there is $ \si >0$ such that $\tilde{\psi}$ satisfies the equation 
\beq
\label{5.74}
 \frac{\p^2 \tilde{\psi}}{\p x_1 ^2 } 
+ \si ^2 \sum_{j=2}^3 \frac{\p ^2\tilde{\psi}}{\p x_j ^2 } + F(|\tilde{\psi}|^2) \tilde{\psi} = 0 
\qquad \mbox{ in } \Do '(\R^3). 
\eeq
Let $ \psi_k ^* = \tilde{\psi}_{n_k} ( \cdot + x_k). $
Since $ \psi _n$ solves (\ref{1.3}) with $ c_n$ instead of $c$,  it is obvious that $ \psi _k^*$ satisfies 
\beq
\label{5.75}
i c_{n_k} \frac{\p \psi _k^*}{\p x_1 } + \frac{\p^2 \psi _k^* }{\p x_1 ^2 } 
+ \si _{n_k} ^2 \sum_{j=2}^3 \frac{\p ^2\psi _k^*}{\p x_j ^2 } + F(|\psi _k^*|^2) \psi_k^* = 0 
\qquad \mbox{ in } \Do '(\R^3). 
\eeq
It is easy to see that $ \psi_k^* \lra \tilde{\psi} $ and $ F(|\psi _k^*|^2) \psi_k^* \lra F(|\tilde{\psi}|^2) \tilde{\psi}$ 
in $ \Do'(\R^3)$.

We show that $ (\si_{n_k})_{k \geq 1}$ is bounded. 
We argue by contradiction and we assume that it contains a subsequence, still denoted the same, that tends to $\infty$. 
Multiplying (\ref{5.75}) by $\frac{1}{ \si _{n_k} ^2}$ and passing to the limit as  $ k \lra \infty$ we get
\beq
\label{5.76}
 \frac{\p ^2\tilde{\psi}}{\p x_2 ^2 } + \frac{\p ^2\tilde{\psi}}{\p x_3 ^2 } = 0 \qquad \mbox{ in } \Do '(\R^3).
\eeq
Since $ \frac{\p ^2 \tilde{\psi}}{\p x_j \p x_k} \in L_{loc}^p(\R^3)$ for any $ p \in [1, \infty)$, we infer that the above 
equality holds in $ L_{loc}^p(\R^3)$ for any $ p \in [1, \infty)$. 
By the Sobolev embedding (see Lemma 7 and Remark 4.2 p. 774-775 in \cite{PG}) we know that there is $ \al \in \C$ such that $ |\al | = 1$ and $ \tilde{\psi} - \al \in L^6(\R^3)$. 
Let $ \chi \in C_c^{\infty} ( \R^3) $ be a cut-off function such that $ \chi = 1 $ on $B(0,1)$ and $supp(\chi) \subset B(0,2)$. 
Taking the scalar product (in $\C$) of (\ref{5.76}) by $ \chi(\frac xn) (\psi - \al)$ and letting $ n\lra \infty$ we find
$\ii_{\R^3} \big| \frac{\p \tilde{\psi}}{\p x_2} \big| ^2 +   \big| \frac{\p \tilde{\psi}}{\p x_3} \big| ^2 \, dx = 0.$
Since $ \tilde{\psi } \in C^{1, \al}(\R^3)$, we conclude that $  \frac{\p \tilde{\psi}}{\p x_2} = \frac{\p \tilde{\psi}}{\p x_3} = 0$, 
hence $ \tilde{\psi }$ depends only on $x_1$. Together with the fact that $ \frac{\p \tilde{\psi}}{\p x_1} \in L^2( \R^3)$ 
this implies that $\tilde{\psi}$ is constant, a contradiction. Thus $ (\si_{n_k})_{k \geq 1}$ is bounded.

If there is a subsequence $( \si_{n_{k_j}})_{j \geq 1}$ such that $  \si_{n_{k_j}} \lra \si _*$  as $ j \lra \infty$,
passing to the limit in (\ref{5.75}) we discover 
$$
 \frac{\p^2 \tilde{\psi}}{\p x_1 ^2 } 
+ \si _* ^2 \sum_{j=2}^3 \frac{\p ^2\tilde{\psi}}{\p x_j ^2 } + F(|\tilde{\psi}|^2) \tilde{\psi} = 0 
\qquad \mbox{ in } \Do '(\R^3). 
$$
If $ \si _* \neq \si$, comparing the above equation to (\ref{5.74}) we find 
$\frac{\p ^2\tilde{\psi}}{\p x_2 ^2 } + \frac{\p ^2\tilde{\psi}}{\p x_3 ^2 } = 0 $ in $ \Do'(\R^3)$ 
and arguing as previously we infer that $\tilde{\psi}$ is constant, a contradiction. 
We conclude that necessarily $\si_{n_k}\lra \si $ as $ k \lra \infty$. 
Denoting $ \psi = \tilde{\psi}_{1, \frac{1}{\si}}$, we easily see that $ \psi $ minimizes $ E$ in $ \Co _0$ and is a ground state of 
(\ref{5.650}).  Then  $(\psi_{n_k})_{k \geq 1}$ and  $ \psi $ satisfy the conclusion of Proposition \ref{P5.14} (i).

\medskip

(ii) By the Sobolev embedding there are  $ \al, \, \al _k \in \C$ of modulus $1$ and $ C_S >0 $ such that 
$$
\| \psi_{n_k} - \al _k \|_{L^{2^*}(\R^N)} \leq C_S \| \nabla \psi_{n_k} \|_{L^{2}(\R^N)}   \qquad \mbox{ and } \qquad 
\| \psi - \al  \|_{L^{2^*}(\R^N)} \leq C_S \| \nabla \psi \|_{L^{2}(\R^N)} .
$$
We may assume that $ \al =1$ for otherwise we multiply $ \psi_{n_k}$ and $ \psi $ by $ \al ^{-1}$.
(In fact we have $\psi = \al \psi_0$, where $ \psi _0$ is real-valued, but we do not need this observation.)

Let $ R > 0 $ be arbitrary, but fixed. By (i) there exists $ k(R) \in \N$ such that for all $ k \geq k(R)$ we have 
$\| \psi_{n_k} ( \cdot + x_k) - \psi \|_{L^{2^*} (B(0,R))} < 1$. Then we find
$$
\| \al _k - 1\| _{L^{2^*} (B(0,R))} \leq \| \psi_{n_k} ( \cdot + x_k) -  \al _k \|_{L^{2^*} (\R^N )}
+ \| \psi_{n_k} ( \cdot + x_k) - \psi \|_{L^{2^*} (B(0,R))}  
+ \| \psi - 1 \| _{L^{2^*} (\R^N )} \leq C 
$$
for any $ k \geq k(R)$, where $C$ does not depend on $k$. This implies that $ \al _k \lra 1$. 

Let $ \psi _k ^* = \al _k^{-1} \psi_{n_k} ( \cdot + x_k) $, so that  
 $ \psi _k ^* - \psi \in L^{2^*} (\R^N)$.  Using (i) and the Sobolev embedding  we get 
\beq
\label{5.77}
\| \psi_k ^* - \psi \|_{L^{2^*} (\R^N )} \leq C_S \| \nabla \psi _k ^* - \nabla \psi \|_{L^{2} (\R^N )} \lra 0 \qquad \mbox{ as } k \lra \infty. 
\eeq
By (i), $ \nabla \psi _{k}^* $ is bounded in $ L^2( \R^N)$ and $ \psi _k^*$ is a traveling wave to (\ref{1.1}) of speed $ c_{n_k}$. 
It follows from   Step 1 in the proof of Lemma \ref{L7.1} below that there is $L>0$, independent of $k$, such that 
$$
\| \nabla \psi _k ^* \| _{L^{\infty}(\R^N)} \leq L \qquad \mbox{ and } \qquad \| \nabla \psi  \| _{L^{\infty}(\R^N)} \leq L.
$$
By interpolation we get 
\beq
\label{5.78}
\| \nabla \psi_k ^* - \nabla \psi \|_{L^{p} (\R^N )} \lra 0 \qquad \mbox{ as } k \lra \infty \quad \mbox{ for any } p \in [2, \infty).
\eeq
Using (\ref{5.77}), (\ref{5.78}) and the Sobolev embedding we infer that 
\beq
\label{5.79}
\| \psi_k ^* - \psi \|_{L^{p} (\R^N )} \lra 0 \qquad \mbox{ as } k \lra \infty \quad \mbox{ for any } p \in [2^*, \infty].
\eeq

We claim that $ \|F(|\psi_k ^*|^2) \psi _k ^*  - F(|\psi |^2) \psi \|_{L^{p} (\R^N )} \lra 0 $ 
as $  k \lra \infty  $ for any $ p \in [2^*, \infty).$
To see this fix $ \de > 0 $ such that $ F$ is $ C^1 $ on $[1 - 2 \de, 1 + 2 \de]$ (such $ \de $ exists by (A1)). 
Since $ \psi - 1 \in L^{2^*}( \R^N)$ and $ \| \nabla \psi  \| _{L^{\infty}(\R^N)} \leq L$
 we have $ \psi \lra 1 $ as $ |x | \lra \infty$, hence there exists $ R(\de) > 0 $ verifying $ \big| \, |\psi | - 1 \big| < \de$
on $ \R^N \setminus B(0, R(\de))$. 
By (\ref{5.79}) there is $ k_{\de} \in \N$ such that $\| \psi_k ^* - \psi \|_{L^{\infty} (\R^N )} < \de $ for $ k \geq k_{\de}$. 
The mapping $ z \longmapsto F(|z|^2) z $ is Lipschitz on $\{ z \in \C \; | \; 1 - 2 \de \leq |z| \leq 1 + 2 \de \}$, 
hence there is $C>0$ such that
\beq
\label{5.80}
\big| F(|\psi_k ^*|^2) \psi _k ^*  - F( |\psi |^2) \psi \big| \leq C | \psi_k^* - \psi | \qquad
\mbox{ on } \R^N \setminus B(0, R(\de)) \mbox{ for all } k \geq k_{\de}.
\eeq
Since $  F(|\psi_k ^*|^2) \psi _k ^*  - F(|\psi |^2) \psi $ is bounded and tends a.e. to zero, 
using  Lebesgue's dominated convergence theorem we get
\beq
\label{5.81}
\|F(|\psi_k ^*|^2) \psi _k ^*  - F(|\psi |^2) \psi \|_{L^{p} (B(0, \de))} \lra 0
\qquad \mbox{ for any } p \in [1, \infty).
\eeq
Now the claim follows from (\ref{5.79}) - (\ref{5.81}). 

Using the equations satisfied by $ \psi_k^*$ and $\psi$ we get 
$$
\Delta ( \psi_k^* - \psi ) = - i c_{n_k} \frac{ \p \psi_k^*}{\p x_1} - \left( F(|\psi_k ^*|^2) \psi _k ^*  - F(|\psi |^2) \psi \right). 
$$
From the above we infer that $\| \Delta ( \psi_k^* - \psi) \|_{L^p(\R^N)} \lra 0 $ for any $ p \in [2^*, \infty)$, then using 
(\ref{5.79}) and the inequality $\| f \|_{W^{2,p}( \R^N ) } \leq C_p \left( \| f \|_{L^p( \R^N ) } + \| \Delta f \|_{L^p( \R^N ) } \right)$
we get the desired conclusion. 
\hfill
$\Box$

\section{Small energy traveling waves} 
\label{sectionsmallenergy}

The aim of this section is to prove Proposition \ref{smallE}.
The next lemma shows that 
the modulus of  traveling waves  
of small energy is close to $1$. 


\begin{Lemma}
\label{L7.1} Let $ N \geq 2$. 
Assume that (A1) and ((A2)  or (A3)) hold. 

(i) For any $ \e >0 $ there exists $ M(\e ) > 0$ such that for any $c \in [0, v _s]$  and for 
any solution $ \psi \in \Eo$ of (\ref{1.3}) with $\| \nabla \psi \|_{L^2( \R^N) } < M(\e) $ we have
\beq
\label{7.1}
\big| \, | \psi (x) | - 1 \big| < \e \qquad \mbox{ for all } x \in \R^N.
\eeq

(ii) Let $ p > N p_0$, where $ p_0 $ is as in (A2) (respectively $ p \geq 1$ if (A3) is satisfied).
 For any $ \e >0 $ there exists $ \ell _p(\e ) > 0$ such that for any 
$c \in [0, v _s]$  and for 
any solution $ \psi \in \Eo$ of (\ref{1.3}) with $\| \, | \psi | - 1  \|_{L^p( \R^N) } < \ell_p(\e) $, 
(\ref{7.1}) holds. 
\end{Lemma}

{\it Proof.} 
Assume first that (A1) and (A2) are satisfied. 

We will prove that there is $L>0$ such that any solution $ \psi \in \Eo $ of (\ref{1.3}) such that 
$\| \nabla \psi \|_{L^2( \R^N) }$ is sufficiently small 
(respectively $\| | \psi | - r_0  \|_{L^2( \R^N) } $ is sufficiently small) 
satisfies 
\beq
\label{7.2}
\| \nabla \psi \|_{L^{\infty}(\R^N)} \leq L.
\eeq

{\it Step 1. } We prove (\ref{7.2}) if $ N \geq 3$ and $\| \nabla \psi \|_{L^2( \R^N) } \leq M$, 
where $M>0$ is fixed. 

Using the Sobolev embedding, for any $ \phi \in \Eo $ such that  $\| \nabla \phi \|_{L^2( \R^N)} \leq M$ 
we get
$$
\| ( |\phi | - 2  )_+ \| _{L^{2^*}( \R^N) } \leq C_S \| \nabla | \phi | \, \|_{L^2( \R^N)} 
\leq C_S \| \nabla \phi  \|_{L^2( \R^N)}.
$$
Since $ |\phi| \leq 2  + ( |\phi | - 2  )_+, $
we see that $ \phi $ is bounded in $ L^{2^*} + L^{\infty}(\R^N)$.
It follows that for any $ R > 0 $ there exists $C_{R, M } > 0 $ such that for any 
$ \phi \in \Eo $ as above we have
$$
\| \phi \|_{H^1( B(x, R )) } \leq C_{R, M} \qquad \mbox{ for all } x \in \R^N.
$$

If $ c \in [0, v_s]$, $ \psi \in \Eo $ is a solution of (\ref{1.3})   and  $\| \nabla \psi \|_{L^2( \R^N) } \leq M$, 
using (\ref{3.11}) and a standard bootstrap argument (which works thanks to (A2)) 
we infer that for any $ p \in [2, \infty)$ 
 there is $ \tilde{C}_p > 0 $ (depending only on $F$, $N$, $p$ and $M$) such that 
$$
\| \psi \|_{W^{2,p}( B(x, 1 ))} \leq \tilde{C}_p  \qquad \mbox{ for all } x \in \R^N.
$$
Then the Sobolev embedding implies that $ \psi \in C^{1, \al }(\R^N) $ for all $ \al \in [0, 1)$ 
and there is $ L>0$ such that (\ref{7.2}) holds.

\medskip

{\it Step 2.} Proof of (i) in the case $N \geq 3$. 

Fix $ \e > 0$.
There is $L>0$ such that any solution $ \psi \in \Eo $ of (\ref{1.3}) 
with $\| \nabla \psi \| _{L^2( \R^N)} \leq 1$ satisfies (\ref{7.2}). 
If $ \psi $ is such a solution and $ \big|\, | \psi (x_0) | - 1 \big| \geq \e $ for some $ x _0 \in \R^N$, 
from (\ref{7.2}) we infer that 
$ \big|\, | \psi (x) | - 1 \big| \geq \frac{\e}{2} $ for any $ x \in B(x_0, \frac{ \e}{2L})$. 
Then using the Sobolev embedding we get 
$$
C_S \| \nabla \psi \| _{L^2( \R^N)} \geq \| \, | \psi | - 1 \|_{L^{2^*}(\R^N)} 
\geq  \| \, | \psi | - 1 \|_{L^{2^*}(B(x_0, \frac{\e}{2L}))} 
\geq \frac{ \e}{2} \left( \left(\frac{ \e}{2 L} \right) ^N \Lo ^N (B(0,1)) \right)^{\frac{1}{2^*}}.
$$
We conclude that if 
$\| \nabla \psi \| _{L^2( \R^N)} < \min \left(\! 1, 
\frac{ \e}{2 C_S} \left( \left(\frac{ \e}{2 L} \right) ^N \Lo ^N (B(0,1)) \right)^{\frac{1}{2^*}} \right)$, 
then 
$ \psi $ satisfies~(\ref{7.1}). 

\medskip

{\it Step 3.} Proof of (\ref{7.2}) if $N=2$ and $\| \nabla \psi \| _{L^2( \R^2)} $ is sufficiently small. 

By (\ref{4.2}) there is $ M_1 >0$ such that for any $ \phi \in \Eo $ 
with $ \| \nabla \phi \|_{L^2( \R^2)}  \leq M_1$ we have
\beq
\label{7.3}
\frac 14 \int_{\R^2} \left( \ph^2( |\phi|) - 1 \right) ^2 \, dx 
\leq \int_{\R^2} V(| \phi |^2) ) \, dx 
\leq \frac 34 \int_{\R^2} \left( \ph^2( |\phi|) - 1 \right) ^2 \, dx .
\eeq
Let $ \psi \in \Eo $ be a solution of (\ref{1.3}). 
By Proposition \ref{P4.12} (ii) we have $ \psi \in W_{loc}^{2, p} ( \R^2) $ and this regularity is enough to 
prove that $ \psi $ satisfies the Pohozaev identity
\beq
\label{7.4}
- \int_{\R^2} \bigg\vert \frac{ \p \psi}{\p x_1} \bigg\vert ^2 \, dx
+ \int_{\R^2} \bigg\vert \frac{ \p \psi}{\p x_2} \bigg\vert ^2 \, dx
+ \int_{\R^2} V(| \psi |^2) ) \, dx  = 0
\eeq
(see Proposition 4.1 p. 1091 in \cite{M8}).
In particular, if $ \| \nabla \psi \|_{L^2( \R^2)}  \leq M_1$ by (\ref{7.3}) and (\ref{7.4}) we get  
\beq
\label{7.5}
  \int_{\R^2} \left( \ph^2( |\psi|) - 1 \right) ^2 \, dx 
\leq 4 \int_{\R^2} V(| \psi |^2) ) \, dx 
\leq 4 \int_{\R^2} \bigg\vert \frac{ \p \psi}{\p x_1} \bigg\vert ^2 \, dx \leq 4 M_1
\eeq
and Corollary \ref{C4.2} implies that there is some $ M_2 > 0 $ (independent of $\psi $) such that 
$\| \, | \psi | - 1 \| _{L^2( \R^2)} \leq M_2$. 
We infer that for any $ R > 0 $ there is $ M_3(R) > 0 $ (independent of $\psi $) such that 
$ \| \psi \| _{H^1 (B(x, R) ) } \leq M_3(R)$ and hence, by the Sobolev embedding, 
$ \| \psi \| _{L^p (B(x, R) ) } \leq C_p (R)$ for all $ x \in \R^2$ and 
$ p \in [2, \infty)$. 
 Using (\ref{3.11}) and an easy bootstrap argument we get 
$ \| \psi \| _{W^{2,p} (B(x, 1) ) } \leq \tilde{C}_p $ for all $ x \in \R^2$ and 
$ p \in [1, \infty)$. 
As in Step 1 we conclude that there is $ L>0$ such that any solution $ \psi \in \Eo $ of (\ref{1.3}) with 
$ \| \nabla \psi \|_{L^2( \R^2)}  \leq M_1$ satisfies (\ref{7.2}). 

\medskip

{\it Step 4.} Proof of (i) if $N=2$. 

Fix $ \e > 0 $. Let $ \eta $ be as in (\ref{3.19}) and $M_1$ as in step 3.  
If $ \psi \in \Eo $ is a solution of (\ref{1.3}) with $\| \nabla \psi \|_{L^2( \R^2)} \leq M_1 $ 
and there is $ x_0 \in \R^2$ such that $ \big|\, | \psi (x_0) | - 1 \big| \geq \e $, 
using (\ref{7.2}) we infer that 
$ \big| \, | \psi (x) | - 1 \big| \geq \frac{\e}{2} $ for any $ x \in B(x_0, \frac{ \e}{2L})$, hence 
$ \left( \ph^2(|\psi |) - 1 \right)^2 \geq \eta( \frac{\e}{2}) $ on $  B(x_0, \frac{ \e}{2L})$
and therefore
$$
\int_{\R^2} 
 \left( \ph^2(|\psi |) - 1 \right)^2  \, dx 
 \geq \int_{B(x_0, \frac{ \e}{ 2L} ) }  \left( \ph^2(|\psi |) - 1 \right)^2 \, dx 
\geq \pi \left( \frac{ \e}{ 2L}  \right)^2 \eta \left( \frac{\e}{2} \right) .
$$
On the other hand, by (\ref{7.5}) we have 
$  \ii_{\R^2} 
 \left( \ph^2(|\psi |) - 1 \right)^2  \, dx  \leq 4 \| \nabla \psi \|_{L^2( \R^2)} ^2$. 
We conclude that necessarily $\big| \, | \psi | - 1 \big| < \e $ on $ \R^2$ if 
$\| \nabla \psi \|_{L^2( \R^2)} ^2 < \frac {\pi}{4}  \left( \frac{ \e}{ 2L}  \right)^2 \eta( \frac{\e}{2}).$

\medskip

{\it Step 5.} Proof of (\ref{7.2}) if $\| \, |\psi | - 1 \|_{L^p(\R^N)} \leq M$  and $ p > Np_0$.

By Proposition \ref{P4.12} (ii) we know that $ \psi $ and $ \nabla \psi $ belong to $ L^{\infty}(\R^N)$. 
We will prove that $\| \psi \|_{L^{\infty} (\R^N)}$ and $\| \nabla \psi \|_{L^{\infty} (\R^N)}$ 
are bounded uniformly with respect to $\psi$. 
The constants $ C_j$ below depend only on $ M, F, p, N$, but not on $\psi$. 

Let $ \phi (x) = e^{\frac{i c x_1}{2}} \psi (x)$, 
so that $ |\phi | = |\psi |$ and $ \phi $ satisfies the equation 
\beq
\label{7.6}
\Delta \phi + \left( \frac{ c^2}{4} + F(|\phi|^2) \right) \phi = 0 \qquad \mbox{ in } \R^N.
\eeq
For all $ x \in \R^N$ we have $ \| \phi \|_{L^p(B(x,2))} \leq C_1$, where $C_1$ depends only 
on $M$. Fix $ r = (\frac{p}{ 2 p_0} )^-$ such that $ N p_0 < 2r p_0 < p$ and 
$ (2 p_0 + 1) r > p$. In particular, we have $ r > \frac N2 \geq 1$. Since 
$ \big\vert \left( \frac{ c^2}{4} + F(|\phi|^2) \right) \phi \big\vert \leq C_2 + C_3 |\phi |^{2 p_0 + 1}$, 
using (\ref{7.6}) we find that for all $ x \in \R^N$ we have
\beq
\label{7.7}
\| \Delta \phi \| _{L^r(B(x,2))} 
\leq C_4 + C_5 \| \phi \|_{L^{\infty} (B(x,2))} ^{2 p_0 + 1 - \frac pr } \| \phi \|_{L^{p} (B(x,2))} ^{ \frac pr } 
\leq C_6 + C_7 \| \phi \|_{L^{\infty} (\R^N)} ^{2 p_0 + 1 - \frac pr }.
\eeq
It is obvious that $ |\phi | \leq C_8 + C_9 \| \phi \|^{ 2 p_0 + 1 - \frac pr }_{L^\infty(\R^N)} 
| \phi |^{ \frac pr }$, hence $ \phi $ satisfies 
$$ 
\| \phi \| _{L^r(B(x,2))} \leq (\mathcal{L}^N (B(0,2) )^{\frac 1r} C_8 
+ C_9 C_1^{ \frac pr } \| \phi \|^{ 2 p_0 + 1 - \frac pr }_{L^\infty(\R^N)}  .
$$ 
Then, using (\ref{3.11}) we infer that for all $ x \in \R^N$, 
$$
\| \phi \|_{W^{2,r}(B(x,1))} \leq C_{10} + C_{11} \| \phi \|_{L^{\infty} (\R^N)} ^{2 p_0 + 1 - \frac pr }.
$$
Since $ r > \frac N2$, the Sobolev embedding implies 
$ \| \phi \|_{L^{\infty} (B(x,1))} \leq C_s \| \phi \|_{W^{2,r}(B(x,1))}.$
Choose $ x _0 \in \R^N$ such that $ \| \phi \|_{L^{\infty} (B(x_0,1))} \geq 
\frac 12 \| \phi \|_{L^{\infty} (\R^N)} .$ We have 
$$
\frac{1}{2 C_s} \| \phi \|_{L^{\infty} (\R^N)} 
\leq \frac{1}{ C_s} \| \phi \|_{L^{\infty} (B(x_0,1))} 
\leq \| \phi \|_{W^{2,r}(B(x_0,1))}
\leq C_{10} + C_{11} \| \phi \|_{L^{\infty} (\R^N)} ^{2 p_0 + 1 - \frac pr }.
$$
Since $ 2 p_0 + 1 - \frac pr < 1 $ by the choice of $ r$, 
the above inequality implies that there is $ C_{12} > 0$ such that $\| \phi \|_{L^{\infty} (\R^N)}   \leq C_{12}$. 
Then using (\ref{7.6}) and (\ref{3.11}) we infer that 
$  \| \phi \|_{W^{2,q}(B(x,1))} \leq C(q)$ for all $ x \in ( \R^N)$ and all $ q \in (1, \infty)$, 
and the Sobolev embedding implies 
$\| \nabla \phi \|_{L^{\infty} (\R^N)}   \leq C_{13}$ for some $ C_{13} > 0$.
Since $ \psi (x)= e^{- \frac{ i c x_1}{2}} \phi (x)$, 
the conclusion follows.

\medskip

{\it Step 6.} Proof of (ii). 

Let $ \psi $ be a solution of (\ref{1.3}) such that $\| \, |\psi | - 1 \|_{L^p( \R^N)} \leq 1$. 
By step 5, there is $ L > 0$ (independent of $ \psi $) such that $(\ref{7.2})$ holds. 
If there is $ x_0 \in \R^N$ such that $\big|\, |\psi( x_0)| - 1 \big| \geq \e$, 
we have $\big|\, |\psi | - 1 \big| \geq \frac{\e }{2} $ on $ B ( x_0, \frac{\e}{2 L})$ and consequently 
$$
\| \, |\psi | - 1 \|_{L^p( \R^N)}  
\geq \| \, |\psi | - 1 \|_{L^p( B(x_0, \frac{\e}{2L}))}  
\geq \frac{ \e}{2} \left( \left(\frac{ \e}{2L} \right)^N \Lo ^N ( B(0,1)) \right) ^{\frac 1p}.
$$
Thus necessarily $\big| \, |\psi( x)| - 1 \big| < \e$ on  $\R^N$ if 
$\| \, |\psi | - 1 \|_{L^p( \R^N)}  < \min \left( 1 , \frac{ \e}{2} \left( \left( \frac{ \e}{2L} \right)^N \Lo ^N ( B(0,1)) \right) ^{\frac 1p} \right).$

\medskip

If (A1) and (A3) hold,  it follows from the proof of Proposition 2.2 (i) p. 1078 in \cite{M8} that 
there is $ L > 0$ such that (\ref{7.2}) holds for any $ c \in [0, v_s]$ and any solution $ \psi  \in \Eo $ of (\ref{1.3}).
Therefore the conclusions of steps 1, 3 and 5 are automatically satisfied. 
The rest of the proof is exactly as  above.
\hfill 
$\Box$ 

\medskip

By (A1) we may fix $ \beta _* > 0 $ such that 
$ \frac 14 ( s - 1) ^2 \leq V( s) \leq \frac 34 ( s - 1) ^2$ 
if $ | \sqrt{ s} - 1 | \leq \beta_*$. 

Let $U\in \Eo $ be a traveling wave to (\ref{1.1}) such that 
 $ 1 - \beta _*  \leq |U | \leq 1 + \beta _* $. It is clear that 
\beq
\label{equival}
\frac 14 ( |U |^2 - 1) ^2 \leq V( |U|^2) \leq \frac 34 ( |U|^2 - 1) ^2 \qquad \mbox{ on } \R^N.
\eeq
It is an easy consequence of
Theorem 3 p. 38 and of Lemma C1 p. 66 in \cite{BBM} that there exists a lifting 
$ U = \rho e^{i \theta} $ on $ \R^N$, where 
$ \rho, \theta \in W_{loc}^{2,p}(\R^N) $ for any $p \in [1, \infty)$. 
Then (\ref{1.3}) can be written in the form 
\beq
\label{syst}
\left\{ 
\begin{array}{l}
\ds \Delta \rho - \rho |\nabla \theta|^2 + \rho F(\rho^2) 
= c \rho \frac{\partial \theta}{\partial x_1},
\\
\\
\ds \mbox{div} ( \rho^2 \nabla \theta) 
= - \frac{c}{2} \frac{\partial }{\partial x_1}(\rho^2 - 1 ) .
\end{array}
\right.
\eeq
Multiplying the first equation in (\ref{syst}) by $\rho$ we get 
\beq
\label{7.9}
\frac 12 \Delta ( \rho ^2 - 1) - |\nabla U|^2 + \rho ^2 F( \rho ^2) - c ( \rho ^2 - 1) \frac{ \p \theta }{\p x _1} = c  \frac{ \p \theta }{\p x _1} .
\eeq
The second equation in (\ref{syst}) can be written as 
\beq
\label{7.10}
\mbox{div} (( \rho ^2 - 1) \nabla \theta ) + \frac{c}{2} \frac{\partial }{\partial x_1}(\rho^2 - 1 )
= -  \Delta \theta. 
\eeq
We set $ \eta = \rho^2 - 1 $ 
and define $ g : [-1 , +\infty ) $ by 
$ g(s) =v_s^2 s + 2 ( 1 + s) F(1 + s )  $, so that 
$ g(s) = \mathcal{O}(s^2) $ for $ s \lra 0 $. 
Taking the Laplacian of (\ref{7.9}) and applying the operator $ c \frac{ \p }{\p x_1}$  to (\ref{7.10}), 
then summing up the resulting equalities  we find  
\beq
\label{7.11}
\left[ \Delta^2 \!- v_s^2 \Delta + c^2 \partial_{x_1}^2 \right] \eta = 
\Delta \left( 2 |\nabla U|^2 \!- g(\eta) + 2 c \eta \partial_{x_1} \theta \right) 
- 2 c \partial_{x_1} (\mbox{div}  ( \eta \nabla \theta)) 
\quad \mbox{ in } \So '( \R^N).
\eeq
Notice that 
 the right-hand side of (\ref{7.11}) contains terms that are 
(at least) quadratic. 
We write (\ref{7.11}) using the Fourier transform as 
\beq
\label{magique}
\hat{\eta} (\xi) = \mathcal{L}_c (\xi) \hat{\Upsilon}(\xi) ,
\eeq
where
\beq
\label{7.14}
 \wh{\Upsilon}(\xi) = - \mathcal{F} ( 2 |\nabla U|^2 - g(\eta) ) 
- 2 c \frac{|\xi|^2-\xi_1^2}{|\xi|^2} \mathcal{F}(\eta \p_{x_1} \phi) 
+ 2c \sum_{j=2}^N \frac{\xi_1 \xi_j}{|\xi|^2} \mathcal{F} ( \eta \p_{x_j} \phi ) 
\eeq
and
\beq
\label{7.15}
 \mathcal{L}_c (\xi) =  \frac{|\xi|^2}{ |\xi|^4 + v_s^2 |\xi|^2 - c^2 \xi_1^2} .
\eeq

On the other hand, we know that $U$ satisfies the Pohozaev identity (\ref{5.19}).     
Using (\ref{lift}) and the Cauchy-Schwarz identity we have
$$
| Q( U) | = \Big\vert \ii_{\R^N} (\rho ^2 - 1) \theta _{x_1}   \, dx \Big\vert
\leq \| \eta \|_{L^2( \R^N)} \| \theta_{x_1}  \|_{L^2( \R^N)} 
\leq \frac{1}{ 1 - \beta_* }  \| \eta \|_{L^2( \R^N)}    \| \nabla U \|_{L^2( \R^N)}  .
$$
Inserting this estimate into 
(\ref{5.19}), using (\ref{equival}) and the fact that 
$|c| \leq v_s $ we get 
\beq
\label{7.16}
(N-2) \| \nabla U \|_{L^2( \R^N )}^2 
- \frac{(N-1)v_s}{ 1 - \beta_* }  \| \eta \|_{L^2( \R^N)}    \| \nabla U \|_{L^2( \R^N)} 
+ \frac{N}{4} \| \eta \|_{L^2( \R^N)}  ^2 \leq 0. 
\eeq

{\bf The case $\boldsymbol{N \geq 3}.$  }
If $ N \geq 3$, let $ a_1 \leq a_2$ be the two roots of the equation 
$(N-2) y^2 - \frac{(N-1)v_s}{ 1 - \beta_* } y + \frac{N}{4}  = 0.$ 
It is obvious that $ a_1$ and $a_2$ are positive and from (\ref{7.16}) we infer that 
\beq
\label{7.17}
a_1  \| \eta \|_{L^2( \R^N)} \leq \| \nabla U \|_{L^2( \R^N )} \leq a_2 \| \eta \|_{L^2( \R^N)}.
\eeq

\medskip

{\it Proof of Proposition  \ref{smallE}  for $N \geq 3$}. 
We use the ideas introduced in \cite{BGS} and \cite{dL}.

In the following $C_j$ and $ K_j$ are positive constants depending only on $N$ and $F$.

Let $ \beta_*$ be as above.
By Lemma \ref{L7.1}, there are $ M_1, \, \ell_1 > 0 $ such that any solution $ U \in \Eo $ to (\ref{1.3}) with 
$\ii_{\R^N} |\nabla U|^2 \, dx \leq M_1$ 
(respectively with $\ii_{\R^N} \left(|U|^2 - 1 \right)^2\, dx \leq \ell _1$ if (A3)
 holds or if (A2) holds and $ p_0 < \frac 2N$) satisfies $ 1 - \beta_* \leq |U| \leq 1 + \beta _*$ and, in addition, (\ref{7.2}) is verified.
Then we have a lifting $U = \rho e^{i \theta } $ and (\ref{equival})-(\ref{7.17}) hold.
Since $ g ( \eta) = \Oo( \eta ^2)$, it follows from (\ref{7.17}) that 
\beq
\begin{array}{rcl}
\label{7.18}
\| 2 |\nabla U|^2 - g ( \eta ) \| _{L^1( \R^N)} 
& \leq &  2 \| \nabla U \|_{L^2( \R^N)} ^2 + C_1 \| \eta \|_{L^2(\R^N)} ^2 
\leq C_2  \| \eta \|_{L^2(\R^N)} ^2  
\\
& \leq & C_3  \| \nabla U \|_{L^2( \R^N)} ^2. 
\end{array}
\eeq
On the other hand, from   $ 1 - \beta_* \leq |U| \leq 1 + \beta_*$  and (\ref{7.2})  we get 
$ \| 2 |\nabla U|^2 - g ( \eta ) \| _{L^{\infty}( \R^N)}  \leq C_4$ and then, by interpolation, 
\beq 
\label{7.19}
\| 2 |\nabla U|^2 - g ( \eta ) \| _{L^p( \R^N)} 
\leq C_1 (p)  \| \eta \|_{L^2(\R^N)} ^{\frac 2p} , 
\eeq
 respectively 
\beq
\label{7.20}
\| 2 |\nabla U|^2 - g ( \eta ) \| _{L^p( \R^N)} 
\leq K_1 (p)  \| \nabla U \|_{L^2(\R^N)} ^{\frac 2p} 
\eeq
for any $ p \in [1, \infty)$.
It is obvious that 
$| \eta \p_{x_j} \theta | \leq \frac{1}{1 - \beta_* } |\eta | \cdot |\nabla U | $ and, as above, we find 
\beq
\label{7.21}
\| \eta \p_{x_j} \theta \| _{L^p( \R^N)} 
\leq C_2 (p)  \| \eta \|_{L^2(\R^N)} ^{\frac 2p} ,
\quad \mbox{ and } \quad 
\| \eta \p_{x_j} \theta \| _{L^p( \R^N)} 
\leq K_2 (p)  \| \nabla U  \|_{L^2(\R^N)} ^{\frac 2p} .
\eeq
By the standard theory of Riesz operators (see, e.g., \cite{stein}), the functions 
$ \xi \longmapsto \frac{ \xi_j \xi_k}{|\xi |^2}$ are Fourier multipliers from $L^p(\R^N)$ to $L^p(\R^N)$, 
$1 < p < \infty$.
Using (\ref{7.14}) and (\ref{7.19})-(\ref{7.21}) we infer that 
$\Upsilon \in L^p(\R^N)$ for $1 < p < \infty$ and
\beq
\label{7.22}
\| \Upsilon \| _{L^p( \R^N)} 
\leq C_3 (p)  \| \eta \|_{L^2(\R^N)} ^{\frac 2p} ,
\quad \mbox{ respectively } \quad 
\| \Upsilon \| _{L^p( \R^N)} 
\leq K_3 (p)  \| \nabla U  \|_{L^2(\R^N)} ^{\frac 2p} .
\eeq

We will use the following result, which is  Lemma 3.3 p. 377 in \cite{dL} with 
$ \alpha = \frac{2}{2N-1} $ and $q=2$. Notice that $ \frac{1}{1-\alpha} = \frac{2N-1}{2N-3} <2 $ 
if $ N \geq 3 $.

\begin{Lemma} (\cite{dL})
\label{multiplier}
Let  $N \geq 3$ and let $p_N = \frac{2(2N-1)}{2N+3} \in (1 , 2) $.
There exists a constant $K_N$, depending only on $N$, 
such that for any $ c \in [ 0 ,v_s ] $ and  any $f \in L^{p_N}(\R^N)$ we have
$$
\| \mathcal{F}^{-1} (\mathcal{L}_c (\xi) \mathcal{F}(f) ) \|_{L^2(\R^N)} 
\leq K _N\| f \|_{L^{p_N}(\R^N)} .
$$
\end{Lemma}

From (\ref{magique}), Lemma \ref{multiplier} and (\ref{7.22}) we get 
\beq
\label{7.23}
\| \eta \|_{L^2(\R^N)} \leq K_N \| \Upsilon \|_{L^{p_N}(\R^N)} \leq K _NC_3(p_N) \| \eta \|_{L^2(\R^N)} ^{\frac{2}{p_N}}.
\eeq
Since $ \frac{2}{p_N} > 1$, (\ref{7.23}) implies that there is $ \ell_* > 0 $ (depending only on $N$ and $F$)
such that  $ \| \eta \|_{L^2(\R^N)}  \geq \ell_*$, or $  \| \eta \|_{L^2(\R^N)}  = 0$.
In the latter case from (\ref{7.17})  we get $ \| \nabla U \|_{L^2(\R^N)} = 0$, hence $U$ is constant.

From (\ref{7.23}) and (\ref{7.17})  we obtain 
$$
 \| \nabla U \|_{L^2(\R^N)}  \leq a_2 \| \eta \|_{L^2(\R^N)}  
\leq 
a_2  K_N C_3(p_N) \| \eta \|_{L^2(\R^N)} ^{\frac{2}{p_N}}
\leq a_1 ^{- \frac{2}{p_N}} a_2  K _NC_3(p_N) \| \nabla U \|_{L^2(\R^N)} ^{\frac{2}{p_N}}.
$$
As above we infer that there is $ k_* > 0$ such that either $ \| \nabla U \|_{L^2(\R^N)}   \geq k_* $, or $U$ is constant.
\hfill
$\Box$

\medskip

{\bf The case $\boldsymbol{N = 2}.$ }
If $N=2$, 
from (\ref{7.16}) we infer that 
$\| \eta\|_{L^2(\R^N)} \leq \frac{2 v_s}{1 - \beta_*} \| \nabla U \|_{L^2(\R^N)} $.
However, 
the Pohozaev identities alone do not imply an estimate of the form 
$ \| \nabla U \|_{L^2(\R^N)} \leq C \| \eta\|_{L^2(\R^N)}$.
To prove this we need the following two  identities, which are valid in any space dimension   and are of independent interest. 

\begin{Lemma}
\label{L7.3}
Let $ U = \rho e^{ i \theta }  \in \Eo $ be a solution of (\ref{1.3}), where $ \inf \rho >0$ and $ \rho $ is bounded.
Then we have 
\beq
\label{grincheux}
2 \int_{\R^N} \rho^2 |\nabla \theta |^2 \ dx = 
- c \int_{\R^N} ( \rho^2 - 1 ) \p_{x_1} \theta \ dx 
\qquad \mbox{ and } 
\eeq
\beq
\label{dormeur}
\int_{\R^N} 2 \rho |\nabla \rho|^2 
+ \rho (\rho^2 -1) |\nabla \theta|^2 
- \rho (\rho^2 -1) F(\rho^2) \ dx 
= - c \int_{\R^N} \rho (\rho^2 -1) \p_{x_1} \theta \ dx .
\eeq
\end{Lemma} 

{\it Proof. } 
Formally, $U$ is a critical point of the functional $E_c = E - cQ$. 
Denoting $U(s) = \rho e^{ i s \theta}$ one would expect that 
$\frac{d}{ds} _{\mid s=1}  (E_c (U(s)) = 0$ and this  is precisely (\ref{grincheux}).

In the case of the Gross-Pitaevskii equation, (\ref{grincheux}) was proven  in \cite{BGS}  
(see Lemma 2.8 p. 594 there) 
by multiplying the second equation in (\ref{syst}) by $\theta$, then integrating by parts. 
The integrations are justified by the particular decay  at infinity of traveling waves for the Gross-Pitaevskii equation. 
Since such decay properties have not been rigorously established for other nonlinearities, we proceed as follows.

For $R > 0 $, we denote 
$ \bar{\theta} = \frac{1}{\mathcal{H}^{N-1} (\p B(0,R) )} \ii_{\p B_R} \theta \ d \mathcal{H}^{N-1} ,$
 we multiply the second  equation  in (\ref{syst})
by $\theta - \bar{\theta}$ and
 integrate by parts over $B(0,R)$. We get
\beq
\label{7.26}
\begin{array}{l}
\ds 2 \ii_{B(0,R)} \rho^2 |\nabla \theta |^2 \ dx 
- 2 \ii_{\p B(0,R)} \rho^2  \frac{\p \theta}{\p \nu}  (\theta - \bar{\theta}) \ d\mathcal{H}^{N-1} \\ 
\\
= \ds - c \ii_{B(0,R)} ( \rho^2 - 1 ) \p_{x_1} \theta \ dx 
+ c \int_{\p B(0,R)} (\rho^2-1 ) (\theta - \bar{\theta}) \nu_1\ d\mathcal{H}^{N-1} ,
\end{array}
\eeq
where $ \nu $ is the outward unit normal to $ \p B(0, R)$.
By the Poincar\'e inequality we have  for some constant $C$ 
independent of $R$,
$$ \| \theta - \bar{\theta} \|_{L^2(\p B(0,R))} \leq 
C R \| \nabla \theta \|_{L^2(\p B(0,R))} . $$
Using the boundedness of $ \rho $ and the Cauchy-Schwarz inequality we have for $ R \geq 1$
\begin{align*}
\Big| 2 \ii_{\p B(0,R)} \rho^2 (\theta - \bar{\theta}) \frac{\p \theta}{\p \nu} 
 \ d\mathcal{H}^{N-1} \Big| 
+ \Big| c \ii_{\p B(0,R)} (\rho^2-1) (\theta - \bar{\theta}) \nu_1\ d\mathcal{H}^{N-1} \Big| 
\\ 
\leq C R \ii_{\p B(0,R)} ( \rho^2 - 1 )^2 
+ | \nabla \theta |^2 \ d \mathcal{H}^{N-1} . 
\end{align*}
Since $ \rho^2 - 1 \in L^2(\R^N) $ and $\nabla \theta \in L^2(\R^N)$, we have
$$ \ii_1^{+\infty} \left( \int_{\p B(0,R)} ( \rho^2 - 1 )^2 
+ |\nabla \theta|^2 \ d\mathcal{H}^{N-1} \right) \ dR = 
\ii_{ \{ |x| \geq 1 \} } ( \rho^2 - 1 )^2 + |\nabla \theta|^2 \ dx < \infty  , $$
hence  there exists a  sequence $ R_j \lra + \infty $ such that
$$ \ii_{\p B (0, R_j)} ( \rho^2 - 1 )^2 + |\nabla \theta|^2 \ d\mathcal{H}^{N-1} 
\leq \frac{1}{R_j \ln R_j} . $$
Writing (\ref{7.26}) for each $j$, then  passing to the limit  as $ j \lra \infty$   
we obtain \eqref{grincheux}.

It is easily seen that $ \rho ^2 - 1 \in H^1( \R^N)$. 
 Multiplying the first equation in 
\eqref{syst} by $\rho^2 -1 $ and  using the standard integration by parts formula for 
$H^1$ functions (cf. \cite{brezis} p. 197) we get \eqref{dormeur}.
\hfill
$\Box$

\medskip

Using \eqref{grincheux}  and the Cauchy-Schwarz inequality we get
$$ 2 \ii_{\R^N} \rho^2 |\nabla \theta |^2 \ dx = 
- c \ii_{\R^N} ( \rho^2 - 1 ) \p_{x_1} \theta \ dx 
\leq C \left( \ii_{\R^N} ( \rho^2 - 1 )^2 \ dx \right)^{1/2} 
\left( \ii_{\R^N} \rho^2 |\nabla \theta |^2 \ dx \right)^{1/2} , $$
from which it comes
$$ \ii_{\R^N} \rho^2 |\nabla \theta |^2 \ dx \leq C 
\int_{\R^N} \eta^2 \ dx . $$

Using  \eqref{dormeur}, the fact that 
 $ 0 < 1 - \beta_* \leq \rho \leq 1 + \beta_*$,
the inequality $ 2 ab \leq a^2 + b^2$ and the above estimate we find
\begin{align*}
2( 1 - \beta_*) \ii_{\R^N} |\nabla \rho|^2 \ dx 
\leq & \ \ii_{\R^N} 2 \rho |\nabla \rho|^2  \ dx \\
= & \ - \int_{\R^N} \rho \eta |\nabla \theta |^2 \ dx 
- c \int_{\R^N} \eta  \rho \p_{x_1} \theta  \ dx 
+ \int_{\R^N} \rho^2 \eta F(\rho^2) \ dx \\ 
\leq & \ C \int_{\R^N} \rho^2 |\nabla \theta |^2 + \eta^2 \ dx \leq 
C \int_{\R^N} \eta^2 \ dx .
\end{align*}

It follows from the above inequalities that in the case $N = 2$, there exist two positive constants $a_1, \, a_2$ such that any solution $ U \in \Eo $ to (\ref{1.3})
with $ 0 \leq c \leq v_s$ and $ 1 - \beta_* \leq |U| \leq 1 + \beta_*$ satisfies (\ref{7.17}).

\medskip 

{\it Proof of Proposition \ref{smallE} if $N=2$,  (A4) holds and $ F''(1) = 3$.}
The  strategy used in the case $N \geq 3$ has to be adapted: small energy traveling waves 
do exist when $N=2$ and $ F''(1) \not = 3 $ (see Theorem \ref{T4.7}, Proposition \ref{P4.12} 
and Theorem \ref{T4.13}). This is related to 
the fact that Lemma \ref{multiplier} does not apply if $N=2$. The  proof  relies on an expansion in the small parameter 
$\eta$ and the observation  that when the energy is small, we must have 
$  \p_{x_1} \phi \simeq - c \eta / 2 $. 
Since $ v_s^2 = 2= - 2  F'(1) $ and 
$  F''(1) = 3 $, by (A4)   the function $g$ has the expansion as $ s \lra 0 $
\begin{align*}
g(s) = & \, v_s^2 s + 2 (1 +s) F (1 +s) 
= v_s^2 s + 2 (1 +s) \Big( s F'(1 ) + \frac12 s^2 F''(1 ) 
+ \mathcal{O} (s^3) \Big) 
\\ 
= & \,  s^2 + \mathcal{O} (s^3) .
\end{align*}

By Lemma \ref{L7.1}, there are $ M_1, \, \ell_1 > 0 $ such that any solution $ U \in \Eo $ to (\ref{1.3}) with $ c \in [0, v_s]$ and 
$\ii_{\R^2} |\nabla U|^2 \, dx \leq M_1$ (respectively $\ii_{\R^2} \left(|U|^2 - 1 \right)^2\, dx \leq \ell _1$ 
if (A3) holds or if (A2) holds and $ p_0 <  1$) satisfies $ 1 - \beta_* \leq |U| \leq 1 + \beta _*$, the estimate  (\ref{7.2}) is verified,
 we have a lifting $U = \rho e^{i \theta } $ and 
all the statements above are valid.

Recalling that $\Upsilon$ is defined by (\ref{7.14}), 
we observe that in the  expression of $ 2 |\nabla U|^2 - g(\eta) $ 
we have the almost cancellation of two quadratic terms: 
$ 2 \rho^2 (\p_{x_1} \theta)^2 -  \eta^2 
\simeq 2 ( ( \p_{x_1} \theta)^2 - \frac{v_s^2}{4} \eta^2 ) $ 
is much smaller than quadratic if $  \p_{x_1} \theta \simeq v_s \eta / 2 $. 
We now quantify this idea and split the proof into 7 steps. We denote
\beq
\label{7.27}
 h =  \p_{x_1} \theta + \frac c2 \eta . 
\eeq
By Lemma \ref{L7.1},  $ \| \eta \|_{L^\infty(\R^2)}$ can be made arbitrarily small 
by taking  $M_1$ (respectively $\ell _1$) sufficiently small.
 Moreover,  using (\ref{7.17}) we get 
\beq
\label{tookiki}
\| \eta \|_{L^4(\R^2)}^4 \leq 
\| \eta \|_{L^\infty(\R^2)}^2 \| \eta \|_{L^2(\R^2)}^2 \leq C  \| \eta \|_{L^2(\R^2)}^2
\leq C  \| \nabla U \|_{L^2(\R^2)}^2 . 
\eeq

\medskip 

{\it Step 1.} There is $ C >0$, depending only on $F$, 
such that if $M_1$ (respectively $ \ell_1$) is small enough,
$$ \ii_{\R^2} h^2 + (\p_{x_2} \theta)^2 + ( v_s ^2 - c^2)  \eta^2 \ dx 
\leq C \| \eta \|^4_{L^4(\R^2)} . $$

The starting point is the integral identity
$$
\ii_{\R^2} \rho^2 |\nabla \theta|^2 + V(\rho^2) + c (\rho^2 - 1) \p_{x_1} \theta
\ d x = 0 , 
$$ which comes from the combination of (\ref{grincheux}) and the 
Pohozaev identity 
$ \ii_{\R^2} 2 V(\rho^2) + c (\rho^2 - 1) \p_{x_1} \theta \ dx = 0 $
(see Proposition 4.1 in \cite{M8}).
From (A4) with $ F''(1) = 3 $  we have the Taylor expansion of the potential 
$$
V (\rho^2 ) = V (1+ \eta ) 
= \frac{v_s^2}{4} \eta^2 - \frac16 F''(1) \eta^3 + \mathcal{O}(\eta^4) 
= \frac{v_s^2}{4} \eta^2 - \frac{v_s^2}{4 } \eta^3 + \mathcal{O}(\eta^4)  .
$$ 
Therefore, the above integral identity gives
$$
 \ii_{\R^2} ( 1 + \eta ) (\p_{x_2} \theta)^2 
+  (\p_{x_1} \theta)^2 + \eta (\p_{x_1} \theta)^2 
+ \frac{v_s^2}{4} \eta^2 - \frac{v_s^2}{4 } \eta^3 + \mathcal{O}(\eta^4) 
+ c \eta \p_{x_1} \theta
\ d x = 0 .
$$
Then the identity $ h^2 =  (\p_{x_1} \theta )^2 + c \eta \p_{x_1} \theta
+ \frac{c^2}{4 } \eta^2 $ gives
$$
 \ii_{\R^2} ( 1 + \eta ) (\p_{x_2} \theta)^2 
+ h^2 + \frac{v_s^2-c^2}{4 } \eta^2 
+ \eta (\p_{x_1} \theta)^2 
- \frac{v_s^2}{4 } \eta^3 
\ d x = - \ii_{\R^2} \mathcal{O}(\eta^4) \ dx ,
$$
hence, rearranging the cubic terms,
\beq
\label{eqrase}
 \ii_{\R^2} ( 1 + \eta ) (\p_{x_2} \theta)^2 + h^2 
+ \frac{v_s^2-c^2}{4 } \eta^2 \left( 1 - \eta \right) \ dx 
= - \ii_{\R^2} \eta  h \left( h - c \eta \right) 
+ \mathcal{O}(\eta^4) \ dx .
\eeq
For the left-hand side, we have $ 1 + \eta \geq \frac 12 $ and 
$ 1 - \eta \geq \frac12 $ if $M_1$ or $ \ell _1$  are sufficiently small (because 
$ \| \eta \|_{L^\infty(\R^2)} $ is small). We now estimate the right-hand side. 
Since $ \| \eta \|_{L^\infty(\R^2)}$ is small, we have 
$ | \ii_{\R^2} \mathcal{O}(\eta^4) \ dx | \leq C \| \eta \|_{L^4(\R^2)}^4 $ 
and by Cauchy-Schwarz and the inequality $ 2 a b \leq a^2 + b^2 $,
\begin{align*}
\Big| \ii_{\R^2} \eta  h \left( h - c \eta \right) \ dx \Big| 
\leq & \, \| \eta \|_{L^\infty(\R^2)}  \ii_{\R^2} h^2 \ dx 
+ c  \Big( \ii_{\R^2} h^2 \ dx \Big)^{\frac12} 
\Big( \ii_{\R^2} \eta^4 \ dx \Big)^{\frac12} \\
\leq & \,  \frac{1}{2} \ii_{\R^2} h^2 \ dx + C \| \eta \|^4_{L^4(\R^2)} ,
\end{align*}
provided that $M_1 $ or $ \ell _1$ are small enough,  where $C$ depends only on $F$. Inserting 
these estimates into (\ref{eqrase}) yields the result.

\medskip

{\it Step 2.} There exists $ C $, depending only on $F$, 
such that for $M_1$ (respectively $\ell _1$)  small enough,
$$ 
\int_{\R^2} | \nabla \rho |^2 + ( v_s^2 - c^2 ) \eta^2 \ dx 
\leq C \| \eta \|^2_{L^4(\R^2)} .
 $$
We start from (\ref{dormeur}), that we write in the form
$$
\int_{\R^2} 2 \rho |\nabla \rho|^2 \ dx 
= - \int_{\R^2} \rho \eta \left( ( \p_{x_1} \theta)^2 + ( \p_{x_2} \theta)^2 \right) 
- \rho \eta F(\rho^2) + c \rho \eta \p_{x_1} \theta \ dx .
$$
Using the expansion $ F(\rho^2) = \eta F'(1) + \mathcal{O}(\eta^2) 
= - \frac{v_s^2 \eta}{2 } + \mathcal{O}(\eta^2) $, this gives
\beq
\label{eqlabousse}
\int_{\R^2} 2 \rho |\nabla \rho|^2 + \frac{v_s^2 - c^2}{2  } \rho \eta^2 \ dx 
= - \int_{\R^2} \rho \eta \left( ( \p_{x_1} \theta)^2 + ( \p_{x_2} \theta)^2 \right) 
+ c \rho \eta h + \mathcal{O} ( |\eta|^3 )\ dx .
\eeq
Note that by the Cauchy-Schwarz inequality,
\beq
\label{eta3}
\| \eta \|_{L^3(\R^2)}^3 \leq \| \eta \|_{L^4(\R^2)}^2 \| \eta \|_{L^2(\R^2)} .
\eeq
Since either $ \|\eta \| _{L^2( \R^2)} ^2\leq \ell_1$ or 
$  \|\nabla U \| _{L^2( \R^2)} ^2\leq M_1$ and then, by (\ref{7.17}), 
 $ \|\eta \| _{L^2( \R^2)} ^2\leq \frac{ M_1}{a_1^2}$, we get
\beq
\label{7.32}
 \Big| \int_{\R^2} \mathcal{O} ( |\eta|^3 )\ dx \Big| 
\leq C \| \eta \|_{L^3(\R^2)}^3 
\leq C \| \eta \|_{L^2(\R^2)} \| \eta \|_{L^4(\R^2)}^2 
\leq C \| \eta \|_{L^4(\R^2)}^2. 
\eeq
Recall that  $1 - \beta_* \leq \rho \leq 1 + \beta_*$ and using step 1 we find
$$ \Big| \int_{\R^2} \rho \eta ( \p_{x_2} \theta)^2 \ dx \Big| 
\leq C  \int_{\R^2} ( \p_{x_2} \theta)^2 \ dx 
\leq C \| \eta \|_{L^4(\R^2)}^4 .$$
 Since $ c\in [0,v_s]$, from step 1 and the Cauchy-Schwarz inequality we obtain
$$ \Big| \int_{\R^2} c \rho \eta h \ dx \Big| 
\leq C \| \eta \|_{L^2(\R^2)} \| h \|_{L^2(\R^2)} 
\leq C \| \eta \|_{L^4(\R^2)}^2 , $$
 Using the definition of $h$, step 1 and (\ref{7.32}) we now  estimate 
\beq
\label{cool}
\begin{array}{l}
\ds \Big| \ii_{\R^2}  \rho \eta ( \p_{x_1} \theta)^2 \ dx \Big| 
 \leq   C \ii_{\R^2} |\eta| \left( h  -  \frac{c}{2} \eta \right)^2 \ dx 
\\
\\
\ds \leq   C \| \eta \|_{L^\infty(\R^2)} \int_{\R^2} h^2 \ dx 
+ C \int_{\R^2} |\eta|^3 \ dx 
\leq C \| \eta \|_{L^4(\R^2)}^2 . 
\end{array}
\eeq
Summing up  the above estimates and using (\ref{eqlabousse})  yields the conclusion.

\medskip 

In steps 1 and  2  we have not used the fact that $ F''(1) = 3 $. 
Since $ g(s) = \frac{v_s^2}{2 } s^2 + \mathcal{O}(s^3)$ when 
$  F''(1) = 3 $, it is natural to write (\ref{magique})  in the form
\begin{align}
\label{magic}
\wh{\eta} (\xi) = & \, 
- \mathcal{L}_c (\xi) \mathcal{F} \left( 
2  (\p_{x_1} \theta)^2 - \frac{v_s^2}{2 } \eta^2 \right)
\nonumber \\
& \, - \mathcal{L}_c (\xi) \mathcal{F} \left( 
2 \eta (\p_{x_1} \theta)^2 + 2 \rho^2 (\p_{x_2} \theta)^2 + 2 | \nabla \rho|^2 
- \left[ g(\eta) - \frac{v_s^2}{2 } \eta^2 \right] \right) 
\\ 
& \, - 2 c \mathcal{L}_c (\xi) \frac{\xi_2^2}{|\xi|^2} \mathcal{F}(\eta \p_{x_1} \theta) 
 + 2c \frac{\xi_1 \xi_2}{|\xi|^2} \mathcal{L}_c (\xi) \mathcal{F} ( \eta \p_{x_2} \theta ) 
. \nonumber 
\end{align}
where we recall that
$ \mathcal{L}_c (\xi)  $ is given by (\ref{7.15}).
We expect the term in the first line of (\ref{magic}) to be much smaller than quadratic. 
By the Riesz-Thorin Theorem we have 
$ \| \eta \|_{L^4(\R^2)} \leq C \| \wh{\eta} \|_{L^{4/3}(\R^2)} $.
 We will estimate the $L^{4/3} $ norm of all the terms in the right-hand side of (\ref{magic}) and we will show that 
they are bounded by  $C \| \eta \|_{L^4(\R^2)} ^2$.

\medskip 

{\it Step 3.} We have, for some constant $C$ depending only 
on $F$,
$$ 
\Big\| 2c \frac{\xi_1 \xi_2}{|\xi|^2} \mathcal{L}_c (\xi) \mathcal{F} 
( \eta \p_{x_2} \theta ) \Big\|_{L^{4/3}(\R^2)} \leq C \| \eta \|_{L^4(\R^2)}^2 .
$$
Indeed,  by the continuity of $ \mathcal{F} : L^1(\R^2) \lra L^\infty(\R^2) $ 
and the  Cauchy-Schwarz inequality one has 
\begin{align*}
\Big\| 2c \frac{\xi_1 \xi_2}{|\xi|^2} \mathcal{L}_c (\xi) 
\mathcal{F} ( \eta \p_{x_2} \theta ) \Big\|_{L^{4/3}(\R^2)} 
\leq & \, C \| \mathcal{F} ( \eta \p_{x_2} \theta ) \|_{L^\infty(\R^2)} 
\Big\| \frac{\xi_1 \xi_2}{|\xi|^2} \mathcal{L}_c (\xi) \Big\|_{L^{4/3}(\R^2)} 
\\ \leq & \, C \| \eta \p_{x_2} \theta \|_{L^1(\R^2)} 
\Big\| \frac{\xi_1 \xi_2}{|\xi|^2} \mathcal{L}_c (\xi) \Big\|_{L^{4/3}(\R^2)}  
\\ \leq & \, C \| \eta \|_{L^2(\R^2)} \| \p_{x_2} \theta \|_{L^2(\R^2)} 
\Big\| \frac{\xi_1 \xi_2}{|\xi|^2} \mathcal{L}_c (\xi) \Big\|_{L^{4/3}(\R^2)} 
\\ \leq & \, C \| \eta \|_{L^4(\R^2)}^2  
\Big\| \frac{\xi_1 \xi_2}{|\xi|^2} \mathcal{L}_c (\xi) \Big\|_{L^{4/3}(\R^2)} ,
\end{align*}
where we have used the estimate 
$\| \p _{x_2} \theta \|_{L^2(\R^2)} \leq C \| \eta \|_{L^4(\R^2)}$ (see Step 1) 
and the fact that $\| \eta \|_{L^2(\R^2)}$ is bounded.
Thus it suffices to prove that $ \Big\| \frac{\xi_1 \xi_2}{|\xi|^2} \mathcal{L}_c (\xi) \Big\|_{L^{4/3}(\R^2)} $ is bounded independently of $c$.

Using polar coordinates, we find for all $ q >1$, 
$$
\begin{array}{l}
\ds \| \mathcal{L}_c (\xi) \|_{L^{q}(\R^2)} ^{q}
=   \ii_{\R^2} \frac{|\xi|^{2q} d \xi }{ ( |\xi|^4 + v_s^2 |\xi|^2 - c^2 \xi_1^2)^{q}} 
= 4 \ii_0^{\pi/2} \int_0^{+\infty} \frac{ r\, d r \, d \vartheta }{ 
( r^2 + v_s^2  - c^2 \cos^2\vartheta ) ^{q}} 
\\ 
\\
=   \frac{2}{q-1}
 {\ds \ii _0^{\pi/2}  \frac{ d \vartheta }{ ( v_s^2  - c^2 \cos^2 \vartheta ) ^{q-1}}  }
\leq \frac{2}{q-1} 
{\ds  \ii_0^{\pi/2}  \frac{ d \vartheta }{ ( v_s^2  - v_s^2 \cos^2 \vartheta ) ^{q-1}} }
= \frac{2}{(q-1)v_s^{2(q-1)}}
 {\ds \ii_0^{\pi/2} \frac{ d \vartheta }{ (\sin \vartheta )^{2(q-1)}} .}
\end{array}
$$
Since the last integral is finite and does not depend on $c$ if $ 2(q-1) < 1$, we get 
\beq
\label{enorme43}
\sup_{0 \leq c \leq v_s} \| \mathcal{L}_c (\xi) \|_{L^{q}(\R^2)} \leq C_q < \infty 
\qquad \mbox{ for any } q \in \left(1, \frac 32\right) .
\eeq
In particular we have 
$ \Big\| \frac{\xi_1 \xi_2}{|\xi|^2} \mathcal{L}_c (\xi) \Big\|_{L^{4/3}(\R^2)} 
\leq \| \mathcal{L}_c (\xi) \|_{L^{4/3}(\R^2)} \leq C_{\frac 43} $
  for $ 0 \leq c \leq v_s $ and this concludes the proof of step 3.

\medskip 

 {\it Step 4.} There holds
$$ 
\Big\| 2c \frac{ \xi_2^2}{|\xi|^2} \mathcal{L}_c (\xi) \mathcal{F} ( \eta \p_{x_1} \theta ) 
\Big\|_{L^{4/3}(\R^2)} \leq C \| \eta \|_{L^4(\R^2)}^2 . 
$$
From the definition of $ h $ we have 
$ \eta \p_{x_1} \theta =  \eta h - \frac{c \eta^2}{2 } $, thus
$$ 
\Big\| 2c \frac{ \xi_2^2}{|\xi|^2} \mathcal{L}_c (\xi) 
\mathcal{F} ( \eta \p_{x_1} \theta ) \Big\|_{L^{4/3}(\R^2)} 
\leq C \Big\| \frac{ \xi_2^2}{|\xi|^2} \mathcal{L}_c (\xi) \mathcal{F} ( \eta^2 ) 
\Big\|_{L^{4/3}(\R^2)} + C \| \mathcal{L}_c (\xi) \mathcal{F} ( \eta h ) \|_{L^{4/3}(\R^2)} . $$
The second term is estimated as in Step 3, using (\ref{enorme43}), step 1 and the fact that $\| \eta \|_{L^2(\R^2)}$ is bounded:
$$ \| \mathcal{L}_c (\xi) \mathcal{F} ( \eta h ) \|_{L^{4/3}(\R^2)} 
\leq \| \mathcal{L}_c (\xi) \|_{L^{4/3}(\R^2)} 
\| \mathcal{F} ( \eta h ) \|_{L^{\infty}(\R^2)}
\leq C \| \eta \|_{L^2(\R^2)} \| h \|_{L^2(\R^2)} \leq C \| \eta \|_{L^4(\R^2)}^2 . $$
For the first term  we first observe that, since $ c^2 \leq v_s^2 $,
$$ \Big| \frac{ \xi_2^2}{|\xi|^2} \mathcal{L}_c (\xi) \Big| = 
\frac{\xi_2^{2} }{ |\xi|^4 + v_s^2 |\xi|^2 - c^2 \xi_1^2} 
\leq \frac{ \xi_2^2 }{ v_s^2 |\xi|^2 - v_s^2 \xi_1^2} = \frac{1}{v_s^2} . $$
Hence, using the estimate $ \| f\|_{L^4} \leq \|f\|_{L^{\infty}}^{\frac 23} \|f\|_{L^{4/3}}^{\frac 13}$, 
we get for $0 \leq c \leq v_s$,
$$
 \Big\| \frac{ \xi_2^2}{|\xi|^2} \mathcal{L}_c (\xi) \Big\|_{L^4(\R^2)} 
\leq \frac{1}{v_s^{4/3}} \Big\| \frac{ \xi_2^2}{|\xi|^2} \mathcal{L}_c (\xi) \Big\|_{L^{4/3}(\R^2)} ^{\frac 13}
\leq \frac{1}{v_s^{4/3}} \| \mathcal{L}_c (\xi) \|_{L^{4/3}(\R^2)} ^{\frac 13}
\leq C.
$$ 
(Warning:  $ \mathcal{L}_c $ is not uniformly bounded in $ L^4(\R^2) $ as $ c \ra v_s $.)
As a consequence, using the generalized H\"older inequality with 
$ \frac{1}{4/3} = \frac{1}{4} + \frac{1}{2} $ and the Plancherel formula,
$$ \Big\| \frac{ \xi_2^2}{|\xi|^2} \mathcal{L}_c (\xi) \mathcal{F} ( \eta^2 ) 
\Big\|_{L^{4/3}(\R^2)} 
\leq \Big\| \frac{ \xi_2^2}{|\xi|^2} \mathcal{L}_c (\xi) \Big\|_{L^4(\R^2)} 
\| \mathcal{F} ( \eta^2 ) \|_{L^{2}(\R^2)} 
\leq 
C \| \eta \|_{L^4(\R^2)}^2 . $$
Combining the above estimates gives the desired   conclusion.

\medskip 

{\it Step 5.} If $F''(1) = 3$  we have
$$ \Big\| \mathcal{L}_c (\xi) \mathcal{F} \left( 
2 \eta (\p_{x_1} \theta)^2 + 2 \rho^2 (\p_{x_2} \theta)^2 + 2 | \nabla \rho|^2 
- \left[ g(\eta) - \frac{v_s^2}{2 } \eta^2 \right] \right) 
\Big\|_{L^{4/3}(\R^2)} \leq C \| \eta \|_{L^4(\R^2)}^2 . $$

By (\ref{enorme43}) and   the inequality
$$ \| \mathcal{L}_c (\xi) \mathcal{F}( H ) \|_{L^{4/3}(\R^2)}  \leq
\| \mathcal{L}_c (\xi) \|_{L^{4/3}(\R^2)}  \|\mathcal{F}( H ) \|_{L^{\infty}(\R^2)} 
\leq  C_{\frac43} \| H \|_{L^1(\R^2)} 
$$ 
it suffices to estimate
$$
 \Big\| 2 \eta (\p_{x_1} \theta)^2 + 2 \rho^2 (\p_{x_2} \theta)^2 
+ 2 | \nabla \rho|^2 - \left[ g(\eta) - \frac{v_s^2}{2 } \eta^2 \right] 
\Big\|_{L^{1}(\R^2)} . 
$$
We estimate each term separately.
We have already seen that $ g(s) = \frac{v_s^2}{2 } s^2 + \mathcal{O}(s^3)$ 
as $s \ra 0$ because $F''(1) = 3$. By    (\ref{eta3}) we obtain
$$ \Big\| g(\eta) - \frac{v_s^2}{2 } \eta^2 \Big\|_{L^{1}(\R^2)} 
\leq C \| \eta \|_{L^3(\R^2)}^3 \leq C \| \eta \|_{L^4(\R^2)}^2 .
$$
 From step 2  we have
$$ \Big\| |\nabla \rho |^2 \Big\|_{L^{1}(\R^2)} 
= \int_{\R^2} | \nabla \rho |^2 \ dx \leq C \| \eta \|_{L^4(\R^2)}^2 $$
and from step 1  we get
$$ \Big\| \rho^2 (\p_{x_2} \theta)^2 \Big\|_{L^{1}(\R^2)} 
\leq C \int_{\R^2} ( \p_{x_2} \theta )^2 \ dx \leq C \| \eta \|_{L^4(\R^2)}^4
\leq C \| \eta \|_{L^4(\R^2)}^2 . $$
Finally,  as in (\ref{cool})    we infer  that
$$ \Big\| \eta (\p_{x_1} \theta )^2 \Big\|_{L^{1}(\R^2)} 
\leq C \| \eta \|_{L^4(\R^2)}^2 . $$
Gathering the above estimates we get the conclusion.

\medskip 

{\it Step 6.} The following estimate holds:
$$ \Big\| \mathcal{L}_c (\xi) \mathcal{F} \left( 
2  (\p_{x_1} \theta)^2 - \frac{v_s^2}{2 } \eta^2 \right) 
\Big\|_{L^{4/3}(\R^2)} \leq C \| \eta \|_{L^4(\R^2)}^2 .
$$
Indeed, arguing as is step 5 and using the definition of $h$, the Cauchy-Schwarz inequality  and step 1 we deduce
$$
\begin{array}{l}
\ds  \Big\| \mathcal{L}_c (\xi) \mathcal{F} \left( 
2  (\p_{x_1} \theta)^2 - \frac{v_s^2}{2 } \eta^2 \right) 
\Big\|_{L^{4/3}(\R^2)} 
\leq C_{\frac43} \Big\|   2 \left( h - \frac{ c}{2} \eta \right)^2 -   \frac{v_s^2}{2 } \eta^2 \Big\|_{L^1(\R^2)} 
\\
\\
\ds
\leq C \| h^2 \|_{L^1(\R^2)} + C \| \eta \|_{L^2(\R^2)} \| h \|_{L^2(\R^2)} + 
\frac{v_s ^2 - c^2}{2 } \ii_{\R^2} \eta ^2 \, dx 
\leq C \| \eta \|_{L^4(\R^2)}^2 .
\end{array}
$$

\medskip 
{\it Step 7.}    Conclusion. 


Using  the Riesz-Thorin 
theorem, we have $ \| \eta \|_{L^{4}(\R^2)} \leq C \| \hat{\eta} \|_{L^{4/3}(\R^2)} $. 
Coming back to (\ref{magic}) and  gathering the estimates in steps 3-6, we deduce
$$ 
\| \eta \|_{L^{4}(\R^2)} \leq C \| \hat{\eta} \|_{L^{4/3}(\R^2)} 
\leq C \| \eta \|_{L^{4}(\R^2)}^2,
 $$
where $C$ depends only on $F$.
Consequently, either $\| \eta \|_{L^{4}(\R^2)} = 0 $, or there is a constant $\kappa > 0$ such that  $\| \eta \|_{L^{4}(\R^2)} \geq \kappa$. 
If $\| \eta \|_{L^{4}(\R^2)} = 0 $ we have $ \eta = 0 $ a.e. and from (\ref{7.17}) we get 
$\| \nabla U \|_{L^{2}(\R^2)} = 0 $, hence $U$ is constant.
If  $\| \eta \|_{L^{4}(\R^2)} \geq \kappa$, (\ref{tookiki})  implies that there are $ \ell_*> 0$ 
and $ k_* > 0$ such that  $\| \eta \|_{L^{2}(\R^2)} \geq \ell_* $ and  $\| \nabla U \|_{L^{2}(\R^2)} \geq  k_*. $ 
The proof of Proposition \ref{smallE} is complete.
\hfill
$\Box$

\bigskip

\noindent
{\bf Acknowledgement. } We acknowledge the support of the French ANR (Agence 
Nationale de la Recherche) under Grant ANR-09-JCJC-0095-01-ArDyPitEq.

\end{document}